%% file: dHascoAss.tex
\documentclass{smf-l}

\usepackage{a4wide,amssymb,amsrefs,bbm,calligra,comment,graphicx,mathrsfs,runumber,txfonts}

\usepackage[2cell,arrow,curve,matrix,web]{xy}
\UseHalfTwocells
\UseTwocells
\newdir{ >}{{}*!/-10pt/\dir{>}}
\newdir{ (}{{}*!/-10pt/\dir^{(}}
\newdir{ )}{{}*!/-10pt/\dir_{(}}

\def\calli#1{\textup{\!\textcalligra{#1}\,}}

\def\A{{\mathscr A}}
\def\bA{{\mathbf A}}
\def\ba{{\mathbf a}}
\def\B{{\mathbf B}}
\def\bb{{\mathbf b}}
\def\aB{{\mathscr B}}
\def\bC{{\mathbf C}}
\def\E{{\mathrm E}}
\def\F{{\mathscr F}}
\def\FF{{\mathbb F}}
\def\G{{\mathbb G}}
\def\GG{{\mathbb G}}

\def\H{{\mathscr H}}
\def\L{{\mathrm L}}
\def\LL{{\mathbb L}}
\def\M{{\mathscr M}}
\def\NN{{\mathbb N}}
\def\P{{\mathrm P}}
\def\QQ{{\mathbb Q}}
\def\rR{{\mathrm R}}
\def\R{{\mathscr R}}
\def\RR{{\mathbb R}}
\def\V{{\mathscr V}}
\def\W{{\mathscr W}}
\def\Z{{\mathscr Z}}
\def\ZZ{{\mathbb Z}}

\def\Topt{{\mathbf{Top}^*}}
\def\Sp{{\mathbf{Spec}}}

\def\op{^{\mathrm{op}}}
\def\ho{_\simeq}
\def\1{^{-1}}
\def\comp{{\scriptstyle\Box}}

\def\Pair{{\calli{Pair}}}
\def\Mod{{\mathbf{Mod}}}
\def\pair{{\mathit{pair}}}
\def\PAR{{\mathrm{PAR}}}
\def\pre{{\mathrm{pre}}}

\DeclareMathOperator{\Aut}{Aut}
\DeclareMathOperator{\Ext}{Ext}
\DeclareMathOperator{\Tor}{Tor}
\DeclareMathOperator{\Cotor}{Cotor}
\DeclareMathOperator{\Der}{Der}
\DeclareMathOperator{\Ider}{Ider}
\DeclareMathOperator{\Hom}{Hom}
\DeclareMathOperator{\im}{im}
\DeclareMathOperator{\coker}{coker}
\DeclareMathOperator{\Sq}{Sq}
\DeclareMathOperator{\Tr}{Tr}
\DeclareMathOperator{\id}{id}
\DeclareMathOperator{\cote}{\Box}

\def\xto#1{\xrightarrow[]{#1}}
\def\xot#1{\xleftarrow[]{#1}}
\def\set#1{{\left\{#1\right\}}}

\def\hog#1{\left\llbracket#1\right\rrbracket}

\def\brk#1{\left\langle#1\right\rangle}
\def\alignbox#1{\begin{aligned}#1\end{aligned}}
\def\smat#1{{\left(\begin{smallmatrix}#1\end{smallmatrix}\right)}}
\def\restr#1#2{\left.#1\right|_{#2}}
\def\dblb#1{{\bar{\bar{#1}}}}
\let\qf\dblb
\def\Alg{{\mathbf{Alg}}}
\def\Coalg{{\mathbf{Coalg}}}
\def\ESL{{\mathrm{ESL}}}

\def\oo{{\mathbbm1}}

\let\ge\geqslant
\let\le\leqslant
\let\d\partial
\let\x\times
\let\ox\otimes
\let\ot\leftarrow
\let\then\Rightarrow
\let\onto\twoheadrightarrow
\let\into\rightarrowtail

\let\incl\hookrightarrow
\let\b b
\let\l\ell
\let\r r
\let\i\iota
\let\ph\varphi
\let\k\varkappa
\let\eps\varepsilon
\let\ro\varrho

\let\bref\eqref

\makeindex

\theodef{Theorem}
\theodef{Lemma}
\theodef{Proposition}
\theodef{Corollary}
\theodef{Conjecture}

\newtheorem{Theocite}[subsection]{Theorem \cite{Baues}}

\theoremstyle{definition}
\theodef{Definition}
\theodef{Definitions}
\theodef{Example}
\theodef{Examples}

\theoremstyle{remark}
\theodef{Remark}

\numberwithin{section}{chapter}

\begin{document}

\frontmatter

\title{Dualization of the Hopf algebra of secondary cohomology operations and
the Adams spectral sequence}

\author{Hans-Joachim Baues}
\address{Max-Planck-Institut f\"ur Mathematik, Vivatsgasse 7, D-53111 Bonn, Germany}
\email{baues@mpim-bonn.mpg.de}

\author{Mamuka Jibladze}
\address{Razmadze Mathematical Institute, Alexidze st. 1, Tbilisi 0193, Georgia}
\curraddr{Max-Planck-Institut f\"ur Mathematik, Vivatsgasse 7, D-53111 Bonn, Germany}
\email{jib@rmi.acnet.ge}

%\thanks{The second author was supported in part by the RTN Network HPRN-CT-2002-00287}

\date{\today}

\begin{abstract}
We describe the dualization of the algebra of secondary cohomology operations in terms of generators extending the Milnor dual of the Steenrod algebra. In
this way we obtain explicit formul\ae\ for the computation of the $\E_3$-term of the Adams spectral sequence converging to the stable homotopy groups of spheres.
\end{abstract}

\maketitle

\tableofcontents

\include{dHascoAss_p}

\mainmatter
\include{dHascoAss_c1}

\include{dHascoAss_c2}

\include{dHascoAss_c3}

\include{dHascoAss_c4}

\include{dHascoAss_c5}

\include{dHascoAss_c6}

\include{dHascoAss_c7}

\include{dHascoAss_c8}

\include{dHascoAss_c9}

\backmatter

\include{dHascoAss_b}
\include{dHascoAss_i}
\end{document}

%% file: dHascoAss_p.tex
\chapter*{Introduction}

Spheres are the most elementary compact spaces, but the simple question of
counting essential maps between spheres turned out to be a landmark
problem. In fact, progress in algebraic topology might be measured by its
impact on this question. Topologists have worked on the problem of describing the
homotopy groups of spheres for around 80 years and there is still no
satisfactory solution in sight. Many approaches have been developed: a
distinguished one is the Adams spectral sequence
$$
\E_2,\E_3,\E_4,...
$$
converging to homotopy groups of spheres. Adams computed the $\E_2$-term and
showed that
$$
\E_2=\Ext_\A(\FF,\FF)
$$
is algebraically determined by $\Ext$-groups associated to the Steenrod
algebra $\A$. Hence $\E_2$ is an upper bound for homotopy groups of spheres
and is given by an algebraic resolution of the prime field $\FF=\FF_p$ over
the algebra $\A$. The Steenrod algebra $\A$ is in fact a Hopf algebra with
wonderful algebraic properties. Milnor showed that the dual algebra
$$
\A_*=\Hom(\A,\FF)
$$
is a polynomial algebra. Topologically the Steenrod algebra is the algebra of
primary cohomology operations. Adams' formula on $E_2$ shows a 
fundamental connection between homotopy groups of spheres and primary 
cohomology operations. Much work in the literature is exploiting 
this connection. However, since $E_2$ is only an upper bound, one cannot 
expect the Steenrod algebra to be sufficient to determine homotopy groups of 
spheres. In fact, for this the ``algebra of all higher cohomology operations"
is needed. The structure of this total algebra is highly unknown; it is not even 
clear what kind of algebra is needed to describe the additive properties of 
higher cohomology operations. The structure of the Adams spectral sequence 
$E_2, E_3, \ldots$ shows that the total algebra can be approximated by
constructing inductively primary, secondary, tertiary \ldots operations. In doing
so one might be able to grasp the total algebra. This is the program of computing 
homotopy groups of spheres via higher cohomology operations. The first step
beyond Adams' result is understanding the algebra of secondary
cohomology operations which, surprisingly, turned out to be a differential
algebra, namely a pair algebra.

In the book \cite{Baues} the pair algebra $\aB$ of secondary cohomology
operations is computed and this enriches the known algebraic structure of the
Steenrod algebra considerably. The pair algebra $\aB$ is given by an exact
sequence
$$
\begin{aligned}
\xymatrix@1{
\Sigma\A\ar@{ >->}[r]&\aB_1\ar[r]^\d&\aB_0\ar@{->>}[r]^q&\A.
}
\end{aligned}
\eqno{(*)}
$$
Here $\aB_0$ is the free associative algebra over $\GG=\ZZ/p^2\ZZ$ generated
by the Steenrod operations which also generate $\A$ and $q$ is the identity on
generators. Moreover there is a multiplication map
$$
m:\aB_0\!\ox\!\aB_1\oplus\aB_1\!\ox\!\aB_0\to\aB_1
$$
and a diagonal map
$$
\Delta:\aB_1\to(\aB_0\!\ox\!\aB_1\oplus\aB_1\!\ox\!\aB_0)/\sim
$$
such that $\aB=(\aB,m,\Delta)$ is a ``secondary Hopf algebra'', see
\cite{Baues}, inducing the Hopf algebra structure of the Steenrod algebra
$\A$. It is proven in \cite{Baues} that the structure of $\aB$ as a secondary 
Hopf algebra together with the explicit invariants $L$ and $S$ determines $\aB$ 
up to isomorphism. The nature of secondary homotopy operations leads
forcibly to this kind of new algebraic object which has wonderful properties
shedding light on the structure of the Steenrod algebra $\A$ as a Hopf algebra.
By a striking result of Milnor, the dual $\A_{\ast}$ of the Hopf algebra $\A$
is a polynomial algebra with a nice diagonal which, for many purposes, 
is easier to deal with than the algebra $\A$ itself which is given by generators, 
the Steenrod squares, and Adem relations. Thus this paper also describes
the dualization $\aB_{\ast}$ of the secondary Hopf algebra $\aB$. We compute the 
invariants dual to $L$ and $S$ by explicit and easy formul\ae. Therefore computations
in terms of $\aB$ can equivalently be carried out in terms of the dual $\aB_{\ast}$ and
often the dual formul{\ae} are easier to handle. In this paper we use the secondary 
Hopf algebra $\aB$ and its dual $\aB_{\ast}$ for computating a 
secondary resolution which determines the differential $d_{(2)}$ on $E_2$ and 
hence $E_3$.

%Since the work of Adams \cites{AdamsH,AdamsS} it is generally believed that
%the $\E_3$-term of the Adams spectral sequence can be computed by the use of
%secondary cohomology operations. This was done by
Adams %for
computed
those special values of the differentials $d_{(2)}$ in $\E_2$ which are related
to the Hopf invariant 1 problem. In the book of Ravenel \cite{Ravenel} one
finds a list of all differentials up to degree 60 which, however, is only
tentative in degrees $\ge46$. Corrections of published differentials in low
degrees were made by Bruner \cite{Brunercorr}. An explicit method for
computing the differential $d_{(2)}$ in general, however, has not been achieved in the
literature. But it is done in the present paper. Our result is thus showing the 
global computable nature of the $E_3$--term of the Adams spectral sequence.
According to Ravenel's observer, ``who looks to the far distant homotopy groups
of spheres through a telescope,'' such a global result on $E_3$ seemed 
impossible for a long time.

We show that the differential $d_{(2)}$ and the $\E_3$-term can be completely
computed by the formula
$$
\E_3=\Ext_\aB(\GG^\Sigma,\GG^\Sigma)
$$
where the secondary $\Ext$-groups $\Ext_\aB$ are given by an algebraic
secondary resolution associated to the pair algebra $\aB$. The computation of
$\E_3$ yields a new algebraic upper bound of homotopy groups of spheres
improving the Adams bound given by $\E_2$.

In order to do explicit computations of the new bound $\E_3$ one has to carry
out two tasks. On the one hand one has to describe the algebraic structure of
the secondary Hopf algebra $\aB$ explicitly by equations which a computer can
deal with in an easy way. On the other hand one has to choose a secondary
resolution associated to $\aB$, by solving inductively a system of explicit
equations determined by $\aB$.

In the first part (chapters \ref{sext}, \ref{secsteen}, \ref{E3}) of this
paper we describe the algebra which yields the secondary resolution associated
to $\aB$ and which determines the differential $d_{(2)}$ on $\E_2$ by the
resolution. In the second part (chapters \ref{Hpa}, \ref{gens}, \ref{LS},
\ref{xi}, \ref{A}) we study the algebraic properties of $\aB$ and of the
dualization of $\aB$. In particular we show that the results of Milnor on the
dual Steenrod algebra $\A_*$ have secondary analogues. For the dualization of
$\aB$ we proceed as follows. The projection $q:\aB_0\onto\A$ in $(*)$ above
admits a factorization
$$
q:\aB_0\onto\F_0\onto\A
$$
where $\F_0=\aB_0\ox\FF$ is the free associative algebra over $\FF=\ZZ/p\ZZ$
generated by the Steenrod operations. Now let
\begin{align*}
R_\aB&=\textrm{kernel}(\aB_0\to\A)\\
R_\F&=\textrm{kernel}(\F_0\to\A).
\end{align*}
Then one has an exact sequence of $\FF$-vector spaces
$$
\A\into R_\aB\ox\FF\onto R_\F
$$
which can be dualized by applying the functor $\Hom(-,\FF)$. Moreover the
exact sequence of $\FF$-vector spaces
$$
\Sigma\A\into\aB_1\ox\FF\onto R_\aB\ox\FF
$$
can be dualized by $\Hom(-,\FF)$. The main results of this work describe in detail the multiplication in $\aB$ and the diagonal in $\aB$ on the level of $\aB_1\ox\FF$ and on the dual $\Hom(\aB_1,\FF)$. In this way we obtain explicit formul\ae\ describing the algebraic structure of $\aB$ and of the dual of $\aB$. Of course the dual of $\aB$ determines $\aB$ and vice versa.

We use these formul\ae\ for computer calculations of the secondary resolution associated to $\aB$ and we derive in this way the differentials $d_{(2)}$ on $\E_2$. In section \ref{d2} we do such computations up to degree 40 in order to confirm the algebraic equations achieved in the book \cite{Baues}. The goal is to compute $\E_3$ up to degree 210 as this was done for $\E_2$ by Nassau \cite{Nassau}. A more effective computer implementation of $\E_3$, which is left to the interested reader, relies on the computation of the dual of $\aB$, see the formul{\ae} in section \ref{cobcomp} below. The functions needed for the implementation are described in the paper by tables of values in low degrees. These tables should be helpful to control the implementation.

%\aufm{Hans-Joachim Baues,\\Mamuka Jibladze}

%% file: dHascoAss_c1.tex
\chapter{Secondary $\Ext$-groups associated to pair algebras}\label{sext}

In this chapter we introduce algebraically secondary $\Ext$-groups $\Ext_B$
over a pair algebra $B$. In \cite{BJ5} we already studied secondary
$\Ext$-groups in an additive track category which yield the $\Ext$-groups
$\Ext_B$ as a special case if one considers the track category of
$B$-modules. In chapter \ref{E3} we shall see thet the E$_3$-term of the
Adams spectral sequence is given by secondary $\Ext$-groups over the pair
algebra $\aB$ of secondary cohomology operations.

\section{Modules over pair algebras}
\label{secodjf}
We here recall from \cite{Baues} the notion of pair modules, pair
algebras, and pair modules over a pair algebra $B$. The category
$B$-$\Mod$ of pair modules over $B$ is an additive track category in which
we consider secondary resolutions as defined in \cite{BJ5}. Using such
secondary resolutions we shall obtain the secondary derived functors
$\Ext_B$ in section \ref{secodif}.

Let $k$ be a commutative ring with unit and let $\Mod$ be the category of
$k$-modules and $k$-linear maps. This is a symmetric monoidal category via
the tensor product $A\!\ox\!B$ over $k$ of $k$-modules $A$, $B$. A \emph{pair}
of modules is a morphism
\begin{equation}
X=\left(X_1\xto\d X_0\right)
\end{equation}
in $\Mod$. We write $\pi_0(X)=\coker\d$ and $\pi_1(X)=\ker\d$. A
\emph{morphism} $f:X\to Y$ of pairs is a commutative diagram
$$
\xymatrix
{
X_1\ar[r]^{f_1}\ar[d]_\d&Y_1\ar[d]^\d\\
X_0\ar[r]^{f_0}&Y_0.
}
$$
Evidently pairs with these morphisms form a category $\Pair(\Mod)$ and
one has functors
$$
\pi_1, \pi_0 : \Pair(\Mod)\to\Mod.
$$
A pair morphism is called a
\emph{weak equivalence} if it induces isomorphisms on $\pi_0$ and $\pi_1$.

Clearly a pair in $\Mod$ coincides with a chain complex concentrated in
degrees 0 and 1. For two pairs $X$ and $Y$ the tensor product of the
complexes corresponding to them is concentrated in degrees in 0, 1 and 2
and is given by
$$
X_1\!\ox\!Y_1\xto{\d_1}X_1\!\ox\!Y_0\oplus
X_0\!\ox\!Y_1\xto{\d_0}X_0\!\ox\!Y_0
$$
with $\d_0=(\d\ox1,1\ox\d)$ and $\d_1=(-1\ox\d,\d\ox1)$. Truncating $X\ox Y$
we get the pair
$$
X\bar\otimes Y=
\left((X\bar\otimes Y)_1=\coker(\d_1)\xto\d X_0\ox Y_0=(X\bar\otimes Y)_0\right)
$$
with $\d$ induced by $\d_0$.

\begin{Remark}\label{trunc}
Note that the full embedding of the category of pairs into the category of
chain complexes induced by the above identification has a left adjoint
$\Tr$ given by truncation: for a chain complex
$$
C=\left(...\to C_2\xto{\d_1}C_1\xto{\d_0}C_0\xto{\d_{-1}}C_{-1}\to...\right),
$$
one has
$$
\Tr(C)=\left(\coker(\d_1)\xto{\bar{\d_0}}C_0\right),
$$
with $\bar{\d_0}$ induced by $\d_0$. Then clearly one has
$$
X\bar\otimes Y=\Tr(X\ox Y).
$$
Using the fact that $\Tr$ is a reflection onto a full subcategory, one
easily checks that the category $\Pair(\Mod)$ together with the tensor
product $\bar\otimes$ and unit $k=(0\to k)$ is a symmetric monoidal category,
and $\Tr$ is a monoidal functor.
\end{Remark}

We define the tensor product $A\ox B$ of two graded modules in the usual
way, i.~e. by
$$
(A\ox B)^n=\bigoplus_{i+j=n}A^i\ox B^j.
$$

A \emph{pair module} is a graded object of $\Pair(\Mod)$, i.~e. a
sequence $X^n=(\d:X_1^n\to X_0^n)$ of pairs in $\Mod$. We identify such a
pair module $X$ with the underlying morphism $\d$ of degree 0 between
graded modules
$$
X=\left(X_1\xto\d X_0\right).
$$
Now the tensor product $X\bar\otimes Y$ of graded pair modules $X$, $Y$ is defined by 
\begin{equation}\label{grapr}
(X\bar\otimes Y)^n=\bigoplus_{i+j=n}X^i\bar\otimes Y^j.
\end{equation}
This defines a monoidal structure on the category of graded pair modules.
Morphisms in this category are of degree 0.

For two morphisms $f,g:X\to Y$ between graded pair modules, a
\emph{homotopy} $H:f\then g$ is a morphism $H:X_0\to Y_1$ of degree 0 as in
the diagram
\begin{equation}\label{homoto}
\alignbox{
\xymatrix
{
X_1\ar@<.5ex>[r]^{f_1}\ar@<-.5ex>[r]_{g_1}\ar[d]_\d&Y_1\ar[d]^\d\\
X_0\ar@<.5ex>[r]^{f_0}\ar@<-.5ex>[r]_{g_0}\ar[ur]|H&Y_0,
}}
\end{equation}
satisfying $f_0-g_0=\d H$ and $f_1-g_1=H\d$.

A \emph{pair algebra} $B$ is a monoid in the monoidal category of graded pair
modules, with multiplication
$$
\mu:B\bar\otimes B\to B.
$$
We assume that $B$ is concentrated in nonnegative degrees, that is $B^n=0$
for $n<0$.

A \emph{left $B$-module} is a graded pair module $M$ together with a left
action
$$
\mu:B\bar\otimes M\to M
$$
of the monoid $B$ on $M$.

More explicitly pair algebras and modules over them can be described as
follows.

\begin{Definition}
A \emph{pair algebra} $B$ is a graded pair
$$
\d:B_1\to B_0
$$
in $\Mod$ with $B_1^n=B_0^n=0$ for $n<0$ such that $B_0$ is a graded
algebra in $\Mod$, $B_1$ is a graded $B_0$-$B_0$-bimodule, and $\d$ is a
bimodule homomorphism. Moreover for $x,y\in B_1$ the equality
$$
\d(x)y=x\d(y)
$$
holds in $B_1$.
\end{Definition}

It is easy to see that there results an exact sequence of graded
$B_0$-$B_0$-bimodules
$$
0\to\pi_1B\to B_1\xto\d B_0\to\pi_0B\to0
$$
where in fact $\pi_0B$ is a $k$-algebra, $\pi_1B$ is a
$\pi_0B$-$\pi_0B$-bimodule, and $B_0\to\pi_0(B)$ is a homomorphism of
algebras.

\begin{Definition}\label{bmod}
A \emph{(left) module} over a pair algebra $B$ is a graded pair
$M=(\d:M_1\to M_0)$ in $\Mod$ such that $M_1$ and $M_0$ are left
$B_0$-modules and $\d$ is $B_0$-linear. Moreover, a $B_0$-linear map
$$
\bar\mu:B_1\!\otimes_{B_0}\!M_0\to M_1
$$
is given fitting in the commutative diagram
$$
\xymatrix{
B_1\otimes_{B_0}M_1\ar[r]^{1\ox\d}\ar[d]_\mu&B_1\otimes_{B_0}M_0\ar[dl]^{\bar\mu}\ar[d]^\mu\\
M_1\ar[r]_\d&M_0,
}
$$
where $\mu(b\ox m)=\d(b)m$ for $b\in B_1$ and $m\in M_1\cup M_0$.

For an indeterminate element $x$ of degree $n=|x|$ let $B[x]$ denote the
$B$-module with $B[x]_i$ consisting of expressions $bx$ with $b\in B_i$,
$i=0,1$, with $bx$ having degree $|b|+n$, and structure maps given by
$\d(bx)=\d(b)x$, $\mu(b'\ox bx)=(b'b)x$ and $\bar\mu(b'\ox bx)=(b'b)x$.

A \emph{free} $B$-module is a direct sum of several copies of modules of
the form $B[x]$, with $x\in I$ for some set $I$ of indeterminates of
possibly different degrees. It will be denoted
$$
B[I]=\bigoplus_{x\in I}B[x].
$$

For a left $B$-module $M$ one has the exact sequence of $B_0$-modules
$$
0\to\pi_1M\to M_1\to M_0\to\pi_0M\to0
$$
where $\pi_0M$ and $\pi_1M$ are actually $\pi_0B$-modules.

Let $B$-$\Mod$ be the category of left modules over the pair algebra $B$.
Morphisms $f=(f_0,f_1):M\to N$ are pair morphisms which are
$B$-equivariant, that is,$f_0$ and $f_1$ are $B_0$-equivariant and 
compatible with $\bar\mu$ above, i.~e. the diagram
$$
\xymatrix{
B_1\ox_{B_0}M_0\ar[r]^-{\bar\mu}\ar[d]_{1\ox f_0}&M_1\ar[d]^{f_1}\\
B_1\ox_{B_0}N_0\ar[r]^-{\bar\mu}&N_1
}
$$
commutes.

For two such maps $f,g:M\to N$ a track $H:f\then g$ is a degree zero map
\begin{equation}\label{track}
H:M_0\to N_1
\end{equation}
satisfying $f_0-g_0=\d H$ and $f_1-g_1=H\d$ such that $H$ is
$B_0$-equivariant. For tracks $H:f\then g$, $K:g\then h$ their composition
$K\comp H:f\then h$ is $K+H$.
\end{Definition}

\begin{Proposition}
For a pair algebra $B$, the category $B$-$\Mod$ with the above track
structure is a well-defined additive track category.
\end{Proposition}

\begin{proof}
For a morphism $f=(f_0,f_1):M\to N$ between $B$-modules, one has
$$
\Aut(f)=\set{H\in\Hom_{B_0}(M_0,N_1)\ |\ \d H=f_0-f_0,H\d=f_1-f_1}\cong\Hom_{\pi_0B}(\pi_0M,\pi_1N).
$$
Since this group is abelian, by \cite{Baues&JibladzeI} we know that $B$-$\Mod$
is a linear track extension of its homotopy category by the bifunctor $D$
with $D(M,N)=\Hom_{\pi_0B}(\pi_0M,\pi_1N)$. It thus remains to show that the
homotopy category is additive and the bifunctor $D$ is biadditive.

By definition the set of morphisms $[M,N]$ between objects $M$, $N$ in the
homotopy category is given by the exact sequence of abelian groups
$$
\Hom_{B_0}(M_0,N_1)\to\Hom_B(M,N)\onto[M,N].
$$
This makes evident the abelian group structure on the hom-sets $[M,N]$. Bilinearity
of composition follows from consideration of the commutative diagram
$$
\xymatrix{
\Hom_{B_0}(M_0,N_1)\!\ox\!\Hom_B(N,P)\oplus
\Hom_B(M,N)\!\ox\!\Hom_{B_0}(N_0,P_1)\ar[d]\ar[r]^-\mu
&\Hom_{B_0}(M_0,P_1)\ar[d]\\
\Hom_B(M,N)\ox\Hom_B(N,P)\ar[r]\ar@{->>}[d]
&\Hom_B(M,P)\ar@{->>}[d]\\
[M,N]\ox[N,P]\ar@{-->}[r]
&[M,P]
}
$$
with exact columns, where $\mu(H\!\ox\!g+f\!\ox\!K)=g_1H+Kf_0$. It also shows
that the functor $B$-$\Mod\to B$-$\Mod\ho$ is linear. Since this functor
is the identity on objects, it follows that the homotopy category is additive.

Now note that both functors $\pi_0$, $\pi_1$ factor to define functors on
$B$-$\Mod\ho$. Since these functors are evidently additive, it follows that
$D=\Hom_{\pi_0B}(\pi_0,\pi_1)$ is a biadditive bifunctor.
\end{proof}

\begin{Lemma}\label{freehom}
If $M$ is a free $B$-module, then the canonical map
$$
[M,N]\to\Hom_{\pi_0B}(\pi_0M,\pi_0N)
$$
is an isomorphism for any $B$-module $N$. 
\end{Lemma}

\begin{proof}
Let $(g_i)_{i\in I}$ be a free generating set for $M$. Given a
$\pi_0(B)$-equivariant homomorphism $f:\pi_0M\to\pi_0N$, define its lifting
$\tilde f$ to $M$ by specifying $\tilde f(g_i)=n_i$, with $n_i$ chosen
arbitrarily from the class $f([g_i])=[n_i]$.

To show monomorphicity, given $f:M\to N$ such that $\pi_0f=0$, this means
that $\im f_0\subset\im\d$, so we can choose $H(g_i)\in N_1$ in such a way
that $\d H(g_i)=f_0(g_i)$. This then extends uniquely to a $B_0$-module
homomorphism $H:M_0\to N_1$ with $\d H=f_0$; moreover any element of $M_1$ is a linear
combination of elements of the form $b_1g_i$ with $b_1\in B_1$, and for
these one has $H\d(b_1g_i)=H(\d(b_1)g_i)=\d(b_1)H(g_i)$. But
$f_1(b_1g_i)=b_1f_0(g_i)=b_1\d H(g_i)=\d(b_1)H(g_i)$ too, so $H\d=f_1$.
This shows that $f$ is nullhomotopic.
\end{proof}

\section{$\Sigma$-structure}

\begin{Definition}
The \emph{suspension} $\Sigma X$ of a graded object $X=(X^n)_{n\in\ZZ}$ is
given by degree shift, $(\Sigma X)^n=X^{n-1}$.
\end{Definition}

Let $\Sigma:X\to\Sigma X$ be the map of degree 1 given by the identity. If
$X$ is a left $A$-module over the graded algebra $A$ then $\Sigma X$ is a
left $A$-module via
\begin{equation}\label{suspact}
a\cdot\Sigma x=(-1)^{|a|}\Sigma(a\cdot x)
\end{equation}
for $a\in A$, $x\in X$. On the other hand if $\Sigma X$ is a right $A$-module then
$(\Sigma x)\cdot a=\Sigma(x\cdot a)$ yields the right $A$-module structure
on $\Sigma X$.

\begin{Definition}\label{sigmas}
A \emph{$\Sigma$-module} is a graded pair module $X=(\d:X_1\to X_0)$
equipped with an isomorphism
$$
\sigma:\pi_1X\cong\Sigma\pi_0X
$$
of graded $k$-modules. We then call $\sigma$ a \emph{$\Sigma$-structure} of
$X$. A $\Sigma$-map between $\Sigma$-modules is a map $f$ between pair
modules such that $\sigma(\pi_1f)=\Sigma(\pi_0f)\sigma$. If $X$ is a pair
algebra then a $\Sigma$-structure is an isomorphism of
$\pi_0X$-$\pi_0X$-bimodules. If $X$ is a left module over a pair algebra $B$
then a $\Sigma$-structure of $X$ is an isomorphism $\sigma$ of left
$\pi_0B$-modules. Let
$$
(B\textrm{-}\Mod)^\Sigma\subset B\textrm{-}\Mod
$$
be the track category of $B$-modules with $\Sigma$-structure and
$\Sigma$-maps.
\end{Definition}

\begin{Lemma}
Suspension of a $B$-module $M$ has by \eqref{suspact} the structure of a
$B$-module and $\Sigma M$ has a $\Sigma$-structure if $M$ has one.
\end{Lemma}

\begin{proof}
Given $\sigma:\pi_1M\cong\Sigma\pi_0M$ one defines a $\Sigma$-structure on
$\Sigma M$ via
$$
\pi_1\Sigma
M=\Sigma\pi_1M\xto{\Sigma\sigma}\Sigma\Sigma\pi_0M=\Sigma\pi_0\Sigma M.
$$
\end{proof}

Hence we get suspension functors between track categories
$$
\xymatrix{
B\textrm{-}\Mod\ar[r]^\Sigma&B\textrm{-}\Mod\\
(B\textrm{-}\Mod)^\Sigma\ar[u]\ar[r]^\Sigma&(B\textrm{-}\Mod)^\Sigma.\ar[u]
}
$$

\begin{Lemma}\label{add}
The track category $(B\mathrm{-}\Mod)^\Sigma$ is $\LL$-additive in the sense of
\cite{BJ5}, with $\LL=\Sigma\1$, as well as $\RR$-additive,
with $\RR=\Sigma$.
\end{Lemma}

\begin{proof}
The statement of the lemma means that the bifunctor
$$
D(M,N)=\Aut(0_{M,N})
$$
is either left- or right-representable, i.~e. there is an endofunctor
$\LL$, respectively $\RR$ of $(B$-$\Mod)^\Sigma$ and a binatural
isomorphism $D(M,N)\cong[\LL M,N]$, resp. $D(M,N)\cong[M,\RR N]$.

Now by \eqref{track}, a track in $\Aut(0_{M,N})$ is a $B_0$-module
homomorphism $H:M_0\to N_1$ with $\d H=H\d=0$; hence
$$
D(M,N)\cong\Hom_{\pi_0B}(\pi_0M,\pi_1N)\cong\Hom_{\pi_0B}(\pi_0\Sigma\1
M,\pi_0N)\cong\Hom_{\pi_0B}(\pi_0M,\pi_0\Sigma N).
$$
\end{proof}

\begin{Lemma}
If $B$ is a pair algebra with $\Sigma$-structure then each free $B$-module
has a $\Sigma$-structure.
\end{Lemma}

\begin{proof}
This is clear from the description of free modules in \ref{bmod}.
\end{proof}

\section{The secondary differential over pair algebras}\label{secodif}

For a pair algebra $B$ with a $\Sigma$-structure, for a $\Sigma$-module
$M$ over $B$, and a module $N$ over $B$ we now define the \emph{secondary differential}
$$
d_{(2)}:\Ext^n_{\pi_0B}(\pi_0M,\pi_0N)\to\Ext^{n+2}_{\pi_0B}(\pi_0M,\pi_1N).
$$
Here $d_{(2)}=d_{(2)}(M,N)$ depends on the $B$-modules $M$ and $N$ and is
natural in $M$ and $N$ with respect to maps in $(B\mathrm{-}\Mod)^\Sigma$. For the
definition of $d_{(2)}$ we consider secondary chain complexes and secondary
resolutions. In \cite{BJ5} such a construction was performed in the
generality of an arbitrary $\LL$-additive track category. We will first
present the construction of $d_{(2)}$ for the track category of pair
modules and then will indicate how this construction is a particular case
of the more general situation discussed in \cite{BJ5}.

\begin{Definition}\label{secs}
For a pair algebra $B$, a \emph{secondary chain complex} $M_\bullet$ in $B$-$\Mod$
is given by a diagram to be
$$
\xymatrix@!{
...\ar[r]
&M_{n+2,1}\ar[r]^{d_{n+1,1}}\ar[d]|\hole^>(.75){\d_{n+2}}
&M_{n+1,1}\ar[r]^{d_{n,1}}\ar[d]|\hole^>(.75){\d_{n+1}}
&M_{n,1}\ar[r]^{d_{n-1,1}}\ar[d]|\hole^>(.75){\d_n}
&M_{n-1,1}\ar[r]\ar[d]|\hole^>(.75){\d_{n-1}}
&...\\
...\ar[r]\ar[urr]^<(.3){H_{n+1}}
&M_{n+2,0}\ar[r]_{d_{n+1,0}}\ar[urr]^<(.3){H_n}
&M_{n+1,0}\ar[r]_{d_{n,0}}\ar[urr]^<(.3){H_{n-1}}
&M_{n,0}\ar[r]_{d_{n-1,0}}\ar[urr]
&M_{n-1,0}\ar[r]
&...\\
}
$$
where each $M_n=(\partial_n:M_{n,1}\to M_{n,0})$ is a $B$-module, each
$d_n=(d_{n,0},d_{n,1})$ is a morphism in $B$-$\Mod$, each $H_n$ is
$B_0$-linear and moreover the identities
\begin{align*}
d_{n,0}d_{n+1,0}&=\d_nH_n\\
d_{n,1}d_{n+1,1}&=H_n\d_{n+2}\\
\intertext{and}
H_nd_{n+2,0}&=d_{n,1}H_{n+1}
\end{align*}
hold for all $n\in\ZZ$. We thus see that in this case a secondary complex is
the same as a graded version of a \emph{multicomplex} (see e.~g. \cite{Meyer}) with only two
nonzero rows.

One then defines the \emph{total complex} Tot$(M_\bullet)$ to be
$$
...\ot
M_{n-1,0}\oplus M_{n-2,1}
\xot{
\left(
\begin{smallmatrix}
d_{n-1,0}&-\d_{n-1}\\
H_{n-2}&-d_{n-2,1}
\end{smallmatrix}
\right)
}
M_{n,0}\oplus M_{n-1,1}
\xot{
\left(
\begin{smallmatrix}
d_{n,0}&-\d_n\\
H_{n-1}&-d_{n-1,1}
\end{smallmatrix}
\right)
}
M_{n+1,0}\oplus M_{n,1}
\ot
...
$$
Cycles and boundaries in this complex will be called secondary cycles,
resp. secondary boundaries of $M_\bullet$. Thus a secondary $n$-cycle in
$M_\bullet$ is a pair $(c,\gamma)$ with $c\in M_{n,0}$, $\gamma\in
M_{n-1,1}$ such that $d_{n-1,0}c=\d_{n-1}\gamma$,
$H_{n-2}c=d_{n-2,1}\gamma$ and such a cycle is a
boundary iff there exist $b\in M_{n+1,0}$ and $\beta\in M_{n,1}$ with
$c=d_{n,0}b+\d_n\beta$ and $\gamma=H_{n-1}b+d_{n-1,1}\beta$. A secondary
complex $M_\bullet$ is called \emph{exact} if its total complex is exact, that is, if
secondary cycles are secondary boundaries.
\end{Definition}

Let us now consider a secondary chain complex $M_\bullet$ in $B$-$\Mod$. It
is clear then that
$$
...\to\pi_0M_{n+2}\xto{\pi_0d_{n+1}}\pi_0M_{n+1}\xto{\pi_0d_n}\pi_0M_n\xto{\pi_0d_{n-1}}\pi_0M_{n-1}\to...
\leqno{\pi_0M_\bullet:}
$$
is a chain complex of $\pi_0B$-modules. The next result corresponds to
\cite{BJ5}*{lemma 3.5}.

\begin{Proposition}\label{exact}
Let $M_\bullet$ be a secondary complex consisting of $\Sigma$-modules and
$\Sigma$-maps between them. If $\pi_0(M_\bullet)$ is an exact complex then
$M_\bullet$ is an exact secondary complex. Conversely, if $\pi_0M_\bullet$ is
bounded below then secondary exactness of $M_\bullet$ implies exactness of
$\pi_0M_\bullet$.
\end{Proposition}

\begin{proof} The proof consists in translating the argument from the
analogous general statement in \cite{BJ5} to our setting.
Suppose first that $\pi_0M_\bullet$ is an exact complex, and consider a
secondary cycle $(c,\gamma)\in M_{n,0}\oplus M_{n-1,1}$, i.~e. one has
$d_{n-1,0}c=\d_{n-1}\gamma$ and $H_{n-2}c=d_{n-2,1}\gamma$. Then in
particular $[c]\in\pi_0M_n$ is a cycle, so there exists
$[b]\in\pi_0M_{n+1}$ with $[c]=\pi_0(d_n)[b]$. Take $b\in[b]$, then
$c-d_{n,0}b=\d_n\beta$ for some $\beta\in M_{n+1,1}$. Consider
$\delta=\gamma-H_{n-1}b-d_{n-1,1}\beta$. One has
$\d_{n-1}\delta=\d_{n-1}\gamma-\d_{n-1}H_{n-1}b-\d_{n-1}d_{n-1,1}\beta=d_{n-1,0}c-d_{n-1,0}d_{n,0}b-d_{n-1,0}\d_n\beta=0$,
so that $\delta$ is an element of $\pi_1M_n$. Moreover
$d_{n-2,1}\delta=d_{n-2,1}\gamma-d_{n-2,1}H_{n-1}b-d_{n-2,1}d_{n-1,1}\beta=H_{n-2}c-H_{n-2}d_{n,0}b-H_{n-2}\d_n\beta=0$,
i.~e. $\delta$ is a cycle in $\pi_1M_\bullet$. Since by assumption
$\pi_0M_\bullet$ is exact, taking into account the $\Sigma$-structure
$\pi_1M_\bullet$ is exact too, so that there exists $\psi\in\pi_1M_n$
with $\delta=d_{n-1,1}\psi$. Define $\tilde\beta=\beta+\psi$. Then
$d_{n,0}b+\d_n\tilde\beta=d_{n,0}b+\d_n\beta=c$ since $\psi\in\ker\d_n$.
Moreover
$d_{n-1,1}\tilde\beta=d_{n-1,1}\beta+d_{n-1,1}\psi=d_{n-1,1}\beta+\delta=\gamma-H_{n-1}b$,
which means that $(c,\gamma)$ is the boundary of $(b,\tilde\beta)$. Thus
$M_\bullet$ is an exact secondary complex.

Conversely suppose $M_\bullet$ is exact, and $\pi_0M_\bullet$ bounded
below. Given a cycle $[c]\in\pi_0(M_n)$, represent it by a $c\in M_{n,0}$.
Then $\pi_0d_{n-1}[c]=0$ implies $d_{n-1,0}c\in\im\d_{n-1}$, so there is a
$\gamma\in M_{n-1,1}$ such that $d_{n-1,0}c=\d_{n-1}\gamma$. Consider
$\omega=d_{n-2,1}\gamma-H_{n-2}c$. One has
$\d_{n-2}\omega=\d_{n-2}d_{n-2,1}\gamma-\d_{n-2}H_{n-2}c=d_{n-2,0}\d_{n-1}\gamma-d_{n-2,0}d_{n-1,0}c=0$,
i.~e. $\omega$ is an element of $\pi_1M_{n-2}$. Moreover
$d_{n-3,1}\omega=d_{n-3,1}d_{n-2,1}\gamma-d_{n-3,1}H_{n-2}c=H_{n-3}\d_{n-1}\gamma-H_{n-3}d_{n,0}c=0$,
so $\omega$ is a $n-2$-dimensional cycle in $\pi_1M_\bullet$. Using the
$\Sigma$-structure, this then gives a $n-3$-dimensional cycle in
$\pi_0M_\bullet$. Now since $\pi_0M_\bullet$ is bounded below, we might
assume by induction that it is exact in dimension $n-3$, so that $\omega$
is a boundary. That is, there exists $\alpha\in\pi_1M_{n-1}$ with
$d_{n-2,1}\alpha=\omega$. Define $\tilde\gamma=\gamma-\alpha$; then one
has
$d_{n-2,1}\tilde\gamma=d_{n-2,1}\gamma-d_{n-2,1}\alpha=d_{n-2,1}\gamma-\omega=H_{n-2}c$.
Moreover $\d_{n-1}\tilde\gamma=\d_{n-1}\gamma=d_{n-1,0}c$ since
$\alpha\in\ker(\d){n-1}$. Thus $(c,\tilde\gamma)$ is a secondary cycle,
and by secondary exactness of $M_\bullet$ there exists a pair $(b,\beta)$
with $c=d_{n,0}b+\d_n\beta$. Then $[c]=\pi_0(d_n)[b]$, i.~e. $c$ is a
boundary.
\end{proof}

\begin{Definition}
Let $B$ be a pair algebra with $\Sigma$-structure.
A \emph{secondary resolution} of a $\Sigma$-module $M=(\d:M_1\to M_0)$ over
$B$ is an exact secondary complex $F_\bullet$ in $(B\mathrm{-}\Mod)^\Sigma$ of the form
$$
\xymatrix@!{
...\ar[r]\ar@{}[d]|\cdots
&F_{31}\ar[r]^{d_{21}}\ar[d]|\hole^>(.75){\d_3}
&F_{21}\ar[r]^{d_{11}}\ar[d]|\hole^>(.75){\d_2}
&F_{11}\ar[r]^{d_{01}}\ar[d]|\hole^>(.75){\d_1}
&F_{01}\ar[r]^{\epsilon_1}\ar[d]|\hole^>(.75){\d_0}
&M_1\ar[d]|\hole^>(.9)\d\ar[r]
&0\ar[r]\ar[d]|\hole
&0\ar[r]\ar[d]
&...\ar@{}[d]|\cdots
\\
...\ar[r]\ar[urr]^<(.3){H_2}
&F_{30}\ar[r]_{d_{20}}\ar[urr]^<(.3){H_1}
&F_{20}\ar[r]_{d_{10}}\ar[urr]^<(.3){H_0}
&F_{10}\ar[r]_{d_{00}}\ar[urr]^<(.3){\hat\epsilon}
&F_{00}\ar[r]_{\epsilon_0}\ar[urr]
&M_0\ar[r]\ar[urr]
&0\ar[r]
&0\ar[r]
&...
}
$$
where each $F_n=(\d_n:F_{n1}\to F_{n0})$ is a free $B$-module.
\end{Definition}

It follows from \ref{exact} that for any secondary resolution $F_\bullet$
of a $B$-module $M$ with $\Sigma$-structure, $\pi_0F_\bullet$ will be a
free resolution of the $\pi_0B$-module $\pi_0M$, so that in particular one
has
$$
\Ext^n_{\pi_0B}(\pi_0M,U)\cong H^n\Hom(\pi_0F_\bullet,U)
$$
for all $n$ and any $\pi_0B$-module $U$.

\begin{Definition}\label{secod}
Given a pair algebra $B$ with $\Sigma$-structure, a $\Sigma$-module $M$
over $B$, a module $N$ over $B$ and a secondary resolution $F_\bullet$ of $M$, we define
the \emph{secondary differential}
$$
d_{(2)}:\Ext^n_{\pi_0B}(\pi_0M,\pi_0N)\to\Ext^{n+2}_{\pi_0B}(\pi_0M,\pi_1N)
$$
in the following way. Suppose given a class
$[c]\in\Ext^n_{\pi_0B}(\pi_0M,\pi_0N)$.
First represent it by some element in $\Hom_{\pi_0B}(\pi_0F_n,\pi_0N)$
which is a cocycle, i.~e. its composite with $\pi_0(d_n)$ is 0. By
\ref{freehom} we know that the natural maps
$$
[F_n,N]\to\Hom_{\pi_0B}(\pi_0F,\pi_0N)
$$
are isomorphisms, hence to any such element corresponds a homotopy class
in $[F_n,N]$ which is also a cocycle, i.~e. value of $[d_n,N]$ on it is
zero. Take a representative map $c:F_n\to N$ from this homotopy class. Then
$cd_n$ is nullhomotopic, so we can find a $B_0$-equivariant map $H:F_{n+1,0}\to N_1$
such that in the diagram
$$
\xymatrix@!{
&F_{n+2,1}\ar[r]^{d_{n+1,1}}\ar[d]_{\d_{n+2}}
&F_{n+1,1}\ar[r]^{d_{n,1}}\ar[d]|\hole^>(.75){\d_{n+1}}
&F_{n,1}\ar[r]^{c_1}\ar[d]|\hole^>(.75){\d_n}
&N_1\ar[d]^\d\\
&F_{n+2,0}\ar[r]_{d_{n+1,0}}\ar[urr]^<(.3){H_n}
&F_{n+1,0}\ar[r]_{d_{n,0}}\ar[urr]^<(.3)H
&F_{n,0}\ar[r]_{c_0}
&N_0.
}
$$
one has $c_0d_{n,0}=\d H$, $c_1d_{n,1}=H\d_{n+1}$ and $\d c_1=c_0\d_n$.
Then taking $\Gamma=c_1H_n-Hd_{n+1,0}$ one has $\d\Gamma=0$,
$\Gamma\d_{n+2}=0$, so $\Gamma$ determines a map
$\bar\Gamma:\coker\d_{n+2}\to\ker\d$, i.~e. from $\pi_0F_{n+2}$ to
$\pi_1N$. Moreover $\bar\Gamma\pi_0(d_{n+2})=0$, so it is a cocycle in
$\Hom(\pi_0(F_\bullet),\pi_1N)$ and we define
$$
d_{(2)}[c]=[\bar\Gamma]\in\Ext^{n+2}_{\pi_0B}(\pi_0M,\pi_1N).
$$
\end{Definition}

\begin{Definition}\label{secodef}
Let $M$ and $N$ be $B$-modules with $\Sigma$-structure. Then also all
the $B$-modules $\Sigma^kM$, $\Sigma^kN$ have $\Sigma$-structures and we
get by \ref{secod} the secondary differential
$$
\xymatrix@!C=16em{
\Ext^n_{\pi_0B}(\pi_0M,\pi_0\Sigma^kN)\ar@{=}[d]\ar[r]^{d_{(2)}(M,\Sigma^kN)}&\Ext^{n+2}_{\pi_0B}(\pi_0M,\pi_1\Sigma^kN)\ar@{=}[d]\\
\Ext^n_{\pi_0B}(\pi_0M,\Sigma^k\pi_0N)\ar[r]^d&\Ext^{n+2}_{\pi_0B}(\pi_0M,\Sigma^{k+1}\pi_0N).
}
$$
In case the composite
$$
\Ext^{n-2}_{\pi_0B}(\pi_0M,\Sigma^{k-1}\pi_0N)\xto d\Ext^n_{\pi_0B}(\pi_0M,\Sigma^k\pi_0N)\xto d\Ext^{n+2}_{\pi_0B}(\pi_0M,\Sigma^{k+1}\pi_0N)
$$
vanishes we define the \emph{secondary $\Ext$-groups} to be the quotient
groups
$$
\Ext^n_B(M,N)^k:=\ker d/\im d.
$$
\end{Definition}

\begin{Theorem}\label{sigmacoinc}
For a $\Sigma$-algebra $B$, a $B$-module $M$ with $\Sigma$-structure and
any $B$-module $N$, the secondary differential $d_{(2)}$ in \ref{secod}
coincides with the secondary differential
$$
d_{(2)}:\Ext^n_\ba(M,N)\to\Ext^{n+2}_\ba(M,N)
$$
from \cite{BJ5}*{Section 4} as constructed for the $\LL$-additive track
category $(B\mathrm{-}\Mod)^\Sigma$ in \ref{add}, relative to the
subcategory $\bb$ of free $B$-modules with $\ba=\bb\ho$.
\end{Theorem}

\begin{proof}
We begin by recalling the appropriate notions from \cite{BJ5}. There
secondary chain complexes $A_\bullet=(A_n,d_n,\delta_n)_{n\in\ZZ}$ are
defined in an arbitrary additive track category $\B$. They consist of objects
$A_n$, morphisms $d_n:A_{n+1}\to A_n$ and tracks
$\delta_n:d_nd_{n+1}\then0_{A_{n+2},A_n}$, $n\in\ZZ$, such that the
equality of tracks
$$
\delta_nd_{n+2}=d_n\delta_{n+1}
$$
holds for all $n$. For an object $X$, an $X$-valued $n$-cycle in a
secondary chain complex $A_\bullet$ is defined to be a pair $(c,\gamma)$
consisting of a morphism $c:X\to A_n$ and a track
$\gamma:d_{n-1}c\then0_{X,A_{n-1}}$ such that the equality of tracks
$$
\delta_{n-2}c=d_{n-2}\gamma
$$
is satisfied. Such a cycle is called a \emph{boundary} if there exists a
map $b:X\to A_{n+1}$ and a track $\beta:c\then d_nb$ such that the equality
$$
\gamma=\delta_{n-1}b\comp d_{n-1}\beta
$$
holds. Here the right hand side is given by track addition. A secondary chain 
complex is called $X$-exact if every $X$-valued
cycle in it is a boundary. Similarly it is called \emph{$\bb$-exact}, if it
is $X$-exact for every object $X$ in $\bb$, where $\bb$ is a track
subcategory of $\B$. A secondary $\bb$-resolution of an object $A$ is a
$\bb$-exact secondary chain complex $A_\bullet$ with $A_n=0$ for $n<-1$,
$A_{-1}=A$, $A_n\in\bb$ for $n\ne-1$; the last differentials will be then
denoted $d_{-1}=\epsilon:A_0\to A$, $\delta_{-1}=\hat\epsilon:\epsilon
d_0\to0_{A_1,A}$ and the pair $(\epsilon,\hat\epsilon)$ will be called the
\emph{augmentation} of the resolution. It is clear that any secondary
chain complex $(A_\bullet,d_\bullet,\delta_\bullet)$ in $\B$ gives rise to
a chain complex $(A_\bullet,[d_\bullet])$, in the ordinary sense, in the
homotopy category $\B\ho$ of $\B$. Moreover if $\B$ is $\Sigma$-additive,
i.~e. there exists a functor $\Sigma$ and isomorphisms
$\Aut(0_{X,Y})\cong[\Sigma X,Y]$, natural in $X$, $Y$, then $\bb$-exactness
of $(A_\bullet,d_\bullet,\delta_\bullet)$ implies $\bb\ho$-exactness of
$(A_\bullet,[d_\bullet])$ in the sense that the chain complex of abelian
groups $[X,(A_\bullet,[d_\bullet])]$ will be exact for each $X\in\bb$. In
\cite{BJ5}, the notion of $\bb\ho$-relative derived functors has been
developed using such $\bb\ho$-resolutions, which we also recall.

For an additive subcategory $\ba=\bb\ho$ of the homotopy category $\B\ho$,
the $\ba$-relative left derived functors $\L^\ba_nF$, $n\ge0$, of a
functor $F:\B\ho\to\A$ from $\B\ho$ to an abelian category $\mathscr
A$ are defined by
$$
(\L^\ba_nF)A=H_n(F(A_\bullet)),
$$
where $A_\bullet$ is given by any $\ba$-resolution of $A$. Similarly,
the $\ba$-relative right derived functors of a contravariant functor
$F:\B\ho\op\to\A$ are given by
$$
(\rR_\ba^nF)A=H^n(F(A_\bullet)).
$$
In particular, for the contravariant functor $F=[\_,B]$ we get the $\ba$-relative $\Ext$-groups
$$
\Ext^n_\ba(A,B):=(\rR_\ba^n[\_,B])A=H^n([A_\bullet,B])
$$
for any $\ba$-exact resolution $A_\bullet$ of $A$.
Similarly, for the contravariant functor $\Aut(0_{\_,B})$ which assigns to
an object $A$ the group $\Aut(0_{A,B})$ of all tracks
$\alpha:0_{A,B}\then0_{A,B}$ from the zero map $A\to*\to B$ to itself, one
gets the groups of $\ba$-derived automorphisms
$$
\Aut^n_\ba(A,B):=(\rR^n_\ba\Aut(0_{\_,B}))(A).
$$

It is proved in \cite{BJ5} that under mild conditions (existence of a
subset of $\ba$ such that every object of $\ba$ is a direct summand of a
direct sum of objects from that subset) every object has an
$\ba$-resolution, and that the resulting groups do not depend on the choice
of a resolution. 

We next recall the construction of the secondary differential from
\cite{BJ5}. This is a map of the form
$$
d_{(2)}:\Ext^n_\ba(A,B)\to\Aut^n_\ba(0_{A,B});
$$
it is constructed from any secondary $\bb$-resolution
$(A_\bullet,d_\bullet,\delta_\bullet,\epsilon,\hat\epsilon)$ of the object
$A$. Given an element $[c]\in\Ext^n_\ba(A,B)$, one first represents it
by an $n$-cocycle in $[(A_\bullet,[d_\bullet]),B]$, i.~e. by a homotopy
class $[c]\in[A_n,B]$ with $[cd_n]=0$. One then chooses an actual
representative $c:A_n\to B$ of it in $\B$ and a track $\gamma:0\then
cd_n$. It can be shown that the composite track
$\Gamma=c\delta_n\comp\gamma d_{n+1}\in\Aut(0_{A_{n+2},B})$ satisfies
$\Gamma d_{n+1}=0$, so it is an $(n+2)$-cocycle in the cochain complex
$\Aut(0_{(A_\bullet,[d_\bullet]),B})\cong[(\Sigma A_\bullet,[\Sigma
d_\bullet]),B]$, so determines a cohomology class
$d{(2)}([c])=[\Gamma]\in\Ext^{n+2}_\ba(\Sigma A,B)$. It is proved in
\cite{BJ5}*{4.2} that the above construction does not indeed depend on
choices.

Now turning to our situation, it is straightforward to verify that a
secondary chain complex in the sense of \cite{BJ5} in the track category
$B$-$\Mod$ is the same as a 2-complex in the sense of \ref{secs}, and
that the two notions of exactness coincide. In particular then the notions
of resolution are also equivalent.

The track subcategory $\bb$ of free modules is generated by coproducts from
a single object, so $\bb\ho$-resolutions of any $B$-module exist. In fact
it follows from \cite{BJ5}*{2.13} that any $B$-module has a secondary
$\bb$-resolution too.

Moreover there are natural isomorphisms
$$
\Aut(0_{M,N})\cong\Hom_{\pi_0B}(\pi_0M,\pi_1N).
$$
Indeed a track from the zero map to itself is a $B_0$-module homomorphism
$H:M_0\to N_1$ with $\d H=0$, $H\d=0$, so $H$ factors through $M_0\onto\pi_0M$
and over $\pi_1N\rightarrowtail N_1$. 

Hence the proof is finished with the following lemma.
\end{proof}

\begin{Lemma}
For any $B$-modules $M$, $N$ there are isomorphisms
$$
\Ext^n_\ba(M,N)\cong\Ext^n_{\pi_0B}(\pi_0M,\pi_0N)
$$
and
$$
(\rR^n_\ba(\Hom_{\pi_0B}(\pi_0(\_),\pi_1N)))(M)\cong\Ext^n_{\pi_0B}(\pi_0M,\pi_1N).
$$
\end{Lemma}

\begin{proof}
By definition the groups $\Ext^*_\ba(M,N)$, respectively
$(\rR^n_\ba(\Hom_{B_0}(\pi_0(\_),\pi_1N)))(M)$, are cohomology groups of the
complex $[F_\bullet,N]$, resp. $\Hom_{\pi_0B}(\pi_0(F_\bullet),\pi_1N)$,
where $F_\bullet$ is some $\ba$-resolution of $M$. We can choose for
$F_\bullet$ some secondary $\bb$-resolution of $M$. Then $\pi_0F_\bullet$
is a free $\pi_0B$-resolution of $\pi_0M$, which makes evident the second
isomorphism. For the first, just note in addition that by \ref{freehom}
$[F_\bullet,N]$ is isomorphic to $\Hom_{B_0}(\pi_0(F_\bullet),\pi_0N)$.
\end{proof}

%% file: dHascoAss_c2.tex
\chapter{The pair algebra $\aB$ of secondary cohomology
operations}\label{secsteen}

The algebra $\aB$ of secondary cohomology operations is a pair algebra with
$\Sigma$-structure which as a Hopf algebra was explicitly computed in
\cite{Baues}. In particular the multiplication map $A$ of $\aB$ was
determined in \cite{Baues} by an algorithm. In this chapter we recall the
topological definition of the pair algebra $\aB$ and the definition of the
multiplication map $A$. The main results of this work will provide methods for
the computation of $A$ or its dual multiplication map $A_*$. 
We express in terms of $A$ the secondary $\Ext$-groups $\Ext_\aB$ over the pair algebra $\aB$.
This yields the computation of the E$_3$-term of the Adams spectral sequence
in the next chapter.

\section{The track category of spectra}

In this section we introduce the notion of stable maps and stable tracks
between spectra. This yields the track category of spectra. See also
\cite{Baues}*{section 2.5}.

\begin{Definition}\label{sp}
A \emph{spectrum} $X$ is a sequence of maps
$$
X_i\xto r\Omega X_{i+1},\ i\in\ZZ
$$
in the category $\Topt$ of pointed spaces. This is an $\Omega$-spectrum if
$r$ is a homotopy equivalence for all $i$.

A \emph{stable homotopy class} $f:X\to Y$ between spectra is a sequence of
homotopy classes $f_i\in[X_i,Y_i]$ such that the squares
$$
\xymatrix{
X_i\ar[r]^{f_i}\ar[d]^r&Y_i\ar[d]^r\\
\Omega X_{i+1}\ar[r]^{\Omega f_{i+1}}&\Omega Y_{i+1}
}
$$
commute in $\Topt\ho$. The category $\Sp$ consists of spectra and stable
homotopy classes as morphisms. Its full subcategory $\Omega$-$\Sp$
consisting of $\Omega$-spectra is equivalent to the homotopy category
of spectra considered as a Quillen model category as in the work on symmetric spectra of 
M. Hovey, B. Shipley and J. Smith \cite{Hoveyetal}. For us the classical notion of a spectrum
as above is sufficient.

A \emph{stable map} $f=(f_i,\tilde f_i)_i:X\to Y$ between spectra is a sequence
of diagrams in the track category $\hog{\Topt}$ $(i\in\ZZ)$
$$
\xymatrix{
X_i\ar[r]^{f_i}\ar[d]_r\drtwocell\omit{^{\tilde f_i\ }}&Y_i\ar[d]^r\\
\Omega X_{i+1}\ar[r]_{\Omega f_{i+1}}&\Omega Y_{i+1}.
}
$$
Obvious composition of such maps yields the category
$$
\hog{\Sp}_0.
$$
It is the underlying category of a track category $\hog{\Sp}$ with tracks
$(H:f\then g)\in\hog{\Sp}_1$ given by sequences
$$
H_i:f_i\then g_i
$$
of tracks in $\Topt$ such that the diagrams
$$
\xymatrix@C=8em@!R=3em{
X_i\ar[r]_{f_i}\ar[d]_r\drtwocell\omit{^{\tilde
f_i\ }}\ruppertwocell^{g_i}{^H_i\ }
&Y_i\ar[d]^r\\
\Omega X_{i+1}\ar[r]^{\Omega f_{i+1}}\rlowertwocell_{\Omega
g_{i+1}}{_\Omega H_{i+1}\hskip4em}
&\Omega Y_{i+1}
}
$$
paste to $\tilde g_i$. This yields a well-defined track category
$\hog{\Sp}$. Moreover
$$
\hog{\Sp}\ho\cong\Sp
$$
is an isomorphism of categories. Let $\hog{X,Y}$ be the groupoid of
morphisms $X\to Y$ in $\hog{\Sp}_0$ and let $\hog{X,Y}_1^0$ be the set of
pairs $(f,H)$ where $f:X\to Y$ is a map and $H:f\then0$ is a track in
$\hog{\Sp}$, i.~e. a stable homotopy class of nullhomotopies for $f$.

For a spectrum $X$ let $\Sigma^kX$ be the \emph{shifted spectrum} with
$(\Sigma^kX)_n=X_{n+k}$ and the commutative diagram
$$
\xymatrix{
(\Sigma^kX)_n\ar@{=}[d]\ar[r]^-r&\Omega(\Sigma^kX)_{n+1}\ar@{=}[d]\\
X_{n+k}\ar[r]^-r&\Omega(X_{n+k+1})
}
$$
defining $r$ for $\Sigma^kX$. A map $f:Y\to\Sigma^kX$ is also called a map
$f$ \emph{of degree $k$} from $Y$ to $X$.
\end{Definition}

\section{The pair algebra $\aB$ and secondary cohomology of spectra as a
$\aB$-module}

The secondary cohomology of a space was introduced in \cite{Baues}*{section
6.3}. We here consider the corresponding notion of secondary cohomology of
a spectrum. 

Let $\FF$ be a prime field $\FF=\ZZ/p\ZZ$ and let $Z$ denote the Eilenberg-Mac
Lane spectrum with 
$$
Z^n=K(\FF,n)
$$
chosen as in \cite{Baues}. Here $Z^n$ is a topological $\FF$-vector space
and the homotopy equivalence $Z^n\to\Omega Z^{n+1}$ is $\FF$-linear. This
shows that for a spectrum $X$ the sets $\hog{X,\Sigma^kZ}_0$ and
$\hog{X,\Sigma^kZ}_1^0$, of stable maps and stable 0-tracks repectively, are
$\FF$-vector spaces.

We now recall the definition of the pair algebra $\aB=(\d:\aB_1\to\aB_0)$
of secondary cohomology operations from \cite{Baues}. Let $\G=\ZZ/p^2\ZZ$ and let
$$
\aB_0=T_\G(E_\A)
$$
be the $\G$-tensor algebra generated by the subset
$$
E_\A=
\begin{cases}
\set{\Sq^1,\Sq^2,...}&\textrm{for $p=2$},\\
\set{\P^1,\P^2,...}\cup\set{\beta,\beta\P^1,\beta\P^2,...}&\textrm{for odd
$p$}
\end{cases}
$$
of the mod $p$ Steenrod algebra $\A$. We define $\aB_1$ by the
pullback diagram of graded abelian groups
\begin{equation}\label{defb1}
\begin{aligned}
\xymatrix{
&\Sigma\A\ar@{ >->}[d]\\
\aB_1\ar[r]\ar[d]_\d\ar@{}[dr]|<{\Bigg{\lrcorner}}&\hog{Z,\Sigma^*Z}_1^0\ar[d]^\d\\
\aB_0\ar[r]^-s&\hog{Z,\Sigma^*Z}_0\ar@{->>}[d]\\
&\A.
}
\end{aligned}
\end{equation}
in which the right hand column is an exact sequence. Here we choose for
$\alpha\in E_\A$ a stable map $s(\alpha):Z\to\Sigma^{|\alpha|}Z$
representing $\alpha$ and we define $s$ to be the $\G$-linear map given on
monomials $a_1\cdots a_n$ in the free monoid Mon$(E_\A)$ generated by
$E_\A$ by the composites
$$
s(a_1\cdots a_n)=s(a_1)\cdots s(a_n).
$$
It is proved in \cite{Baues}*{5.2.3} that $s$ defines a pseudofunctor, that is,
there is a well-defined track
$$
\Gamma:s(a\cdot b)\then s(a)\circ s(b)
$$
for $a,b\in\aB_0$ such that for any $a$, $b$, $c$ pasting of tracks in the
diagram
$$
\includegraphics{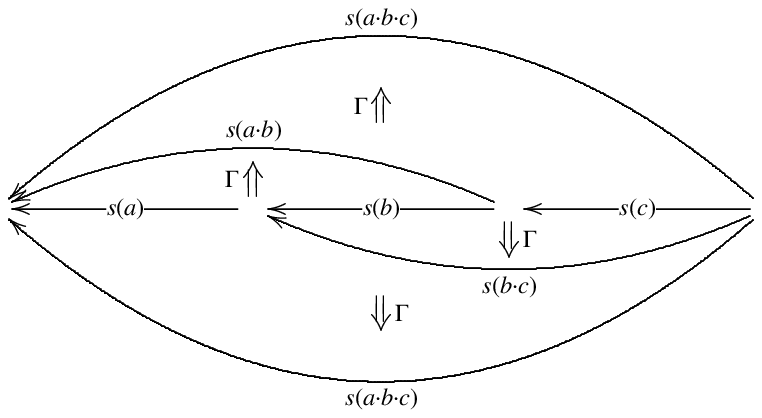}
$$
yields the identity track. Now $\aB_1$ is a $\aB_0$-$\aB_0$-bimodule by
defining
$$
a(b,z)=(a\cdot b,a\bullet z)
$$
with $a\bullet z$ given by pasting $s(a)z$ and $\Gamma$. Similarly
$$
(b,z)a=(b\cdot a,z\bullet a)
$$
where $z\bullet a$ is obtained by pasting $zs(a)$ and $\Gamma$. Then it is
shown in \cite{Baues} that $\aB=(\d:\aB_1\to\aB_0)$ is a well-defined pair
algebra with $\pi_0\aB=\A$ and $\Sigma$-structure $\pi_1\aB=\Sigma\A$.

For a spectrum $X$ let
$$
\H(X)_0=\aB_0\hog{X,\Sigma^*Z}_0
$$
be the free $\aB_0$-module generated by the graded set
$\hog{X,\Sigma^*Z}_0$. We define $\H(X)_1$ by the pullback diagram
$$
\xymatrix{
&\Sigma H^*X\ar@{ >->}[d]\\
\H(X)_1\ar[r]\ar[d]_\d\ar@{}[dr]|<{\Bigg{\lrcorner}}&\hog{X,\Sigma^*Z}_1^0\ar[d]^\d\\
\H(X)_0\ar[r]^-s&\hog{X,\Sigma^*Z}_0\ar@{->>}[d]\\
&H^*X
}
$$
where $s$ is the $\G$-linear map which is the identity on generators and is
defined on words $a_1\cdots a_n\cdot u$ by the composite $s(a_1)\cdots
s(a_n)s(u)$ for $a_i$ as above and $u\in\hog{X,\Sigma^*Z}_0$. Again $s$ is
a pseudofunctor and with actions $\bullet$ defined as above we see that the
graded pair module
$$
\H(X)=\left(\H(X)_1\xto\d\H(X)_0\right)
$$
is a $\aB$-module. We call $\H(X)$ the \emph{secondary cohomology} of the
spectrum $X$. Of course $\H(X)$ has a $\Sigma$-structure in the sense of
\ref{sigmas} above.

\begin{Example}
Let $\G^\Sigma$ be the $\aB$-module given by the augmentation $\aB\to\G^
\Sigma$ in \cite{Baues}. Recall that $\G^\Sigma$ is the pair
$$
\G^\Sigma=\left(\FF\oplus\Sigma\FF\xto\d\G\right)
$$
with $\restr\d\FF$ the inclusion nad $\restr\d{\Sigma\FF}=0$. Then the
sphere spectrum $S^0$ admits a weak equivalence of $\aB$-modules
$$
\H(S^0)\xto\sim\G^\Sigma.
$$
Compare \cite{Baues}*{12.1.5}.
\end{Example}

%% file: dHascoAss_c3.tex
\chapter{Computation of the $\E_3$-term of the Adams spectral
sequence as a secondary $\Ext$-group}\label{E3}

We show that the E$_3$-term of the Adams spectral sequence (computing stable
maps in $\set{Y,X}^*_p$) is given by the secondary $\Ext$-groups
$$
\E_3(Y,X)=\Ext_\aB(\H X,\H Y).
$$
Here $\H X$ is the secondary cohomology of the spectrum $X$ which is the
$\aB$-module $\G^\Sigma$ if $X$ is the sphere spectrum $S^0$. This leads to
an algorithm for the computation of the group
$$
\E_3(S^0,S^0)=\Ext_\aB(\G^\Sigma,\G^\Sigma)
$$
which is a new explicit approximation of stable homotopy groups of spheres
improving the Adams approximation
$$
\E_2(S^0,S^0)=\Ext_\A(\FF,\FF).
$$
An implementation of our algorithm computed $\E_3(S^0,S^0)$ by now up
to degree 40. In this range our results confirm the known results in the
literature, see for example the book of Ravenel \cite{Ravenel}.

\section{The ${\mathrm E}_3$-term of the Adams spectral sequence}

We now are ready to formulate the algebraic equivalent of the E$_3$-term of
the Adams spectral sequence. Let $X$ be a spectrum of finite type and $Y$ a
finite dimensional spectrum. Then for each prime $p$ there is a spectral
sequence $\E_*=\E_*(Y,X)$ with
\begin{align*}
\E_*&\Longrightarrow[Y,\Sigma^*X]_p\\
\E_2&=\Ext_\A(H^*X,H^*Y).
\end{align*}

\begin{Theorem}\label{e3der}
The $\E_3$-term $\E_3=\E_3(Y,X)$ of the Adams spectral sequence is given by
the secondary $\Ext$ group defined in \ref{secodef}
$$
\E_3=\Ext_\aB(\H^*X,\H^*Y).
$$
\end{Theorem}

\begin{Corollary}
If $X$ and $Y$ are both the sphere spectrum we get
$$
\E_3(S^0,S^0)=\Ext_\aB(\G^\Sigma,\G^\Sigma).
$$
\end{Corollary}

Since the pair algebra $\aB$ is computed in \cite{Baues} completely we see
that $\E_3(S^0,S^0)$ is algebraically determined. This leads to the
algorithm below computing $\E_3(S^0,S^0)$.

The proof of \ref{e3der} is based on the following result in \cite{Baues}.
Consider the track categories
\begin{align*}
\bb&\subset\hog{\Sp}\\
\bb'&\subset(\aB\mathrm{-}\Mod)^\Sigma
\end{align*}
where $\hog{\Sp}$ is the track category of spectra in \ref{sp} and
$(\aB\mathrm{-}\Mod)^\Sigma$ is the track category of $\aB$-modules with
$\Sigma$-structure in \ref{sigmas} with the pair algebra $\aB$ defined by
\eqref{defb1}. Let $\bb$ be the full track subcategory of $\hog{\Sp}$
consisting of finite products of shifted Eilenberg-Mac Lane
spectra $\Sigma^kZ^*$. Moreover let $\bb'$ be the full track subcategory of
$(\aB\mathrm{-}\Mod)^\Sigma$ consisting of finitely generated free
$\aB$-modules. As in \cite{BJ5}*{4.3} we obtain for spectra $X$, $Y$ in
\ref{e3der} the track categories
\begin{align*}
\set{Y,X}\bb&\subset\hog{\Sp}\\
\bb'\set{\H X,\H Y}&\subset(\aB\mathrm{-}\Mod)^\Sigma
\end{align*}
with $\set{Y,X}\bb$ obtained by adding to $\bb$ the objects $X$, $Y$ and all
morphisms and tracks from $\hog{X,Z}$, $\hog{Y,Z}$ for all objects $Z$ in
$\bb$. It is proved in \cite{Baues}*{5.5.6} that the following result holds
which shows that we can apply \cite{BJ5}*{5.1}.

\begin{Theocite}\label{strictif}
There is a strict track equivalence
$$
(\set{Y,X}\bb)\op\xto\sim\bb'\set{\H X,\H Y}.
$$
\end{Theocite}
\qed

\begin{proof}[Proof of \ref{e3der}]
By the main result 7.3 in \cite{BJ5} we have a description of the
differential $d_{(2)}$ in the Adams spectral sequence by the following
commutative diagram
$$
\xymatrix{
\Ext^n_{\ba\op}(X,Y)^m\ar[r]^{d_{(2)}}\ar[d]^\cong
&\Ext^{n+2}_{\ba\op}(X,Y)^{m+1}\ar[d]^\cong\\
\Ext_{\mathscr A}^n(H^*X,H^*Y)^m\ar[r]^{d_{(2)}}
&\Ext_{\mathscr A}^{n+2}(H^*X,H^*Y)^{m+1}
}
$$
where $\ba=\bb\ho$. On the other hand the differential
$d_{(2)}$ defining the secondary $\Ext$-group $\Ext_\aB(\H X,\H Y)$ is by \ref{sigmacoinc}
given by the commutative diagram
$$
\xymatrix{
\Ext^n_{\ba'}(\H X,\H Y)^m\ar[r]\ar@{=}[d]
&\Ext^{n+2}_{\ba'}(\H X,\H Y)^{m+1}\ar@{=}[d]\\
\Ext_{\mathscr A}^n(H^*X,H^*Y)^m\ar[r]
&\Ext_{\mathscr A}^{n+2}(H^*X,H^*Y)^{m+1}
}
$$
where $\ba'=\bb'\ho$. Now \cite{BJ5}*{5.1} shows by \ref{strictif} that the top rows of these diagrams coincide.
\end{proof}

\section{The algorithm for the computation of $d_{(2)}$ on
$\Ext_\A(\FF,\FF)$ in terms of the multiplication maps}\label{d2}

Suppose now given some projective resolution of the left $\A$-module $\FF$. For definiteness, we will work with the minimal
resolution
\begin{equation}\label{minireso}
\FF\ot\A\brk{g_0^0}\ot\A\brk{g_1^{2^n}\mid
n\ge0}\ot\A\brk{g_2^{2^i+2^j}\mid
|i-j|\ne1}\ot...,
\end{equation}
where $g_m^d$, $d\ge m$, is a generator of the $m$-th resolving module in degree $d$.
Sometimes there are more than one generators with the same $m$ and $d$, in which
case the further ones will be denoted by $'g_m^d$, $''g_m^d, \cdots$.

These generators and values of the differential on them can be computed
effectively; for example, $d(g_1^{2^n})=\Sq^{2^n}g_0^0$ and
$d(g_m^m)=\Sq^1g_{m-1}^{m-1}$; moreover e.~g. an
algorithm from \cite{Bruner} gives
$$
\alignbox{
d(g_2^4)&=\Sq^3g_1^1+\Sq^2g_1^2\\
d(g_2^5)&=\Sq^4g_1^1+\Sq^2\Sq^1g_1^2+\Sq^1g_1^4\\
d(g_2^8)&=\Sq^6g_1^2+(\Sq^4+\Sq^3\Sq^1)g_1^4\\
d(g_2^9)&=\Sq^8g_1^1+(\Sq^5+\Sq^4\Sq^1)g_1^4+\Sq^1g_1^8\\
d(g_2^{10})&=(\Sq^8+\Sq^5\Sq^2\Sq^1)g_1^2+(\Sq^5\Sq^1+\Sq^4\Sq^2)g_1^4+\Sq^2g_1^8\\
d(g_2^{16})&=(\Sq^{12}+\Sq^9\Sq^2\Sq^1+\Sq^8\Sq^3\Sq^1)g_1^4+(\Sq^8+\Sq^7\Sq^1+\Sq^6\Sq^2)g_1^8\\
\cdots,\\
d(g_3^6)&=\Sq^4g_2^2+\Sq^2g_2^4+\Sq^1g_2^5\\
d(g_3^{10})&=\Sq^8g_2^2+(\Sq^5+\Sq^4\Sq^1)g_2^5+\Sq^1g_2^9\\
d(g_3^{11})&=(\Sq^7+\Sq^4\Sq^2\Sq^1)g_2^4+\Sq^6g_2^5+\Sq^2\Sq^1g_2^8\\
d(g_3^{12})&=\Sq^8g_2^4+(\Sq^6\Sq^1+\Sq^5\Sq^2)g_2^5+(\Sq^4+\Sq^3\Sq^1)g_2^8+\Sq^3g_2^9+\Sq^2g_2^{10}\\
\cdots,}
$$
$$
\alignbox{
d(g_4^{11})&=\Sq^8g_3^3+(\Sq^5+\Sq^4\Sq^1)g_3^6+\Sq^1g_3^{10}\\
d(g_4^{13})&=\Sq^8\Sq^2g_3^3+(\Sq^7+\Sq^4\Sq^2\Sq^1)g_3^6+\Sq^2\Sq^1g_3^{10}+\Sq^2g_3^{11}\\
\cdots,\\
d(g_5^{14})&=\Sq^{10}g_4^4+\Sq^2\Sq^1g_4^{11}\\
d(g_5^{16})&=\Sq^{12}g_4^4+\Sq^4\Sq^1g_4^{11}+\Sq^3g_4^{13}\\
\cdots,
\\
d(g_6^{16})&=\Sq^{11}g_5^5+\Sq^2g_5^{14}\\
\cdots,
}
$$
etc.

By understanding the above formul\ae\ as matrices (i.~e. by applying
$\chi$ degreewise to them), each such resolution gives rise to a sequence
of $\aB$-module homomorphisms
\begin{equation}\label{precomplex}
\G^\Sigma\ot\aB\brk{g_0^0}\ot\aB\brk{g_1^{2^n}\mid
n\ge0}\ot\aB\brk{g_2^{2^i+2^j}\mid
|i-j|\ne1}\ot...,
\end{equation}
which is far from being exact --- in fact even the composites of consecutive
maps are not zero. In more detail, one has commutative diagrams
$$
\xymatrix{
2\G\ar[d]
&R_\aB^0g_0^0\ar[l]_{\epsilon_0}\ar[d]
&0\ar[l]\ar[d]
&...\ar[l]\\
\G
&\aB_0^0g_0^0\ar[l]_{\epsilon_0}
&0\ar[l]
&...\ar[l]
}
$$
in degree 0,
$$
\xymatrix{
\FF\ar[d]
&R_\aB^1g_0^0\oplus\A^0g_0^0\ar[l]_-{(0,\epsilon)}\ar[d]
&R_\aB^0g_1^1\ar[l]_-{\binom d0}\ar[d]
&0\ar[l]\ar[d]
&...\ar[l]\\
0
&\aB_0^1g_0^0\ar[l]
&\aB_0^0g_1^1\ar[l]_d
&0\ar[l]
&...\ar[l]
}
$$
in degree 1,
$$
\xymatrix{
0\ar[d]
&R_\aB^2g_0^0\oplus\A^1g_0^0\ar[l]\ar[d]
&\left(R_\aB^1g_1^1\oplus R_\aB^0g_1^2\right)\oplus{\mathscr
A}^0g_1^1\ar[l]_-{\smat{d&0\\0&d}}\ar[d]
&R_\aB^0g_2^2\ar[l]_-{\binom d0}\ar[d]
&0\ar[l]\ar[d]
&...\ar[l]\\
0
&\aB_0^2g_0^0\ar[l]
&\aB_0^1g_1^1\oplus\aB_0^0g_1^2\ar[l]_-d
&\aB_0^0g_2^2\ar[l]_-d
&0\ar[l]
&...\ar[l]
}
$$
in degree 2, ...
$$
\xymatrix{
0\ar[d]
&R_\aB^ng_0^0\oplus\A^{n-1}g_0^0\ar[l]\ar[d]
&\bigoplus_{2^i\le n}R_\aB^{n-2^i}g_1^{2^i}
\oplus\bigoplus_{2^i\le n-1}\A^{n-1-2^i}g_1^{2^i}\ar[l]_-{\smat{d&0\\0&d}}\ar[d]
&...\ar[l]\\
0
&\aB_0^ng_0^0\ar[l]
&\bigoplus_{2^i\le n}\aB_0^{n-2^i}g_1^{2^i}\ar[l]_-d
&...\ar[l]
}
$$
in degree $n$, etc.

Our task is then to complete these diagrams into an exact secondary complex
via certain (degree preserving) maps
$$
\delta_m=\binom{\delta^R_m}{\delta^\A_m}:\aB_0\brk{g_{m+2}^n\mid n}\to(R_\aB\oplus\Sigma\A)\brk{g_m^n\mid
n}.
$$

Now for these maps to form a secondary complex, according to
\ref{secs}.1 one must have $\d\delta=d_0d_0$,
$\delta\d=d_1d_1$, and $d_1\delta=\delta d_0$. One sees easily that
these equations together with the requirement that $\delta$ be left
$\aB_0$-module homomorphism are equivalent to
\begin{align}
\delta^R&=dd,\\
\label{deltaeqs}\delta^\A(bg)&=\pi(b)\delta^\A(g)+A(\pi(b),dd(g)),\\
d\delta^\A&=\delta^\A d,
\end{align}
for $b\in\aB_0$, $g$ one of the $g_m^n$, and $A(a,rg):=A(a,r)g$ for
$a\in\A$, $r\in R_\aB$. Hence $\delta$ is completely determined by the elements
\begin{equation}\label{deltagen}
\delta^\A_m(g_{m+2}^n)\in\bigoplus_k\A^{n-k-1}\brk{g_m^k}
\end{equation}
which, to form a secondary complex, are only required to satisfy
$$
d\delta_m^\A(g_{m+2}^n)=\delta_{m-1}^\A d(g_{m+2}^n),
$$
where on the right $\delta_{m-1}^\A$ is extended to
$\aB_0\brk{g_{m+1}^*}$ via \ref{deltaeqs}. In addition secondary
exactness must hold, which by \ref{secs} means that the (ordinary) complex
$$
\ot
\aB_0\brk{g_{m-1}^*}\oplus(R_\aB\oplus\Sigma\A)\brk{g_{m-2}^*}
\ot
\aB_0\brk{g_m^*}\oplus(R_\aB\oplus\Sigma\A)\brk{g_{m-1}^*}
\ot
\aB_0\brk{g_{m+1}^*}\oplus(R_\aB\oplus\Sigma\A)\brk{g_m^*}
\ot
$$
with differentials
$$
\smat{d_{m+1}&i_{m+1}&0\\d_md_{m+1}&d_m&0\\\delta_m^\A&0&d_m}:
\aB_0\brk{g_{m+2}^*}\oplus R_\aB\brk{g_{m+1}^*}\oplus\Sigma\A\brk{g_{m+1}^*}
\to
\aB_0\brk{g_{m+1}^*}\oplus R_\aB\brk{g_m^*}\oplus\Sigma\A\brk{g_m^*}
$$
is exact. Then straightforward checking shows that one can eliminate $R_\aB$ from this complex altogether, so that its exactness is equivalent to the exactness of a smaller complex
$$
\ot
\aB_0\brk{g_{m-1}^*}\oplus\Sigma\A\brk{g_{m-2}^*}
\ot
\aB_0\brk{g_m^*}\oplus\Sigma\A\brk{g_{m-1}^*}
\ot
\aB_0\brk{g_{m+1}^*}\oplus\Sigma\A\brk{g_m^*}
\ot
$$
with differentials
$$
\smat{d_{m+1}&0\\\delta_m^\A&d_m}:
\aB_0\brk{g_{m+2}^*}\oplus\Sigma\A\brk{g_{m+1}^*}
\to
\aB_0\brk{g_{m+1}^*}\oplus\Sigma\A\brk{g_m^*}.
$$
Note also that by \ref{deltaeqs} $\delta^\A$ factors through
$\pi$ to give
$$
\bar\delta_m:\A\brk{g^*_{m+2}}\to\Sigma\A\brk{g^*_m}.
$$
It follows that secondary exactness of the resulting complex is equivalent
to exactness of the \emph{mapping cone} of this $\bar\delta$, i.~e. to the
requirement that $\bar\delta$ is a quasiisomorphism. On the other hand, the
complex $(\A\brk{g^*_*},d_*)$ is acyclic by construction, so any
of its self-maps is a quasiisomorphism. We thus obtain

\begin{Theorem}\label{delta}
Completions of the diagram \ref{precomplex}
to an exact secondary complex are in one-to-one correspondence with
maps $\delta_m:\A\brk{g_{m+2}^*}\to\Sigma\A\brk{g^*_m}$
satisfying
\begin{equation}\label{maineq}
d\delta g=\delta dg,
\end{equation}
with $\delta(ag)$ for $a\in\A$ defined by
$$
\delta(ag)=a\delta(g)+A(a,ddg)
$$
where $A(a,rg)$ for $r\in R_\aB$ is interpreted as $A(a,r)g$.
\end{Theorem}
\qed

Later in chapter \ref{diff} we will need to dualize the map $\delta$. For this purpose it is more convenient to reformulate the conditions in \ref{delta} above in terms of commutative diagrams.

Let
$$
W_p=\bigoplus_{q\ge0}W_p^q
$$
denote the free graded $\GG$-module spanned by the generators $g_p^q$, so that we can write
$$
\aB_0\brk{g_p^q\ \mid\ q\ge0}=\aB_0\ox W_p.
$$
The differential in the $\aB$-lifting of \eqref{minireso}, being $\aB$-equivariant, is then given by the composite
$$
\aB_0\ox W_{p+1}\xto{1\ox d}\aB_0\ox\aB_0\ox W_p\xto{m\ox1}\aB_0\ox W_p,
$$
where
$$
d:W_{p+1}\to\aB_0\ox W_p
$$
is the restriction of this differential to the generators. As a linear
operator, this $d$ is given by the same matrix as the one giving the operator
of the same name in \eqref{minireso}, i.~e. it is obtained by applying the map
$\chi$ componentwise to the latter.

Moreover let us denote
$$
V_p=W_p\ox\FF,
$$
so that similarly to the above the differential of \eqref{minireso} itself can be given by the same formul\ae, with $\A$ in place of $\aB_0$ and $\V_p$ in place of $\W_p$. Then by \ref{delta} the whole map $\delta$ is determined by its restriction
$$
\delta^\A:V_{p+2}\to\Sigma\A\ox V_p
$$
(cf. \eqref{deltagen}). Indeed \ref{delta} implies that $\delta$ is given by the sum of the two composites in the diagram
\begin{equation}\label{deltasum}
\alignbox{
\xymatrix{
&\A\ox\Sigma\A\ox V_p\ar[dr]^{m\ox1}\\
\A\ox V_{p+2}\ar[ur]^{1\ox\delta^\A}\ar[dr]_{1\ox\ph}&&\Sigma\A\ox V_p.\\
&\A\ox R_\aB\ox V_p\ar[ur]_{A\ox1}
}
}
\end{equation}
Here we set $\ph=dd\ox\FF$, where the map $dd$ is the composite
$$
W_{p+2}\xto d\aB_0\ox W_{p+1}\xto{1\ox d}\aB_0\ox\aB_0\ox W_p\xto{m\ox1}\aB_0\ox W_p
$$
whose image, as we know, lies in
$$
R_\aB\ox W_p\subset\aB_0\ox W_p.
$$
In other words, there is a commutative diagram
\begin{equation}\label{dd}
\alignbox{
\xymatrix{
&\aB_0\ox W_{p+1}\ar[rr]^{1\ox d}&&\aB_0\ox\aB_0\ox W_p\ar[dr]^{m\ox1}\\
W_{p+2}\ar@{->>}[d]\ar[ur]^d\ar@{-->}[drr]_{dd}&&&&\aB_0\ox W_p\\
V_{p+2}\ar@{-->}[drr]_\ph&&R_\aB\ox W_p\ar@{ >->}[urr]\ar@{->>}[d]\\
&&R_\aB\ox V_p
}
}
\end{equation}

Then in terms of the above diagrams of $\FF$-vector spaces, the condition of \ref{delta} can be expressed as follows:

\begin{Corollary}
Completions of \ref{precomplex} to a secondary resolution are in one-to-one correspondence with sequences of maps
$$
\delta^\A_p:V_{p+2}\to\Sigma\A\ox V_p,\ \ p\ge0
$$
making the diagrams below commute, with $\ph$ defined by \eqref{dd}.
\begin{equation}\label{diadelta}
\alignbox{
\xymatrix{
&\Sigma\A\ox V_{p+1}\ar[r]^-{1\ox d}&\Sigma\A\ox\A\ox V_p\ar[dr]^{m\ox1}\\
V_{p+3}\ar[ur]^{\delta^\A_{p+1}}\ar[dr]_d&&\A\ox R_\aB\ox V_p\ar[r]^{A\ox1}&\Sigma\A\ox V_p\\
&\A\ox V_{p+2}\ar[ur]^{1\ox\ph}\ar@{}[urr]|{+}\ar[r]_-{1\ox\delta^\A_p}&\A\ox\Sigma\A\ox V_p\ar[ur]_{m\ox1}
}
}
\end{equation}
\end{Corollary}\qed

We can use this to construct the secondary resolution inductively.
Just start by introducing values of $\delta$ on the generators as
expressions with indeterminate coefficients; the equation \eqref{maineq} will
impose linear conditions on these coefficients. These are then solved
degree by degree. For example, in degree 2 one may have
$$
\delta(g_2^2)=\eta_2^2(\Sq^1)\Sq^1g_0^0
$$
for some $\eta_2^2(\Sq^1)\in\F$. Similarly in degree 3 one may have
$$
\delta(g_3^3)=\eta_3^3(\Sq^1)\Sq^1g_1^1+\eta_3^3(1)g_1^2.
$$
Then one will get
$$
d\delta(g_3^3)=\eta_3^3(\Sq^1)\Sq^1d(g_1^1)+\eta_3^3(1)d(g_1^2)=\eta_3^3(\Sq^1)\Sq^1\Sq^1g_0^0+\eta_3^3(1)\Sq^2g_0^0=\eta_3^3(1)\Sq^2g_0^0
$$
and
\begin{multline*}
\delta d(g_3^3)=\delta(\Sq^1g_2^2)\\
=\Sq^1\delta(g_2^2)+A(\Sq^1,dd(g_2^2))=\eta_2^2(\Sq^1)\Sq^1\Sq^1g_0^0+A(\Sq^1,d(\Sq^1g_1^1))=A(\Sq^1,\Sq^1\Sq^1g_0^0)=0;
\end{multline*}
thus \eqref{maineq} forces $\eta^3_3(1)=0$.

Similarly one puts $\delta(g_m^d)=\sum_{m-2\le d'\le d-1}\sum_a\eta_m^d(a)ag_{m-2}^{d'}$,
with $a$ running over a basis in $\A^{d-1-d'}$, and then substituting this in
\eqref{maineq} gives linear equations on the numbers $\eta_m^d(a)$. Solving
these equations and choosing the remaining $\eta$'s arbitrarily then gives
values of the differential $\delta$ in the secondary resolution.

Then finally to obtain the secondary differential
$$
d_{(2)}:\Ext^n_\A(\FF,\FF)^m\to\Ext^{n+2}_\A(\FF,\FF)^{m+1}
$$
from this $\delta$, one just applies the functor $\Hom_\A(\_,\FF)$ to the
initial minimal resolution and calculates the map induced by $\delta$ on
cohomology of the resulting cochain complex, i.~e. on $\Ext^*_\A(\FF,\FF)$.
In fact since \eqref{minireso} is a minimal resolution, the value of
$\Hom_\A(\_,\FF)$ on it coincides with its own cohomology and is the
$\FF$-vector space of those linear maps $\A\brk{g_*^*}\to\FF$ which vanish
on all elements of the form $ag_*^*$ with $a$ of positive degree.

Let us then identify $\Ext^*_\A(\FF,\FF)$ with this space and choose a basis
in it consisting of elements $\hat g_m^d$ defined as the maps sending the
generator $g_m^d$ to 1 and all other generators to 0. One then has
$$
(d_{(2)}(\hat g_m^d))(g_{m'}^{d'})=\hat g_m^d\delta(g_{m'}^{d'}).
$$
The right hand side is nonzero precisely when $g_m^d$ appears in
$\delta(g_{m'}^{d'})$ with coefficient 1, i.~e. one has
\begin{equation}\label{d2gen}
d_{(2)}(\hat g_m^d)=\sum_{\textrm{$g_m^d$ appears in
$\delta(g_{m+2}^{d+1})$}}\hat g_{m+2}^{d+1}.
\end{equation}

For example, looking at the table of values of $\delta$ below we see that
the first instance of a $g_m^d$ appearing with coefficient 1 in a value of
$\delta$ on a generator is
$$
\delta(g_3^{17})=g_1^{16}+ \Sq^{12} g_1^4+\Sq^{10}\Sq^4
g_1^2+(\Sq^9\Sq^4\Sq^2+\Sq^{10}\Sq^5+\Sq^{11}\Sq^4)g_1^1.
$$
This means
$$
d_{(2)}(\hat g_1^{16})=\hat g_3^{17}
$$
and moreover $d_{(2)}(\hat g_m^d)=0$ for all $g_m^d$ with $d<17$ (one can
check all cases for each given $d$ since the number of generators $g_m^d$
for each given $d$ is finite).

Treating similarly the rest of the table below we find that the only
nonzero values of $d_{(2)}$ on generators of degree $<40$ are as follows:
$$
\begin{array}{rl}
d_{(2)}(\hat g_1^{16})&=\hat g_3^{17}\\
d_{(2)}(\hat g_4^{21})&=\hat g_6^{22}\\
d_{(2)}(\hat g_4^{22})&=\hat g_6^{23}\\
d_{(2)}(\hat g_5^{23})&=\hat g_7^{24}\\
d_{(2)}(\hat g_7^{30})&=\hat g_9^{31}\\
d_{(2)}(\hat g_8^{31})&=\hat g_{10}^{32}\\
d_{(2)}(\hat g_1^{32})&=\hat g_3^{33}\\
d_{(2)}(\hat g_2^{33})&=\hat g_4^{34}\\
d_{(2)}(\hat g_7^{33})&=\hat g_9^{34}\\
d_{(2)}(\hat g_8^{33})&=\hat g_{10}^{34}\\
d_{(2)}({}'\hat g_3^{34})&=\hat g_5^{35}\\
d_{(2)}(\hat g_8^{34})&=\hat g_{10}^{35}\\
d_{(2)}(\hat g_7^{36})&=\hat g_9^{37}\\
d_{(2)}(\hat g_8^{37})&=\hat g_{10}^{38}.
\end{array}
$$
These data can be summarized in the following picture, thus confirming
calculations presented in Ravenel's book \cite{Ravenel}.

\

\newcount\s
\def\latticebody{
\s=\latticeA
\advance\s by\latticeB
\ifnum\s<40\drop{.}\else\fi}
\ \hfill
\xy
*\xybox{
0;<.3cm,0cm>:
,0,{\xylattice0{40}0{15}}
,(0,0)*{\bullet}
,(0,1)*{\bullet}
,(0,2)*{\bullet}
,(0,3)*{\bullet}
,(0,4)*{\bullet}
,(0,5)*{\bullet}
,(0,6)*{\bullet}
,(0,7)*{\bullet}
,(0,8)*{\bullet}
,(0,9)*{\bullet}
,(0,10)*{\bullet}
,(0,11)*{\bullet}
,(0,12)*{\bullet}
,(0,13)*{\bullet}
,(0,14)*{\bullet}
,(0,15)*{\bullet}
,(1,1)*{\bullet}
,(3,1)*{\bullet}
,(7,1)*{\bullet}
,(15,1)*{\circ}="1.16"
,(31,1)*{\circ}="1.32"
,(2,2)*{\bullet}
,(3,2)*{\bullet}
,(6,2)*{\bullet}
,(7,2)*{\bullet}
,(8,2)*{\bullet}
,(14,2)*{\bullet}
,(15,2)*{\bullet}
,(16,2)*{\bullet}
,(18,2)*{\bullet}
,(30,2)*{\bullet}
,(31,2)*{\circ}="2.33"
,(32,2)*{\bullet}
,(34,2)*{\bullet}
,(3,3)*{\bullet}
,(7,3)*{\bullet}
,(8,3)*{\bullet}
,(9,3)*{\bullet}
,(14,3)*{\circ};"1.16"**\dir{-}
,(15,3)*{\bullet}
,(17,3)*{\bullet}
,(18,3)*{\bullet}
,(19,3)*{\bullet}
,(21,3)*{\bullet}
,(30,3)*{\circ};"1.32"**\dir{-}
,(31.15,2.85)*{\bullet}
,(30.85,3.15)*{\circ}="3.34"
,(33,3)*{\bullet}
,(34,3)*{\bullet}
,(7,4)*{\bullet}
,(9,4)*{\bullet}
,(14,4)*{\bullet}
,(15,4)*{\bullet}
,(17,4)*{\circ}="4.21"
,(18.15,3.85)*{\bullet}
,(17.85,4.15)*{\circ}="4.22"
,(20,4)*{\bullet}
,(22,4)*{\bullet}
,(23,4)*{\bullet}
,(30,4)*{\circ};"2.33"**\dir{-}
,(31,4)*{\bullet}
,(32,4)*{\bullet}
,(33,4)*{\bullet}
,(34,4)*{\bullet}
,(9,5)*{\bullet}
,(11,5)*{\bullet}
,(14,5)*{\bullet}
,(14.85,5.15)*{\bullet}
,(15.15,4.85)*{\bullet}
,(17,5)*{\bullet}
,(18,5)*{\circ}="5.23"
,(20,5)*{\bullet}
,(21,5)*{\bullet}
,(23,5)*{\bullet}
,(24,5)*{\bullet}
,(30,5)*{\circ};"3.34"**\dir{-}
,(30.85,5.15)*{\bullet}
,(31.15,4.85)*{\bullet}
,(33,5)*{\bullet}
,(10,6)*{\bullet}
,(11,6)*{\bullet}
,(14,6)*{\bullet}
,(15,6)*{\bullet}
,(16,6)*{\circ};"4.21"**\dir{-}
,(17,6)*{\circ};"4.22"**\dir{-}
,(20,6)*{\bullet}
,(23,6)*{\bullet}
,(26,6)*{\bullet}
,(30,6)*{\bullet}
,(31,6)*{\bullet}
,(32,6)*{\bullet}
,(11,7)*{\bullet}
,(15,7)*{\bullet}
,(16,7)*{\bullet}
,(17,7)*{\circ};"5.23"**\dir{-}
,(23,7)*{\circ}="7.30"
,(26,7)*{\circ}="7.33"
,(29,7)*{\circ}="7.36"
,(30,7)*{\bullet}
,(31,7)*{\bullet}
,(32,7)*{\bullet}
,(15,8)*{\bullet}
,(17,8)*{\bullet}
,(22,8)*{\bullet}
,(23,8)*{\circ}="8.31"
,(25,8)*{\circ}="8.33"
,(26,8)*{\circ}="8.34"
,(28,8)*{\bullet}
,(29,8)*{\circ}="8.37"
,(30,8)*{\bullet}
,(30.85,8.15)*{\bullet}
,(31.15,7.85)*{\bullet}
,(17,9)*{\bullet}
,(19,9)*{\bullet}
,(22,9)*{\circ};"7.30"**\dir{-}
,(22.85,9.15)*{\bullet}
,(23.15,8.85)*{\bullet}
,(25,9)*{\circ};"7.33"**\dir{-}
,(26,9)*{\circ}="9.35"
,(28,9)*{\circ};"7.36"**\dir{-}
,(29,9)*{\bullet}
,(30,9)*{\bullet}
,(18,10)*{\bullet}
,(19,10)*{\bullet}
,(22,10)*{\circ};"8.31"**\dir{-}
,(23,10)*{\bullet}
,(24,10)*{\circ};"8.33"**\dir{-}
,(25,10)*{\circ};"8.34"**\dir{-}
,(28,10)*{\circ};"8.37"**\dir{-}
,(19,11)*{\bullet}
,(23,11)*{\bullet}
,(24,11)*{\bullet}
,(25,11)*{\circ};"9.35"**\dir{-}
,(23,12)*{\bullet}
,(25,12)*{\bullet}
,(25,13)*{\bullet}
,(0,17)*{\ }
,(42,0)*{\ }
}="O"
,{"O"+L \ar "O"+R*+!LD{d-m}}
,{"O"+D \ar "O"+U*+!RD{m}}
\endxy
\hfill\

\section{The table of values of the differential $\delta$ in the secondary
resolution for $\G^\Sigma$}

The following table presents results of computer calculations of the
differential $\delta$. Note that it does not have invariant meaning since
it depends on the choices involved in determination of the multiplication
map $A$, of the resolution and of those indeterminate coefficients
$\eta_m^d(a)$ which remain undetermined after the conditions
\eqref{maineq} are satisfied. The resulting secondary differential
$d_{(2)}$ however does not depend on these choices and is canonically
determined.

$$
\begin{array}{rl}
\delta(g_2^2) &= 0\\
\ \\
\delta(g_3^3) &= 0\\
\ \\
\delta(g_2^4) &= 0\\
\delta(g_4^4) &= 0\\
\ \\
\delta(g_2^5) &= 0\\
\delta(g_5^5) &= 0\\
\ \\
\delta(g_3^6) &= \Sq^4 g_1^1\\
\delta(g_6^6) &= 0\\
\ \\
\delta(g_7^7) &= 0\\
\ \\
\delta(g_2^8) &= 0\\
\delta(g_8^8) &= 0\\
\ \\
\delta(g_2^9) &= 0\\
\delta(g_9^9) &= 0\\
\ \\
\delta(g_2^{10}) &= 0\\
\delta(g_3^{10}) &= (\Sq^4\Sq^2\Sq^1 + \Sq^7) g_1^2\\
 &+ \Sq^8 g_1^1\\
\delta(g_{10}^{10}) &= 0\\
\end{array}
$$
\hfill

$$
\begin{array}{rl}
\delta(g_3^{11}) &= (\Sq^7\Sq^1 + \Sq^8) g_1^2\\
 &+ \Sq^6\Sq^3 g_1^1\\
\delta(g_4^{11}) &= \Sq^5 g_2^5\\
 &+\Sq^4\Sq^2 g_2^4\\
\delta(g_{11}^{11}) &= 0\\
\ \\
\delta(g_3^{12}) &= \Sq^7\Sq^3 g_1^1\\
\delta(g_{12}^{12}) &= 0\\
\ \\
\delta(g_4^{13}) &= \Sq^4 g_2^8\\
 &+ (\Sq^7 + \Sq^5\Sq^2) g_2^5\\
 &+ (\Sq^8 + \Sq^6\Sq^2) g_2^4\\
 &+ (\Sq^7\Sq^3 + \Sq^8\Sq^2 + \Sq^{10}) g_2^2\\
\delta(g_{13}^{13}) &= 0\\
\ \\
\delta(g_5^{14}) &= \Sq^4\Sq^2\Sq^1 g_3^6\\
 &+ (\Sq^7\Sq^3 + \Sq^8\Sq^2) g_3^3\\
\delta(g_{14}^{14}) &= 0\\
\ \\
\delta(g_2^{16}) &= 0\\
\delta(g_5^{16}) &=\Sq^3 g_3^{12}\\
 &+ \Sq^4 g_3^{11}\\
 &+ \Sq^5 g_3^{10}\\
 &+ \Sq^{10}\Sq^2 g_3^3\\
\delta(g_6^{16}) &= 0\\
\ \\
\delta(g_2^{17}) &= 0\\
\delta(g_3^{17}) &= g_1^{16}\\
 &+ \Sq^{12} g_1^4\\
 &+\Sq^{10}\Sq^4 g_1^2\\
 &+ (\Sq^9\Sq^4\Sq^2 + \Sq^{10}\Sq^5 + \Sq^{11}\Sq^4) g_1^1\\
\delta(g_6^{17}) &= (\Sq^5 + \Sq^4\Sq^1) g_4^{11}\\
 &+ (\Sq^{12} +\Sq^{10}\Sq^2) g_4^4\\
\ \\
\delta(g_2^{18}) &= 0\\
\delta(g_3^{18}) &=(\Sq^{11}\Sq^4 + \Sq^8\Sq^4\Sq^2\Sq^1) g_1^2\\
 &+(\Sq^{10}\Sq^4\Sq^2 + \Sq^{11}\Sq^5 + \Sq^{12}\Sq^4 +\Sq^{14}\Sq^2 + \Sq^{16}) g_1^1\\
\delta(g_4^{18}) &=(\Sq^6\Sq^1 + \Sq^7) g_2^{10}\\
 &+(\Sq^6\Sq^3 + \Sq^7\Sq^2 + \Sq^9) g_2^8\\
 &+\Sq^8\Sq^4g_2^5\\
 &+(\Sq^{10}\Sq^2\Sq^1 + \Sq^{13} + \Sq^{11}\Sq^2 +\Sq^{12}\Sq^1)g_2^4\\
 &+(\Sq^9\Sq^4\Sq^2 + \Sq^{15} + \Sq^{12}\Sq^3 +\Sq^{10}\Sq^5)g_2^2\\
\delta(g_7^{18}) &= \Sq^2\Sq^1 g_5^{14}\\
\ \\
\delta(g_4^{19}) &=\Sq^9 g_2^9\\
 &+(\Sq^{10} + \Sq^8\Sq^2) g_2^8\\
 &+ \Sq^{11}\Sq^2g_2^5\\
 &+ ( \Sq^{11}\Sq^2\Sq^1 + \Sq^{13}\Sq^1 + \Sq^8\Sq^4\Sq^2 + \Sq^{10}\Sq^3\Sq^1) g_2^4\\
 &+(\Sq^{14}\Sq^2 +\Sq^{10}\Sq^4\Sq^2 + \Sq^{12}\Sq^4) g_2^2\\
\delta(g_5^{19}) &=\Sq^1 g_3^{17}\\
 &+ \Sq^4\Sq^2 g_3^{12}\\
 &+ \Sq^4\Sq^2\Sq^1 g_3^{11}\\
 &+ (\Sq^6\Sq^2 + \Sq^8) g_3^{10}\\
 &+ (\Sq^8\Sq^4 +\Sq^{11}\Sq^1) g_3^6\\
 &+ (\Sq^{13}\Sq^2 + \Sq^{10}\Sq^5 + \Sq^{15} + \Sq^{11}\Sq^4) g_3^3
 \end{array}
$$

$$
\begin{array}{rl}
\delta(g_2^{20}) &= 0\\
\delta(g_3^{20}) &= (\Sq^{15} + \Sq^9\Sq^4\Sq^2) g_1^4\\
 &+ (\Sq^{12}\Sq^5 + \Sq^{13}\Sq^4 + \Sq^{16}\Sq^1) g_1^2\\
 &+ (\Sq^{11}\Sq^5\Sq^2 + \Sq^{15}\Sq^3 + \Sq^{18} + \Sq^{12}\Sq^6)g_1^1\\
\delta(g_5^{20}) &= \Sq^4\Sq^2\Sq^1 g_3^{12}\\
 &+ (\Sq^7\Sq^1 + \Sq^8)g_3^{11}\\
 &+ (\Sq^{10}\Sq^3 + \Sq^8\Sq^4\Sq^1 + \Sq^{13} + \Sq^{11}\Sq^2)g_3^6\\
 &+ (\Sq^{13}\Sq^3 + \Sq^{10}\Sq^4\Sq^2 + \Sq^{11}\Sq^5 + \Sq^{12}\Sq^4)g_3^3\\
\delta({}'g_5^{20}) &= \Sq^5\Sq^2 g_3^{12}\\
 &+ \Sq^7\Sq^2 g_3^{10}\\
 &+ (\Sq^{12}\Sq^1 + \Sq^{10}\Sq^3 + \Sq^8\Sq^4\Sq^1 + \Sq^{10}\Sq^2\Sq^1
 + \Sq^{11}\Sq^2) g_3^6\\
 &+ (\Sq^{14}\Sq^2 + \Sq^{13}\Sq^3 + \Sq^{11}\Sq^5 + \Sq^{16} + \Sq^{12}\Sq^4)g_3^3\\
\delta(g_6^{20}) &= (\Sq^6\Sq^2 + \Sq^8) g_4^{11}\\
 &+ (\Sq^{13}\Sq^2 + \Sq^{15} + \Sq^{11}\Sq^4) g_4^4\\
\ \\
\delta(g_3^{21}) &= (\Sq^{15}\Sq^2\Sq^1 + \Sq^{17}\Sq^1 +
 \Sq^{12}\Sq^6) g_1^2\\
 &+ ( \Sq^{13}\Sq^4\Sq^2 + \Sq^{15}\Sq^4 + \Sq^{16}\Sq^3 + \Sq^{17}\Sq^2 + \Sq^{19}) g_1^1\\
\delta(g_4^{21}) &=\Sq^3 g_2^{17}\\
 &+ (\Sq^{10} +\Sq^9\Sq^1) g_2^{10}\\
 &+ (\Sq^9\Sq^3 + \Sq^{11}\Sq^1) g_2^8\\
 &+ (\Sq^{15} +\Sq^{13}\Sq^2 + \Sq^{10}\Sq^5) g_2^5\\
 &+ ( \Sq^{13}\Sq^2\Sq^1 + \Sq^{12}\Sq^3\Sq^1 + \Sq^{12}\Sq^4 + \Sq^9\Sq^4\Sq^2\Sq^1 + \Sq^{10}\Sq^4\Sq^2) g_2^4\\
 &+ (\Sq^{16}\Sq^2 + \Sq^{12}\Sq^6 + \Sq^{15}\Sq^3) g_2^2\\
\delta(g_6^{21}) &=  (\Sq^7 + \Sq^6\Sq^1) g_4^{13}\\
 &+ (\Sq^9 + \Sq^8\Sq^1) g_4^{11}\\
 &+ \Sq^{11}\Sq^5 g_4^4\\
\ \\
\delta(g_3^{22}) &=\Sq^{17} g_1^4\\
 &+ (\Sq^{16}\Sq^2\Sq^1 + \Sq^{13}\Sq^6 +\Sq^{12}\Sq^4\Sq^2\Sq^1 + \Sq^{12}\Sq^6\Sq^1) g_1^2\\
 &+ ( \Sq^{13}\Sq^5\Sq^2 + \Sq^{17}\Sq^3 + \Sq^{18}\Sq^2 +\Sq^{14}\Sq^4\Sq^2)g_1^1\\
\delta(g_4^{22}) &=\Sq^4 g_2^{17}\\
 &+ \Sq^{11} g_2^{10}\\
 &+ (\Sq^{12} + \Sq^9\Sq^3) g_2^9\\
 &+ (\Sq^9\Sq^4 + \Sq^{13} + \Sq^8\Sq^4\Sq^1)g_2^8\\
 &+ \Sq^{12}\Sq^4g_2^5\\
 &+ \Sq^{15}\Sq^2 g_2^4\\
 &+ (\Sq^{13}\Sq^4\Sq^2 + \Sq^{19} + \Sq^{13}\Sq^6 + \Sq^{14}\Sq^5)g_2^2\\
\delta({}'g_4^{22}) &= \Sq^2\Sq^1 g_2^{18}\\
 &+ (\Sq^8\Sq^4 + \Sq^{12}) g_2^9\\
 &+ (\Sq^9\Sq^4 + \Sq^{13} + \Sq^{12}\Sq^1) g_2^8\\
 &+ (\Sq^{16} + \Sq^{13}\Sq^3) g_2^5\\
 &+ (\Sq^{15}\Sq^2 + \Sq^{16}\Sq^1 + \Sq^{13}\Sq^4 + \Sq^{11}\Sq^4\Sq^2)g_2^4\\
 &+ (\Sq^{14}\Sq^5 + \Sq^{19} + \Sq^{17}\Sq^2) g_2^2\\
\delta(g_5^{22}) &= (\Sq^7\Sq^2 + \Sq^6\Sq^2\Sq^1 + \Sq^6\Sq^3)g_3^{12}\\
 &+ \Sq^{10} g_3^{11} + (\Sq^9\Sq^2 + \Sq^8\Sq^3 + \Sq^{11}) g_3^{10}\\
 &+(\Sq^{14}\Sq^1 + \Sq^{11}\Sq^3\Sq^1 + \Sq^{12}\Sq^3 +\Sq^{13}\Sq^2)g_3^6\\
 &+ \Sq^{13}\Sq^5 g_3^3\\
\delta(g_6^{22}) &= g_4^{21}\\
 &+ (\Sq^6\Sq^2 + \Sq^8 + \Sq^7\Sq^1) g_4^{13}\\
 &+ \Sq^{10} g_4^{11}\\
 &+ (\Sq^{13}\Sq^4 + \Sq^{15}\Sq^2 + \Sq^{17}) g_4^4\\
\delta(g_7^{22}) &= (\Sq^{13}\Sq^3 + \Sq^{14}\Sq^2 + \Sq^{16}) g_5^5\\
\end{array}
$$

%% file: dHascoAss_c4.tex
\chapter[Hopf pair algebras and Hopf pair coalgebras]{Hopf pair algebras and Hopf pair coalgebras representing the algebra
of secondary cohomology operations}\label{Hpa}

We describe a modification $\aB^\FF$ of the algebra $\aB$ of secondary
cohomology operations in chapter \ref{secsteen} which is suitable for
dualization. The resulting object $\aB^\FF$ and the dual object $\aB_\FF$
will be used to give an alternative description of the multiplication map $A$
and the dual multiplication map $A_*$. All triple Massey products in the
Steenrod algebra can be deduced from $\aB^\FF$ or $\aB_\FF$ and from $A$ and
$A_*$.

We first recall the notions of pair modules and pair algebras from chapter
\ref{sext} and give the corresponding dual notions. Next we define the
concept of $M$-algebras and $N$-coalgebras, where $M$ is a folding system and
$N$ an unfolding system. An $M$-algebra is a variation on the notion of a
$[p]$-algebra from \cite{Baues}. We show that the algebra $\aB$ of secondary
cohomology operations gives rise to a comonoid $\aB^\FF$ in the monoidal
category of $M$-algebras, and we describe the dual object $\aB_\FF$, which is
a monoid in the monoidal category of $N$-coalgebras.

In chapter \ref{LS} we study the algebraic objects $\aB^\FF$ and $\aB_\FF$ in terms
of generators. This way we obtain explicit descriptions which can be used
for computations. In particular we characterize algebraically
multiplication maps $A_\phi$ and comultiplication maps $A^\psi$ which
determine $\aB^\FF$ and $\aB_\FF$ completely, see sections \ref{rcomp},
\ref{bcomp}, \ref{cobcomp}. For the dual object $\aB_\FF$ the inclusion of
polynomial algebras $\A_*\subset\F_*$ will be crucial. Here $\A_*$ is the
Milnor dual of the Steenrod algebra and $\F_*$ is the dual of a free
associative algebra.

\section{Pair modules and pair algebras}\label{pairs}

We here recall from \ref{secodjf} the following notation in order to prepare the 
reader for the dualization of this notation in the next section.
Let $k$ be a commutative ring (usually it will be actually a prime field
$\FF=\FF_p=\ZZ/p\ZZ$ for some prime $p$) and let $\Mod$ be the category of
finite dimensional $k$-modules (i.~e. $k$-vector spaces) and $k$-linear
maps. 
A \emph{pair module} is a
homomorphism
\begin{equation}\label{pair}
X=\left(X_1\xto\d X_0\right)
\end{equation}
in $\Mod$.  We write $\pi_0(X)=\coker\d$ and $\pi_1(X)=\ker\d$.

For two pair modules $X$ and $Y$ the tensor product of the
complexes corresponding to them is concentrated in degrees in 0, 1 and 2
and is given by
\begin{equation}\label{tens}
X_1\!\ox\!Y_1\xto{\d_1}X_1\!\ox\!Y_0\oplus
X_0\!\ox\!Y_1\xto{\d_0}X_0\!\ox\!Y_0
\end{equation}
with $\d_0=(\d\ox1,1\ox\d)$ and $\d_1=\binom{-1\ox\d}{\d\ox1}$. Truncating this
chain complex we get the pair module
$$
X\bar\ox Y=
\left((X\bar\ox Y)_1=\coker(\d_1)\xto\d X_0\ox Y_0
=(X\bar\ox Y)_0\right)
$$
with $\d$ induced by $\d_0$. Clearly one has $\pi_0(X\bar\ox Y)\cong\pi_0(X)\ox\pi_0(Y)$
and 

\begin{equation}\label{tox}
\pi_1(X\bar\ox Y)\cong\pi_1(X)\!\ox\!\pi_0(Y)\oplus\pi_0(X)\!\ox\!\pi_1(Y).
\end{equation}

We next consider the category $\Mod^\cdot$ of \emph{graded} modules, i.~e.
graded objects in $\Mod$ (graded $k$-vector spaces
$A^\cdot=(A^n)_{n\in\ZZ}$ with upper indices, which in each degree have
finite dimension). For graded modules $A^\cdot$, $B^\cdot$ we define their
graded tensor product $A^\cdot\ox B^\cdot$ in the usual way with
an interchange
\begin{equation}\label{symm}
T_{A^\cdot,B^\cdot}:A^\cdot\ox B^\cdot\xto\cong B^\cdot\ox A^\cdot .
\end{equation}

A \emph{graded pair module} is a graded object of $\Mod_*$, i.~e. a
sequence $X^n=(\d^n:X_1^n\to X_0^n)$ with $n\in\ZZ$ of pair modules. 
The tensor product $X^\cdot\bar\ox Y^\cdot$ of graded pair modules
$X^\cdot$, $Y^\cdot$ is defined by 
\begin{equation}\label{grpr}
(X^\cdot\bar\ox Y^\cdot)^n=\bigoplus_{i+j=n}X^i\bar\ox Y^j.
\end{equation}

For two morphisms $f,g:X^\cdot\to Y^\cdot$ between graded pair modules, a
\emph{homotopy} $H:f\then g$ is a morphism $H:X^\cdot_0\to Y^\cdot_1$ of degree 0 satisfying
\begin{equation}
f_0-g_0=\d H \textrm{ and } f_1-g_1=H\d .
\end{equation}

\begin{Definition}
A \emph{pair algebra} $B^\cdot$ is a graded pair module, i.~e. an object
$$
\d^\cdot:B^\cdot_1\to B^\cdot_0
$$
in $\Mod^\cdot_*$ with $B_1^n=B_0^n=0$ for $n<0$ such that $B^\cdot_0$ is a graded
algebra in $\Mod^\cdot$, $B^\cdot_1$ is a graded $B^\cdot_0$-$B^\cdot_0$-bimodule, and
$\d^\cdot$ is a bimodule homomorphism. Moreover for $x,y\in B^\cdot_1$ the equality
\begin{equation}\label{paireq}
\d(x)y=x\d(y)
\end{equation}
holds in $B^\cdot_1$.
\end{Definition}

It is easy to see that a graded pair algebra $B^\cdot$ yields an exact
sequence of graded $B^\cdot_0$-$B^\cdot_0$-bimodules
\begin{equation}\label{piseq}
0\to\pi_1B^\cdot\to B^\cdot_1\xto\d B^\cdot_0\to\pi_0B^\cdot\to0
\end{equation}
where in fact $\pi_0B^\cdot$ is a graded $k$-algebra, $\pi_1B^\cdot$ is a graded
$\pi_0B^\cdot$-$\pi_0B^\cdot$-bimodule, and $B^\cdot_0\to\pi_0B^\cdot$ is a homomorphism of
graded $k$-algebras.

The tensor product of pair algebras has a natural pair algebra structure,
as it happens in any symmetric monoidal category.

We are mainly interested in two examples of pair algebras defined below in
sections \ref{grel} and \ref{sec} respectively: the \emph{$\GG$-relation
pair algebra $\R$ of the Steenrod algebra $\A$} and the \emph{pair algebra
$\aB$ of secondary cohomology operations} deduced from
\cite{Baues}*{5.5.2}.

By the work of Milnor \cite{Milnor} it is well known that the dual of the
Steenrod algebra $\A$ is a polynomial algebra and this fact yields
important algebraic properties of $\A$. For this reason we also consider
the dual of the $\GG$-relation pair algebra $\R$ of $\A$ and the dual of
the pair algebra $\aB$ of secondary cohomology operations. The duality
functor $D$ is studied in the next section.

\section{Pair comodules and pair coalgebras}\label{coalg}

This section is exactly dual to the previous one. There is a contravariant
self-equivalence of categories
$$
D=\Hom_k(\_,k):\Mod\op\to\Mod
$$
which carries a vector space $V$ in $\Mod$ to its dual
$$
DV=\Hom_k(V,k).
$$
We also denote the dual of $V$ by $V_*=DV$, for example, the dual of the
Steenrod algebra $\A$ is $\A_*=D(\A)$. We can apply the functor
$\Hom_k(\_,k)$ to dualize straightforwardly all notions of section
\ref{pairs}. Explicitly, one gets:

A \emph{pair comodule} is a homomorphism
\begin{equation}\label{copair}
X=\left(X^1\xot dX^0\right)
\end{equation}
in $\Mod$. We write $\pi^0(X)=\ker d$ and $\pi^1(X)=\coker d$.
The dual of a pair module $X$ is a pair comodule
\begin{align*}
DX&=\Hom_k(X,k)\\
&=(D\d:DX_0\to DX_1)
\end{align*}
with $(DX)^i=D(X_i)$. 
A \emph{morphism} $f:X\to Y$ of pair comodules is a commutative diagram
$$
\xymatrix
{
X^1\ar[r]^{f^1}&Y^1\\
X^0\ar[u]^d\ar[r]^{f_0}&Y^0\ar[u]^d.
}
$$
Evidently pair comodules with these morphisms form a category $\Mod^*$ and
one has functors
$$
\pi^0, \pi^1 : \Mod^*\to\Mod.
$$
which are compatible with the duality functor $D$, that is, for any pair
module $X$ one has
$$
\pi_i(DX)=D(\pi_iX)\textrm{ for }i=0,1.
$$
A morphism of pair comodules is called a \emph{weak equivalence} if it
induces isomorphisms on $\pi^0$ and $\pi^1$.

Clearly a pair comodule is the same as a cochain complex concentrated in
degrees 0 and 1. For two pair comodules $X$ and $Y$ the tensor product
of the cochain complexes is concentrated in degrees in 0, 1 and 2 and is
given by
$$
X^1\!\ox\!Y^1\xot{d^1}X^1\!\ox\!Y^0\oplus
X^0\!\ox\!Y^1\xot{d^0}X^0\!\ox\!Y^0
$$
with $d^0=\binom{d\ox1}{1\ox d}$ and $d^1=(-1\ox d,d\ox1)$. Cotruncating this
cochain complex we get the pair comodule
$$
X\dblb\ox Y=
\left((X\dblb\ox Y)^1=\ker(d^1)\xot dX^0\ox Y^0=(X\dblb\ox Y)^0\right)
$$
with $d$ induced by $d_0$. One readily checks the natural isomorphism
\begin{equation}
D(X\bar\ox Y)\cong DX\dblb\ox DY.
\end{equation}

\begin{Remark}[compare \ref{trunc}]
Note that the full embedding of the category of pair comodules into the category of
cochain complexes induced by the above identification has a right adjoint
$\Tr^*$ given by cotruncation: for a cochain complex
$$
C^*=\left(...\ot C^2\xot{d^1}C^1\xot{d^0}C^0\xot{d^{-1}}C^{-1}\ot...\right),
$$
one has
$$
\Tr^*(C^*)=\left(\ker(d^1)\xot{{\bar d}^0}C^0\right),
$$
with ${\bar d}^0$ induced by $d^0$. Then clearly one has
$$
X\dblb\ox Y=\Tr^*(X\ox Y).
$$
Using the fact that $\Tr^*$ is a coreflection onto a full subcategory, one
easily checks that the category $\Mod^*$ together with the tensor
product $\dblb\ox$ and unit $k^*=(0\ot k)$ is a symmetric monoidal category,
and $\Tr^*$ is a monoidal functor.
\end{Remark}

We next consider the category $\Mod_\bullet$ of \emph{graded} modules, i.~e.
graded objects in $\Mod$ (graded $k$-vector spaces
$A_\cdot=(A_n)_{n\in\ZZ}$ with lower indices which in each degree have
finite dimension). For graded modules $A_\cdot$, $B_\cdot$ we define their
graded tensor product $A_\cdot\ox B_\cdot$ again in the usual way, i.~e.
by
$$
(A_\cdot\ox B_\cdot)_n=\bigoplus_{i+j=n}A_i\ox B_j.
$$

A \emph{graded pair comodule} is a graded object of $\Mod^*$, i.~e. a
sequence $X_n=(d_n:X^0_n\to X^1_n)$ of pair comodules. We can also identify such a
graded pair comodule $X_\cdot$ with the underlying morphism $d$ of degree 0 between
graded modules
$$
X_\cdot=\left(X_\cdot^1\xot{d_\cdot}X_\cdot^0\right).
$$
Now the tensor product $X_\cdot\dblb\ox Y_\cdot$ of graded pair comodules
$X_\cdot$, $Y_\cdot$ is defined by 
\begin{equation}\label{cogrpr}
(X_\cdot\dblb\ox Y_\cdot)_n=\bigoplus_{i+j=n}X_i\dblb\ox Y_j.
\end{equation}
This defines a monoidal structure on the category $\Mod_\bullet$ of graded
pair comodules. Morphisms in this category are of degree 0.

For two morphisms $f,g:X_\cdot\to Y_\cdot$ between graded pair comodules, a
\emph{homotopy} $H:f\then g$ is a morphism $H:X_\cdot^1\to Y_\cdot^0$ of degree 0 as in
the diagram
\begin{equation}\label{cohomot}
\alignbox{
\xymatrix
{
X_\cdot^1\ar[dr]|H\ar@<.5ex>[r]^{f^1}\ar@<-.5ex>[r]_{g^1}&Y_\cdot^1\\
X_\cdot^0\ar[u]_d\ar@<.5ex>[r]^{f^0}\ar@<-.5ex>[r]_{g^0}&Y_\cdot^0\ar[u]^d,
}}
\end{equation}
satisfying $f^0-g^0=Hd$ and $f^1-g^1=dH$.

A \emph{pair coalgebra} $B_\cdot$ is a comonoid in the monoidal category of graded pair
comodules, with the diagonal
$$
\delta:B_\cdot\to B_\cdot\dblb\ox B_\cdot.
$$
We assume that $B_\cdot$ is concentrated in nonnegative degrees, that is
$B_n=0$ for $n<0$.

Of course the duality functor $D$ yields a duality functor
$$
D:(\Mod^\cdot_*)\op\to\Mod_\cdot^*
$$
which is compatible with the monoidal structure, i.~e.
$$
D(X^\cdot\bar\ox Y^\cdot)\cong(DX^\cdot)\dblb\ox(DY^\cdot).
$$
We also write $D(X^\cdot)=X_\cdot$.

More explicitly pair coalgebras can be described as follows.

\begin{Definition}
A \emph{pair coalgebra} $B_\cdot$ is a graded pair comodule, i.~e. an object
$$
d_\cdot:B_\cdot^0\to B_\cdot^1
$$
in $\Mod_\bullet^*$ with $B^1_n=B^0_n=0$ for $n<0$ such that $B_\cdot^0$ is
a graded coalgebra in $\Mod_\bullet$, $B_\cdot^1$ is a graded
$B_\cdot^0$-$B_\cdot^0$-bicomodule, and $d_\cdot$ is abullebullet
homomorphism. Moreover the diagram 
$$
\xymatrix{
B_\cdot^1\ar[r]^\lambda\ar[d]_\rho&
B_\cdot^0\ox B_\cdot^1\ar[d]_{d_\cdot\ox1}\\
B_\cdot^1\ox B_\cdot^0\ar[r]^{1\ox d_\cdot}
&B_\cdot^1\ox B_\cdot^1
}
$$
commutes, where $\lambda$, resp. $\rho$ is the left, resp. right coaction.
\end{Definition}

It is easy to see that there results an exact sequence of graded
$B_\cdot^0$-$B_\cdot^0$-bicomodules dual to \eqref{piseq}
\begin{equation}\label{copiseq}
0\ot\pi^1B_\cdot\ot B_\cdot^1\xot{d_\cdot}B_\cdot^0\ot\pi^0B_\cdot\ot0
\end{equation}
where in fact $\pi^0B_\cdot$ is a graded $k$-coalgebra, $\pi^1B_\cdot$ is a graded
$\pi^0B_\cdot$-$\pi^0B_\cdot$-bicomodule, and $B_\cdot^0\ot\pi^0B_\cdot$ is a homomorphism of
graded $k$-coalgebras.

One sees easily that the notions in this section correspond to those in the
previous section under the duality functor $D=\Hom_k(\_,k)$. In particular,
$D$ carries (graded) pair algebras to (graded) pair coalgebras.

\section{Folding systems}

In this section we associate to a ``right module system'' $M$ a category
of $M$-algebras $\Alg^\r_M$ which is a monoidal category if $M$ is a
``folding system''. Our main examples given by the $\GG$-relation pair
algebra $\R$ of the Steenrod algebra $\A$ and by the pair algebra $\aB$ of
secondary cohomology operations are in fact comonoids in monoidal
categories of such type, see sections \ref{grel} and \ref{sec}. This
generalizes the well known fact that the Steenrod algebra $\A$ is a Hopf
algebra, i.~e. a comonoid in the category of algebras.

\begin{Definition}
Let $\bA$ be a subcategory of the category of graded $k$-algebras. A
\emph{right module system} $M$ over $\bA$ is an assignment, to each
$A\in\bA$, of a right $A$-module $M(A)$, and, to each homomorphism
$f:A\to A'$ in $\bA$, of a homomorphism $f_*:M(A)\to M(A')$ which is
$f$-equivariant, i.~e.
$$
f_*(xa)=f_*(x)f(a)
$$
for any $a\in A$, $x\in M(A)$. The assignment must be functorial, i.~e. one
must have $(\id_A)_*=\id_{M(A)}$ for all $A$ and $(fg)_*=f_*g_*$ for all
composable $f$, $g$.

There are the obvious similar notions of a \emph{left module system} and a
\emph{bimodule system} on a category of graded $k$-algebras $\bA$. Clearly
any bimodule system can be considered as a left module system and a right
module system by forgetting part of the structure.
\end{Definition}

\begin{Examples}\label{sysex}
One obvious example is the bimodule system $\oo$ given by $\oo(A)=A$,
$f_*=f$ for all $A$ and $f$. Another example is the bimodule system
$\Sigma$ given by the suspension. That is, $\Sigma A$ is given by the shift
$$
\Sigma:A^{n-1}=(\Sigma A)^n
$$
($n\in\ZZ$) which is the identity map denoted by $\Sigma$. The bimodule
structure for $a,m\in A$ is given by
\begin{align*}
a(\Sigma m)&=(-1)^{\deg(a)}\Sigma(am),\\
(\Sigma m)a&=\Sigma(ma).
\end{align*}
We shall need the \emph{interchange} of $\Sigma$ which for graded modules
$U$, $V$, $W$ is the isomorphism
\begin{equation}\label{sigma}
\sigma_{U,V,W}:U\ox(\Sigma V)\ox W\xto\cong\Sigma(U\ox V\ox W)
\end{equation}
which carries $u\ox\Sigma v\ox w$ to $(-1)^{\deg(u)}\Sigma(u\ox v\ox w)$.

Clearly a direct sum of module systems is again a module system of the
same kind, so that in particular we get a bimodule system $\oo\oplus\Sigma$
with $(\oo\oplus\Sigma)(A)=A\oplus\Sigma A$. 
\end{Examples}

We are mainly interested in the bimodule system $\oo$ and the bimodule
system $\oo\oplus\Sigma$ which are in fact both folding systems, see
\bref{foldex} below.

\begin{Definition}\label{malg}
For a right module system $M$ on the category of algebras $\bA$ and an algebra
$A$ from $\bA$, an \emph{$M$-algebra of type $A$} is a pair $D_*=(\d:D_1\to D_0)$
with $\pi_0(D_*)=A$ and $\pi_1(D_*)=M(A)$, such that $D_0$ is a
$k$-algebra, the quotient homomorphism $D_0\onto\pi_0D_*=A$ is a
homomorphism of algebras, $D_1$ is a right $D_0$-module, $\d$ is a
homomorphism of right $D_0$-modules, and the induced structure of a right
$\pi_0(D_*)$-module on $\pi_1(D_*)$ conicides with the original right
$A$-module structure on $M$. For $A$, $A'$ in $\bA$, an $M$-algebra $D_*$ of
type $A$, and another one $D'_*$ of type $A'$, a morphism $D_*\to D'_*$ of
$M$-pair algebras is defined to be a commutative diagram of the form
$$
\xymatrix{
0\ar[r]
&M(A)\ar[r]\ar[d]_{f_*}
&D_1\ar[r]^\d\ar[d]^{f_1}
&D_0\ar[r]\ar[d]^{f_0}
&A\ar[r]\ar[d]^f
&0\\
0\ar[r]
&M(A')\ar[r]
&D'_1\ar[r]_{\d'}
&D'_0\ar[r]
&A'\ar[r]
&0
}
$$
where $f_0$ is a homomorphism of algebras and $f_1$ is a right
$f_0$-equivariant $k$-linear map. It is clear how to compose such
morphisms, so that $M$-algebras form a category which we denote
$\Alg^\r_M$.

With obvious modifications, we also get notions of $M$-algebra of type $A$
when $M$ is a left module system or a bimodule system; the corresonding
categories of algebras will be denoted by $\Alg^\l_M$ and $\Alg^\b_M$,
respectively. Moreover, for a bimodule system $M$ there is also a further
full subcategory
$$
\Alg^\pair_M\subset\Alg^\b_M
$$
whose objects, called \emph{$M$-pair algebras} are those $M$-algebras which
satisfy the pair algebra equation $(\d x)y=x\d y$ for all $x,y\in D_1$. 
\end{Definition}

\begin{Remark}\label{initial}
Note that if $\bA$ contains $k$, then $\Alg^?_M$ has an initial object
given by the $M$-algebra $I=(0:M(k)\to k)$ of type $k$. Moreover if $\bA$
contains the trivial algebra $0$, then $\Alg^?_M$ also has a terminal
object --- the $M$-algebra $0=M(0)\to0$ of type $0$. Here ? stands for
$\l$, $\r$ or $\b$ if $M$ is a left-, right-, or bimodule system,
respectively.
\end{Remark}

\begin{Definition}\label{defold}
Let $\bA$ be a category of graded algebras as above which in addition is closed
under tensor product, i.~e. $k$ belongs to $\bA$ and for any $A$, $A'$ from
$\bA$ the algebra $A\ox_kA'$ also belongs to $\bA$. A \emph{right folding system}
on $\bA$ is then defined to be a right module system $M$ on $\bA$ together with
the collection of right $A\ox_kA'$-module homomorphisms
\begin{align*}
\lambda_{A,A'}&:A\ox_kM(A')\to M(A\ox_kA'),\\
\rho_{A,A'}&:M(A)\ox_kA'\to M(A\ox_kA')
\end{align*}
for all $A$, $A'$ in $\bA$ which are natural in the sense that for any
homomorphisms $f:A\to A_1$, $f':A'\to A'_1$ in $\bA$ the diagrams
\begin{equation}\label{natfold}
\alignbox{
\xymatrix{
A\ox_kM(A')\ar[r]^{\lambda_{A,A'}}\ar[d]_{f\ox f'_*}
&M(A\ox_kA')\ar[d]^{(f\ox f')_*}\\
A_1\ox_kM(A'_1)\ar[r]^{\lambda_{A_1,A'_1}}
&M(A_1\ox_kA'_1)
}}
,\ \ \ \ \ \ \ \ \ \ \ \ \ 
\alignbox{
\xymatrix{
M(A)\ox_kA'\ar[r]^{\rho_{A,A'}}\ar[d]_{f_*\ox f'}
&M(A\ox_kA')\ar[d]^{(f\ox f')_*}\\
M(A_1)\ox_kA'_1\ar[r]^{\rho_{A_1,A'_1}}
&M(A_1\ox_kA'_1)
}
}
\end{equation}
commute. Moreover the homomorphisms
\begin{equation}\label{unitfold}
\alignbox{
\lambda_{k,A}&:k\ox_kM(A)\to M(k\ox_kA),\\
\rho_{A,k}&:M(A)\ox_kk\to M(A\ox_kk)
}
\end{equation}
must coincide with the obvious isomorphisms and the diagrams
\begin{gather}\label{leftfold}
\alignbox{\xymatrix{
&A\ox_kM(A'\ox_kA'')\ar[dr]^{\lambda_{A,A'\ox_kA''}}\\
A\ox_kA'\ox_kM(A'')\ar[rr]^{\lambda_{A\ox_kA',A''}}\ar[ur]^{1\ox\lambda_{A',A''}}
&&M(A\ox_kA'\ox_kA''),
}}\\\label{rightfold}
\alignbox{\xymatrix{
&M(A\ox_kA')\ox_kA''\ar[dr]^{\rho_{A\ox_kA',A''}}\\
M(A)\ox_kA'\ox_kA''\ar[rr]^{\rho_{A,A'\ox_kA''}}\ar[ur]^{\rho_{A,A'}A\ox1}
&&M(A\ox_kA'\ox_kA''),
}}\\\label{midfold}
\alignbox{\xymatrix{
&M(A\ox_kA')\ox_kA''\ar[dr]^{\rho_{A\ox_kA',A''}}\\
A\ox_kM(A')\ox_kA''\ar[ur]^{\lambda_{A,A'}\ox1}\ar[dr]_{1\ox\rho_{A',A''}}
&&M(A\ox_kA'\ox_kA'')\\
&A\ox_kM(A'\ox_kA'')\ar[ur]_{\lambda_{A,A'\ox_kA''}}
}}
\end{gather}
must commute for all $A$, $A'$, $A''$ in $\bA$. A folding system is called
\emph{symmetric} if in addition the diagrams
$$
\xymatrix{
A\ox_kM(A')\ar[r]^{\lambda_{A,A'}}\ar[d]_{T_{A,M(A')}}&M(A\ox_kA')\ar[d]^{M(T_{A,A'})}\\
M(A')\ox_kA\ar[r]^{\rho_{A',A}}&M(A'\ox_kA)
}
$$
commute for all $A$, $A'$, where $T$ is the graded
interchange operator given in \eqref{symm}.

Once again, we have the corresponding obvious notions of a left folding
system and a bifolding system.
\end{Definition}

For a right folding system $M$, the category $\Alg^\r_M$ has a monoidal
structure given by the \emph{folding product} $\hat\ox$ below. Given an
$M$-algebra $D$ of type $A$ and another one, $D'$ of type $A'$, we define
an $M$-pair algebra $D\hat\ox D'$ of type $A\ox A'$ as the lower row in the
diagram
\begin{equation}\label{fopr}
\alignbox{
\xymatrix{
0\ar[r]
&A\!\ox\!M(A')\oplus M(A)\!\ox\!A'
\ar[r]\ar[d]_{(\lambda_{A,A'},\rho_{A,A'})}\ar@{}[dr]|{\textrm{push}}
&(D\bar\ox D')_1\ar[r]^{\d_{\bar\ox}}\ar[d]
&(D\bar\ox D')_0\ar[r]\ar@{=}[d]
&A\ox A'\ar[r]\ar@{=}[d]
&0\\
0\ar[r]
&M(A\ox A')\ar[r]
&(D\hat\ox D')_1\ar[r]_{\d_{\hat\ox}}
&D_0\ox D'_0\ar[r]
&A\ox A'\ar[r]
&0.
}
}
\end{equation}
Here the leftmost square is required to be pushout, and the upper row is
exact by \bref{tox}.

\begin{Proposition}\label{monofold}
For any right (resp. left, bi-) folding system $M$, the folding product
defines a monoidal structure on $\Alg^\r_M$ (resp. $\Alg^\l_M$,
$\Alg^\b_M$, $\Alg^\pair_M$), with unit object $I=(0:M(k)\to
k)$. If moreover the folding system is symmetric, then this monoidal
structure is symmetric.
\end{Proposition}

We only will use the monoidal categories $\Alg^\r_{\oo\oplus\Sigma}$ and
$\Alg^\pair_\oo$.

\begin{proof}
To begin with, let us show that $\hat\ox$ is functorial, i.~e. for any
morphisms $f:D\to E$, $f':D'\to E'$ in $\Alg_M$, let us define a morphism
$f\hat\ox f':D\hat\ox E\to D'\hat\ox E'$ in a way compatible with
identities and composition. We put $(f\hat\ox f')_0=f_0\hat\ox f'_0$, and
define $(f\hat\ox f')_1$ as the unique homomorphism making the following
diagram commute:
$$
\xymatrix@!C=6em{
B\!\ox\!M(B')\oplus M(B)\!\ox\!B'\ar[rrr]\ar[ddd]_{(\lambda_{B,B'},\rho_{B,B'})}
&&&(E\bar\ox E')_1\ar[ddd]\\
&A\!\ox\!M(A')\oplus M(A)\!\ox\!A'
\ar[r]\ar[d]_{(\lambda_{A,A'},\rho_{A,A'})}
\ar[ul]_{f\!\ox\!f'_*\oplus f_*\!\ox\!f'}
&(D\bar\ox D')_1\ar[d]\ar[ur]^{(f\bar\ox f')_1}\\
&M(A\ox A')\ar[r]\ar[ld]^{(f\ox f')_*}
&(D\hat\ox D')_1\ar@{-->}[dr]_{(f\hat\ox f')_1}\\
M(B\ox B')\ar[rrr]
&&&(E\hat\ox E')_1
}
$$
where the left hand trapezoid commutes by \eqref{natfold}. Using the universal
property of pushout it is clear that right equivariance of $f_1$ and $f_1'$ iplies
that of $(f\hat\ox f')_1$ so that this indeed defines a morphism in
$\Alg_M$. The same universality implies compatibility with composition.

Next to show that $I=(0:M(k)\to k)$ is a unit object first note that for an
$M$-algebra $D$ by \bref{trunc} one has
$$
I\bar\ox D=\Tr_*\left(M(k)\ox D_1\xto{\binom0{1\ox\d}}D_1\oplus
M(k)\!\ox\!D_0\xto{(\d,0)}D_0\right)\cong\left(D_1\oplus
M(k)\!\ox\!A\xto{(\d,0)}D_0\right).
$$
From this using \eqref{unitfold} it is easy to see that $(I\hat\ox D)_1$ is given by the pushout
$$
\xymatrix{
M(A)\oplus
M(k)\!\ox\!A\ar[r]^{\mathrm{incl}\oplus1}\ar[d]_{\mathrm{proj}}
&D_1\oplus M(k)\!\ox\!A\ar[d]\\
M(A)\ar[r]&(I\hat\ox D)_1
}
$$
so that there is a canonical isomorphism $(I\hat\ox D)_1\cong D_1$ compatible
with the canonical isomorphism $k\ox D_0\cong D_0$. Symmetrically, one
constructs the isomorphism $D\hat\ox I\cong D$.

Turning now to associativity, first note that the tensor product
\eqref{tens} can be equivalently stated as defining $(D\bar\ox D')_1$ by
the requirement that the diagram
$$
\xymatrix@!=1.5em{
&D_1\ox D'_1\ar[dr]\ar[dl]\ar@{}[dd]|{\textrm{push}}\\
D_0\ox D'_1\ar[dr]&&D_1\ox D'_0\ar[dl]\\
&(D\bar\ox D')_1
}
$$
be pushout. Then combining diagrams we see that $(D\hat\ox D')_1$ can be
equivalently defined as the colimit of the following diagram:
\begin{equation}\label{altox}
\alignbox{
\xymatrix@!=2.5em{
D_0\ox M(A')\ar[d]\ar[dr]|\hole
&D_1\ox D'_1\ar[dr]\ar[dl]
&M(A)\ox D'_0\ar[d]\ar[dl]|\hole\\
D_0\ox D'_1
&M(A\ox A')
&D_1\ox D'_0
}
}
\end{equation}
where the map $D_0\ox M(A')\to M(A\ox A')$ is the composite $D_0\ox
M(A')\to A\ox M(A')\to M(A\ox A')$ and similarly for $M(A)\ox D'_0\to
M(A\ox A')$. Hence $((D\hat\ox D')\hat\ox D'')_1$ is given by the colimit of the diagram
$$
\xymatrix@!=5em{
D_0\ox D'_0\ox M(A'')\ar[d]\ar[dr]|\hole
&(D\hat\ox D')_1\ox D''_1\ar[dr]\ar[dl]
&M(A\ox A')\ox D''_0\ar[d]\ar[dl]|\hole\\
D_0\ox D'_0\ox D''_1
&M(A\ox A'\ox A'')
&(D\hat\ox D')_1\ox D''_0.
}
$$
Substituting here the diagram for $(D\hat\ox D')_1$ we obtain that this is
the same as the colimit of a diagram of the form
$$
\xymatrix@!C=6em{
&&&&D_0\ox D'_1\ox D''_0\\
&&D_0\ox D'_1\ox D''_1\ar[dll]\ar[urr]
&D_0\ox M(A')\ox D''_0\ar[ur]\ar[dl]\\
D_0\ox D'_0\ox D''_1
&D_0\ox D'_1\ox M(A'')\ar[l]\ar[r]
&M(A\ox A'\ox A'')
&&D_1\ox D'_1\ox D''_0\ar[uu]\ar[dd]\\
&&D_1\ox D'_0\ox D''_1\ar[ull]\ar[drr]
&M(A)\ox D'_0\ox D''_0\ar[ul]\ar[dr]\\
&&&&D_1\ox D'_0\ox D''_0.
}
$$
Treating now $(D\hat\ox(D'\hat\ox D''))_1$ in the same way we obtain that
it is colimit of a diagram with same objects; then, using
\eqref{leftfold}, \eqref{midfold}, and \eqref{rightfold}, one can see that
also morphisms in these diagrams are the same.

Finally, suppose that $M$ is a symmetric folding system. Then for any $M$-algebras
$D$, $D'$ of type $A$, $A'$ respectively, there is a commutative diagram
$$
\xymatrix@!=.2em{
&&&&M(A\ox A')\ar[ddd]\\
\\
\\
&&&&M(A'\ox A)\\
&&D_0\ox M(A')\ar[uuuurr]\ar[ddddll]\ar[dr]
&&&&M(A)\ox D'_0\ar[uuuull]\ar[ddddrr]\ar[dl]\\
&&&M(A')\ox D_0\ar[uur]\ar[ddl]
&&D'_0\ox M(A)\ar[uul]\ar[ddr]\\
\\
&&D'_1\ox D_0
&&D'_1\ox D_1\ar[ll]\ar[rr]
&&D'_0\ox D_1\\
D_0\ox D'_1\ar[urr]
&&&&D_1\ox D'_1\ar[u]\ar[llll]\ar[rrrr]
&&&&D_1\ox D'_0\ar[ull]
}
$$
which induces a map from the colimit of the outer triangle to that of the
inner one, i.~e. by \eqref{altox} a map $(D\hat\ox D')_1\to(D'\hat\ox
D)_1$. It is then straightforward to check that this defines an interchange for
the monoidal structure.
\end{proof}

\begin{Examples}\label{foldex}
The bimodule system $\oo$ above clearly has the structure
of a folding system, with $\lambda$ and $\rho$ both identity maps. Also the
bimodule system $\oo\oplus\Sigma$ is a folding system via the obvious
isomorphisms
\begin{align}
\label{lambda+}
\lambda_{A,A'}&:A\ox(A'\oplus\Sigma
A')\cong A\!\ox\!A'\oplus A\!\ox\!\Sigma A'\xto{1\oplus\sigma}
A\!\ox\!A'\oplus\Sigma(A\!\ox\!A'),\\
\label{rho+}
\rho_{A,A'}&:(A\oplus\Sigma A)\ox A'\cong A\!\ox\!A'\oplus(\Sigma A)\!\ox\!A'
\cong A\!\ox\!A'\oplus\Sigma(A\!\ox\!A')
\end{align}
where in \eqref{lambda+}, the interchange \eqref{sigma} for $\Sigma$ is used.
\end{Examples}

\begin{Lemma}\label{sigmafold}
The isomorphisms \eqref{lambda+}, \eqref{rho+} give the bimodule system
$\oo\oplus\Sigma$ with the structure of a symmetric folding system on any
category $\bA$ of algebras closed under tensor products. 
\end{Lemma}

\begin{proof}
It is obvious that $\oo$ with the identity maps is a folding system, and
that a direct sum of folding systems is a folding system again, so it
suffices to show that $\Sigma$ is a folding system.

The right diagram in \eqref{natfold} is trivially commutative, while
commutativity of the left one follows from
\begin{multline*}
\sigma_{A_1,M(A'_1)}(f(a)\ox\Sigma f'(a'))
=(-1)^{\deg(a)}\Sigma(f(a)\ox f'(a'))\\
=\Sigma(f\ox f')((-1)^{\deg(a)}\Sigma(a\ox a'))
=\Sigma(f\ox f')\sigma_{A,M(A')}(a\ox\Sigma a')
\end{multline*}
for any $a\in A$, $a'\in A'$, $f:A\to A_1$, $f':A'\to A'_1$. Next, the
diagrams \eqref{unitfold} commute since $k$ is concentrated in degree 0.

The diagrams \eqref{rightfold} commute trivially, as only right actions are
involved. Commutativity of \eqref{leftfold} follows from the obvious
equality
$$
(-1)^{\deg(a)}\Sigma(a\ox(-1)^{\deg(a')}a'\ox a'')=(-1)^{\deg(a\ox
a')}\Sigma(a\ox a'\ox a'')
$$
and that of \eqref{midfold} is also obvious from
$$
\xymatrix@M=1em@!C=3em{
&(-1)^{\deg(a)}\Sigma(a\ox a')\ox a''\ar@{|->}[dr]\\
a\ox\Sigma(a')\ox a''\ar@{|->}[ur]\ar@{|->}[dr]
&&(-1)^{\deg(a)}\Sigma(a\ox a'\ox a'')\\
&a\ox\Sigma(a'\ox a'')\ar@{|->}[ur]
}
$$
\end{proof}

Thus by \bref{monofold} the folding system $\oo\oplus\Sigma$ yields a
well-defined monoidal category $\Alg^\r_{\oo\oplus\Sigma}$ of
\emph{$\oo\oplus\Sigma$-algebras} as in \bref{malg}. The initial object and
at the same time the unit for the monoidal structure of
$\Alg^\r_{\oo\oplus\Sigma}$ is by \bref{initial} and \bref{monofold}
$$
I_{\oo\oplus\Sigma}=\left(\FF\oplus\Sigma\FF\xto0\FF\right).
$$
For $\Alg^\r_\oo$ it is
$$
I_\oo=\left(\FF\xto0\FF\right).
$$
The projections $q:A\oplus\Sigma A\to A$ can be used to construct a
monoidal functor
\begin{equation}\label{k1}
q:\Alg^\r_{\oo\oplus\Sigma}\to\Alg^\r_\oo
\end{equation}
carrying an object $D$ in $\Alg^\r_{\oo\oplus\Sigma}$ to the pushout in
the following diagram
$$
\xymatrix{
A\oplus\Sigma A\ar@{ >->}[r]\ar[d]_q\ar@{}[dr]|{\textrm{push}}
&D_1\ar[d]\ar[r]
&D_0\ar@{=}[d]\ar@{->>}[r]
&A\ar@{=}[d]\\
A\ar@{ >->}[r]
&q(D)_1\ar[r]
&q(D)_0\ar@{->>}[r]
&A.
}
$$
Evidently $q(I_{\oo\oplus\Sigma})=I_\oo$.

\section{Unfolding systems}\label{unfold}

It is clear how to dualize the constructions from the previous section
along the lines of section \ref{coalg}. We will not give detailed definitions but
only briefly indicate the underlying structures.

We thus consider a category $\bC$ of graded $k$-coalgebras, and define a
right comodule system $N$ on $\bC$ as an assignment, to each coalgebra $C$
in $\bC$, of a $C$-comodule $N(C)$, and to each homomorphism $f:C\to C'$ of
coalgebras of an $f$-equivariant homomorphism $f_*:N(C)\to N(C')$, i.~e.
the diagram
$$
\xymatrix{
N(C)\ar[r]^-{\textrm{coaction}}\ar[d]_{f_*}&N(C)\ox C\ar[d]^{f_*\ox f}\\
N(C')\ar[r]^-{\textrm{coaction}}&N(C')\ox C'
}
$$
is required to commute. Similarly one defines left comodule systems and
bicomodule systems. As before, we have a bicomodule system $\oo$ given by
$\oo(C)=C$ and also $\Sigma$, $\oo\oplus\Sigma$ defined dually to
\bref{sysex}.

Then further for a right comodule system $N$ on $\bC$ and for a coalgebra $C$
from $\bC$ one defines an $N$-coalgebra of type $C$ by dualizing
\bref{malg}. It is thus a pair $D^*=(d:D^0\to D^1)$ where $D^0$ is a
coalgebra, $D^1$ is a right $D^0$-comodule and $d$ is a comodule
homomorphism. Moreover one must have $\pi^0(D^*)=C$, $\pi^1(D^*)=N(C)$, and
the $C$-comodule structure on $N(C)$ induced by this must be the one coming
from the comodule system $N$. With morphisms defined dually to \bref{malg},
the $N$-coalgebras form a category $\Coalg^\r_N$. Similarly one
defines categories $\Coalg^\l_N$ and
$\Coalg^\pair_N\subset\Coalg^\b_N$ for a left,
resp. bicomodule system $N$. These categories have the initial object
$0:0\to N(0)$ and the terminal object $0:k\to N(k)$.
 
Also dually to \bref{defold} one defines \emph{unfolding systems} as
comodule systems $N$ equipped with $C\ox C'$-comodule homomorphisms
\begin{align*}
l^{C,C'}&:N(C\ox C')\to C\ox N(C')\\
r^{C,C'}&:N(C\ox C')\to N(C)\ox C'
\end{align*}
for all $C,C'\in\bC$ required to satisfy obvious duals to the diagrams
\eqref{natfold} -- \eqref{midfold}. Also there is an obvious notion of a
symmetric unfolding system.

Then for an unfolding system $N$ we can dualize \eqref{fopr} to obtain
definition of the \emph{unfolding product} $D\check\ox D'$ of $N$-coalgebras
via the upper row in the diagram
$$
\xymatrix{
0\ar[r]
&C\ox C'\ar[r]\ar@{=}[d]
&D^0\ox{D'}^0\ar[r]^{d^{\check\ox}}\ar@{=}[d]
&(D\check\ox D')^1\ar[d]\ar[r]\ar@{}[dr]|{\mathrm{pull}}
&N(C\ox C')\ar[d]^{\binom{l^{C,C'}}{r^{C,C'}}}\ar[r]
&0\\
0\ar[r]
&C\ox C'\ar[r]
&(D\dblb\ox D')^0\ar[r]^{d^{\dblb\ox}}
&(D\dblb\ox D')^1\ar[r]
&C\!\ox\!N(C')\oplus N(C)\!\ox\!C'\ar[r]
&0
}
$$
where now the rightmost square is required to be pullback and the lower row
is exact by the dual of \bref{tox}.

It is then straightforward to dualize \bref{monofold}, so we conclude that
for any unfolding system $N$ the unfolding product equips the category
$\Coalg_N^?$ with the structure of a monoidal category, symmetric
if $N$ is symmetric. Here, ``?'' stands for ``r'', ``l'', ``b'' or
``pair'', according to the type of $N$. Obviously also the dual of
\bref{sigmafold} holds, so that the categories
$\Coalg_\oo^\pair$ and
$\Coalg_{\oo\oplus\Sigma}^\r$ have monoidal
structures given by the unfolding product.

\section{The $\GG$-relation pair algebra of the Steenrod
algebra}\label{grel}

Fix a prime $p$, and let $\GG=\ZZ/p^2\ZZ$ be the ring of integers mod $p^2$,
with the quotient map $\GG\onto\FF=\FF_p=\ZZ/p\ZZ$. Let $\A$ be the mod $p$
Steenrod algebra and let
$$
\E_\A=
\begin{cases}
\set{\Sq^1,\Sq^2,...}&\textrm{for $p=2$},\\
\set{\P^1,\P^2,...}\cup\set{\beta,\beta\P^1,\beta\P^2,...}&\textrm{for odd
$p$}
\end{cases}
$$
be the set of generators of the algebra $\A$. We consider the following
algebras and homomorphisms
\begin{equation}\label{BF}
\alignbox{
\xymatrix{
**[l]q:\aB_0\ar@{=}[d]\ar@{->>}[r]
&\F_0\ar@{=}[d]\ar@{->>}[r]^{q_\F}
&\A\\
T_\GG(\E_\A)
&T_\FF(\E_\A)&.
}
}
\end{equation}
For a commutative ring $k$, $T_k(S)$ denotes the free
associative $k$-algebra with unit generated by the set $S$, i.~e. the tensor
algebra of the free $k$-module on $S$. The map $q_\F$ is the algebra
homomorphism which is the identity on $E_\A$. For $f\in\F_0$ we denote the
element $q_\F(f)\in\A$ by
$$
\qf f=q_\F(f).
$$
Let $R_\aB$ denote the kernel of $q$, i.~e. there is a short exact sequence
$$
\xymatrix{
R_\aB\ar@{ >->}[r]&\aB_0\ar@{->>}[r]^q&\A.
}
$$
This short exact sequence gives rise to a long exact sequence
$$
\xymatrix{
\Tor(R_\aB,\FF)\ar@{ >->}[r]
&\Tor(\aB_0,\FF)\ar[r]
&\Tor(\A,\FF)\ar[r]^i
&R_\aB\ox\FF\ar[r]
&\aB_0\ox\FF\ar@{->>}[r]
&\A\ox\FF.
}
$$
Here $A\ox\FF\cong A/pA$ and $\Tor(A,\FF)$  is just the $p$-torsion part of $A$
for an abelian group $A$, so the connecting homomorphism $i$ sends
$a=q(b)+p\aB_0$ to $pb+pR_\aB$. It follows that the second homomorphism in
the above sequence is zero. Moreover clearly we can identify
$\aB_0\ox\FF=\F_0$ and $\Tor(\A,\FF)=\A$, so that there is an exact
sequence
\begin{equation}\label{prel}
\alignbox{
\xymatrix{
\A\ar@{ >->}[r]^i
&\R^\FF_1\ar[r]^\d\ar@{=}[d]
&\R^\FF_0\ar@{->>}[r]\ar@{=}[d]
&\A\\
&R_\aB\ox\FF
&\F_0
}
}
\end{equation}
One has
\begin{Lemma}\label{parel}
The pair $\R^\FF=(\d:\R^\FF_1\to\R^\FF_0)$ above has a pair algebra
structure compatible with the standard bimodule structure of $\A$ on
itself, so that $\R^\FF$ yields an object in $\Alg^\pair_\oo$,
see \bref{malg}.
\end{Lemma}

\begin{proof}
Clearly mod $p$ reduction of any pair algebra over $\GG$ is a pair algebra
over $\FF$. Then let $\R^\FF$ be the mod $p$ reduction of the pair algebra
$R_\aB\into\aB_0$. Thus the $\F_0$-$\F_0$-bimodule structure on
$\R^\FF_1=R_\aB/pR_\aB$ is just the mod $p$ reduction of the
$\aB_0$-$\aB_0$-bimodule structure on $R_\aB$, i.~e. $b'+p\aB_0\in\R^\FF_0=\aB_0/p\aB_0$
acts on $r+pR_\aB\in\R^\FF_1=R_\aB/pR_\aB$ via
$$
(b'+p\aB_0)(r+pR_\aB)=b'r+pR_\aB.
$$
Moreover the above inclusion $\A\into R_\aB/pR_\aB$ sends an element $q(b)$
to $pb+pR_\aB$. Then the action of $a'=q(b')\in\A$ on $i(a)=pb+pR_\aB\in
i(\A)=\ker\d$ induced by this pair algebra is given as follows:
$$
a'i(a)=q_\F(b'+p\aB_0)(pb+pR_\aB)=pb'b+pR_\aB=iq(b'b)=i(a'a)
$$
and similarly for the right action.
\end{proof}

We call the object $\R^\FF$ of the category $\Alg^\pair_\oo$ 
the \emph{$\GG$-relation pair algebra of $\A$}.

\begin{Theorem}\label{relcom}
The $\oo$-pair algebra $\R^\FF$ has a structure of a cocommutative comonoid
in the symmetric monoidal category $\Alg^\pair_\oo$.
\end{Theorem}

\begin{proof}
For $n\ge0$, let $R_\aB^{(n)}$ denote the kernel of the map $q^{\ox n}$, so
that there is a short exact sequence
$$
\xymatrix{
R_\aB^{(n)}\ar@{ >->}[r]&\aB_0^{\ox n}\ar@{->>}[r]^{q^{\ox n}}&\A^{\ox n}
}
$$
and similarly to \bref{parel} there is a pair algebra of the form
$$
\xymatrix{
\A^{\ox n}\ar@{ >->}[r]
&R_\aB^{(n)}\ox\FF\ar[r]
&\F_0^{\ox n}\ar@{->>}[r]^{q_\F^{\ox n}}
&\A^{\ox n}
}
$$
determining an object $\R^{(n)}$ in $\Alg^\pair_\oo$. Then one
has the following lemma which yields natural examples of folding products
in $\Alg^\pair_\oo$.
\begin{Lemma}\label{rn}
There is a canonical isomorphism $\R^{(n)}\cong(\R^\FF)^{\hat\ox n}$ in
$\Alg_\oo^\pair$.
\end{Lemma}

\begin{proof}
Using induction, we will assume given an isomorphism
$\alpha_n:(\R^\FF)^{\hat\ox n}\cong\R^{(n)}$ and construct $\alpha_{n+1}$
in a canonical way. To do this it clearly suffices to construct a
canonical isomorphism $\R^\FF\hat\ox\R^{(n)}\cong\R^{(n+1)}$ as then its
composite with $\R^\FF\hat\ox\alpha_n$ will give $\alpha_{n+1}$.

To construct a map $(\R^\FF\hat\ox\R^{(n)})_1\to\R^{(n+1)}_1$ means by
\eqref{altox} the same as to find three dashed arrows making the diagram
$$
\xymatrix@!=1em{
&&R_\aB\ox R_\aB^{(n)}\ox\FF\ar[ddl]\ar[ddr]\\
\\
&\F_0\ox R_\aB^{(n)}\ar@{-->}[dr]&&R_\aB\ox\F_0^{\ox n}\ar@{-->}[dl]\\
&&R_\aB^{(n+1)}\ox\FF\\
\F_0\ox\A^{\ox n}\ar[rr]\ar[uur]&&\A^{\ox(n+1)}\ar@{-->}[u]&&\A\ox\F_0^{\ox
n}\ar[ll]\ar[uul]
}
$$
commute. For this we use the commutative diagram
$$
\xymatrix@!=1em{
&&R_\aB\ox R_\aB^{(n)}\ar[ddl]\ar[ddr]\\
\\
&\aB_0\ox R_\aB^{(n)}\ar[dr]&&R_\aB\ox\aB_0^{\ox n}\ar[dl]\\
&&R_\aB^{(n+1)}\\
\aB_0\ox\A^{\ox n}\ar[rr]\ar[uur]&&\A^{\ox(n+1)}\ar[u]&&\A\ox\aB_0^{\ox
n};\ar[ll]\ar[uul]
}
$$
This diagram has a commutative subdiagram
$$
\xymatrix@!=1em{
&&p(R_\aB\ox R_\aB^{(n)})\ar[ddl]\ar[ddr]\\
\\
&p\aB_0\ox R_\aB^{(n)}\ar[dr]&&R_\aB\ox p\aB_0^{\ox n}\ar[dl]\\
&&pR_\aB^{(n+1)}\\
p\aB_0\ox\A^{\ox n}\ar[rr]\ar[uur]&&0\ar[u]&&\A\ox p\aB_0^{\ox
n};\ar[ll]\ar[uul]
}
$$
It is obvious that taking the quotient by this subdiagram gives us a diagram of the
kind we need.

We thus obtain a map $(\R^\FF\hat\ox\R^{(n)})_1\to R_\aB^{(n+1)}\ox\FF$.
Moreover by its construction this map fits into the commutative diagram
$$
\xymatrix{
\A^{\ox(n+1)}\ar@{=}[d]\ar@{ >->}[r]
&(\R^\FF\hat\ox\R^{(n)})_1\ar[d]\ar[r]
&\F_0^{\ox(n+1)}\ar@{=}[d]\ar@{->>}[r]
&\A^{\ox(n+1)}\ar@{=}[d]\\
\A^{\ox(n+1)}\ar@{ >->}[r]
&R_\aB^{(n+1)}\ox\FF\ar[r]
&\F_0^{\ox(n+1)}\ar@{->>}[r]
&\A^{\ox(n+1)}
}
$$
with exact rows, hence by the five lemma it is an isomorphism.
\end{proof}

Using the lemma, we next construct the diagonal of $\R^\FF$ given by
$$
\xymatrix{
R_\aB\ox\FF\ar@{=}[r]\ar[d]_{\Delta^\GG\ox1}
&\R^\FF_1\ar[r]^\d\ar[d]_\Delta
&\R^\FF_0\ar@{=}[r]\ar[d]_\Delta
&\F_0\ar[d]^\Delta\\
R_\aB^{(2)}\ox\FF\ar[r]^-\cong
&(\R^\FF\hat\ox\R^\FF)_1\ar[r]^{\d_{\hat\ox}}
&\R^\FF_0\ox\R^\FF_0\ar@{=}[r]
&\F_0\ox\F_0.
}
$$
Here $\Delta^\GG$ is defined by the commutative diagram
\begin{equation}\label{deltag}
\alignbox{
\xymatrix{
R_\aB\ar[d]_{\Delta^\GG}\ar@{ (->}[r]
&\aB_0\ar[d]_{\Delta^\GG}\\
R_\aB^{(2)}\ar@{ (->}[r]
&\aB_0\ox\aB_0,
}
}
\end{equation}
where the diagonal $\Delta^\GG$ on $\aB_0$ is defined on generators by
$$
\begin{aligned}
\Delta^\GG(\Sq^n)=
\sum_{i=0}^n\Sq^i\ox\Sq^{n-i}\hskip4.5em &\textrm{for $p=2$,}\\
\left.
\begin{aligned}
\Delta^\GG(\beta)&=\beta\ox1+1\ox\beta,\\
\Delta^\GG(\P^n)&=\sum_{i+j=n}\P^i\ox\P^j,\\
\Delta^\GG(\P^n_\beta)&=\sum_{i+j=n}(\P^i_\beta\ox\P^j+\P^i\ox\P^j_\beta)
\end{aligned}
\right\}&\textrm{for odd $p$}
\end{aligned}
$$
(with $\Sq^0=1$, $\P^0=1$ as usual) and extended to the whole $\aB_0$ as
the unique algebra homomorphism with respect to the algebra structure on
$\aB_0\ox\aB_0$ given by the nonstandard interchange formula
$$
\xymatrix{
&\aB_0\ox\aB_0\ox\aB_0\ox\aB_0\ar[dr]^{\mu\ox\mu}\\
\aB_0\ox\aB_0\ox\aB_0\ox\aB_0\ar[ur]^{1\ox T^\GG\ox1}\ar[rr]^{\mu_\ox}
&&\aB_0\ox\aB_0
}
$$
with
\begin{align*}
&T^\GG:\aB_0\ox\aB_0\xto\cong\aB_0\ox\aB_0\\
&T^\GG(x\ox y)=(-1)^{p\deg(x)\deg(y)}y\ox x.
\end{align*}
In particular, clearly for all $p$ one has $T^\GG\Delta^\GG=\Delta^\GG$, i.~e. the
coalgebra structure on $\aB_0$ is cocommutative.

The counit for $\R^\FF$ is given by the diagram
\begin{equation}\label{rcounit}
\alignbox{
\xymatrix{
\A\ar@{ >->}[r]\ar[d]^\epsilon
&R_\aB\ox\FF\ar[r]\ar@{-->}[d]
&\F_0\ar@{->>}[r]\ar[d]^\epsilon
&\A\ar[d]^\epsilon\\
\FF\ar@{=}[r]
&\FF\ar[r]^0
&\FF\ar@{=}[r]
&\FF
}
}
\end{equation}
where the map $R_\aB\ox\FF\to\FF$ sends the generator $(p1_{\aB_0})\ox1$ in degree
$0$ to 1 and all elements in higher degrees to zero. It is then clear from
the formula for $\Delta^\GG$ that this indeed gives a counit for this
diagonal.

Finally, to prove coassociativity, by the lemma it suffices to consider the
diagram
$$
\xymatrix@!=2em{
&&R_\aB\ar@{_(->}[d]\ar[ddll]_{\Delta^\GG}\ar[ddrr]^{\Delta^\GG}\\
&&\aB_0\ar[dl]_{\Delta^\GG}\ar[dr]^{\Delta^\GG}\\
R^{(2)}_\aB\ar@{ (->}[r]\ar[ddrr]&\aB_0^{\ox2}\ar[dr]_{1\ox\Delta^\GG}
&&\aB_0^{\ox2}\ar[dl]^{\Delta^\GG\ox1}&R^{(2)}_\aB\ar[ddll]\ar@{_(->}[l]\\
&&\aB_0^{\ox3}\\
&&R^{(3)}_\aB\ar@{ (->}[u]
}
$$
\end{proof}

\section{The algebra of secondary cohomology operations}\label{sec}

Let us next consider a derivation of degree 0 of the form
$$
\k:\A\to\Sigma\A,
$$
uniquely determined by
\begin{equation}\label{kappa}
\begin{aligned}
\k\Sq^n=\Sigma\Sq^{n-1}\hskip.7em &\textrm{for $p=2$,}\\
\left.
\begin{aligned}
\k\beta&=\Sigma1,\\
\k(\P^i)&=0,i\ge0
\end{aligned}
\right\}&\textrm{for odd $p$.}
\end{aligned}
\end{equation}
We will use $\k$ to define an $\A$-$\A$-bimodule
$$
\A\oplus_\k\Sigma\A
$$
as follows. The right $\A$-module structure is the same as on
$\A\oplus\Sigma\A$ above, i.~e. one has $(x,\Sigma y)a=(xa,\Sigma ya)$. As
for the left $\A$-module structure, it is given by
$$
a(x,\Sigma y)=(ax,(-1)^{\deg(a)}\Sigma ay+\k(a)x).
$$
There is a short exact sequence of $\A$-$\A$-bimodules
$$
0\to\Sigma\A\to\A\oplus_\k\Sigma\A\to\A\to0
$$
given by the standard inclusion and projection.

\begin{Remark}\label{hoch}
The above short exact sequence of bimodules and the derivation $\k$
correspond to each other under the well known description of the first
Hochschild cohomology group in terms of bimodule extensions and derivations,
respectively. Indeed, more generally recall that for a graded $k$-algebra
$A$ and an $A$-$A$-bimodule $M$, one of the possible definitions of the
Hochschild cohomology of $A$ with coefficients in $M$ is
$$
HH^n(A;M)=\Ext^n_{A\ox_kA^\circ}(A,M).
$$
On the other hand, $HH^1(A;M)$ can be also described in terms of
derivations. Recall that an $M$-valued derivation on $A$ is a $k$-linear
map $\k:A\to M$ of degree 0 satisfying
$$
\k(xy)=\k(x)y+(-1)^{\deg(x)}x\k(y)
$$
for any $x,y\in A$. Such derivations form a $k$-vector space $\Der(A;M)$.
A derivation $\k=\i_m$ is called inner if there is an $m\in M$ such that
$$
\k(x)=mx-(-1)^{\deg(x)}xm=\i_m(x)
$$
for all $x\in A$. These form a subspace $\Ider(A;M)\subset\Der(A;M)$ and
one has an isomorphism $HH^1(A;M)\cong\Der(A;M)/\Ider(A;M)$. Moreover
there is an exact sequence
$$
0\to HH^0(A;M)\to M\xto{\i_{\_}}\Der(A;M)\to HH^1(A;M)\to0.
$$

Explicitly, the isomorphism
$$
\Der(A;M)/\Ider(A;M)\cong\Ext^1_{A\ox A^\circ}(A,M),
$$
can be described by assigning to a class of a derivation $\k:A\to M$ the class of the
extension
$$
0\to M\to A\oplus_\k M\to A\to0
$$
where as a vector space, $A\oplus_\k M=A\oplus M$, the maps are the canonical
inclusion and projection and the bimodule structure is given by
\begin{align*}
a(x,m)&=(ax,am+\k(a)x),\\
(x,m)a&=(xa,ma).
\end{align*}

Obviously $\A\oplus_\k\Sigma\A$ above is an example of this
construction.
\end{Remark}

\begin{Definition}\label{hpa}
A \emph{Hopf pair algebra} $\V$ (associated to $\A$) is a pair algebra
$\d:\V_1\to\V_0$ over $\FF$ together with the following commutative diagram in the
category of $\F_0$-$\F_0$-bimodules
\begin{equation}\label{hpad}
\alignbox{
\xymatrix{
\Sigma\A\ar@{=}[r]\ar@{ >->}[d]
&\Sigma\A\ar@{ >->}[d]\\
\A\oplus_\k\Sigma\A\ar@{ >->}[r]\ar@{->>}[d]_q
&\V_1\ar[r]^\d\ar@{->>}[d]_q
&\V_0\ar@{=}[d]\ar@{->>}[r]
&\A\ar@{=}[d]\\
\A\ar@{ >->}[r]
&\R^\FF_1\ar[r]
&\R^\FF_0\ar@{->>}[r]
&\A
}
}
\end{equation}
with exact rows and columns. The pair morphism $q:\V\to\R^\FF$ will be
called the \emph{$\GG$-structure} of $\V$. Moreover $\V$ has a structure
of a comonoid in $\Alg^\r_{\oo\oplus\Sigma}$ and $q$ is compatible with
the $\Alg^\pair_\oo$-comonoid structure on $\R^\FF$ in
\bref{relcom}, in  the sense that the diagrams
\begin{equation}\label{diacomp}
\alignbox{
\xymatrix{
\V_1\ar[r]^-{\Delta_\V}\ar[d]_q
&(\V\hat\ox\V)_1\ar[d]^{q\hat\ox q}\\
\R^\FF_1\ar[r]^-{\Delta_\R}
&(\R^\FF\hat\ox\R^\FF)_1
}
}
\end{equation}
and
\begin{equation}\label{coucomp}
\alignbox{
\xymatrix{
\V_1\ar[d]_q\ar[r]^-{\epsilon_\V}
&\FF\oplus\Sigma\FF\ar[d]\\
\R^\FF_1\ar[r]^-{\epsilon_\R}
&\FF
}
}
\end{equation}
commute.
\end{Definition}

We next observe that the following diagrams commute:
$$
\xymatrix{
\A\ar[rr]^\k\ar[d]_\delta
&&\Sigma\A\ar[d]^{\Sigma\delta}\\
\A\ox\A\ar[r]^{\k\ox1}
&\Sigma\A\ox\A\ar@{=}[r]
&\Sigma(\A\!\ox\!\A),
}\ \ \ \ \ \ \ 
\xymatrix{
\A\ar[rr]^\k\ar[d]_\delta
&&\Sigma\A\ar[d]^{\Sigma\delta}\\
\A\ox\A\ar[r]^{1\ox\k}
&\A\ox\Sigma\A\ar[r]^\sigma
&\Sigma(\A\!\ox\!\A)
}
$$
where $\sigma$ is the interchange for $\Sigma$ in \eqref{sigma}. Or, on
elements,
\begin{equation}\label{kapid}
\sum\k(a_\l)\ox a_\r=\sum\k(a)_\l\ox\k(a)_\r=\sum\sigma(a_\l\ox\k(a_\r)),
\end{equation}
where we use the Sweedler notation for the diagonal
$$
\delta(x)=\sum x_\l\ox x_\r.
$$

\begin{Remark}
The above identities have a simple explanation using dualization. We will
see in \bref{dkappa} below that the map dual to $\k$ is the map
$\Sigma\A_*\to\A_*$ given, for $p=2$, by multiplication with the degree 1 generator
$\zeta_1\in\A_*$ and for odd $p$ by the degree 1 generator $\tau_0$.
Then the duals of \eqref{kapid} are the obvious identities for any $x,y\in\A_*$
$$
(\zeta_1x)y=\zeta_1(xy)=x(\zeta_1y)
$$
for $p=2$ and
$$
(\tau_0x)y=\tau_0(xy)=(-1)^{\deg(x)}x(\tau_0y)
$$
for odd $p$ (recall that $\A_*$ is graded commutative).
\end{Remark}

Using \eqref{kapid} we prove:

\begin{Lemma}\label{vleft}
For a Hopf pair algebra $\V$ there is a unique left action of $\F_0$ on
$(\V\hat\ox\V)_1$ such that the quotient map
$$
(\V\bar\ox\V)_1\onto(\V\hat\ox\V)_1
$$
is $\F_0$-equivariant. Here we use the pair algebra structure on
$\V\bar\ox\V$ to equip $(\V\bar\ox\V)_1$ with an $\F_0\ox\F_0$-bimodule
structure and then turn it into a left $\F_0$-module via restriction of
scalars along $\Delta:\F_0\to\F_0\ox\F_0$. 
\end{Lemma}

\begin{proof}
Uniqueness is clear as the module structure on the quotient of any module
$M$ by a submodule is clearly uniquely determined by the module structure
on $M$.

For the existence, consider the diagram 
\begin{equation}\label{vcolim}
\alignbox{
\xymatrix{
\F_0\ox(\A\oplus_\k\Sigma\A)\ar[d]\ar[dr]|\hole
&\V_1\ox\V_1\ar[dr]\ar[dl]
&(\A\oplus_\k\Sigma\A)\ox\F_0\ar[d]\ar[dl]|\hole\\
\F_0\ox\V_1
&\A\!\ox\!\A\oplus_{\k\ox\!1}\Sigma(\A\!\ox\!\A)
&\V_1\ox\F_0.
}
}
\end{equation}
whose colimit, by \eqref{altox}, is $(\V\hat\ox\V)_1$, with the right
$\F_0\ox\F_0$-module structure coming from the category
$\Alg^\r_{\oo\oplus\Sigma}$. It then suffices to show that all maps in
this diagram are also left $\F_0$-equivariant, if one uses the left
$\F_0$-module structure by restricting scalars along the diagonal
$\F_0\to\F_0\ox\F_0$.

This is trivial except possibly for two of the maps involved. For the map
$$
\Phi:\F_0\ox(\A\oplus_\k\Sigma\A)\to\A\!\ox\!\A\oplus_{\k\!\ox\!\A}\Sigma(\A\!\ox\!\A)
$$
given by
$$
\Phi(f'\ox(x,\Sigma y))=(\qf{f'}\ox x,(-1)^{\deg(f')}\Sigma\qf{f'}\ox y),
$$
this amounts to checking that for any $f,f'\in\F_0$ and $x,y\in\A$ one must have
$$
\sum(f_\l\ox f_\r)(\qf{f'}\ox x,(-1)^{\deg(f')}\Sigma\qf{f'}\ox y)
=\Phi((-1)^{\deg(f_\r)\deg(f')}\sum f_\l f'\ox(\qf{f_\r}x,(-1)^{\deg(f_\r)}\Sigma
\qf{f_\r}y+\k(\qf{f_\r})x)),
$$
where again the above Sweedler notation
$$
\Delta(f)=\sum f_\l\ox f_\r,
$$
is used for the diagonal of $\F_0$ too, and $\qf{f'}$ denotes $q_\F(f')$ by
the notation in \eqref{BF}.

The left hand side expression then expands as
\begin{multline*}
\sum((-1)^{\deg(f_\r)\deg(f')}\qf{f_\l}\qf{f'}\ox\qf{f_\r}x,\\
(-1)^{\deg(f_\r)\deg(f')}(-1)^{\deg(f)}(-1)^{\deg(f')}\Sigma\qf{f_\l}
\qf{f'}\ox\qf{f_\r}y
+(-1)^{\deg(f_\r)\deg(f')}\k(\qf{f_\l})\qf{f'}\ox\qf{f_\r}x)
\end{multline*}
and the right hand side expands as
$$
(-1)^{\deg(f_\r)\deg(f')}
\sum(\qf{f_\l}\qf{f'}\ox \qf{f_\r}x,
(-1)^{\deg(f_\l f')}((-1)^{\deg(f_\r)}\Sigma\qf{f_\l}\qf{f'}\ox\qf{f_\r}y
+\qf{f_\l}\qf{f'}\ox\k(\qf{f_\r})x)).
$$
Thus left equivariance of $\Phi$ is equivalent to the equality
$$
\sum\k(\qf{f_\l})\qf{f'}\ox\qf{f_\r}x=\sum(-1)^{\deg(f_\l f')}\qf{f_\l}
\qf{f'}\ox\k(\qf{f_\r})x.
$$
This is easily deduced from
$$
\sum\k(\qf{f_\l})\ox\qf{f_\r}=\sum(-1)^{\deg(f_\l)}\qf{f_\l}\ox\k(\qf{f_\r}),
$$
which is an instance of \eqref{kapid}.

For another map
$$
\Psi:(\A\oplus_\k\Sigma\A)\ox\F_0\to\A\!\ox\!\A\oplus_{\k\!\ox\!\A}\Sigma(\A\!\ox\!\A)
$$
given by
$$
\Psi((x,\Sigma y)\ox f')=(x\ox\qf{f'},\Sigma y\ox\qf{f'})
$$
the equality to check is
$$
\sum(f_\l\ox f_\r)(x\ox\qf{f'},\Sigma y\ox\qf{f'})\\
=\Psi((-1)^{\deg(f_\r)\deg(x,\Sigma y)}\sum(\qf{f_\l}x,(-1)^{\deg(f_\l)}\Sigma
\qf{f_\l}y+\k(\qf{f_\l})x)\ox f_\r f').
$$
Here the left hand side expands as
$$
\sum((-1)^{\deg(f_\r)\deg(x)}\qf{f_\l}x\ox\qf{f_\r}\qf{f'},\\
(-1)^{\deg(f_\r)\deg(\Sigma y)}(-1)^{\deg(f_\l)}\qf{f_\l}y\ox\qf{f_\r}\qf{f'}
+(-1)^{\deg(f_\r)\deg(x)}\k(\qf{f_\l})x\ox\qf{f_\r}\qf{f'})
$$
and the right hand side expands as
$$
(-1)^{\deg(f_\r)\deg(x,\Sigma y)}\sum(\qf{f_\l}x\ox\qf{f_\r}
\qf{f'},(-1)^{\deg(f_\l)}\Sigma\qf{f_\l}y\ox
\qf{f_\r}\qf{f'}+\k(\qf{f_\l})x\ox\qf{f_\r}\qf{f'});
$$
these two expressions are visibly the same.
\end{proof}

Given this left module structure on $(\V\hat\ox\V)_1$, one can measure the
deviation from left equivariance of the diagonal
$\Delta_\V:\V_1\to(\V\hat\ox\V)_1$. For that, consider the map $\hat
L:\V_0\ox\V_1\to(\V\hat\ox\V)_1$ given by
$$
\hat L(f\ox x):=\Delta_\V(fx)-f\cdot\Delta_\V(x),
$$
for any $f\in\F_0=\V_0$, $x\in\V_1$, where $\cdot$ denotes the left
$\F_0$-module action defined in \bref{vleft}. Since the diagonal
$\Delta_\R$ of $\R^\FF$ is left equivariant, it follows from
\eqref{diacomp} that the image of $\hat L$ lies in the kernel of the map
$q\hat\ox q$, i.~e. in $\Sigma\A\ox\A$. Moreover if $f=\d v_1$ for some
$v_1\in\V_1$, then one has
$$
\Delta_\V(\d(v_1)x)=\Delta_\V(v_1\d x)=\Delta_\V(v_1)\Delta_\F(\d
x)=\Delta_\V(v_1)\d_{\hat\ox}\Delta_\V(x)=\d_{\hat\ox}\Delta_\V(v_1)\Delta_\V(x)
=\Delta_\F(\d v_1)\Delta_\V(x),
$$
so that the image of $\d\ox\V_1:\V_1\ox\V_1\to\V_0\ox\V_1$ lies in the
kernel of $\hat L$. Similarly commutativity of
\begin{equation}\label{diacomm}
\alignbox{
\xymatrix{
\V_1\ar[r]^-{\Delta_\V}\ar[d]_\d
&(\V\hat\ox\V)_1\ar[d]^{\d_{\hat\ox}}\\
\V_0\ar[r]^-{\Delta_\F}
&\V_0\ox\V_0
}
}
\end{equation}
implies that $\V_0\ox\ker\d$ is in the kernel of $\hat L$. It then follows that
$L$ factors uniquely through a map
$$
\A\ox R_\F=(\V_0/\im\d)\ox(\V_1/\ker\d)\to\ker(q\hat\ox q)=\Sigma\A\ox\A.
$$

\begin{Definition}\label{lv}
The map
$$
L_\V:\A\ox R_\F\to\Sigma\A\ox\A
$$
given by the unique factorization of the map $\hat L$ above is characterized by
the deviation of the diagonal $\Delta_\V$ of the Hopf pair algebra $\V$
from left equivariance. That is, one has
$$
\Delta_\V(fx)=f\cdot\Delta_\V(x)+L_\V(\qf f\ox\d x)
$$
for any $f\in\F_0=\V_0$, $x\in\V_1$ and the action $\cdot$ from \bref{vleft}.
\end{Definition}

Similarly one can measure the deviation of $\Delta_\V:\V_1\to(\V\hat\ox\V)_1$
from cocommutativity by means of the map $\hat S:\V_1\to(\V\hat\ox\V)_1$ given
by
$$
\hat S(x):=\Delta_\V(x)-T\Delta_\V(x),
$$
where $T:(\V\hat\ox\V)_1\to(\V\hat\ox\V)_1$ is the interchange operator for
$\Alg^\r_{\oo\oplus\Sigma}$ as constructed in \bref{monofold}. Then
similarly to $\hat L$ above, $\hat S$ admits a factorization in the following way.
First, by commutativity of \eqref{diacomp} one has
$$
(q\hat\ox q)T\Delta_\V=T(q\hat\ox q)\Delta_\V=T\Delta_\R q=\Delta_\R
q=(q\hat\ox q)\Delta_\V,
$$
since the $\Alg^\pair_\oo$-comonoid $\R^\FF$ is cocommutative.
Thus the image of $\hat S$ is contained in $\ker(q\hat\ox q)=\Sigma\A\ox\A$.
Next, commutativity of \eqref{diacomm} implies that $\ker\d$ is contained
in the kernel of $\hat S$. Hence $\hat S$ factors uniquely as follows
$$
R_\F=\V_1/\ker\d\to\ker(q\hat\ox q)=\Sigma\A\ox\A.
$$

\begin{Definition}\label{S}
The map
$$
S_\V:R_\F\to\Sigma\A\ox\A
$$
given by the unique factorization of the map $\hat S$ above is characterized by
the deviation of the diagonal $\Delta_\V$ of the Hopf pair algebra $\V$
from cocommutativity. That is, one has
$$
T\Delta_\V(x)=\Delta_\V(x)+S_\V(\d x)
$$
for any $x\in\V_1$.
\end{Definition}

It is clear from these definitions that $L_\V$ and $S_\V$ are well defined
maps by the Hopf pair algebra $\V$. Below in \bref{lao} we define the left
action operator $L:\A\ox R_\F\to\Sigma\A\ox\A$ and the symmetry operator
$S:R_\F\to\Sigma\A\ox\A$ with $L=0$ and $S=0$ if $p$ is odd. For $p=2$
these operators are quite intricate but explicitly given. We also will
study the dualization of $S$ and $L$. 

The next two results are essentially reformulations of the main results in
the book \cite{Baues}.

\begin{Theorem}[Existence]\label{exist}
There exists a Hopf pair algebra $\V$ with $L_\V=L$ and $S_\V=S$.
\end{Theorem}

\begin{Theorem}[Uniqueness]\label{unique}
The Hopf pair algebra $\V$ satisfying $L_\V=L$ and $S_\V=S$ is unique up
to an isomorphism over the $\GG$-structure $\V\to\R^\FF$ and under the
kernel $\A\oplus_\k\Sigma\A\into\V$.
\end{Theorem}

The Hopf pair algebra appearing in these theorems is the \emph{algebra of
secondary cohomology operations} over $\FF$, denoted by
$\aB^\FF=(\aB_1^\FF\to\aB_0^\FF)=\aB\ox\FF$. The algebra $\aB$ has been defined
over $\GG$ in \cite{Baues}.

\begin{proof}[Proof of \bref{exist}]
Recall that in \cite{Baues}*{12.1.8} a folding product $\hat\ox$ is defined
for pair $\GG$-algebras in such a way that $\aB$ has a comonoid structure with
respect to it, i.~e. a \emph{secondary Hopf algebra} structure. Let
$$
\Delta_1:\aB_1\to(\aB\hat\ox\aB)_1
$$
be the corresponding \emph{secondary diagonal} from
\cite{Baues}*{(12.2.2)}. It is proved in \cite{Baues}*{14.4} that the left
action operator $L$ satisfies
$$
\Delta_1(bx)=b\Delta_1(x)+L(q(b)\ox(\d x\ox1))
$$
for $b\in\aB_0$, $x\in\aB_1$, $\d x\ox1\in R_\aB\ox\FF=\R^\FF_1$.
Also in \cite{Baues}*{14.5} it is proved that the symmetry operator $S$
satisfies
$$
T\Delta_1(x)=\Delta_1(x)+S(\d x\ox1)
$$
for $x\in\aB_1$. Moreover it is proved in \cite{Baues}*{15.3.13} that the
secondary Hopf algebra $\aB$ is determined uniquely up to isomorphism by the
maps $\k$, $L$ and $S$.

Consider now the diagram
$$
\xymatrix{
\Sigma\A\ar@{=}[r]\ar@{ >->}[d]
&\Sigma\A\ar@{ >->}[d]\\
\A\oplus_\k\Sigma\A\ar@{ >->}[r]^{i_\k}\ar@{->>}[d]_q
&\aB_1\ox\FF\ar[r]^{\d\ox1}\ar@{->>}[d]_{q=\d\ox1}
&\aB_0\ox\FF\ar@{=}[d]\ar@{->>}[r]
&\A\ar@{=}[d]\\
\A\ar@{ >->}[r]
&R_\aB\ox\FF\ar[r]
&\F_0\ar@{->>}[r]
&\A.
}
$$
Here the inclusion $i_\k:\A\oplus_\k\Sigma\A\into\aB_1\ox\FF$ is given by
the inclusion $\Sigma\A\subset\aB_1$ and by the map
$$
\A\to\aB_1\ox\FF
$$
which assigns to an element $q(b)\in\A$, for $b\in\aB_0$, the element
$[p]\cdot b\ox1$. Then it is clear that $i_\k$ is a right $\A$-module
homomorphism. Moreover it is also a left $\A$-module homomorphism since for
$b\in\aB_0$ the following identity holds in $\aB_1$:
$$
b\cdot[p]-[p]\cdot b=\k(b).
$$
Compare \cite{Baues}*{A20 in the introduction}. Now one can check that the
properties of $\aB$ established in \cite{Baues} yield the result.
\end{proof}

\begin{Remark}\label{tmp}
For elements $\alpha,\beta,\gamma\in\A$ with $\alpha\beta=0$ and
$\beta\gamma=0$ the \emph{triple Massey product}
$$
\brk{\alpha,\beta,\gamma}\in\A/(\alpha\A+\A\gamma)
$$
is defined. Here the degree of elements in $\brk{\alpha,\beta,\gamma}$ is
$\deg(\alpha)+\deg(\beta)+\deg(\gamma)-1$. We can compute
$\brk{\alpha,\beta,\gamma}$ by use of the Hopf pair algebra $\aB^\FF$ above
as follows. For this we consider the maps
$$
\xymatrix{
\A&\aB_0\supset R_\aB\ar@{->>}[l]_{q_\aB}\ar@{->>}[r]^{q_R}&R_\aB\ox\FF.
}
$$
We choose elements $\bar\alpha,\bar\beta,\bar\gamma\in\aB_0$ which $q_\aB$
carries to $\alpha,\beta,\gamma$ respectively. Then we know that the
products $\bar\alpha\bar\beta$, $\bar\beta\bar\gamma$ are elements in
$R_\aB$ for which we can choose elements $x,y\in\aB_1\ox\FF$ with
\begin{align*}
q(x)&=q_R(\bar\alpha\bar\beta),\\
q(y)&=q_R(\bar\beta\bar\gamma).
\end{align*}
Then the bimodule structure of $\aB_1\ox\FF$ yields the element $\bar\alpha
y-x\bar\gamma$ in the kernel $\Sigma\A$ of $q:\aB_1\ox\FF\to R_\aB\ox\FF$.
Now $\bar\alpha y-x\bar\gamma\in\Sigma\A$ represents
$\brk{\alpha,\beta,\gamma}$, see \cite{Baues}.
\end{Remark}

\section{The dual of the $\GG$-relation pair algebra}

We next turn to the dualization of the $\GG$-relation pair algebra of the
Steenrod algebra from section \ref{grel}.

For this we just apply the duality functor $D$ to \eqref{prel}. There results an exact
sequence
$$
\xymatrix{
\A_*\ar@{ >->}[r]
&\R^0_\FF\ar[r]^d
&\R^1_\FF\ar@{->>}[r]
&\A_*,
}
$$
i.~e. the sequence
\begin{equation}\label{drel}
\alignbox{
\xymatrix{
\A_*\ar@{ >->}[r]
&D(\R^\FF_0)\ar@{=}[d]\ar[r]^{D(\d)}&D(\R^\FF_1)\ar@{=}[d]\ar@{->>}[r]
&\A_*\\
&\Hom(\F_0,\FF)&\Hom(R_\aB,\FF).
}
}
\end{equation}
In particular, by the dual of \bref{parel} one has
\begin{Lemma}\label{dparel}
The pair $\R_\FF=(d:\R^0_\FF\to\R^1_\FF)$ has a pair coalgebra structure
compatible with the standard bicomodule structure of $\A_*$ over itself, so
that $\R_\FF$ yields an object in $\Coalg^\pair_\oo$, see section \ref{unfold}.
\end{Lemma}

Moreover the dual of \bref{relcom} takes place, i.~e. one has
\begin{Theorem}\label{drelcom}
The pair coalgebra $\R_\FF$ has a structure of a commutative monoid in the
category $\Coalg_\oo^\pair$ with respect to the
unfolding product $\check\ox$.
\end{Theorem}\qed

The proof uses the duals of the pair algebras $\R^{(n)}$, $n\ge0$, from
\bref{relcom}. Namely, applying to the short exact sequence
$$
\xymatrix{
R_\aB^{(n)}\ar@{ >->}[r]&\aB_0^{\ox n}\ar@{->>}[r]^{q^{\ox n}}&\A^{\ox n}
}
$$
the functor $D=\Hom(\_,\FF)$ gives, similarly to \bref{dparel}, a pair
coalgebra
$$
\R_*^{(n)}=\left(
\xymatrix{
\A_*^{\ox n}\ar@{ >->}[r]&\F_*^{\ox
n}\ar[r]&{R^{(n)}_\aB}_*\ar@{->>}[r]&\A_*^{\ox n}
}
\right)
$$
such that the following dual of \bref{rn} holds:
\begin{Lemma}
There is a canonical isomorphism $\R_*^{(n)}\cong(\R_\FF)^{\check\ox n}$ in
$\Coalg^\pair_\oo$.
\end{Lemma}\qed

Using this lemma one constructs the $\check\ox$-monoid structure on
$\R_\FF$ by the diagram
$$
\xymatrix{
\F_*\ox\F_*\ar@{=}[r]\ar[d]^{\Delta_*}
&\R^0_\FF\ox\R^0_\FF\ar[r]^{d_{\check\ox}}\ar[d]^\mu
&(\R_\FF\check\ox\R_\FF)^1\ar[r]^-\cong\ar[d]^\mu
&{R_\aB^{(2)}}_*\ar[d]^{\Delta^\GG_*}\\
\F_*\ar@{=}[r]
&\R^0_\FF\ar[r]^d
&\R_\FF^1\ar@{=}[r]
&{R_\aB}_*
}
$$
with $\Delta^\GG$ as in \bref{deltag}.

Moreover the unit of $\R_\FF$ is given by the dual of \bref{rcounit}, i.~e.
by the diagram
$$
\xymatrix{
\FF\ar@{=}[r]\ar[d]^1
&\FF\ar[r]^0\ar[d]^1
&\FF\ar@{=}[r]\ar@{-->}[d]
&\FF\ar[d]^1\\
\A_*\ar@{ >->}[r]
&\F_*\ar[r]
&{R_\aB}_*\ar@{->>}[r]
&\A_*
}
$$
so that the unit element of ${R_\aB}_*$ is the map $R_\aB\to\FF$ sending the
generator $p1_{\aB_0}$ in degree 0 to $1$ and all elements in higher degrees
to zero.

\section{Hopf pair coalgebras}

We next turn to the dualization of the notion of a Hopf pair algebra from
\bref{hpa}, using the dual $\R_\FF$ of $\R^\FF$ from the previous section.

\begin{Definition}\label{hpc}
A \emph{Hopf pair coalgebra} $\W$ (associated to $\A_*$) is a pair
coalgebra $d:\W^0\to\W^1$ over $\FF$ together with the following
commutative diagram in the category of $\F_*$-$\F_*$-bicomodules
$$
\xymatrix{
\A_*\ar@{ >->}[r]\ar@{=}[d]
&\R^0_\FF\ar[r]\ar@{=}[d]
&\R^1_\FF\ar@{ >->}[d]^i\ar@{->>}[r]
&\A_*\ar@{ >->}[d]^i\\
\A_*\ar@{ >->}[r]
&\W^0\ar[r]^d
&\W^1\ar@{->>}[d]^{\pi_\Sigma}\ar@{->>}[r]^-{\binom\pi{\pi_\Sigma}}
&\A_*\oplus_{\k_*}\Sigma\A_*\ar@{->>}[d]\\
&&\Sigma\A_*\ar@{=}[r]
&\Sigma\A_*
}
$$
with exact rows and columns. The pair morphism $i:\R_\FF\to\W$ will be
called the \emph{$\GG$-structure of $\W$}. Moreover $\W$ must be equipped
with a structure of a monoid $(m_\W,1_\W)$ in $\Coalg^\r_{\oo\oplus\Sigma}$ such
that $i$ is compatible with the $\Coalg^\pair_\oo$-monoid
structure on $\R_\FF$ from \bref{drelcom}, i.~e. diagrams dual to
\eqref{diacomp} and \eqref{coucomp}
$$
\xymatrix{
(\R_\FF\check\ox\R_\FF)^1\ar[r]^-{m_\R}\ar[d]_{i\check\ox i}
&\R_\FF^1\ar[d]^i\\
(\W\check\ox\W)^1\ar[r]^-{m_\W}
&\W^1,
}\ \ \ \ 
\xymatrix{
\FF\ar[r]^{1_\R}\ar[d]
&\R_\FF^1\ar[d]^i\\
\FF\oplus\Sigma\FF\ar[r]^{1_\W}
&\W^1
}
$$
commute.
\end{Definition}

We next note that the dual of \bref{vleft} holds; more precisely, one has
\begin{Lemma}
For a Hopf pair coalgebra $\W$ the subspace
$$
(\W\check\ox\W)^1\subset(\W\dblb\ox\W)^1
$$
is closed under the left coaction of the coalgebra $\F_*$ on
$(\W\dblb\ox\W)^1$ given by the corestriction of scalars along the
multiplication $m_*:\F_*\ox\F_*\to\F_*$ of the left
$\F_*\ox\F_*=(\W\dblb\ox\W)^0$-comodule structure given by the pair
coalgebra $\W\dblb\ox\W$. In other words, there is a unique map
$m^\l:(\W\check\ox\W)^1\to\F_*\ox(\W\check\ox\W)^1$ making the diagram
$$
\xymatrix{
(\W\check\ox\W)^1\ar@{ >->}[d]\ar@{-->}[rr]^{m^\l}
&&\F_*\ox(\W\check\ox\W)^1\ar@{ >->}[d]\\
(\W\dblb\ox\W)^1\ar[r]
&\F_*\ox\F_*\ox(\W\dblb\ox\W)^1\ar[r]^-{m_*\ox1}
&\F_*\ox(\W\dblb\ox\W)^1
}
$$
commute.
\end{Lemma}\qed

Given this left coaction, one can define the dual of the left action
operator in \bref{lv} by measuring deviation of the multiplication
$(\W\check\ox\W)^1\to\W^1$ from being a left comodule homomorphism. For
that, one first observes that the map $\check L:(\W\check\ox\W)^1\to\F_*\ox\W^1$
is given by the difference of two composites in the diagram
$$
\xymatrix{
(\W\check\ox\W)^1\ar[r]^-{m^\l}\ar[d]_{m_\W}
&\F_*\ox(\W\check\ox\W)^1\ar[d]^{1\ox m_\W}\\
\W^1\ar[r]^-{m^\l}
&\F_*\ox\W^1.
}
$$

Then by the argument dual to that before \bref{lv} one sees that the map
$\check L$ factors uniquely through $\coker(i\check\ox
i)=\left((\W\check\ox\W)^1\onto\Sigma\A_*\ox\A_*\right)$ and into
$\ker(d)\ox\im(d)=\left(\A_*\ox{R_\F}_*\into\W^0\ox\W^1\right)$ to yield a map
$\Sigma\A_*\ox\A_*\to\A_*\ox{R_\F}_*$. We thus can make, dually to \bref{lv},
the following

\begin{Definition}
The map
$$
L_\W:\Sigma\A_*\ox\A_*\to\A_*\ox{R_\F}_*
$$
given by the unique factorization of the map $\check L$ above is
characterized by the deviation of the multiplication $m_\W$ of the Hopf
pair coalgebra $\W$ from being a left $\F_*$-comodule homomorphism. That
is, for any $t\in(\W\check\ox\W)^1$ one has
$$
(1\ox m_\W)m^\l(t)=m^\l
m_\W(t)+L_\W(\pi_\Sigma\check\ox\pi_\Sigma)(t).
$$
\end{Definition}

Next, we define a map $S_\W$ in a manner dual to \bref{S}, measuring
noncommutativity of the $\Coalg^\r_{\oo\oplus\Sigma}$-monoid structure on
$\W$. For that, we first consider the map $\check
S:(\W\check\ox\W)^1\to\W^1$ given by
$$
\check S(t)=m_\W T(t)-m_\W(t)
$$
for $t\in(\W\check\ox\W)^1$ and then observe that, dually to \bref{S}, this
map factors uniquely through $\coker(i\check\ox
i)=\left((\W\check\ox\W)^1\onto\Sigma\A_*\ox\A_*\right)$ and into
$\im(d)=\left({R_\F}_*\into\W^1\right)$ so we have

\begin{Definition}
The map
$$
S_\W:\Sigma\A_*\ox\A_*\to{R_\F}_*
$$
given by the unique factorization of the map $\check S$ above is
characterized by being the graded commutator map with respect to the
$\check\ox$-monoid structure on the Hopf pair coalgebra $\W$. That is, for
any $t\in(\W\check\ox\W)^1$ one has
$$
m_\W T(t)=m_\W(t)+S_\W(\pi_\Sigma\check\ox\pi_\Sigma)(t).
$$
\end{Definition}

We now dualize the left action operator \bref{lao} and the symmetry
operator \bref{so}.

\begin{Definition}
The \emph{left coaction operator}
$$
L_*:\A_*\ox\A_*\to\A_*\ox{R_\F}_*
$$
of degree $+1$ is the graded dual of the left action operator \bref{lao}.
\end{Definition}

\begin{Definition}\label{cosym}
The \emph{cosymmetry operator}
$$
S_*:\A_*\ox\A_*\to{R_\F}_*
$$
of degree $+1$ is the graded dual of the symmetry operator \bref{so}.
\end{Definition}

It is clear that the duals of \bref{exist} and \bref{unique} hold. Let us
state these explicitly.

\begin{Theorem}[Existence]
There exists a Hopf pair coalgebra $\W$ with $L_\W=L_*$ and $S_\W=S_*$.
\end{Theorem}

\begin{Theorem}[Uniqueness]
The Hopf pair coalgebra $\W$ satisfying $L_\W=L_*$ and $S_\W=S_*$ is unique
up to an isomorphism over $\W\onto\A_*\oplus_{\k_*}\Sigma\A_*$ and under $\R_\FF\into\W$.
\end{Theorem}

The Hopf pair coalgebra appearing in these theorems will be denoted by
$\aB_\FF=(\aB^0_\FF\to\aB^1_\FF)=D(\aB^\FF)$.

%% file: dHascoAss_c5.tex
\chapter{Generators of $\aB_\FF$ and dual generators of $\aB^\FF$}\label{gens}

In this chapter we describe polynomial generators in the dual Steenrod algebra
$\A_*$ and in the dual of the free tensor algebra $T_\FF(E_\A)$ with the
Cartan diagonal. We use these results to obtain generators in the dual of the
relation module $R_\F$.

\section{The Milnor dual of the Steenrod algebra}

Here we recall the needed facts from \cite{Milnor}. The graded dual of the
Hopf algebra $\A$ is the Milnor Hopf algebra $\A_*=\Hom(\A,\FF)=D(\A)$. It is proved in
\cite{Milnor} that for odd $p$ as an algebra $\A_*$ is a graded polynomial
algebra, i.~e. it is isomorphic to a tensor product of an exterior algebra
on generators of odd degree and a polynomial algebra on generators of even
degree; for $p=2$ the algebra $\A_*$ is a polynomial algebra. Moreover, in
\cite{Milnor}, explicit generators are given in terms of the admissible
basis.

First recall that the admissible basis for $\A$ is given by the following
monomials: for odd $p$ they are of the form
$$
M=\beta^{\epsilon_0}\P^{s_1}\beta^{\epsilon_1}\P^{s_2}\cdots\P^{s_n}\beta^{\epsilon_n}
$$
where $\epsilon_k\in\{0,1\}$ and
$$
s_1\ge\epsilon_1+ps_2,
s_2\ge\epsilon_2+ps_3,\dots,s_{n-1}\ge\epsilon_{n-1}+ps_n,s_n\ge1.
$$
Then let $\xi_k\in\A_{2(p^k-1)}=\Hom(\A^{2(p^k-1)},\FF)$, $k\ge1$ and
$\tau_k\in\A_{2p^k-1}=\Hom(\A^{2p^k-1},\FF)$, $k\ge0$ be given on this
basis by
\begin{equation}
\xi_k(M)=
\begin{cases}
1,&M=\P^{p^{k-1}}\P^{p^{k-2}}\cdots\P^p\P^1,\\
0&\textrm{otherwise}
\end{cases}
\end{equation}
and
\begin{equation}
\tau_k(M)=
\begin{cases}
1,&M=\P^{p^{k-1}}\P^{p^{k-2}}\cdots\P^p\P^1\beta,\\
0&\textrm{otherwise}.
\end{cases}
\end{equation}
As proved in \cite{Milnor}, $\A_*$ is a graded polynomial algebra on these
elements, i.~e. it is generated by the elements $\xi_k$ and $\tau_k$ with
the defining relations
\begin{align*}
\xi_i\xi_j&=\xi_j\xi_i,\\
\xi_i\tau_j&=\tau_j\xi_i,\\
\tau_i\tau_j&=-\tau_j\tau_i
\end{align*}
only.

For $p=2$, the admissible basis for $\A$ is given by the monomials
$$
M=\Sq^{s_1}\Sq^{s_2}\cdots\Sq^{s_n}
$$
with
$$
s_1\ge2s_2,s_2\ge2s_3,\dots,s_{n-1}\ge2s_n,s_n\ge1
$$
and the polynomial generators of $\A_*$ are elements
$\zeta_k\in\A_{2^k-1}=\Hom(\A^{2^k-1},\FF)$ given by
\begin{equation}\label{milgen}
\zeta_k(M)=
\begin{cases}
1,&M=\Sq^{2^{k-1}}\Sq^{2^{k-2}}\cdots\Sq^2\Sq^1,\\
0&\textrm{otherwise}.
\end{cases}
\end{equation}

In terms of these generators, likewise, the coalgebra structure
$m_*:\A_*\to\A_*\ox\A_*$ dual to the multiplication $m$ of $\A$
is determined in \cite{Milnor}. Namely, for odd $p$ one has
\begin{equation}
\alignbox{
m_*(\xi_k)&=\xi_k\ox1+\xi_{k-1}^p\ox\xi_1+\xi_{k-2}^{p^2}\ox\xi_2+\dots+\xi_1^{p^{k-1}}\ox\xi_{k-1}+1\ox\xi_k,\\
m_*(\tau_k)&=\xi_k\ox\tau_0+\xi_{k-1}^p\ox\tau_1+\xi_{k-2}^{p^2}\ox\tau_2+\dots+\xi_1^{p^{k-1}}\ox\tau_{k-1}+1\ox\tau_k+\tau_k\ox1.
}
\end{equation}
For $p=2$ one has
\begin{equation}\label{mildiag}
m_*(\zeta_k)=\zeta_k\ox1+\zeta_{k-1}^2\ox\zeta_1+\zeta_{k-2}^4\ox\zeta_2+\dots+\zeta_1^{2^{k-1}}\ox\zeta_{k-1}+1\ox\zeta_k.
\end{equation}

We will need an expression for the dual $\Sq^1_*:\A_*\to\Sigma\A_*$ to the map
$\Sq^1\cdot:\Sigma\A\to\A$ given by multiplication with $\Sq^1$ from the
left.

\begin{Lemma}\label{sqd}
The map $\Sq^1_*$ is equal to $\frac\partial{\partial\zeta_1}$. That is, on
the monomial basis it is given by
$$
\Sq^1_*(\zeta_1^{n_1}\zeta_2^{n_2}\cdots)=
\begin{cases}
\zeta_1^{n_1-1}\zeta_2^{n_2}\cdots,&n_1\equiv1\mod2\\
0,&n_1\equiv0\mod2.
\end{cases}
$$
\end{Lemma}

\begin{proof}
Note that $\Sq^1_*$ is a derivation, since $\Sq^1\cdot$ is a coderivation,
i.e. the diagram
$$
\xymatrix{
\Sigma\A\ar[rr]^{\Sq^1\cdot}\ar[d]_{\Sigma\delta}
&&\A\ar[d]^\delta\\
\Sigma\A\ox\A\ar[r]^-{\binom{\Sq^1\cdot\ox1}{1\ox\Sq^1\cdot}}
&\A\!\ox\!\A\x\A\!\ox\!\A\ar[r]^-{+}
&\A\ox\A
}
$$
commutes: indeed for any $x\in\A$ one has
$$
\delta(\Sq^1x)=\delta(\Sq^1)\delta(x)=(\Sq^1\ox1+1\ox\Sq^1)\delta(x)
=(\Sq^1\ox1)\delta(x)+(1\ox\Sq^1)\delta(x).
$$
On the other hand, the derivation on the Milnor generators $\Sq^1_*$ acts
as follows:
$$
\Sq^1_*(\zeta_n)(x)=\zeta_n(\Sq^1x)=
\begin{cases}
1,&\Sq^1x=\Sq^{2^{n-1}}\Sq^{2^{n-2}}\cdots\Sq^1,\\
0,&\Sq^1x\ne\Sq^{2^{n-1}}\Sq^{2^{n-2}}\cdots\Sq^1.
\end{cases}
$$
It follows that $\Sq^1_*(\zeta_1)=1$; on the other hand for $n>1$ the
equation $\Sq^1x=\Sq^{2^{n-1}}\Sq^{2^{n-2}}\cdots\Sq^1$ has no solutions,
since it would imply
$\Sq^1\Sq^{2^{n-1}}\Sq^{2^{n-2}}\cdots\Sq^1=\Sq^1\Sq^1x=0$, whereas
actually
$$
\Sq^1\Sq^{2^{n-1}}\Sq^{2^{n-2}}\cdots\Sq^1=\Sq^{1+2^{n-1}}\Sq^{2^{n-2}}\cdots\Sq^1\ne0.
$$

But $\frac\partial{\partial\zeta_1}$ is the unique derivation sending
$\zeta_1$ to 1 and all other $\zeta_n$'s to 0.
\end{proof}

We will also need expression of the dual $\k_*$ of the
derivation $\k$ from \eqref{kappa} in terms of the above generators.

\begin{Lemma}\label{dkappa}
The map $\k_*:\Sigma\A_*\to\A_*$ is equal to the left multiplication by $\tau_0$
for odd $p$ and by $\zeta_1$ for $p=2$.
\end{Lemma}

\begin{proof}
For any linear map $\phi:\A^n\to\FF$ the map
$\k_*(\phi):\A_{n+1}\to\FF$ is the composite of $\phi$ with
$\k:\A_{n+1}\to\A_n$. Thus for $p$ odd one has
\begin{equation}\label{kodd}
\k_*(\phi)(\beta^{\epsilon_0}\P^{s_1}\beta^{\epsilon_1}\P^{s_2}\cdots\P^{s_n}\beta^{\epsilon_n})
=\sum_{\epsilon_k=1}(-1)^{\epsilon_0+\epsilon_1+\dots+\epsilon_{k-1}}
\phi(\beta^{\epsilon_0}\P^{s_1}\beta^{\epsilon_1}\cdots\beta^{\epsilon_{k-1}}\P^{s_k}\P^{s_{k+1}}\beta^{\epsilon_{k+1}}\cdots\P^{s_n}\beta^{\epsilon_n}).
\end{equation}
On the other hand, one has for $M$ as above
$$
(\tau_0\phi)(M)=\sum\tau_0(M_\l)\phi(M_\r)=\sum_{\substack{M_\l=c\beta\\0\ne
c\in\FF}}c\phi(M_\r),
$$
if
$$
\delta(M)=\sum M_\l\ox M_\r.
$$
On the other hand one evidently has
$$
\delta(\beta^{\epsilon_0}\P^{s_1}\beta^{\epsilon_1}\P^{s_2}\cdots\P^{s_n}\beta^{\epsilon_n})
=\sum_{\substack{0\le\iota_0\le\epsilon_0\\0\le i_1\le
s_1\\0\le\iota_1\le\epsilon_1\\\cdots\\0\le i_n\le
s_n\\0\le\iota_n\le\epsilon_n}}
(-1)^{\sum_{0\le\mu<\nu\le n}(\epsilon_\mu-\iota_\mu)\iota_\nu}
\beta^{\iota_0}\P^{i_1}\beta^{\iota_1}\cdots\P^{i_n}\beta^{\iota_n}
\ox
\beta^{\epsilon_0-\iota_0}\P^{s_1-i_1}\beta^{\epsilon_1-\iota_1}\cdots\P^{s_n-i_n}\beta^{\epsilon_n-\iota_n}
$$
so that for
$M=\beta^{\epsilon_0}\P^{s_1}\beta^{\epsilon_1}\cdots\P^{s_n}\beta^{\epsilon_n}$
one has
\begin{multline*}
\sum_{\substack{M_\l=c\beta\\0\ne c\in\FF}}c\phi(M_\r)
=\sum_{\epsilon_k=1}\sum_{\substack{\iota_0=0\\i_1=0\\\cdots\\i_k=0\\\iota_k=1\\i_{k+1}=0\\\cdots\\i_n=0\\\iota_n=0}}
(-1)^{\sum_{0\le\mu<\nu\le n}(\epsilon_\mu-\iota_\mu)\iota_\nu}
\phi(\beta^{\epsilon_0-\iota_0}\P^{s_1-i_1}\beta^{\epsilon_1-\iota_1}\cdots\P^{s_n-i_n}\beta^{\epsilon_n-\iota_n})\\
=\sum_{\epsilon_k=1}
(-1)^{\sum_{0\le\mu<k}\epsilon_\mu}
\phi(\beta^{\epsilon_0}\P^{s_1}\beta^{\epsilon_1}\cdots\P^{s_k}\P^{s_{k+1}}\beta^{\epsilon_{k+1}}\cdots\P^{s_n}\beta^{\epsilon_n})
\end{multline*}
which is the same as \eqref{kodd} above.

Similarly for $p=2$ the map $\k_*(\phi)$ is given by
\begin{equation}\label{k2}
\k_*(\phi)(\Sq^{s_1}\cdots\Sq^{s_n})=\phi(\k(\Sq^{s_1}\cdots\Sq^{s_n}))=\sum_{k=1}^n\phi(\Sq^{s_1}\cdots\Sq^{s_k-1}\cdots\Sq^{s_n})
\end{equation}
and the map $\zeta_1\phi$ is given by
$$
(\zeta_1\phi)(M)=\sum\zeta_1(M_\l)\phi(M_\r)=\sum_{M_\l=\Sq^1}\phi(M_\r).
$$
On the other hand one has
$$
\delta(\Sq^{s_1}\cdots\Sq^{s_n})=\sum_{\substack{0\le i_1\le
s_1\\\cdots\\0\le i_n\le
s_n}}\Sq^{i_1}\cdots\Sq^{i_n}\ox\Sq^{s_1-i_1}\cdots\Sq^{s_n-i_n},
$$
so that for $M=\Sq^{s_1}\cdots\Sq^{s_n}$ one has
$$
\sum_{M_\l=\Sq^1}\phi(M_\r)=\sum_{k=1}^n\sum_{\substack{i_1=0\\\cdots\\i_{k-1}=0\\i_k=1\\i_{k+1}=0\\\cdots\\i_n=0}}
\phi(\Sq^{s_1-i_1}\cdots\Sq^{s_n-i_n})
$$
which is equal to \eqref{k2}.
\end{proof}

It is clear that with respect to the coalgebra structure on $\A_*$ the map
$\k_*$ is a coderivation, i.~e. the diagram
$$
\xymatrix@!C{
\Sigma\A_*\ar[rr]^{\k_*}\ar[d]_{\Sigma m_*}
&&\A_*\ar[d]^{m_*}\\
\Sigma(\A_*\!\ox\!\A_*)\ar[r]^-{\binom1{\sigma}}
&\Sigma\A_*\!\ox\!\A_*\oplus\A_*\!\ox\!\Sigma\A_*\ar[r]^-{(\k_*\ox1,1\ox\k_*)}
&\A_*\ox\A_*
}
$$
is commutative. Here $\sigma$ is the interchange of $\Sigma$ as in
\eqref{sigma}. Then using dual of the construction mentioned in \bref{hoch}
one may equip the vector space $\A_*\oplus\Sigma\A_*$ with a structure of
an $\A_*$-$\A_*$-bicomodule, in such a way that one has a short exact
sequence of $\A_*$-$\A_*$-bicomodules
\begin{equation}
0\to\A_*\to\A_*\oplus_{\k_*}\Sigma\A_*\to\Sigma\A_*\to0.
\end{equation}
Explicitly, one defines the right coaction of $\A_*$ on
$\A_*\oplus_{\k_*}\Sigma\A_*$ as the direct sum of standard coactions on
$\A_*$ and on $\Sigma\A_*$, whereas the left coaction is given by the
composite
$$
\A_*\oplus\Sigma\A_*\xto{m_*\oplus\Sigma m_*}
\A_*\!\ox\!\A_*\oplus\Sigma\A_*\!\ox\!\A_*
\xto{\left(\begin{smallmatrix}1&\k_*\ox1\\0&\sigma\end{smallmatrix}\right)}
\A_*\!\ox\!\A_*\oplus\A_*\!\ox\!\Sigma\A_*\cong\A_*\ox(\A_*\oplus\Sigma\A_*).
$$

\section{The dual of the tensor algebra $\F_0=T_\FF(E_\A)$ for $p=2$}

We begin by recalling the constructions from \cite{Hazewinkel} relevant to our case.

The \emph{Leibniz-Hopf algebra} is the free graded associative ring with
unit $1=Z_0$
\begin{equation}
\Z=T_\ZZ\{Z_1,Z_2,...\}
\end{equation}
on generators $Z_n$, one for each degree $n\ge1$. Here we use notation as
in \eqref{BF}. $\Z$ is a cocommutative Hopf algebra with respect to the
diagonal
$$
\Delta(Z_n)=\sum_{i=0}^nZ_i\ox Z_{n-i}.
$$
Of course for $p=2$ we have $\Z\ox\FF=\F_0=T_\ZZ(E_\A)$ by identifying
$Z_i=\Sq^i$, and moreover the diagonal $\Delta$ corresponds to
$\Delta^\GG\ox\FF$ in \eqref{deltag}. The graded dual of $\Z$ over the
integers is denoted by $\M$; it is proved in \cite{Hazewinkel} that it is
a polynomial algebra. There also a certain set of elements of $\M$ is given; it is
still a conjecture (first formulated by Ditters) that these elements form
a set of polynomial generators for $\M$. If, however, one localizes at any
prime $p$, then there is another set of elements, defined using the so
called \emph{$p$-elementary words}, which, as proved in \cite{Hazewinkel},
is a set of polynomial generators for the localized algebra $\M$. This in
particular gives a polynomial generating set for
$\F_*=\Hom(\F_0,\FF_2)\cong\M/2\M$. Moreover it
turns out that the embedding $\A_*\into\F_*$ given by $\Hom(\A,\FF_2)\into\Hom(\F_0,\FF_2)$ (dual
to the quotient map $\F_0\onto\A$) carries the Milnor generators of $\A_*$
to a subset of these generators.

Choose a basis in $\M$ which is dual to the (noncommutative) monomial
basis in $\Z$: for any sequence $\alpha=(d_1,...,d_n)$ of positive
integers, let $M_\alpha=M_{d_1,...,d_n}$ be the element of the free
abelian group $\M^{d_1+...+d_n}=\Hom(\Z^{d_1+...+d_n},\ZZ)$ determined by
$$
M_{d_1,...,d_n}(Z_{k_1}\cdots Z_{k_m})=
\begin{cases}
1,&(k_1,...,k_m)=(d_1,...,d_n),\\
0&\textrm{otherwise.}
\end{cases}
$$

Since $\Z$ is a free algebra, dually $\M$ is a cofree coalgebra, i.~e. the
diagonal is given by deconcatenation:
\begin{equation}\label{deconc}
\Delta(M_{d_1,...,d_n})=\sum_{i=0}^nM_{d_1,...,d_i}\ox M_{d_{i+1},...,d_n}.
\end{equation}

It is noted in \cite{Hazewinkel} (and easy to check) that in this basis the
multiplication in $\M$ is given by the so called \emph{overlapping shuffle
product}. Rather than defining this rigorously, we will give some examples.
\begin{align*}
M_5M_{2,4,1,9}&=M_{5,2,4,1,9}+M_{7,4,1,9}+M_{2,5,4,1,9}+M_{2,9,1,9}+M_{2,4,5,1,9}+M_{2,4,6,9}\\
&+M_{2,4,1,5,9}+M_{2,4,1,14}+M_{2,4,1,9,5};\\
M_{8,5}M_{1,2}&=M_{8,5,1,2}+M_{8,6,2}+M_{8,1,5,2}+M_{9,5,2}+M_{8,1,7}+M_{9,7}+M_{1,8,5,2}\\
&+M_{1,8,7}+M_{1,8,2,5}+M_{9,2,5}+M_{1,2,8,5}+M_{1,10,5}+M_{8,1,2,5}
\end{align*}
Thus in general, whereas the ordinary shuffle product of the elements, say,
$M_{a_1,a_2,a_3}$ and $M_{b_1,b_2,b_3,b_4,b_5}$ contains all possible
summands like $M_{b_1,a_1,a_2,b_2,b_3,a_3,b_4,b_5}$, the overlapping
shuffle product contains together with each such summand also in addition the summands of the
form
$M_{b_1+a_1,a_2,b_2,b_3,a_3,b_4,b_5}$,
$M_{b_1,a_1,a_2+b_2,b_3,a_3,b_4,b_5}$,
$M_{b_1,a_1,a_2,b_2,b_3+a_3,b_4,b_5}$,
$M_{b_1,a_1,a_2,b_2,b_3,a_3+b_4,b_5}$,
$M_{b_1+a_1,a_2+b_2,b_3,a_3,b_4,b_5}$ and so on, obtained by replacing
an $a_i$ and a $b_j$ standing one next to other with their sum, in all
possible positions.

Note that the algebra of ordinary shuffles is also a polynomial algebra,
but over rationals; it is \emph{not} a polynomial algebra until at least
one prime number remains uninverted. On the other hand, over rationals $\M$
becomes isomorphic to the algebra of ordinary shuffles.

To define a polynomial generating set for $\M$, we need some definitions.
To conform with the admissible basis in the Steenrod algebra, which
consists of monomials with decreasing indices, we will reverse the order of
indices in the definitions from \cite{Hazewinkel}, where the indices go in the
increasing order. Thus in our case statements about some $M_{d_1,...,d_n}$
will be equivalent to the corresponding ones in \cite{Hazewinkel} about
$M_{d_n,...,d_1}$.

\begin{Definitions}
The \emph{lexicographic order} on the basis $M_{d_1,...,d_n}$ of $\M$ is
defined by declaring $M_{d_1,...,d_n}>M_{e_1,...,e_m}$ if either there is an $i$
with $1\le i\le\min(n,m)$ and $d_i>e_i$, $d_n=e_m$, $d_{n-1}=e_{m-1}$, ...,
$d_{n-i+1}=e_{m-i+1}$, $d_{n-i}>e_{m-i}$
or $n>m$ and $d_{n-m+1}=e_1$, $d_{n-m+2}=e_2$, ..., $d_n=e_m$.

A basis element $M_{d_1,...,d_n}$ is \emph{Lyndon} if with respect to this
ordering one has $M_{d_1,...,d_n}<M_{d_1,...,d_i}$ for all $1<i\le n$. For
example, $M_{3,2,3,2,2}$ and $M_{2,2,1,2,1,2,1}$ are Lyndon but
$M_{3,2,2,3,2}$ and $M_{2,1,2,1,2,1}$ are not.

A basis element $M_{d_1,...,d_n}$ is \emph{$\ZZ$-elementary} if no number
$>1$ divides all of the $d_i$, i.~e. $\gcd(d_1,...,d_n)=1$. The set
$\ESL(\ZZ)$ is the set of elementary basis elements of the form
$M_{d_1,...,d_n,d_1,...,d_n,....,d_1,...,d_n}$ (i.~e. $d_1,...,d_n$
repeated any number of times), where $M_{d_1,...,d_n}$ is a 
Lyndon element.

For a prime $p$, a basis element $M_{d_1,...,d_n}$ is called
\emph{$p$-elementary} if there is a $d_i$ not divisible by $p$, i.~e.
$p\nmid\gcd(d_1,...,d_n)$. The set $\ESL(p)$ is defined as the set of
$p$-elementary basis elements of the form
$$
M_{\underbrace{\scriptstyle{d_1,...,d_n,d_1,...,d_n,...,d_1,...,d_n}}_{p^r\textrm{
times}}}
$$
with $d_1,...,d_n$ repeated $p^r$ times for some $r$, where
$M_{d_1,...,d_n}$ is required to be Lyndon.
\end{Definitions}

For example, $M_{15,6,15,6,15,6,15,6}$ is in $\ESL(2)$ but not in $\ESL(\ZZ)$
or in $\ESL(p)$ for any other $p$, whereas $M_{30,6,6}$ is in $\ESL(p)$ for
any $p\ne2,3$ but not in $\ESL(2)$, not in $\ESL(3)$ and not in $\ESL(\ZZ)$.

One then has

\begin{Theorem}[\cite{Hazewinkel}]
The algebra $\M$ is a polynomial algebra.
\end{Theorem}

\begin{Conjecture}[Ditters, \cite{Hazewinkel}]
The set $\ESL(\ZZ)$ is the set of polynomial generators for $\M$.
\end{Conjecture}

\begin{Theorem}[\cite{Hazewinkel}]\label{eslp}
For each prime $p$, the set $\ESL(p)$ is a set of polynomial
generators for $\M_{(p)}=\M\ox\ZZ_{(p)}$, i.~e. if one inverts all primes
except $p$.
\end{Theorem}

In particular, it follows that $\ESL(p)$ is a set of polynomial generators
for $\M/p^n$ over $\ZZ/p^n$ for all $n$.

Here are the polynomial generators in low degrees, over $\ZZ$ and over few
first primes. Note that the numbers of generators in each degree are the
same (as it should be since all these algebras become isomorphic over
$\QQ$). 

$$
\begin{array}{c||c|c|c|c|c}
&1&2&3&4&5\\
\hline
\hline
&&&&&\\
\ZZ&M_1&M_{1,1}&M_{2,1},M_{1,1,1}&M_{3,1},M_{2,1,1},M_{1,1,1,1}&M_{4,1},M_{3,2},M_{3,1,1},
M_{2,2,1},M_{2,1,1,1},M_{1,1,1,1,1}\\
&&&&&\\
\hline
&&&&&\\
p=2&M_1&M_{1,1}&M_3,M_{2,1}&M_{3,1},M_{2,1,1},M_{1,1,1,1}&M_5,M_{4,1},M_{3,2},
M_{3,1,1},M_{2,2,1},M_{2,1,1,1}\\
&&&&&\\
\hline
&&&&&\\
p=3&M_1&M_2&M_{2,1},M_{1,1,1}&M_4,M_{3,1},M_{2,1,1}&M_5,M_{4,1},M_{3,2},
M_{3,1,1},M_{2,2,1},M_{2,1,1,1}\\
&&&&&\\
\hline
&&&&&\\
p=5&M_1&M_2&M_3,M_{2,1}&M_4,M_{3,1},M_{2,1,1}&M_{4,1},M_{3,2},M_{3,1,1},
M_{2,2,1},M_{2,1,1,1},M_{1,1,1,1,1}\\
&&&&&
\end{array}
$$

It is easy to calculate the numbers of polynomial generators in each
degree. Let these numbers be $m_1$, $m_2$, $\cdots$. Then the Poincar\'e
series for the algebra $\M$ (or $\Z$, or $\F$, or $\F_*$, it does not
matter) is
$$
\sum_n\dim(\M_n)t^n=(1-t)^{-m_1}(1-t^2)^{-m_2}(1-t^3)^{-m_3}\cdots;
$$
on the other hand, we know that it is a tensor coalgebra with one generator
in each degree $n\ge1$; this implies that $\dim(\M_n)=2^{n-1}$ for $n\ge1$
(and $\dim(M_0)=1$). Thus we have equality of power series
$$
\prod_{k=1}^\infty(1-t^k)^{-m_k}=1+t+2t^2+4t^3+8t^4+\cdots
=1+t(1+2t+(2t)^2+(2t)^3+\cdots)=1+t\frac1{1-2t}=\frac{1-t}{1-2t}.
$$
Then taking logarithmic derivatives one obtains
$$
\sum_{k=1}^\infty\frac{km_kt^k}{1-t^k}=\frac{2t}{1-2t}-\frac
t{1-t}=t+3t^2+7t^3+\cdots+(2^n-1)t^n+\cdots.
$$
It follows that for all $n$ one has
$$
\sum_{d|n}dm_d=2^n-1,
$$
which by the M\"obius inversion formula gives
$$
m_n=\frac1n\sum_{d|n}\mu(d)(2^{\frac nd}-1).
$$
The latter expression is well known in the literature on combinatorics; it
equals the number of aperiodic bicolored necklaces consisting of $n$
beads, and also the dimension of the $n$th homogeneous component of the
free Lie algebra on two generators. See e.~g. \cite{S}.

\section{The dual of the relation module $R_\F$}

We now turn to the algebra $\F_*=\Hom(\F,\FF_2)\cong\M/2$. By the above,
we know that it, as well as $\M_{(2)}$, is a polynomial algebra on the set
of generators $\ESL(2)$. As an illustration, we will give some expressions of
the $M$-basis elements in terms of sums of overlapping shuffle products of
elements from $\ESL(2)$. We will give these in $\M_{(2)}$ and then their
images in $\F_*$.

\begin{align*}
M_2&=M_1^2-2M_{1,1}\\
&\equiv M_1^2\mod2\\
M_{1,2}&=M_1^3-M_3-M_{2,1}-2M_1M_{1,1}\\
&\equiv M_1^3+M_3+M_{2,1}\mod2\\
M_{1,1,1}&=M_1M_{1,1}-\frac13M_1^3+\frac13M_3\\
&\equiv M_1M_{1,1}+M_1^3+M_3\mod2\\
M_4&=\frac43M_1M_3-\frac13M_1^4+2M_{1,1}^2-4M_{1,1,1,1}\\
&\equiv M_1^4\mod2\\
M_{2,2}&=M_{1,1}^2-2M_1^2M_{1,1}-\frac23M_1M_3+\frac23M_1^4+2M_{1,1,1,1}\\
&\equiv M_{1,1}^2\mod2\\
M_{1,3}&=\frac13M_1^4-\frac13M_1M_3-2M_{1,1}^2-M_{3,1}+4M_{1,1,1,1}\\
&\equiv M_1^4+M_1M_3+M_{3,1}\mod2\\
M_{1,2,1}&=M_1M_{2,1}-M_{3,1}-M_{1,1}^2+2M_1^2M_{1,1}+\frac23M_1M_3-\frac23M_1^4-2M_{1,1,1,1}-2M_{2,1,1}\\
&\equiv M_1M_{2,1}+M_{3,1}+M_{1,1}^2\mod2\\
M_{1,1,2}&=M_{1,1}^2-M_1^2M_{1,1}-\frac13M_1M_3+\frac13M_1^4-2M_{1,1,1,1}+M_{3,1}-M_1M_{2,1}+M_{2,1,1}\\
&\equiv M_{1,1}^2+M_1^2M_{1,1}+M_1M_3+M_1^4+M_{3,1}+M_1M_{2,1}+M_{2,1,1}\mod2
\end{align*}

Moreover it is straightforward to calculate the diagonal in terms of these
generators. For example, in $\F_*$ one has
\begin{align*}
\Delta(M_1)&=1\ox M_1+M_1\ox1,\\
\Delta(M_{1,1})&=1\ox M_{1,1}+M_1\ox M_1+M_{1,1}\ox1,\\
\Delta(M_3)&=1\ox M_3+M_3\ox1\\
\Delta(M_{2,1})&=1\ox M_{2,1}+M_1^2\ox M_1+M_{2,1}\ox1\\
\Delta(M_{3,1})&=1\ox M_{3,1}+M_3\ox M_1+M_{3,1}\ox1\\
\Delta(M_{2,1,1})&=1\ox M_{2,1,1}+M_1^2\ox M_{1,1}+M_{2,1}\ox
M_1+M_{2,1,1}\ox1\\
\Delta(M_{1,1,1,1})&=1\ox M_{1,1,1,1}+M_1\ox M_1M_{1,1}+M_1\ox M_1^3+M_1\ox
M_3+M_{1,1}\ox M_{1,1}+M_1M_{1,1}\ox M_1\\
&+M_1^3\ox M_1+M_3\ox M_1+M_{1,1,1,1}\ox1\\
\Delta(M_{4,1})&=1\ox M_{4,1}+M_1^4\ox M_1+M_{4,1}\ox1\\
\Delta(M_{3,2})&=1\ox M_{3,2}+M_3\ox M_1^2+M_{3,2}\ox1\\
\Delta(M_{2,1,1,1})&=1\ox M_{2,1,1,1}+M_1^2\ox M_1M_{1,1}+M_1^2\ox
M_1^3+M_1^2\ox M_3+M_{2,1}\ox M_{1,1}+M_{2,1,1}\ox M_1\\
&+M_{2,1,1,1}\ox1\\
\Delta(M_5)&=1\ox M_5+M_5\ox1\\
\Delta(M_{3,1,1})&=1\ox M_{3,1,1}+M_3\ox M_{1,1}+M_{3,1}\ox
M_1+M_{3,1,1}\ox1\\
\Delta(M_{2,2,1})&=1\ox M_{2,2,1}+M_1^2\ox M_{2,1}+M_{1,1}^2\ox
M_1+M_{2,2,1}\ox1.
\end{align*}

Also it follows from the results in \cite{Hazewinkel} that one has

\begin{Lemma}\label{pthpow}
For any prime $p$, in $\M_{(p)}$ one has
$$
M_{pd_1,...,pd_n}\equiv M_{d_1,...,d_n}^p\mod p.
$$
\end{Lemma}

To identify the elements to which the Milnor generators $\zeta_k$ of $\A_*$
go under the isomorphism $\F_*\cong\M/2$, we first identify $\A_*$ with
the graded dual of $\A$; then $\zeta_k$ corresponds to a linear form
$\A_{2^k-1}\to\FF$ given by \eqref{milgen}. 

\begin{Proposition}\label{imxi}
Under the embedding $\A_*\into\M/2$, the Milnor generator $\zeta_k$ maps to
the generator $M_{2^{k-1},2^{k-2},...,2,1}$. In particular, this generator is in
$\ESL(2)$, i.~e. is one of the polynomial generators of $\F_*$.
\end{Proposition}

Note that this together with \eqref{deconc} and \bref{pthpow} implies the
Milnor formula \eqref{mildiag} for the diagonal in $\A_*$. Identifying
$\zeta_k$ with its image in $\M/2$ by \bref{imxi}, one obtains
\begin{equation}\label{dizeta}
\begin{aligned}
m_*(\zeta_k)=\Delta(M_{2^{k-1},2^{k-2},...,2,1})=\sum_{i=0}^kM_{2^{k-1},2^{k-2},...,2^i}\ox
M_{2^{i-1},...,2,1}&=\sum_{i=0}^kM_{2^{k-1-i},2^{k-2-i},...,2,1}^{2^i}\ox M_{2^{i-1},...,2,1}\\
&=\sum_{i=0}^k\zeta_{k-i}^{2^i}\ox\zeta_i.
\end{aligned}
\end{equation}

Thus the set $\{\zeta_1,\zeta_2,...\}$ of polynomial generators for $\A_*$ can be
identified with the subset
$$
Q=\{M_1,M_{2,1},M_{4,2,1},M_{8,4,2,1},...\}
$$
of the set of polynomial generators $\ESL(2)$ for $\M/2\cong\F_*$. This in
particular gives an explicit basis for ${R_\F}_*$: it is in one-to-one
correspondence with those monomials in the generators $M_{d_1,...,d_n}$
from $\ESL(2)$ not all of whose variables belong to $Q$. For example, in the first
few dimensions this basis contains the following monomials:
\begin{align*}
&M_{1,1},\\
&M_1M_{1,1},M_3,\\
&M_1^2M_{1,1},M_1M_3,M_{1,1}^2,M_{3,1},M_{2,1,1},M_{1,1,1,1},\\
&M_1^3M_{1,1},M_1^2M_3,M_1M_{1,1}^2,M_1M_{3,1},M_1M_{2,1,1},
M_1M_{1,1,1,1},M_{1,1}M_3,M_{1,1}M_{2,1},M_5,M_{4,1},M_{3,2},M_{3,1,1},M_{2,2,1},\\
&\ \ \ \ \ M_{2,1,1,1}.
\end{align*}
We next note that obviously the embedding $\A_*\into\F_*$ identifies $\F_*$
with a polynomial algebra over $\A_*$, namely one has a canonical
isomorphism
\begin{equation}
\F_*\cong\A_*[\ESL(2)\setminus Q].
\end{equation}
In particular, as an $\A_*$-module $\F_*$ is free on the generating set
$\NN^{(\ESL(2)\setminus Q)}$ (= the free commutative monoid on
$\ESL(2)\setminus Q$). Then obviously the quotient module ${R_\F}_*$ is a
free $\A_*$-module with the generating set $\NN^{(\ESL(2)\setminus
Q)}\setminus\set{1}$.

We will need the dual $\F_*^{\le2}$ of the subspace $\F_0^{\le2}\subset\F_0$ spanned by
the monomials of length $\le2$ in the generators $\Sq^i$. Observe that
$\F_0^{\le2}$ is a subcoalgebra of $\F_0$, so that dually
$\F_*\onto\F_*^{\le2}$ is a quotient algebra. We have
\begin{Proposition}
The algebra $\F_*^{\le2}$ is a quotient of the polynomial algebra on three
generators $M_1$, $M_{1,1}$, $M_{2,1}$ by a single relation
$$
M_1M_{1,1}M_{2,1}+M_{1,1}^3+M_{2,1}^2=0.
$$
\end{Proposition}

\begin{proof}
First of all, it is straightforward to calculate in $\F_*$ the sum of the overlapping
shuffle products
\begin{multline*}
M_1M_{1,1}M_{2,1}+M_{1,1}^3+M_{2,1}^2=\\M_{1,4,1}+M_{2,2,2}+M_{2,3,1}+M_{3,1,2}+M_{3,2,1}+M_{2,1,2,1}+M_{3,1,1,1}+M_{1,1,3,1}+M_{1,2,1,1,1}+M_{1,2,1,2}+M_{1,1,1,2,1}
\end{multline*}
so that indeed this gives zero in $\F^{\le2}_*$. 
Let
$$
X=\FF[x_1,x_2,x_3]/(x_1x_2x_3+x_2^3+x_3^2)
$$
be the graded algebra with deg$(x_i)=i$, $i=1,2,3$, so that there is a
homomorphism of algebras $f:X\to\F_*^{\le2}$ sending $x_1\mapsto M_1$,
$x_2\mapsto M_{1,1}$, $x_3\mapsto M_{2,1}$. It is straightforward
to calculate the Hilbert function of $X$, i.~e. the formal power
series
$$
\sum_n\dim(X_n)t^n;
$$
it is equal to
$$
\frac{1-t^6}{(1-t)(1-t^2)(1-t^3)}.
$$
On the other hand $\F_*^{\le2}$ is dual to $\F_0^{\le2}$ and it is
straightforward also to calculate dimensions of homogeneous components of
this space. One then simply checks that these dimensions coincide for $X$
and for $\F_*^{\le2}$. Thus it suffices to show that $f$ is surjective,
i.~e. that $\F_*^{\le2}$ is generated by (the images of) $M_1$, $M_{1,1}$
and $M_{2,1}$.

We will show by induction on degree that every $M_n$ and $M_{i,j}$ can be obtained
as a polynomial in these three elements. In degree 1, $M_1$ is the only
nonzero element. In degree 2, besides $M_{1,1}$ we have $M_2$ which is
equal to $M_1^2$ by \bref{pthpow}. In degree 3, we have
$$
M_1M_{1,1}=M_{1,2}+M_{2,1}+M_{1,1,1}\equiv M_{1,2}+M_{2,1}\mod\F^{>2}_*
$$
and
$$
M_1^3=M_3+M_{1,2}+M_{2,1},
$$
so that in $\F_*^{\le2}$ we may solve
$$
M_{1,2}\equiv M_1M_{1,1}+M_{2,1}
$$
and
$$
M_3\equiv M_1^3+M_1M_{1,1}.
$$
Given now any degree $n>3$, we can obtain any element $M_{i,j}$ with $i>1$,
$j>1$, $i+j=n$ from elements of lower degree since
$$
M_{i,j}\equiv M_{1,1}M_{i-1,j-1}.
$$
Next we also can obtain the element $M_{n-1,1}$ from
$$
M_{n-1,1}+M_{2,n-2}\equiv M_{2,1}M_{n-3}.
$$
Then we can obtain $M_{1,n-1}$ from
$$
M_{1,n-1}+M_{n-1,1}\equiv M_{1,1}M_{n-2},
$$
and finally we can obtain $M_n$ from
$$
M_n+M_{1,n-1}+M_{n-1,1}\equiv M_1M_{n-1}.
$$
\end{proof}

Let us also identify the dual of the product map
$$
\F_0^{\le1}\ox\F_0^{\le1}\to\F_0^{\le2}
$$
in terms of the above generators. By dualizing it is clear that this dual
is the unique factorization in the diagram
$$
\xymatrix{
\F_*\ar[r]^-{m_*}\ar@{->>}[d]
&\F_*\ox\F_*\ar@{->>}[d]\\
\F_*^{\le2}\ar@{-->}[r]
&\F_*^{\le1}\ox\F_*^{\le1}.
}
$$
In particular, it is an algebra homomorphism. Moreover the algebra
$\F_*^{\le1}$ may be identified with the polynomial algebra on a single
generator $M_1=\zeta_1$, with the quotient map $\F_*\to\F_*^{\le1}$ given
by sending $M_1$ to itself and all other polynomial generators from
$\ESL(2)$ to zero. From this it is straightforward to identify the map
$\F_*^{\le2}\to\F_*^{\le1}\ox\F_*^{\le1}$ with the algebra homomorphism
$$
\FF[x_1,x_2,x_3]/(x_1x_2x_3+x_2^3+x_3^2)\to\FF[y_1,z_1]
$$
given by
\begin{equation}\label{mstar}
\begin{aligned}
x_1&\mapsto y_1+z_1\\
x_2&\mapsto y_1z_1\\
x_3&\mapsto y_1^2z_1.
\end{aligned}
\end{equation}

Let us identify in these terms the map $\F_*^{\le2}\onto{R_\F}^{\le2}_*$. One
clearly has
$$
R_\F^{\le2}=R_\F\cap\F_0^{\le2}
$$
in $\F_0$, so that dually one has that the diagram
$$
\xymatrix{
\F_*\ar@{->>}[r]\ar@{->>}[d]&{R_\F}_*\ar@{->>}[d]\\
\F_*^{\le2}\ar@{->>}[r]&{R_\F^{\le2}}_*
}
$$
is pushout. Thus ${R_\F^{\le2}}_*$ is isomorphic to the quotient of
$\F^{\le2}_*$ by the image of the composite
$\A_*\into\F_*\onto\F_*^{\le2}$. That image is clearly the subalgebra
generated by $M_1$ and $M_{2,1}$.

We can alternatively describe ${R_\F^{\le2}}_*$ in terms of linear forms on
$R_\F^{\le2}\subset\F_0^{\le2}$. It is clear that the latter subspace is
spanned by all Adem relations $[n,m]$, $n<2m$. The map
$\pi:\F^{\le2}_*\onto{R_\F^{\le2}}_*$ assigns to a linear form on
$\F_0^{\le2}$ its restriction to $R_\F^{\le2}$. One then clearly has
\begin{equation}
\pi(M_1^k)=\pi(M_{2,1}^k)=0
\end{equation}
for all $k\ge0$; moreover $\pi(M_{1,1})$ is dual to $[1,1]$ in the basis
given by the elements $[n,m]$, i.~e. $M_{1,1}([1,1])=1$ and
$M_{1,1}([n,m])=0$ for all other $n$, $m$. Moreover for
$x,y\in\F^{\le2}_*$ we have
\begin{equation}
(xy)([n,m])=\sum x([n,m]_\l)y([n,m]_\r)
\end{equation}
in the Sweedler notation
$$
\Delta([n,m])=\sum[n,m]_\l\ox[n,m]_\r.
$$
For example, we have
$$
\Delta([1,2])=(1+T)(1\ox[1,2]+\Sq^1\ox[1,1])
$$
which implies that $M_1M_{1,1}$ is dual to $[1,2]$ in this basis, i.~e.
$(M_1M_{1,1})[1,2]=1$ and $(M_1M_{1,1})[n,m]=0$ for all other $n$, $m$.
Similarly
$$
\Delta([1,3])=(1+T)(1\ox[1,3]+\Sq^1\ox[1,2]+\Sq^2\ox[1,1])
$$
and
$$
\Delta([2,2])=(1+T)(1\ox[2,2]+\Sq^1\ox[1,2]+\Sq^2\ox[1,1])+[1,1]\ox[1,1]
$$
imply that $M_{1,1}^2$ is dual to $[2,2]$ whereas
$(M_1^2M_{1,1})[1,3]=(M_1^2M_{1,1})[2,2]=1$, so that dual to $[1,3]$ is
$M_1^2M_{1,1}+M_{1,1}^2$.

We will also need a description of the dual $\bar R_*$ of $\bar
R=R_\F/(R_\F\cdot R_\F)$. For this first note that similarly to the above
$\F_*\ox\F_*$ is a free $\A_*\ox\A_*$-module on
$\NN^{(\ESL(2)\setminus Q)}\x\NN^{(\ESL(2)\setminus Q)}$ and ${R_\F}_*\ox{R_\F}_*$ is a
free $\A_*\ox\A_*$-module on
$\left(\NN^{(\ESL(2)\setminus
Q)}\setminus\set{1}\right)\x\left(\NN^{(\ESL(2)\setminus
Q)}\setminus\set{1}\right)$. Moreover the diagonal
$\Delta_\F:\F_*\to\F_*\ox\F_*$ and its factorization
$\Delta_R:{R_\F}_*\to{R_\F}_*\ox{R_\F}_*$ through the quotient maps
$\F_*\onto{R_\F}_*$, $\F_*\ox\F_*\onto{R_\F}_*\ox{R_\F}_*$ are obviously both
equivariant with respect to the diagonal $\delta:\A_*\to\A_*\ox\A_*$, i.~e.
one has
\begin{equation}
\begin{aligned}
\Delta_\F(af)&=\delta(a)\Delta_\F(f),\\
\Delta_R(ar)&=\delta(a)\Delta_R(r)
\end{aligned}
\end{equation}
for any $a\in\A_*$, $f\in\F_*$, $r\in{R_\F}_*$.

%% file: dHascoAss_c6.tex
\chapter{The invariants $L$ and $S$ and the dual invariants $L_*$ and $S_*$ in
terms of generators}\label{LS}

As proved in \cite{Baues} there are invariants $L$ and $S$ of the
Steenrod algebra which determine the algebra $\aB$ of secondary cohomology
operations up to isomorphism. Therefore $L$ and $S$ and the dual invariants
$L_*$ and $S_*$ also determine $\aB^\FF$ and $\aB_\FF$ respectively. In this
chapter we recall the definition of $L$ and $S$ and we discuss algebraic
properties of $L_*$ and $S_*$.

\section{The left action operator $L$ and its dual}\label{L*}

We recall constructions of the maps $L$ and $S$  from
\cite{Baues}*{14.4,14.5} of the same kind as the operators in \bref{lv}
and \bref{S} respectively. For that, we first introduce the following
notation:
\begin{equation}\label{rbar}
\bar R:=R_\F/(R_\F\cdot R_\F),
\end{equation}
with the quotient map $R_\F\onto\bar R$ denoted by $r\mapsto\bar r$. There is a
well-defined $\A$-$\A$-bimodule structure on $\bar R$ given by
$$
\qf f\bar r=\overline{fr},\ \ \bar r\qf f=\overline{rf}
$$
for $f\in\F_0$, $r\in R_\F$. As we show below $\bar R$ is free both as a
left and as a right $\A$-module (but not as a bimodule). A basis for $\bar
R$ as a right $\A$-module can be found using the set $\PAR\subset R_\F$ of
\emph{preadmissible relations} as defined in \cite{Baues}*{16.5}. These
are the elements of $R_\F$ of the form
$$
\Sq^{n_1}\cdots\Sq^{n_k}[n,m]
$$
where $[n,m]$, $n<2m$, is an Adem relation, the monomial
$\Sq^{n_1}\cdots\Sq^{n_k}$ is admissible (i.~e. $n_1\ge2n_2$,
$n_2\ge2n_3$, ..., $n_{k-1}\ge2n_k$), and moreover $n_k\ge2n$. It is then
proved in \cite{Baues}*{16.5.2} that $\PAR$ is a basis of $R_\F$ as a free
right $\F_0$-module.

It is equally true that $R_\F$ is a free left $\F_0$-module. An explicit
basis $\PAR'$ of $R_\F$ as a left $\F_0$-module consists of \emph{left
preadmissible} relations --- elements of the form
$$
[n,m]\Sq^{m_1}\cdots\Sq^{m_k}
$$
where $[n,m]$, $n<2m$, is an Adem relation, the monomial
$\Sq^{m_1}\cdots\Sq^{m_k}$ is admissible, and moreover $m\ge2m_1$.

Using this, one also has

\begin{Lemma}\label{barr}
$\bar R$ is free both as a right $\A$-module and as a left $\A$-module.
Moreover, the images $\bar\rho$ of the preadmissible relations $\rho\in\PAR$ under the quotient map
$R_\F\onto\bar R$ form a basis of this free right $\A$-module, and the images
of left preadmissible relations form its basis as a left $\A$-module.
\end{Lemma}

\begin{proof}
This is clear from the obvious isomorphisms
$$
\A\ox_{\F_0}R_\F\cong\bar R\cong R_\F\ox_{\F_0}\A
$$
of left, resp. right $\A$-modules.
\end{proof}

In particular we see that every element of $R_\F$ can be written 
uniquely in the form
\begin{equation}\label{norf}
\rho^{(2)}+\sum_i\alpha_i[n_i,m_i]\beta_i
\end{equation}
with $\rho^{(2)}\in R_\F\cdot R_\F$, $\alpha_i[n_i,m_i]\in\PAR$ and
$\beta_i$ an admissible monomial. Moreover it can be also uniquely written in
the form
\begin{equation}\label{norfl}
\ro^{(2)}+\sum_i\alpha'_i[n'_i,m'_i]\beta'_i
\end{equation}
with $\ro^{(2)}\in R_\F\cdot R_\F$, admissible monomials $\alpha'_i$  and
$[n'_i,m'_i]\beta'_i\in\PAR'$.

\begin{Definition}\label{lao}
The \emph{left action operator}
$$
L:\A\ox R_\F\to\A\ox\A
$$
of degree $-1$ is defined as follows. For odd $p$ let $L$ be the zero map. For $p=2$, let
first the additive map $L_\F:\F_0^{\le2}\to\A\ox\A$ be given by the formula
$$
L_\F(\Sq^n\Sq^m)=\sum_{\substack{n_1+n_2=n\\m_1+m_2=m\\\textrm{$m_1$, $n_2$ odd}}}
\Sq^{n_1}\Sq^{m_1}\ox\Sq^{n_2}\Sq^{m_2}
$$
($n,m\ge0$; remember that $\Sq^0=1$). Equivalently, using the algebra
structure on $\A\ox\A$ one may write
$$
L_\F(\Sq^n\Sq^m)=(1\ox\Sq^1)\delta(\Sq^{n-1})(\Sq^1\ox1)\delta(\Sq^{m-1}).
$$
Restricting this map to $R_\F^{\le2}\subset\F_0^{\le2}$ gives a map
$L_R:R_\F^{\le2}\to\A\ox\A$. It is thus an additive map given on the Adem
relations $[n,m]$, for $0<n<2m$, by
$$
L_R[n,m]=L_\F(\Sq^n\Sq^m)+\sum_{k=\max\{0,n-m+1\}}^{\min\{n/2,m-1\}}\binom{m-k-1}{n-2k}L_\F(\Sq^{n+m-k}\Sq^k).
$$
Next we define the map
$$
\bar L:\A\ox\bar R\to\A\ox\A
$$
as the right $\A$-module homomorphism which satisfies
\begin{equation}\label{lbar}
\bar L(a\ox\overline{\alpha[n,m]})=\delta(\k(a)\qf\alpha)L_R[n,m]
\end{equation}
with $\alpha[n,m]\in\PAR$; by \bref{barr} such a homomorphism exists
and is unique.

Finally, $\bar L$ yields a unique linear map $L:\A\ox R_\F\to\A\ox\A$ by
composing $\bar L$ with the quotient map $\A\ox R_\F\onto\A\ox\bar R$. Thus
one has
$$
L\left(\A\ox\left(R_\F\cdot R_\F\right)\right)=0.
$$
\end{Definition}

The map $L$ is the left action operator in \cite{Baues}*{14.4} where the
following lemma is proved (see \cite{Baues}*{14.4.3}):

\begin{Lemma}\label{lprop}
The map $\bar L$ satisfies the equalities
\begin{align*}
\bar L(a\ox[n,m])&=\k(a)L_R[n,m]\\
\bar L(a\ox br)&=\bar L(ab\ox r)+\delta(a)\bar L(b\ox r)\\
\bar L(a\ox rb)&=\bar L(a\ox r)\delta(b)
\end{align*}
for any $a,b\in\A$, $r\in\bar R$.
\end{Lemma}\qed

We observe that $L$ can be alternatively constructed as follows. Let
$$
\tilde L:\bar R\to\A\ox\A
$$
be the map given by
$$
\tilde L(\bar r)=\bar L(\Sq^1\ox\bar r).
$$
Then one has
\begin{Proposition}\label{ltil}
For any $a\in\A$, $r\in R_\F$ one has
$$
L(a\ox r)=\delta(\k(a))\tilde L(\bar r);
$$
moreover $\tilde L$ is a homomorphism of $\A$-$\A$-bimodules, hence
uniquely determined by its values on the Adem relations, which are
$$
\tilde L([n,m])=L_R[n,m].
$$
\end{Proposition}

\begin{proof}
For any $a\in\A$, $\alpha[n,m]\in\PAR$ and $\beta$ admissible we have
\begin{align*}
L(a\ox\alpha[n,m]\beta)
&=\bar L(a\ox\overline{\alpha[n,m]})\qf\beta\\
&=\delta(\k(a)\qf\alpha)L_R[n,m]\delta\qf\beta\\
&=\delta\k(a)\delta\qf\alpha L_R[n,m]\delta\qf\beta\\
&=\delta\k(a)\delta(\k(\Sq^1)\qf\alpha)L_R[n,m]\delta\qf\beta\\
&=\delta\k(a)L(\Sq^1\ox\alpha[n,m]\beta)\\
&=\delta\k(a)\tilde L(\alpha[n,m]\beta).
\end{align*}
Then using \eqref{norf} we see that the same identity holds for $L(a\ox r)$
with any $r\in R_\F$.

Next for any $a\in\A$, $r\in R_\F$ we have by \bref{lprop} and
$\k\Sq^1=\Sq^0=1$,
$$
\begin{aligned}
\tilde L(a\bar r)
&=\bar L(\Sq^1\ox a\bar r)\\
&=\bar L(\Sq^1a\ox\bar r)+\delta(\Sq^1)\bar L(a\ox\bar r)\\
&=\delta(\k(\Sq^1a))\tilde L(\bar r)+\delta(\Sq^1\k(a))\tilde L(\bar r)\\
&=\delta(\k(\Sq^1a)+\Sq^1\k(a))\tilde L(\bar r)\\
&=\delta(a)\tilde L(\bar r).
\end{aligned}
$$
Thus $\tilde L$ is a left $\A$-module homomorphism. It is also clearly a
right $\A$-module homomorphism since $\bar L$ is.

Finally by \eqref{lbar} we have
$$
L_R[n,m]=\delta(\k(\Sq^1))L_R[n,m]=\bar L(\Sq^1\ox\overline{[n,m]})=\tilde L([n,m]).
$$
\end{proof}

Explicit calculation of the left coaction operator $L_*$ is as follows.
For odd $p$ it is the zero map, and for $p=2$ we first define the additive
map ${L_R}_*:\A_*\ox\A_*\to{R_\F}_*^{\le2}$. It is dual to the composite map
$R_\F^{\le2}\to\A\ox\A$ in the diagram
\begin{equation}
\alignbox{
\xymatrix{
&R_\F^{\le2}\ar@{ >->}[r]\ar@{
>->}[d]\ar@{}[dr]|{\textrm{pull}}%\ar[dddddr]^{L_R}|<(.14)\hole
&R_\F\ar@{ >->}[d]
\\
\F_0^{\le1}\ox\F_0^{\le1}\ar[d]_{\Delta\ox\Delta}\ar@{->>}[r]^m
&\F_0^{\le2}\ar@{ >->}[r]\ar@{-->}[ddd]^{L_\F}
&\F_0
\\
\F_0^{\le1}\ox\F_0^{\le1}\ox\F_0^{\le1}\ox\F_0^{\le1}\ar[d]_{1\ox\Phi\ox\Phi\ox1}
\\
\F_0^{\le1}\ox\F_0^{\le1}\ox\F_0^{\le1}\ox\F_0^{\le1}\ar[d]_{1\ox T\ox1}
\\
\F_0^{\le1}\ox\F_0^{\le1}\ox\F_0^{\le1}\ox\F_0^{\le1}\ar@{->>}[r]^-{m\ox m}
&\F_0^{\le2}\ox\F_0^{\le2}\ar@{ >->}[r]
&\F_0\ox\F_0\ar@{->>}[d]
\\
&&\A\ox\A
}
}
\end{equation}
where $\Phi$ is restriction $\F_0^{\le1}\to\F_0^{\le1}$ of the map $\F_0\to\F_0$ given
by
$$
\Phi(x)=\Sq^1\k(x),
$$
so that one has
$$
\Phi(\Sq^n)=
\begin{cases}
\Sq^n,&n\equiv1\mod2\\
0,&n\equiv0\mod2.
\end{cases}
$$
Indeed by \bref{lao} we have
$$
L_\F(\Sq^n\Sq^m)=(1\ox\Sq^1)\Delta(\Sq^{n-1})(\Sq^1\ox1)\Delta(\Sq^{m-1})
=(1\ox\Sq^1)\Delta\k(\Sq^n)(\Sq^1\ox1)\Delta\k(\Sq^m);
$$
on the other hand we saw in \eqref{kapid} that
$$
\Delta\k=(\k\ox1)\Delta=(1\ox\k)\Delta,
$$
so that we can write
$$
L_\F(\Sq^n\Sq^m)=(1\ox\Sq^1\k)\Delta(\Sq^n)(\Sq^1\k\ox1)\Delta(\Sq^m)
=(1\ox\Phi)\Delta(\Sq^n)(\Phi\ox1)\Delta(\Sq^m).
$$

Therefore, the map dual of $\Phi$ is the map $\Phi_*:\FF[\zeta_1]\to\FF[\zeta_1]$
given by factorization through $\A_*\onto\FF[\zeta_1]$ of the map
$\Phi_*:\A_*\to\A_*$ given on the monomial basis by
$$
\Phi_*(\zeta_1^{n_1}\zeta_2^{n_2}\cdots)=
\begin{cases}
\zeta_1^{n_1}\zeta_2^{n_2}\cdots,&n_1\equiv1\mod2\\
0,&n_1\equiv0\mod2.
\end{cases}
$$
Equivalently, by \eqref{sqd} and \eqref{dkappa}, $\Phi_*=\k_*\Sq^1_*$ is
the map $\zeta_1\frac\partial{\partial\zeta_1}$.

Thus the map ${L_R}_*$ is the composite $\A_*\ox\A_*\to{R_\F^{\le2}}_*$ in the
diagram
\begin{equation}
\alignbox{
\xymatrix{
\A_*\ox\A_*\ar@{ >->}[d]\\
\F_*\ox\F_*\ar@{->>}[r]
&\F^{\le2}_*\ox\F^{\le2}_*\ar@{ >->}[r]^-{m_*\ox
m_*}\ar@{-->}[ddd]_{{L_\F}_*}
&\F^{\le1}_*\ox\F^{\le1}_*\ox\F^{\le1}_*\ox\F^{\le1}_*\ar[d]^{1\ox T\ox1}\\
&&\F^{\le1}_*\ox\F^{\le1}_*\ox\F^{\le1}_*\ox\F^{\le1}_*\ar[d]^{1\ox\Phi_*\ox\Phi_*\ox1}\\
&&\F^{\le1}_*\ox\F^{\le1}_*\ox\F^{\le1}_*\ox\F^{\le1}_*\ar[d]^{\Delta_*\ox\Delta_*}\\
\F_*\ar@{->>}[r]\ar@{->>}[d]\ar@{}[dr]|{\textrm{push}}
&\F^{\le2}_*\ar@{ >->}[r]^{m_*}\ar@{->>}[d]
&\F^{\le1}_*\ox\F^{\le1}_*\\
{R_\F}_*\ar@{->>}[r]
&{R_\F^{\le2}}_*
}
}
\end{equation}

Now by \bref{ltil} we know that $\tilde L$ is a bimodule homomorphism, and
moreover $\bar R$ is generated by $R^{\le2}_\F\cong\bar R^{\le2}\subset\bar
R$ as an $\A$-$\A$-bimodule, so knowledge of $L_R$ (actually already of
$L_\F$ whose restriction it is) determines $\tilde L$ and, by \bref{ltil},
also $L$. Dually, one can reconstruct $\tilde L_*$ and then $L_*$ from
${L_\F}_*$ via the diagram
$$
\xymatrix{
\A_*\ox\A_*\ar@{-->}[r]^{\tilde L_*}\ar[d]_{\textrm{bicoaction}}
&\bar R_*\ar@{ >->}[d]\ar[r]^-{\textrm{bicoaction}}
&\A_*\ox\bar R_*\ox\A_*\ar@{ >->}[d]\\
\A_*\ox(\A_*\ox\A_*)\ox\A_*\ar[r]^-{1\ox L_R\ox1}
&\A_*\ox{R_\F^{\le2}}_*\ox\A_*
&\A_*\ox{R_\F}_*\ox\A_*.\ar@{->>}[l]
}
$$
Here the bicoaction $\A_*\ox\A_*\to\A_*\ox(\A_*\ox\A_*)\ox\A_*$ is the
composite
$$
\xymatrix{
\A_*\ox\A_*\ar[r]^-{m^{(2)}_*\ox m^{(2)}_*}
&(\A_*\ox\A_*\ox\A_*)\ox(\A_*\ox\A_*\ox\A_*)\ar[d]^{(142536)}\\
&(\A_*\ox\A_*)\ox(\A_*\ox\A_*)\ox(\A_*\ox\A_*)\ar[r]^-{\delta_*\ox1\ox1\ox\delta_*}
&\A_*\ox(\A_*\ox\A_*)\ox\A_*
}
$$

We next note the following
\begin{Lemma}\label{bider}
The map $\tilde L_*$ is a biderivation, i.~e.
\begin{align*}
\tilde L_*(x_1x_2,y)&=x_1\tilde L_*(x_2,y)+x_2\tilde L_*(x_1,y),\\
\tilde L_*(x,y_1y_2)&=y_1\tilde L_*(x,y_2)+y_2\tilde L_*(x,y_1)
\end{align*}
for any $x,x_1,x_2,y,y_1,y_2\in\A_*$.
\end{Lemma}

It thus follows that $\tilde L_*$ is fully determined by its values
$\tilde L_*(\zeta_n\ox\zeta_{n'})$ on the Milnor generators.  To calculate
the bicoaction on these, first note that we have
$$
m^{(2)}_*(\zeta_n)=(1\ox m_*)m_*(\zeta_n)
=\sum_{i+i'=n}\zeta_i^{2^{i'}}\ox m_*(\zeta_{i'})
=\sum_{i+j+k=n}\zeta_i^{2^{j+k}}\ox\zeta_j^{2^k}\ox\zeta_k,
$$
where as always $\zeta_0=1$. For the coaction on $\zeta_n\ox\zeta_{n'}$
this then gives in succession
\begin{align*}
\zeta_n\ox\zeta_{n'}
&\mapsto\sum_{\substack{i+j+k=n\\i'+j'+k'=n'}}
\zeta_i^{2^{j+k}}\ox\zeta_j^{2^k}\ox\zeta_k\ox\zeta_{i'}^{2^{j'+k'}}\ox\zeta_{j'}^{2^{k'}}\ox\zeta_{k'}\\
&\mapsto\sum_{\substack{i+j+k=n\\i'+j'+k'=n'}}
\zeta_i^{2^{j+k}}\ox\zeta_{i'}^{2^{j'+k'}}\ox\zeta_j^{2^k}\ox\zeta_{j'}^{2^{k'}}\ox\zeta_k\ox\zeta_{k'}\\
&\mapsto\sum_{\substack{i+j+k=n\\i'+j'+k'=n'}}
\zeta_i^{2^{j+k}}\zeta_{i'}^{2^{j'+k'}}\ox\zeta_j^{2^k}\ox\zeta_{j'}^{2^{k'}}\ox\zeta_k\zeta_{k'},
\end{align*}
so that for the values of $\tilde L_*$ we have the equation
$$
\iota\tilde L_*(\zeta_n\ox\zeta_{n'})=\sum_{\substack{i+j+k=n\\i'+j'+k'=n'}}
\zeta_i^{2^{j+k}}\zeta_{i'}^{2^{j'+k'}}\ox{L_R}_*(\zeta_j^{2^k}\ox\zeta_{j'}^{2^{k'}})\ox\zeta_k\zeta_{k'}
$$
where $\iota$ is the above embedding $\bar
R_*\into\A_*\ox{R_\F^{\le2}}_*\ox\A_*$. Thus we only have to know the
values of ${L_\F}_*$ on the elements of the form
$\zeta_j^{2^k}\ox\zeta_{j'}^{2^{k'}}$ for $j\ge0$, $k\ge0$. Obviously
these values are zero for $j>2$ or $j'>2$. They are also zero for $j=0$ or
$j'=0$ since $\Phi_*(1)=0$. There thus remain four cases $j=j'=1$,
$j=j'=2$, $j=1$, $j'=2$, and $j=2$, $j'=1$. We then have under ${L_\F}_*$
$$
\xymatrixrowsep{.5em}
\xymatrix@!C=3em{
**[l]\zeta_1^{2^k}\ox\zeta_1^{2^{k'}}
\ar@{|->}[r]^{m_*\ox m_*}
&**[r](\zeta_1\ox1+1\ox\zeta_1)^{2^k}\ox(\zeta_1\ox1+1\ox\zeta_1)^{2^{k'}}=\\
**[l]{}&**[r]
\zeta_1^{2^k}\ox1\ox\zeta_1^{2^{k'}}\ox1
+\zeta_1^{2^k}\ox1\ox1\ox\zeta_1^{2^{k'}}
+1\ox\zeta_1^{2^k}\ox\zeta_1^{2^{k'}}\ox1
+1\ox\zeta_1^{2^k}\ox1\ox\zeta_1^{2^{k'}}\\
**[l]{}\ar@{|->}[r]^{1\ox T\ox1}
&**[r]\zeta_1^{2^k}\ox\zeta_1^{2^{k'}}\ox1\ox1
+\zeta_1^{2^k}\ox1\ox1\ox\zeta_1^{2^{k'}}
+1\ox\zeta_1^{2^{k'}}\ox\zeta_1^{2^k}\ox1
+1\ox1\ox\zeta_1^{2^k}\ox\zeta_1^{2^{k'}}\\
**[l]{}\ar@{|->}[r]^{1\ox\Phi_*\ox\Phi_*\ox1}
&**[r]0+0+1\ox\Phi_*\zeta_1^{2^{k'}}\ox\Phi_*\zeta_1^{2^k}\ox1+0\\
**[l]{}\ar@{|->}[r]^{\Delta_*\ox\Delta_*}
&**[r]\Phi_*\zeta_1^{2^{k'}}\ox\Phi_*\zeta_1^{2^k}.
}
$$
We thus have
$$
{L_\F}_*(\zeta_1^{2^k}\ox\zeta_1^{2^{k'}})=
\begin{cases}
M_{1,1},&k=k'=0\\
0&\textrm{otherwise.}
\end{cases}
$$
We next take $j=j'=2$; then
$$
\xymatrixrowsep{.5em}
\xymatrix@!C=3em{
**[l]\zeta_2^{2^k}\ox\zeta_2^{2^{k'}}
\ar@{|->}[r]^{m_*\ox m_*}
&**[r](\zeta_1^2\ox\zeta_1)^{2^k}\ox(\zeta_1^2\ox\zeta_1)^{2^{k'}}=
\zeta_1^{2^{k+1}}\ox\zeta_1^{2^k}\ox\zeta_1^{2^{k'+1}}\ox\zeta_1^{2^{k'}}\\
**[l]{}\ar@{|->}[r]^{1\ox T\ox1}
&**[r]\zeta_1^{2^{k+1}}\ox\zeta_1^{2^{k'+1}}\ox\zeta_1^{2^k}\ox\zeta_1^{2^{k'}}\\
**[l]{}\ar@{|->}[r]^{1\ox\Phi_*\ox\Phi_*\ox1}
&**[r]\zeta_1^{2^{k+1}}\ox\Phi_*\zeta_1^{2^{k'+1}}\ox\Phi_*\zeta_1^{2^k}\ox\zeta_1^{2^{k'}}=0\\
**[l]{}\ar@{|->}[r]^{\Delta_*\ox\Delta_*}
&**[r]0,
}
$$
so that
$$
{L_\F}_*(\zeta_2^{2^k}\ox\zeta_2^{2^{k'}})=0
$$
for all $k$ and $k'$. Next for $j=2$, $j'=1$ we have
$$
\xymatrixrowsep{.5em}
\xymatrix@!C=3em{
**[l]\zeta_2^{2^k}\ox\zeta_1^{2^{k'}}
\ar@{|->}[r]^{m_*\ox m_*}
&**[r](\zeta_1^2\ox\zeta_1)^{2^k}\ox(\zeta_1\ox1+1\ox\zeta_1)^{2^{k'}}=
\zeta_1^{2^{k+1}}\ox\zeta_1^{2^k}\ox\zeta_1^{2^{k'}}\ox1+\zeta_1^{2^{k+1}}\ox\zeta_1^{2^k}\ox1\ox\zeta_1^{2^{k'}}\\
**[l]{}\ar@{|->}[r]^{1\ox T\ox1}
&**[r]\zeta_1^{2^{k+1}}\ox\zeta_1^{2^{k'}}\ox\zeta_1^{2^k}\ox1+\zeta_1^{2^{k+1}}\ox1\ox\zeta_1^{2^k}\ox\zeta_1^{2^{k'}}\\
**[l]{}\ar@{|->}[r]^{1\ox\Phi_*\ox\Phi_*\ox1}
&**[r]\zeta_1^{2^{k+1}}\ox\Phi_*\zeta_1^{2^{k'}}\ox\Phi_*\zeta_1^{2^k}\ox1+0\\
**[l]{}\ar@{|->}[r]^{\Delta_*\ox\Delta_*}
&**[r]\zeta_1^{2^{k+1}}\Phi_*\zeta_1^{2^{k'}}\ox\Phi_*\zeta_1^{2^k},
}
$$
hence
$$
{L_\F}_*(\zeta_2^{2^k}\ox\zeta_1^{2^{k'}})=
\begin{cases}
M_{1,1}^2+M_1M_{2,1},&k=k'=0\\
0&\textrm{otherwise.}
\end{cases}
$$
Finally for $j=1$, $j'=2$ we get
$$
\xymatrixrowsep{.5em}
\xymatrix@!C=3em{
**[l]\zeta_1^{2^k}\ox\zeta_2^{2^{k'}}
\ar@{|->}[r]^{m_*\ox m_*}
&**[r](\zeta_1\ox1+1\ox\zeta_1)^{2^k}\ox(\zeta_1^2\ox\zeta_1)^{2^{k'}}=
\zeta_1^{2^k}\ox1\ox\zeta_1^{2^{k'+1}}\ox\zeta_1^{2^{k'}}+1\ox\zeta_1^{2^k}\ox\zeta_1^{2^{k'+1}}\ox\zeta_1^{2^{k'}}\\
**[l]{}\ar@{|->}[r]^{1\ox T\ox1}
&**[r]
\zeta_1^{2^k}\ox\zeta_1^{2^{k'+1}}\ox1\ox\zeta_1^{2^{k'}}+1\ox\zeta_1^{2^{k'+1}}\ox\zeta_1^{2^k}\ox\zeta_1^{2^{k'}}\\
**[l]{}\ar@{|->}[r]^{1\ox\Phi_*\ox\Phi_*\ox1}
&**[r]0+0\\
**[l]{}\ar@{|->}[r]^{\Delta_*\ox\Delta_*}
&**[r]0,
}
$$
so that
$$
{L_\F}_*(\zeta_1^{2^k}\ox\zeta_2^{2^{k'}})=0
$$
for all $k$ and $k'$.

To pass to ${L_R}_*$ from these values means just omitting all monomials
which do not contain $M_{1,1}$; we thus obtain
\begin{align*}
{L_R}_*(\zeta_1\ox\zeta_1)&=M_{1,1},\\
{L_R}_*(\zeta_2\ox\zeta_1)&=M_{1,1}^2,
\end{align*}
and ${L_R}_*(\zeta_j^{2^k}\ox\zeta_{j'}^{2^{k'}})=0$ in all other cases.

From this we easily obtain
\begin{Proposition}
$\iota\tilde L_*(\zeta_n\ox\zeta_{n'})=\zeta_{n-1}^2\zeta_{n'-1}^2\ox
M_{1,1}\ox1+\zeta_{n-2}^4\zeta_{n'-1}^2\ox M_{1,1}^2\ox1
$
\end{Proposition}\qed

\noindent where now $\zeta_{n-2}=0$ for $n=1$ is understood. Solving $\tilde
L_*(\zeta_n,\zeta_{n'})$ from these equations is then straightforward. In
this way we obtain
\begin{align*}
\tilde L_*(\zeta_1,\zeta_1)&=M_{1,1}\\
\tilde L_*(\zeta_1,\zeta_2)&=M_{2,1,1}\\
\tilde L_*(\zeta_2,\zeta_1)&=M_{2,1,1}+M_{1,1}^2\\
\tilde L_*(\zeta_2,\zeta_2)&=M_{4,1,1}+M_{2,3,1}+M_{2,1,2,1}\\
\tilde L_*(\zeta_1,\zeta_3)&=M_{4,2,1,1}\\
\tilde L_*(\zeta_3,\zeta_1)&=M_{4,2,1,1}+M_{2,1,1}^2\\
\tilde L_*(\zeta_2,\zeta_3)&=M_{6,2,1,1}+M_{4,4,1,1}+M_{4,2,3,1}+M_{4,2,1,2,1}+M_{2,4,2,1,1}\\
\tilde L_*(\zeta_3,\zeta_2)&=M_{6,2,1,1}+M_{4,4,1,1}+M_{4,2,3,1}+M_{4,2,1,2,1}+M_{2,4,2,1,1}\\
                           &+M_5^2+M_{4,1}^2+M_{3,2}^2+M_{2,1,1,1}^2+M_1^2M_{2,1,1}^2+M_1^4M_3^2\\
\tilde L_*(\zeta_3,\zeta_3)&=M_{8,4,1,1}+M_{8,2,3,1}+M_{8,2,1,2,1}+M_{4,6,3,1}+M_{4,6,1,2,1}+M_{4,2,5,2,1}+M_{4,2,4,3,1}+M_{4,2,4,1,2,1}\\
                           &+M_{4,2,1,4,2,1}\\
\tilde L_*(\zeta_1,\zeta_4)&=M_{8,4,2,1,1}\\
\tilde L_*(\zeta_4,\zeta_1)&=M_{8,4,2,1,1}+M_{4,2,1,1}^2\\
\tilde L_*(\zeta_2,\zeta_4)&=M_{10,4,2,1,1}+M_{8,6,2,1,1}+M_{8,4,4,1,1}+M_{8,4,2,3,1}+M_{8,4,2,1,2,1}+M_{8,2,4,2,1,1}+M_{2,8,4,2,1,1}\\
\tilde L_*(\zeta_4,\zeta_2)&=M_{10,4,2,1,1}+M_{8,6,2,1,1}+M_{8,4,4,1,1}+M_{8,4,2,3,1}+M_{8,4,2,1,2,1}+M_{8,2,4,2,1,1}+M_{2,8,4,2,1,1}\\
                           &+M_9^2+M_{7,2}^2+M_{5,4}^2+M_{6,2,1}^2+M_{4,4,1}^2+M_{4,3,2}^2+M_{4,2,1,1,1}^2+M_{3,4,2}^2+M_{2,4,2,1}^2\\
                           &+M_1^2M_{4,2,1,1}^2+M_1^8M_5^2+M_{2,1}^4M_3^2\\
\tilde L_*(\zeta_3,\zeta_4)&=M_{12,6,2,1,1}+M_{12,4,4,1,1}+M_{12,4,2,3,1}+M_{12,4,2,1,2,1}+M_{12,2,4,2,1,1}+M_{8,8,4,1,1}+M_{8,8,2,3,1}\\
                           &+M_{8,8,2,1,2,1}+M_{8,4,6,3,1}+M_{8,4,6,1,2,1}+M_{8,4,2,5,2,1}+M_{8,4,2,4,3,1}+M_{8,4,2,4,1,2,1}+M_{8,4,2,1,4,2,1}\\
                           &+M_{4,10,4,2,1,1}+M_{4,8,6,2,1,1}+M_{4,8,4,4,1,1}+M_{4,8,4,2,3,1}+M_{4,8,4,2,1,2,1}+M_{4,8,2,4,2,1,1}+M_{4,2,8,4,2,1,1},
\end{align*}
etc.

Having $\tilde L_*$ we then can obtain $L_*$ by the dual of \bref{ltil} as
\begin{equation}\label{L}
L_*(x,y)=\sum\zeta_1x_\l y_{\l'}\ox\tilde L_*(x_\r,y_{\r'})
\end{equation}
for $x,y\in\A_*$, with
$$
m_*(x)=\sum x_\l\ox x_\r,\ m_*(y)=\sum y_{\l'}\ox y_{\r'}.
$$

\section{The symmetry operator $S$ and its dual}\label{S*}

\begin{Definition}\label{so}
The \emph{symmetry operator}
$$
S:R_\F\to\A\ox\A
$$
of degree $-1$ is defined as follows. For odd $p$, let $S$ be the zero map.
For $p=2$ let the elements $S_n\in\A\ox\A$, $n\ge0$, be given by
$$
S_n=\sum_{\substack{n_1+n_2=n-1\\\textrm{$n_1$, $n_2$
odd}}}\Sq^{n_1}\ox\Sq^{n_2}=(\Sq^1\ox\Sq^1)\delta(\Sq^{n-3}),
$$
i.~e.
\begin{align*}
S_{2k}&=0,\\
S_{2k+1}&=\sum_{0\le i<k}\Sq^{2i+1}\ox\Sq^{2(k-i)-1},
\end{align*}
$k\ge0$. Then let the linear map $S_\F:\F_0^{\le2}\to\A\ox\A$ be given by
\begin{align*}
S_\F(\Sq^n\Sq^m)&=S_n\delta(\Sq^m)+\delta(\Sq^n)S_m+\delta(\Sq^{n-1})S_{m+1}\\
&=(\Sq^1\ox\Sq^1)\delta(\Sq^{n-3}\Sq^m)+\delta(\Sq^n)(\Sq^1\ox\Sq^1)\delta(\Sq^{m-3})+\delta(\Sq^{n-1})(\Sq^1\ox\Sq^1)\delta(\Sq^{m-2}),
\end{align*}
$n,m\ge0$. Next define the map $S_R:R_\F^{\le2}\to\A\ox\A$ by restriction to
$R_\F^{\le2}\subset\F_0^{\le2}$. Thus on the Adem relations this map is
given by
\begin{equation}\label{sr}
S_R[n,m]=S_\F(\Sq^n\Sq^m)+\sum_{k=\max\{0,n-m+1\}}^{\min\{n/2,m-1\}}\binom{m-k-1}{n-2k}S_\F(\Sq^{n+m-k}\Sq^k).
\end{equation}
Now let us define the map
$$
\bar S:\bar R\to\A\ox\A
$$
as a unique right $\A$-module homomorphism satisfying
$$
\bar S(\overline{\alpha[n,m]})=\delta(\alpha)S_R[n,m]+(1+T)\bar
L(\alpha\ox\overline{[n,m]})
$$
for $\alpha[n,m]\in\PAR$. Then finally this determines a unique linear map $S:R_\F\to\A\ox\A$ by
composing with the quotient map $R_\F\onto\bar R$.
\end{Definition}

The map $S$ is the symmetry operator in \cite{Baues}*{14.5.2} where the
following lemma is proved.

\begin{Lemma}
The map $\bar S$ satisfies the equations
\begin{align*}
\bar S([n,m])&=S_R[n,m]\\
\bar S(ar)&=\delta(a)\bar S(r)+(1+T)\bar L(a\ox r)\\
\bar S(ra)&=\bar S(r)\delta(a)
\end{align*}
for any $0<n<2m$, $a\in\A$ and $r\in\bar R$.
\end{Lemma}

We now turn to the dual $S_*:\A_*\ox\A_*\to{R_\F}_*$ of $S$ (dually to the
above, the image of this operator actually lies in $\bar
R_*\subset{R_\F}_*$ and so defines the operator $\bar
S_*:\A_*\ox\A_*\to\bar R_*$). Since we know that $S_*$ is a biderivation,
it suffices to compute the values $S_*(\zeta_n\ox\zeta_{n'})$. Now
dually to the equation
$$
S(a[n,m]b)
=\delta(a)S_R([n,m])\delta(b)+(1+T)L(a\ox[n,m]b)
=\delta(a)S_R([n,m])\delta(b)+(1+T)\left(\delta\k(a)L_R([n,m])\delta(b)\right)
$$
we have
\begin{multline*}
\iota S_*(\zeta_n\ox\zeta_{n'})=\\
\sum_{\substack{i+j+k=n\\i'+j'+k'=n'}}
\left(
\zeta_i^{2^{j+k}}\zeta_{i'}^{2^{j'+k'}}
\ox{S_R}_*(\zeta_j^{2^k}\ox\zeta_{j'}^{2^{k'}})\ox\zeta_k\zeta_{k'}
+\zeta_1\zeta_i^{2^{j+k}}\zeta_{i'}^{2^{j'+k'}}
\ox\left({L_R}_*(\zeta_j^{2^k}\ox\zeta_{j'}^{2^{k'}})
+{L_R}_*(\zeta_{j'}^{2^{k'}}\ox\zeta_j^{2^k})\right)\ox\zeta_k\zeta_{k'}\right)\\
=
\sum_{\substack{i+j+k=n\\i'+j'+k'=n'}}
\zeta_i^{2^{j+k}}\zeta_{i'}^{2^{j'+k'}}
\ox{S_R}_*(\zeta_j^{2^k}\ox\zeta_{j'}^{2^{k'}})\ox\zeta_k\zeta_{k'}
+
\zeta_1\zeta_{n-2}^4\zeta_{n'-1}^2
\ox M_{1,1}^2\ox1
+
\zeta_1\zeta_{n-1}^2\zeta_{n'-2}^4
\ox M_{1,1}^2\ox1,
\end{multline*}
with $\zeta_0=1$ and $\zeta_n=0$ for $n<0$, as before.

It thus remains to find the values
${S_R}_*(\zeta_j^{2^k}\ox\zeta_{j'}^{2^{k'}})$ --- which in turn are
images of the corresponding values of ${S_\F}_*$ under the map
$\F_*\onto{R_\F}_*$. To find the latter, let us first define another
intermediate operator
$$
S^1:\F_0^{\le1}\to\A\ox\A
$$
by the equation
$$
S^1(\Sq^n)=S_{n+1}=(\Sq^1\ox\Sq^1)\delta\k\k(\Sq^n)=\sum_{\substack{n_1+n_2=n\\\textrm{$n_1$, $n_2$
odd}}}\Sq^{n_1}\ox\Sq^{n_2},
$$
so that we have
$$
S_\F m(\Sq^n\ox\Sq^m)=S_\F(\Sq^n\Sq^m)=S^1\k(\Sq^n)\delta(\Sq^m)+\delta(\Sq^n)S^1\k(\Sq^m)+\delta\k(\Sq^n)S^1(\Sq^m).
$$
We have the dual operator
$$
S^1_*:\A_*\ox\A_*\to\F_*^{\le1}
$$
such that dual
$$
{S_\F}_*:\A_*\ox\A_*\to\F_*^{\le2}
$$
of $S_\F$ is given by
\begin{equation}\label{msf}
\alignbox{
&m_*{S_\F}_*(x\ox y)=\\
&\sum\left(\zeta_1S^1_*(x_\l\ox y_{\l'})\ox(x_\r y_{\r'})^{\le1}
+(x_\l y_{\l'})^{\le1}\ox\zeta_1S^1_*(x_\r\ox y_{\r'})
+(\zeta_1x_\l y_{\l'})^{\le1}\ox S^1_*(x_\r\ox y_{\r'})\right)
}
\end{equation}
where as before we use the Sweedler notation
$$
m_*(x)=\sum x_\l\ox x_\r,\ \ m_*(y)=\sum y_{\l'}\ox y_{\r'}
$$
and
$$
(\_)^{\le1}:\A_*\to\F_*^{\le1}
$$
sends $\zeta_1$ to $M_1$ and all other Milnor generators to 0. Thus we have
\begin{multline*}
m_*{S_\F}_*(\zeta_j^{2^k}\ox\zeta_{j'}^{2^{k'}})=\\
\sum_{\substack{\l+\r=j\\\l'+\r'=j'}}
\zeta_1S^1_*(\zeta_\l^{2^{\r+k}}\ox\zeta_{\l'}^{2^{\r'+k'}})\ox(\zeta_\r^{2^k}\zeta_{\r'}^{2^{k'}})^{\le1}
+(\zeta_\l^{2^{\r+k}}\zeta_{\l'}^{2^{\r'+k'}})^{\le1}\ox\zeta_1S^1_*(\zeta_\r^{2^k}\ox\zeta_{\r'}^{2^{k'}})
+(\zeta_1\zeta_\l^{2^{\r+k}}\zeta_{\l'}^{2^{\r'+k'}})^{\le1}\ox
S^1_*(\zeta_\r^{2^k}\ox\zeta_{\r'}^{2^{k'}})
\end{multline*}

Now the operator $S^1_*$ is obviously given by
\begin{equation}\label{s1}
S^1_*(x\ox y)=
\begin{cases}
xy,&\textrm{$x=\zeta_1^{n_1}$, $y=\zeta_1^{n_2}$, $n_1$, $n_2$ odd},\\
0&\textrm{otherwise,}
\end{cases}
\end{equation}
so that ${S_\F}_*(\zeta_j^{2^k}\ox\zeta_{j'}^{2^{k'}})=0$ whenever $k>0$ or
$k'>0$. And among the remaining values ${S_\F}_*(\zeta_j\ox\zeta_{j'})$ the
only nonzero ones are given by
\begin{align*}
{S_\F}_*(\zeta_1\ox\zeta_1)&=M_3+M_{1,2}=M_1^3+M_{2,1},\\
{S_\F}_*(\zeta_1\ox\zeta_2)={S_\F}_*(\zeta_2\ox\zeta_1)&=M_{2,3}+M_{3,2}=M_1M_{1,1}^2,\\
{S_\F}_*(\zeta_2\ox\zeta_2)&=M_{5,2}+M_{4,3}=M_1M_{2,1}^2.
\end{align*}
Then further passing to ${S_R}_*$ means, as before, removing the monomials
not containing $M_{1,1}$, so that the only nonzero values of the form
${S_R}_*(\zeta_j^{2^k}\ox\zeta_{j'}^{2^{k'}})$ are
$$
{S_R}_*(\zeta_1\ox\zeta_2)={S_R}_*(\zeta_2\ox\zeta_1)=M_1M_{1,1}^2.
$$
Hence we obtain
\begin{Proposition}
$$
\iota S_*(\zeta_n\ox\zeta_{n'})=
\zeta_{n-1}^2\zeta_{n'-2}^4
\ox M_1M_{1,1}^2\ox1
+
\zeta_{n-2}^4\zeta_{n'-1}^2
\ox M_1M_{1,1}^2\ox1
+
\zeta_1\zeta_{n-2}^4\zeta_{n'-1}^2
\ox M_{1,1}^2\ox1
+
\zeta_1\zeta_{n-1}^2\zeta_{n'-2}^4
\ox M_{1,1}^2\ox1.
$$
\end{Proposition}\qed

As for $\tilde L_*$ above, we then solve these equations obtaining e.~g.
\begin{align*}
                     S_*(\zeta_1,\zeta_1)&=0,\\
S_*(\zeta_1,\zeta_2)=S_*(\zeta_2,\zeta_1)&=M_{2,2,1}+M_1M_{1,1}^2,\\
                     S_*(\zeta_2,\zeta_2)&=0,\\
S_*(\zeta_1,\zeta_3)=S_*(\zeta_3,\zeta_1)&=M_{4,2,2,1}+M_1M_{2,1,1}^2,\\
S_*(\zeta_2,\zeta_3)=S_*(\zeta_3,\zeta_2)&=M_{6,2,2,1}+M_{4,4,2,1}+M_{2,4,2,2,1}\\
                                         &+M_1M_5^2+M_1M_{4,1}^2+M_1M_{3,2}^2+M_1M_{2,1,1,1}^2+M_1^3M_{2,1,1}^2+M_1^5M_3^2,\\
                     S_*(\zeta_3,\zeta_3)&=0,\\
S_*(\zeta_1,\zeta_4)=S_*(\zeta_4,\zeta_1)&=M_{8,4,2,2,1}+M_1M_{4,2,1,1}^2,\\
S_*(\zeta_2,\zeta_4)=S_*(\zeta_4,\zeta_2)&=M_{10,4,2,2,1}+M_{8,6,2,2,1}+M_{8,4,4,2,1}+M_{8,2,4,2,2,1}+M_{2,8,4,2,2,1}\\                                        
                                         &+M_1M_9^2+M_1M_{7,2}^2+M_1M_{6,2,1}^2+M_1M_{5,4}^2+M_1M_{4,4,1}^2+M_1M_{4,3,2}^2+M_1M_{4,2,1,1,1}^2\\
                                         &+M_1M_{3,4,2}^2+M_1M_{2,4,2,1}^2+M_1M_{2,1}^2M_3^2+M_1^3M_{4,2,1,1}^2+M_1^9M_5^2,
\end{align*}
etc.

%% file: dHascoAss_c7.tex
\chapter{The extended Steenrod algebra and its cocycle}\label{xi}

We show that the dual invariant $S_*$ determines a singular extension of the
Hopf algebra structure of the Steenrod algebra. We also give a formula for a
cocycle representing the extension. Then we show that $S_*$ is related to a
formula which describes the main result of Kristensen on secondary cohomology
operations. A proof of this formula has not appeared in the literature yet.

\section{Singular extensions of Hopf algebras}

In this section we introduce a singular extension $\hat\A$ of the Steenrod
algebra $\A$ which is determined by the symmetry operator $S$.

\begin{Definition}
A \emph{singular extension} of a Hopf algebra $A$ is a direct sum diagram
$$
\xymatrix{
R\ar@{ >->}@<.5ex>[r]^-i
&\hat A\ar@{->>}@<.5ex>[r]^-p\ar@{->>}@<.5ex>[l]^-q
&A,\ar@{ >->}@<.5ex>[l]^-s
}
$$
i.~e. one has $ps=\id_A$, $qi=\id_R$ and $sp+iq=\id_{\hat A}$, such that
$\hat A$ is an algebra with multiplication $\mu:\hat A\ox\hat A\to\hat A$
and $\hat A$ is also a coalgebra with diagonal $\hat\delta:\hat A\to\hat
A\ox\hat A$. (Here we do not assume that $\hat\delta$ is a homomorphism of
algebras, or equivalently that $\mu$ is a homomorphism of coalgebras, so
that in general $\hat A$ is not a Hopf algebra). In addition $p$ is an algebra
homomorphism, and $s$ is a coalgebra homomorphism. Moreover $(i,p)$ must
be a singular extension of algebras and $(q,s)$ must be a singular
extension of coalgebras. This means that the ideal $R=\ker i$ of the
algebra $\hat A$ is a square zero ideal, i.~e. $xy=0$ for any $x,y\in
R$, and the coideal $R=\coker s$ of the coalgebra $\hat A$ is a
square zero coideal, i.~e. the composite
$$
\hat A\xto{\hat\delta}\hat A\ox\hat A\xto{q\ox q}R\ox R
$$
is zero.
\end{Definition}

It follows that the $\hat A$-$\hat A$-bimodule and $\hat A$-$\hat
A$-bicomodule structures on $R$ descend to $A$-$A$-bimodule and
$A$-$A$-bicomodule structures respectively.

Our basic example of a singular Hopf algebra extension is as follows. We
have seen that $\bar R$ from \eqref{rbar} has an $\A$-$\A$-bimodule
structure. Now it also has an $\A$-$\A$-bicomodule structure as follows. On
the one hand, there is a diagonal $\Delta_R:R_\F\to R^{(2)}_\F=\ker(q_\F\ox
q_\F)$ induced in the commutative diagram
$$
\xymatrix{
R_\F\ar@{ >->}[r]\ar@{-->}[d]_\Delta
&\F_0\ar@{->>}[r]^{q_\F}\ar[d]^\Delta
&\A\ar[d]^\delta\\
R^{(2)}_\F\ar@{ >->}[r]
&\F_0\ox\F_0\ar@{->>}[r]^{q_\F\ox q_\F}
&\A\ox\A
}
$$
with short exact rows. Moreover there is a short exact sequence
$$
\xymatrix{
R_\F\ox R_\F\ar@{ >->}[r]^-{i^{(2)}=\binom{i_\F\ox1}{-1\ox i_\F}}
&\F_0\!\ox\!R_\F\oplus R_\F\!\ox\!\F_0\ar@{->>}[r]
&R_\F^{(2)},
}
$$
where $i_\F:R_\F\incl\F_0$ is the inclusion. Since the composite of the
quotient map
$$
\F_0\!\ox\!R_\F\oplus R_\F\!\ox\!\F_0\onto\A\!\ox\!\bar R\oplus\bar R\!\ox\!\A
$$
with $i^{(2)}$ is obviously zero, we get the induced map
$$
R_\F^{(2)}\to\A\!\ox\!\bar R\oplus\bar R\!\ox\!\A.
$$
Moreover the diagonal of $\F_0$ factors through this map as follows
\begin{equation}\label{bicomod}
\alignbox{
\xymatrix{
\bar R\ar@{-->}[d]^{\binom{\Delta_\l}{\Delta_\r}}
&R_\F\ar@{->>}[l]\ar@{ (->}[r]^{i_\F}\ar@{-->}[d]^\Delta
&\F_0\ar@{->>}[r]^{q_\F}\ar[d]^\Delta
&\A\ar[d]^\delta\\
\A\!\ox\!\bar R\oplus\bar R\!\ox\!\A
&R^{(2)}\ar@{->>}[l]\ar@{ (->}[r]
&\F_0\ox\F_0\ar@{->>}[r]^{q_\F\ox q_\F}
&\A\ox\A
}
}
\end{equation}
giving the left, resp. right coaction $\Delta_\l$, resp. $\Delta_\r$ of the
desired $\A$-$\A$-bicomodule structure on $\bar R$.

Note that the above construction is actually precisely dual to the
standard procedure for equipping the kernel of a singular extension with a
structure of a bimodule over a base. In particular we could use the dual
diagram
\begin{equation}\label{bimod}
\alignbox{
\xymatrix{
\A\!\ox\!\bar R\oplus\bar R\!\ox\!\A\ar@{-->}[d]^{\binom{m_\l}{m_\r}}
&R^{(2)}\ar@{->>}[l]\ar@{ (->}[r]\ar@{-->}[d]^m
&\F_0\ox\F_0\ar@{->>}[r]^{q_\F\ox q_\F}\ar[d]^m
&\A\ox\A\ar[d]^m\\
\bar R
&R_\F\ar@{->>}[l]\ar@{ (->}[r]^{i_\F}
&\F_0\ar@{->>}[r]^{q_\F}
&\A
}
}
\end{equation}
to give $\bar R$ via $m_\l$ and $m_\r$ the structure of $\A$-$\A$-bimodule.

\begin{Theorem}\label{usplit}
There is a unique singular extension of Hopf algebras
$$
\xymatrix{
\Sigma^{-1}\bar R\ar@{ >->}@<.5ex>[r]^-i
&\hat\A\ar@{->>}@<.5ex>[r]^-p\ar@{->>}@<.5ex>[l]^-q
&\A,\ar@{ >->}@<.5ex>[l]^-s
}
$$
where $\hat\A$ is the split singular extension of algebras, that is, as an
algebra
$$
\hat\A=\A\oplus\Sigma^{-1}\bar R
$$
is the semidirect product with multiplication
$$
(a,r)(a',r')=(aa',ar'+ra')
$$
and the following conditions are satisfied.

The induced $\A$-$\A$-bimodule and $\A$-$\A$-bicomodule structures on
$\Sigma^{-1}\bar R$ are given by the ones
indicated in \eqref{bicomod} above, and the diagonal $\hat\delta$ of
the coalgebra $\hat\A$ fits into the commutative diagram
\begin{equation}\label{coext}
\alignbox{
\xymatrix{
\hat\A\ar[rr]^{\hat\delta}\ar@{->>}[d]
&&\hat\A\ox\hat\A\ar[d]^{1+T}\\
\Sigma^{-1}\bar R\ar[r]^-S
&\A\ox\A\ar@{ (->}[r]
&\hat\A\ox\hat\A}
}
\end{equation}
where $S$ is the symmetry operator in \bref{so}.
\end{Theorem}

We will prove this theorem together with the dual statement. Note that
clearly the dual of a singular extension of any Hopf algebra $A$ is a
singular extension of the dual Hopf algebra $A_*$. Clearly then the above
theorem is equivalent to

\begin{Theorem}\label{dusplit}
There is a unique singular extension of Hopf algebras
$$
\xymatrix{
\Sigma^{-1}\bar R_*\ar@{ >->}@<.5ex>[r]^-{q_*}
&\hat\A_*\ar@{->>}@<.5ex>[r]^-{s_*}\ar@{->>}@<.5ex>[l]^-{i_*}
&\A_*,\ar@{ >->}@<.5ex>[l]^-{p_*}
}
$$
where $\hat\A_*$ is the split singular extension of coalgebras, that is, as
a coalgebra
$$
\hat\A_*=\A_*\oplus\Sigma^{-1}\bar R_*
$$
with diagonal
$$
\A_*\oplus\Sigma^{-1}\bar R_*\xto{
\left(\begin{smallmatrix}
m_*&0\\
0&{m_\l}_*\\
0&{m_\r}_*\\
0&0
\end{smallmatrix}\right)
}
\A_*\!\ox\!\A_*\oplus
\A_*\!\ox\!\Sigma^{-1}\bar R_*\oplus
\Sigma^{-1}\bar R_*\!\ox\!\A_*\oplus
\Sigma^{-1}\bar R_*\!\ox\!\Sigma^{-1}\bar R_*
$$
where the diagonal $m_*$ is dual to the multiplication $m:\A\ox\A\to\A$
and ${m_\l}_*$, ${m_\r}_*$ are the $\A_*$-$\A_*$-bicomodule structure maps
dual to the $\A$-$\A$-bimodule structure maps $m_\l:\A\ox\Sigma^{-1}\bar
R\to\Sigma^{-1}\bar R$, $m_\r:\Sigma^{-1}\bar R\ox\A\to\Sigma^{-1}\bar R$
in \eqref{bimod}, where the induced $\A_*$-$\A_*$-bimodule structure on
$\bar R_*$ is dual to the $\A$-$\A$-bicomodule structure indicated in
\eqref{bicomod} above, and where the multiplication $\hat\delta_*$ of the
algebra $\hat\A_*$ satisfies the commutation rule
$$
p_*(y)p_*(x)=p_*(x)p_*(y)+S_*(x\ox y)
$$
for any $x,y\in\A_*$, where
$$
S_*:\A_*\ox\A_*\to\Sigma^{-1}\bar R_*
$$
is the cosymmetry operator from \bref{cosym}.
\end{Theorem}

\begin{proof}[Proof of \bref{usplit} and \bref{dusplit}]
The diagonal $\hat\delta$ can be written as follows
$$
\A\oplus\Sigma^{-1}\bar R\xto{
\left(\begin{smallmatrix}
\phi_{11}&\phi_{12}\\
\phi_{21}&\phi_{22}\\
\phi_{31}&\phi_{32}\\
\phi_{41}&\phi_{42}
\end{smallmatrix}\right)
}
\A\!\ox\!\A\oplus
\A\!\ox\!\Sigma^{-1}\bar R\oplus
\Sigma^{-1}\bar R\!\ox\!\A\oplus
\Sigma^{-1}\bar R\!\ox\!\Sigma^{-1}\bar R.
$$
Then the condition that $s:\A\into\A\oplus\Sigma^{-1}\bar R$ is a coalgebra
homomorphism implies $\phi_{11}=\delta$ and $\phi_{21}=0$, $\phi_{31}=0$,
$\phi_{41}=0$. Moreover the condition that the $\A$-$\A$-bicomodule
structure induced on $\Sigma^{-1}\bar R$ coincides with the one given in
\eqref{bicomod} implies $\phi_{22}=\Delta_\l$, $\phi_{32}=\Delta_\r$. Next the
condition that $(s,q)$ is a singular extension of coalgebras, i.~e. the
coideal $\bar R$ has zero comultiplication, implies $\phi_{42}=0$. Finally,
let us look at the diagram \eqref{coext}. The lower composite in this
diagram sends $(a,r)\in\A\oplus\Sigma^{-1}\bar R$ to
$$
(S(r),0,0,0)\in\A\!\ox\!\A\oplus
\A\!\ox\!\Sigma^{-1}\bar R\oplus
\Sigma^{-1}\bar R\!\ox\!\A\oplus
\Sigma^{-1}\bar R\!\ox\!\Sigma^{-1}\bar R.
$$
The upper composite sends it to
\begin{align*}
(1+T)\hat\delta(a,r)&=(1+T)(\delta(a)+\phi_{12}(r),\Delta_\l(r),\Delta_\r(r),0)\\
&=((1+T)\delta(a)+(1+T)\phi_{12}(r),\Delta_\l(r)+T\Delta_\r(r),\Delta_\r(r)+T\Delta_\l(r),0).
\end{align*}
Since $\delta$ is cocommutative, one has $(1+T)\delta=0$. Moreover
cocommutativity of $\Delta:\F_0\to\F_0\ox\F_0$ implies $T\Delta_\l=\Delta_\r$,
$T\Delta_\r=\Delta_\l$. Thus commutativity of \eqref{coext} is equivalent to the
condition
\begin{equation}\label{xis}
(1+T)\phi_{12}=S:\Sigma^{-1}\bar R\to\A\ox\A.
\end{equation}
Equivalently, passing to the dual we see that the dual map
$\xi_*={\phi_{12}}_*:\A_*\ox\A_*\to\Sigma^{-1}\bar R_*$ must satisfy
$$
\xi_*(1+T)=S_*.
$$
Now it is easy to see that $\xi_*$ is in fact the algebra cocycle
determining the algebra extension
$$
\xymatrix{
\bar R_*\ar@{ >->}[r]^{q_*}&\hat\A_*\ar@{->>}[r]^{s_*}&\A_*,
}
$$
that is, in $\hat\A_*=\A_*\oplus\Sigma^{-1}\bar R_*$ one has
$$
(\alpha,\beta)(\alpha',\beta')=(\alpha\alpha',\alpha\beta'+\beta\alpha'+\xi_*(\alpha\ox\alpha')).
$$
Hence by \eqref{xis} one has
$$
(\alpha,\beta)(\alpha',\beta')-(\alpha',\beta')(\alpha,\beta)=(0,S_*(\alpha\ox\alpha')).
$$
Now recall that $\A_*$ is actually a polynomial algebra. Using this fact it
has been shown in \cite{Baues}*{16.2} that the algebra structure of any of
its singular extensions such as $\hat\A_*$ above is completely determined
by its commutator map, i.~e. by $S_*$. Thus ${\phi_{12}}_*$ and hence the
whole $\phi_{ij}$ matrix is uniquely determined. It is then straightforward
to check that indeed this matrix yields a coalgebra structure on $\hat\A$
with desired properties.
\end{proof}

It follows immediately from \bref{dusplit} (and actually this was also
deduced during its proof) that one has

\begin{Corollary}\label{cosxi}
For the cosymmetry operator $S_*$ from \bref{cosym} there exists a map
$$
\xi_*:\A_*\ox\A_*\to\Sigma\1\bar R_*
$$
which is a 2-cocycle, i.~e. for any $x,y,z\in\A_*$ one has
$$
x\xi_*(y,z)+\xi_*(x,yz)=z\xi_*(x,y)+\xi_*(xy,z)
$$
and such that its symmetrization is equal to $S_*$, i.~e. for any
$x,y\in\A_*$ one has
$$
\xi_*(x,y)+\xi_*(y,x)=S_*(x,y).
$$
\end{Corollary}

\begin{proof}
This follows since any extension
$$
\xymatrix{
M\ar@{ >->}[r]^i&A'\ar@{->>}[r]^p&A
}
$$
of a commutative algebra $A$ by a symmetric $A$-module $M$ is determined by a 2-cocycle
$c:A\ox A\to M$ such that for any $x,y\in A'$ one has
$$
xy-yx=i(c(px,py)-c(py,px)),
$$
i.~e. the commutator map for $A'$ is given by the antisymmetrization of
$c$. Of course for $p=2$ there is no difference between symmetrization and
antisymmetrization.
\end{proof}

\begin{Remark}
The above corollary is easily seen to be exactly dual to
\cite{Baues}*{Theorem 16.1.5}. 
\end{Remark}

Using the extended Steenrod algebra we can next compute the deviation of
the cocycle $\xi_*$ from being an $\A_*$-comodule homomorphism. Namely, let
$$
{\nabla_\xi}_*:\A_*\ox\A_*\to\A_*\ox\Sigma\1\bar R_*
$$
be the difference between the upper and lower composites in the diagram
\begin{equation}\label{nabladef}
\begin{aligned}
\xymatrix{
\A_*\ox\A_*\ar[d]_{\xi_*}\ar[r]^-{\textrm{coaction}}&\A_*\ox\A_*\ox\A_*\ar[d]^{1\ox\xi_*}\\
\Sigma\1\bar R_*\ar[r]^-{\textrm{coaction}}&\A_*\ox\Sigma\1\bar R_*.
}
\end{aligned}
\end{equation}
Thus on elements we have
\begin{equation}\label{nablaelts}
{\nabla_\xi}_*(x,y)=\sum\xi_*(x,y)_\A\ox\xi_*(x,y)_R-\sum x_\l
y_{\l'}\ox\xi_*(x_\r,y_{\r'}),
\end{equation}
where again the Sweedler notation is used,
$$
m_*(x)=\sum x_\l\ox x_\r
$$
for the diagonal
$$
m_*:\A_*\to\A_*\ox\A_*
$$
and
$$
a_*(x)=\sum x_\A\ox x_C
$$
for the coaction
$$
a_*:C\to\A_*\ox C
$$
of a left $\A_*$-comodule $C$.

Let us also denote by ${\nabla_S}_*$ the similar operator but with $S_*$ in
place of $\xi_*$. That is, we define
$$
{\nabla_S}_*(x,y)=\sum S_*(x,y)_\A\ox S_*(x,y)_R-\sum x_\l
y_{\l'}\ox S_*(x_\r,y_{\r'}).
$$

We then obviously have
\begin{equation}\label{nablas}
{\nabla_\xi}_*(x,y)+{\nabla_\xi}_*(y,x)={\nabla_S}_*(x,y)
\end{equation}
for any $x,y\in\A_*$.

\begin{Lemma}
The map ${\nabla_\xi}_*$ above is a 2-cocycle, i.~e. for any $x,y,z\in\A_*$ one has
$$
m_*(x){\nabla_\xi}_*(y,z)+{\nabla_\xi}_*(x,yz)={\nabla_\xi}_*(x,y)m_*(z)+{\nabla_\xi}_*(xy,z).
$$
\end{Lemma}

\begin{proof}
First note that the diagram
$$
\xymatrix{
\A_*\ox\bar R_*\ar[r]^-{m_*\ox\textrm{coaction}}\ar[dd]_{\textrm{action}}
&\A_*\ox\A_*\ox\A_*\ox\bar R_*\ar[dr]^{1\ox T\ox1}\\
&&\A_*\ox\A_*\ox\A_*\ox\bar R_*\ar[d]^{\delta_*\ox\textrm{action}}\\
\bar R_*\ar[rr]^{\textrm{coaction}}&&\A_*\ox\bar R_*
}
$$
commutes --- this follows from the fact that the action and coaction of
$\A_*$ on $\bar R_*$ are induced from the multiplication and
comultiplication in $\F_*$ which is a Hopf algebra.

We thus conclude that the coaction map
$$
\bar R_*\to\A_*\ox\bar R_*
$$
is a homomorphism of $\A_*$-modules, so that its composite with the cocycle
$\xi_*$ is a cocycle. It thus remains to show that the composite
$$
\A_*\ox\A_*\to\A_*\ox\A_*\ox\A_*\to\A_*\ox\bar R_*
$$
in the diagram \eqref{nabladef} is also a cocycle. Let us denote this
composite by $\phi$.

Observe that the Hopf algebra diagram for $\A_*$ expressing interchange of the
multiplication and diagonal can be written on elements as follows:
$$
\sum(xy)_\l\ox(xy)_\r=\sum x_\l y_{\l'}\ox x_\r y_{\r'}.
$$
Using this identity we then have for any $x,y,z\in\A_*$
\begin{align*}
m_*(x)\phi(y,z)&=\sum x_\l y_{\l'}z_{\l''}\ox x_\r\xi_*(y_{\r'},z_{\r''});\\
\phi(x,yz)&=\sum
x_\l(yz)_{\l'}\ox\xi_*(x_\r,(yz)_{\r'})
=(\delta_*\ox\xi_*)\left(\sum x_\l\ox(yz)_{\l'}\ox
x_\r\ox(yz)_{\r'}\right)\\
&=(\delta_*\ox\xi_*)\left(\sum x_\l\ox y_{\l'}z_{\l''}\ox x_\r\ox
y_{\r'}z_{\r''}\right)=\sum x_\l
y_{\l'}z_{\l''}\ox\xi_*(x_\r,y_{\r'}z_{\r''});\\
\phi(xy,z)&=\sum(xy)_\l
z_{\l'}\ox\xi_*((xy)_\r,z_{\r'})=(\delta_*\ox\xi_*)\left(\sum(xy)_\l\ox
z_{\l'}\ox(xy)_\r\ox z_{\r'}\right)\\
&=(\delta_*\ox\xi_*)\left(\sum x_\l y_{\l'}\ox z_{\l''}\ox x_\r y_{\r'}\ox
z_{\r''}\right)=\sum x_\l y_{\l'}z_{\l''}\ox\xi_*(x_\r y_\r',z_{\r''});\\
\phi(x,y)m_*(z)&=\sum x_\l y_{\l'}z_{\l''}\ox\xi_*(x_\r,y_{\r'})z_{\r''}.
\end{align*}
These indentities readily imply that $\phi$ is a cocycle as required.
\end{proof}

We next use the fact the cocycle ${\nabla_\xi}_*$ is defined on a polynomial
algebra and hence can be expressed by its values on generators and by its
(anti)symmetrization ${\nabla_S}_*$. Indeed the proof of \cite{Baues}*{16.2.3}
works in this generality, i.~e. one has

\begin{Proposition}
Let $P=k[\zeta_1,\zeta_2,...]$ be a polynomial algebra over a commutative ring $k$,
let $M$ be a $P$-module, let
$$
\gamma:P\ox P\to M
$$
be a Hochschild 2-cocycle, i.~e. one has
$$
x\gamma(y,z)-\gamma(xy,z)+\gamma(x,yz)-z\gamma(x,y)=0
$$
for all $x,y,z\in P$, and let $\sigma$ be the antisymmetrization of
$\gamma$, i.~e.
$$
\sigma(x,y)=\gamma(x,y)-\gamma(y,x).
$$
Then, up to coboundaries, $\gamma$ can be recovered from $\sigma$, i.~e.
there is a cocycle $\gamma_\sigma$ cohomologous to $\gamma$ which depends
only on $\sigma$.
\end{Proposition}

\begin{proof}
To $\gamma$ corresponds a singular extension of $k$-algebras
$$
\xymatrix{
M\ar@{ >->}[r]^i&E\ar@{->>}[r]^p&P
}
$$
whose isomorphism class uniquely determines the cohomology class of
$\gamma$. Let us choose for each polynomial generator $\zeta_n\in P$ an
element $s(\zeta_n)\in E$ with $ps(\zeta_n)=\zeta_n$. Furthermore let us
choose an ordering on the polynomial generators of $P$,
$\zeta_1<\zeta_2<...$; these data determine uniquely a $k$-linear
section of $p$, by the formula
$$
s(\zeta_{n_1}\zeta_{n_2}\cdots)=s(\zeta_{n_1})s(\zeta_{n_2})\cdots
$$ 
for any finite sequence $n_1\le n_2\le\cdots$ of positive integers. Then we
can use $s$ to construct a cocycle $\gamma_\sigma$ cohomologous to
$\gamma$ determined by
$$
s(xy)=s(x)s(y)+i\gamma_\sigma(x,y).
$$
But if $x$ and $y$ are monomials, then $s(xy)$ and $s(x)s(y)$ differ only
by the order of terms, so that $i\gamma_\sigma(x,y)$ is contained in the
ideal generated by commutators
$$
\gamma_\sigma(\zeta_i,\zeta_j)=s(\zeta_i)s(\zeta_j)-s(\zeta_j)s(\zeta_i)=\sigma(\zeta_i,\zeta_j)
$$
for $i>j$. So in fact one can express each $\gamma_\sigma(x,y)$ by a linear
combination of elements of $M$ of the form $z\sigma(\zeta_i,\zeta_j)$ for
$z\in P$.
\end{proof}

\begin{Remark}
Obviously the above proof actually contains an algorithm for expressing the
cocycle $\gamma_\sigma$ in terms of $\sigma$. For $x=\zeta_{n_1}\zeta_{n_2}\cdots\zeta_{n_k}$ and
$y=\zeta_{m_1}\zeta_{m_2}\cdots\zeta_{m_l}$, with $n_1\le n_2\le\cdots\le
n_k$, $m_1\le m_2\le\cdots\le m_l$, either one has $n_k\le m_1$, in which
case $\gamma_\sigma(x,y)=0$ since $s(x)s(y)=s(xy)$, or one has $n_k>m_1$,
in which case one can write
$$
s(x)s(y)=s(\zeta_{n_1})\cdots
s(\zeta_{n_{k-1}})s(\zeta_{m_1})s(\zeta_{n_k})s(\zeta_{m_2})\cdots
s(\zeta_{m_l})
+\zeta_{n_1}\cdots\zeta_{n_{k-1}}\zeta_{m_2}\cdots\zeta_{m_l}\sigma(\zeta_{m_1},\zeta_{n_k}).
$$
Applying the same procedure again several times one finally arrives at
$s(xy)+$(a sum of elements of the form $z\sigma(\zeta_i,\zeta_j)$). In fact
it is easy to see that one has
$$
\gamma_\sigma(\zeta_{n_1}\zeta_{n_2}\cdots\zeta_{n_k},\zeta_{m_1}\zeta_{m_2}\cdots\zeta_{m_l})=
\sum_{n_i>m_j}
\zeta_{n_1}\cdots\zeta_{n_{i-1}}\zeta_{n_{i+1}}\cdots\zeta_{n_k}
\zeta_{m_1}\cdots\zeta_{m_{j-1}}\zeta_{m_{j+1}}\cdots\zeta_{m_l}\sigma(\zeta_{m_j},\zeta_{n_i}).
$$

In the characteristic $p>0$ case further obvious simplifications occur.
In particular we can choose the cocycle $\xi_*$ in \bref{cosxi} in such a way that the formula
\begin{equation}\label{xifromS}
\begin{aligned}
&\xi_*(\zeta_1^{d_1}\zeta_2^{d_2}\cdots,\zeta_1^{e_1}\zeta_2^{e_2}\cdots)=\\
&\sum_{\substack{i<j\\\textrm{$e_i$, $d_j$ odd}}}
\zeta_1^{d_1+e_1}\cdots\zeta_{i-1}^{d_{i-1}+e_{i-1}}\zeta_i^{d_i+e_i-1}\zeta_{i+1}^{d_{i+1}+e_{i+1}}
\cdots\zeta_{j-1}^{d_{j-1}+e_{j-1}}\zeta_j^{d_j+e_j-1}\zeta_{j+1}^{d_{j+1}+e_{j+1}}\cdots
S_*(\zeta_i,\zeta_j)
\end{aligned}
\end{equation}
holds%, and similarly for $\xi_*$ and $S_*$ replaced by ${\nabla_\xi}_*$ and ${\nabla_S}_*$ respectively:
%\begin{equation}\label{nablaxi}
%\begin{aligned}
%&{\nabla_\xi}_*(\zeta_1^{d_1}\zeta_2^{d_2}\cdots,\zeta_1^{e_1}\zeta_2^{e_2}\cdots)=\\
%&\sum_{\substack{i<j\\\textrm{$e_i$, $d_j$ odd}}}
%\zeta_1^{d_1+e_1}\cdots\zeta_{i-1}^{d_{i-1}+e_{i-1}}\zeta_i^{d_i+e_i-1}\zeta_{i+1}^{d_{i+1}+e_{i+1}}
%\cdots\zeta_{j-1}^{d_{j-1}+e_{j-1}}\zeta_j^{d_j+e_j-1}\zeta_{j+1}^{d_{j+1}+e_{j+1}}\cdots
%{\nabla_S}_*(\zeta_i,\zeta_j).
%\end{aligned}
%\end{equation}
\end{Remark}

The operator ${\nabla_S}_*$ is readily computable. It is a symmetric
biderivation, with ${\nabla_S}_*(x,x)=0$ for all $x$, thus uniquely
determined by its values of the form ${\nabla_S}_*(\zeta_n,\zeta_m)$ for
$n<m$, which are expressed easily from the corresponding values of $S_*$.
For example, one has
\begin{align*}
{\nabla_S}_*(\zeta_1,\zeta_2)&
=\zeta_1\ox M_{1,1}^2,\\
{\nabla_S}_*(\zeta_1,\zeta_3)&
=\zeta_1^5\ox M_{1,1}^2+\zeta_1\ox M_{2,1,1}^2,\\
{\nabla_S}_*(\zeta_2,\zeta_3)&
=\left(\zeta_1^7+\zeta_1\zeta_2^2\right)\ox M_{1,1}^2
+\zeta_1^3\ox M_{2,1,1}^2
+\zeta_1\ox\left(M_1^2M_3+M_1M_{2,1,1}+M_5+M_{4,1}+M_{3,2}+M_{2,1,1,1}\right)^2,\\
{\nabla_S}_*(\zeta_1,\zeta_4)&
=\zeta_1\zeta_2^4\ox M_{1,1}^2
+\zeta_1^9\ox M_{2,1,1}^2
+\zeta_1\ox M_{4,2,1,1}^2,\\
{\nabla_S}_*(\zeta_2,\zeta_4)&
=\left(\zeta_1^3\zeta_2^4+\zeta_1\zeta_3^2\right)\ox M_{1,1}^2
+\zeta_1^{11}\ox M_{2,1,1}^2\\
&+\zeta_1^9\ox\left(M_1^2M_3+M_1M_{2,1,1}+M_5+M_{4,1}+M_{3,2}+M_{2,1,1,1}\right)^2
 +\zeta_1^3\ox M_{4,2,1,1}^2\\
&+\zeta_1\ox\left(
M_1^4M_5
+M_{2,1}^2M_3
+M_1^2M_{4,2,1,1}
+M_9
+M_{7,2}
+M_{6,2,1}
+M_{5,4}
+M_{3,4,2}
+M_{4,3,2}
+M_{4,4,1}
\right.\\
&\phantom{+\zeta_1\ox\ }\left.
+M_{2,4,2,1}
+M_{4,2,1,1,1}
\right)^2,
\end{align*}
etc.

\section{The formula of Kristensen}

We will next use certain elements defined in \cite{Kristensen}*{Theorem 3.3} to
derive more explicit expressions for $\xi_*$, hence for $S_*$,
${\nabla_S}_*$ and ${\nabla_\xi}_*$. We recall that Kristensen defines
$$
A[a,b]=(\Sq^1\ox\Sq^{0,1})\delta\left(\Sq^{a-3}\Sq^{b-2}+\Sq^{a-2}\Sq^{b-3}
+\sum_j\binom{b-1-j}{a-2j}(\Sq^{a+b-j-3}\Sq^{j-2}+\Sq^{a+b-j-2}\Sq^{j-3})\right),
$$
for natural numbers $a,b$. Obviously one has
$$
A[a,b]=(\Sq^1\ox\Sq^{0,1})\delta k([a,b]),
$$
where $k$ is the operator determined by
$$
k(xy)=\k(\k\k(x)\k\k(y))
$$
for $x,y\in\F_0^{\le1}$. We then interpret $A[a,b]$ as an $\FF$-linear
operator of the form
$$
K:\F_0^{\le1}\ox\F_0^{\le1}\to\A\ox\A
$$
given by
$$
K(x\ox y)=(\Sq^1\ox\Sq^{0,1})\delta\k(\k\k(x)\k\k(y))
$$
which is factored through $\F_0^{\le1}\ox\F_0^{\le1}\onto\F_0^{\le2}$ and
then restricted to $R_\F^{\le2}\into\F_0^{\le2}$. We then can dualize $K$
to get

\begin{Definition}
We define an $\FF$-linear operator
$$
K_*:\A_*\ox\A_*\to{R^{\le2}_\F}_*
$$
as composite with the quotient map $\F^{\le2}_*\onto{R^{\le2}_\F}_*$ of the
dual of $K$ above (whose image lies in that of
$m_*:\F^{\le2}_*\into\F^{\le1}_*\ox\F^{\le1}_*$.
\end{Definition}

Thus explicitly, $K_*$ is the composite
$$
\A_*\ox\A_*\xto{\Sq^1\cdot_*\ox\Sq^{0,1}\cdot_*}\A_*\ox\A_*\xto{\delta_*}
\A_*\xto{\zeta_1}\A_*\into\F_*\xto{m_*}\F_*\ox\F_*\onto\F^{\le1}_*\ox\F^{\le1}_*
\xto{M_1^2\ox M_1^2}\F^{\le1}_*\ox\F^{\le1}_*
$$
landing in $\F^{\le2}_*\into\F^{\le1}_*\ox\F^{\le1}_*$ and precomposed with
$\F^{\le2}_*\onto{R_\F^{\le2}}_*$. Or on elements,
$$
K_*(x\ox y)=(M_1^2\ox M_1^2)\left(m_*(\zeta_1\frac{\d x}{\d\zeta_1}\frac{\d
y}{\d\zeta_2})\right)^{\le1}=m_*(M_{1,1}^2M_1\frac{\d x}{\d\zeta_1}\frac{\d
y}{\d\zeta_2})^{\le1}.
$$
One thus has
\begin{equation}\label{kstar}
K_*(\zeta_1^{n_1}\zeta_2^{n_2}\cdots\ox\zeta_1^{m_1}\zeta_2^{m_2}\cdots)=
\begin{cases}
M_1^{n_1+m_1}M_{2,1}^{n_2+m_2-1}M_{1,1}^2,&\textrm{$n_1$, $m_2$ odd,
$n_i=m_i=0$ for $i>2$},\\
0&\textrm{otherwise.}
\end{cases}
\end{equation}

We have
\begin{Proposition}
Symmetrization of the operator $K_*$ dual to the operator $S_R$ in
\eqref{sr}, i.~e. is given by precomposing ${S_\F}_*$ given in \eqref{msf}
with the restriction map $\F_*^{\le2}\onto{R_\F}_*^{\le2}$.
\end{Proposition}

\begin{proof}
From the above formula \eqref{kstar}, for monomials $x=\zeta_1^{n_1}\zeta_2^{n_2}\zeta_3^{n_3}\cdots$ and
$y=\zeta_1^{m_1}\zeta_2^{m_2}\zeta_3^{m_3}\cdots$ we have
$$
K_*(x\ox y)+K_*(y\ox x)=
\begin{cases}
M_1^{n_1+m_1}M_{2,1}^{n_2+m_2-1}M_{1,1}^2,&\textrm{$n_1m_2+m_1n_2$ odd and
$n_i=m_i=0$ for $i>2$,}\\
0&\textrm{otherwise.}
\end{cases}
$$
On the other hand, using the explicit expression \eqref{msf} and the expression for the
operator $S^1_*$ in \eqref{s1} we can write
$$
m_*{S_\F}_*(x\ox y)=\sum_{\substack{x_\l=\zeta_1^{2n-1}\\y_{\l'}=\zeta_1^{2n'-1}}}
\zeta_1^{2(n+n')-1}\ox(x_\r y_{\r'})^{\le1}
+(\zeta_1\ox1+1\ox\zeta_1)\sum_{\substack{x_\r=\zeta_1^{2n-1}\\y_{\r'}=\zeta_1^{2n'-1}}}
(x_\l y_{\l'})^{\le1}\ox\zeta_1^{2(n+n'-1)}.
$$
From the expression \eqref{mildiag} for the Milnor diagonal we thus see
that for monomials $x=\zeta_1^{n_1}\zeta_2^{n_2}\zeta_3^{n_3}\cdots$ and
$y=\zeta_1^{m_1}\zeta_2^{m_2}\zeta_3^{m_3}\cdots$ one has ${S_\F}_*(x\ox
y)=0$ unless $n_i=m_i=0$ for $i>2$, whereas in the remaining cases one
has
\begin{multline*}
m_*{S_\F}_*(\zeta_1^{n_1}\zeta_2^{n_2}\ox\zeta_1^{m_1}\zeta_2^{m_2})=
\sum_{\substack{0\le i\le n_1\\0\le j\le m_1\\\textrm{$i$, $j$ odd}}}
\binom{n_1}i\binom{m_1}j
\zeta_1^{i+j+2(n_2+m_2)+1}\ox\zeta_1^{n_1+m_1-i-j+n_2+m_2}\\
+(\zeta_1\ox1+1\ox\zeta_1)
\sum_{\substack{0\le i\le n_1\\0\le j\le m_1\\\textrm{$n_1-i+n_2$, $m_1-j+m_2$ odd}}}
\binom{n_1}i\binom{m_1}j\zeta_1^{i+j+2(n_2+m_2)}\ox\zeta_1^{n_1+m_1-i-j+n_2+m_2}.
\end{multline*}

Let us now turn back to the symmetrization of $K_*$. We compute its image
under the map $m_*$; by \eqref{mstar} it sends $M_1$ to
$\zeta_1\ox1+1\ox\zeta_1$, $M_{1,1}$ to $\zeta_1\ox\zeta_1$ and $M_{2,1}$
to $\zeta_1^2\ox\zeta_1$. Thus the nonzero values of this image are, for
$n_1m_2+m_1n_2$ odd,
$$
m_*K_*(1+T)(\zeta_1^{n_1}\zeta_2^{n_2}\ox\zeta_1^{m_1}\zeta_2^{m_2})=
(\zeta_1\ox1+1\ox\zeta_1)^{n_1+m_1}(\zeta_1^2\ox\zeta_1)^{n_2+m_2-1}(\zeta_1^2\ox\zeta_1^2).
$$
Then expanding
$(\zeta_1\ox1+1\ox\zeta_1)^{n_1+m_1}=(\zeta_1\ox1+1\ox\zeta_1)^{n_1}(\zeta_1\ox1+1\ox\zeta_1)^{m_1}$
via binomials we obtain
$$
m_*K_*(1+T)(\zeta_1^{n_1}\zeta_2^{n_2}\ox\zeta_1^{m_1}\zeta_2^{m_2})=
\sum_{\substack{0\le i\le n_1\\0\le j\le
m_1}}\binom{n_1}i\binom{m_1}j\zeta_1^{i+j+2(n_2+m_2)}\ox\zeta_1^{n_1+m_1-i-j+n_2+m_2+1}.
$$
It follows that nonzero values of the difference $m_*({S_\F}_*-K_*(1+T))$
on monomials in Milnor generators are equal to
\begin{multline*}
\sum_{\substack{0\le i\le n_1\\0\le j\le m_1\\\textrm{$i$, $j$ odd}}}
\binom{n_1}i\binom{m_1}j
\zeta_1^{i+j+2(n_2+m_2)+1}\ox\zeta_1^{n_1+m_1-i-j+n_2+m_2}\\
+\sum_{\substack{0\le i\le n_1\\0\le j\le m_1\\\textrm{$n_1-i+n_2$, $m_1-j+m_2$ odd}}}
\binom{n_1}i\binom{m_1}j\zeta_1^{i+j+2(n_2+m_2)+1}\ox\zeta_1^{n_1+m_1-i-j+n_2+m_2}\\
+\sum_{\substack{0\le i\le n_1\\0\le j\le m_1\\\textrm{$n_1-i+n_2$, $m_1-j+m_2$ even}}}
\binom{n_1}i\binom{m_1}j\zeta_1^{i+j+2(n_2+m_2)}\ox\zeta_1^{n_1+m_1-i-j+n_2+m_2+1}
\end{multline*}
for $n_1m_2+m_1n_2$ odd and
$m_*{S_\F}_*(\zeta_1^{n_1}\zeta_2^{n_2}\ox\zeta_1^{m_1}\zeta_2^{m_2})$ for
$n_1m_2+m_1n_2$ even.

The first expression can be rewritten as
\begin{multline*}
(\zeta_1^2\ox\zeta_1)^{n_2+m_2}\sum_k\zeta_1^{k+1}\ox\zeta_1^{n_1+m_1-k}\\
\left(
\sum_{\substack{0\le i\le n_1\\0\le k-i\le m_1\\\textrm{$i$, $k-i$ odd}}}
\binom{n_1}i\binom{m_1}{k-i}
+\sum_{\substack{0\le i\le n_1\\0\le k-i\le m_1\\\textrm{$n_1-i+n_2$,
$m_1-k+i+m_2$ odd}}}\binom{n_1}i\binom{m_1}{k-i}
+\sum_{\substack{0\le i\le n_1\\0\le k+1-i\le m_1\\\textrm{$n_1-i+n_2$,
$m_1-k-1+i+m_2$
even}}}\binom{n_1}i\binom{m_1}{k+1-i}
\right)
\end{multline*}
and in the second case we may write
\begin{multline*}
m_*{S_\F}_*(\zeta_1^{n_1}\zeta_2^{n_2}\ox\zeta_1^{m_1}\zeta_2^{m_2})=(\zeta_1^2\ox\zeta_1)^{n_2+m_2}
\sum_k\zeta_1^{k+1}\ox\zeta_1^{n_1+m_1-k}\\\left(\sum_{\substack{0\le i\le n_1\\0\le k-i\le m_1\\\textrm{$i$, $k-i$ odd}}}
\binom{n_1}i\binom{m_1}{k-i}
+\sum_{\substack{0\le i\le n_1\\0\le k-i\le m_1\\\textrm{$n_1-i+n_2$,
$m_1-k+i+m_2$ odd}}}\binom{n_1}i\binom{m_1}{k-i}
+\sum_{\substack{0\le i\le n_1\\0\le k+1-i\le m_1\\\textrm{$n_1-i+n_2$,
$m_1-k-1+i+m_2$
odd}}}\binom{n_1}i\binom{m_1}{k+1-i}\right).
\end{multline*}
One then shows that these expressions lie in the subalgebra of
$\F_*^{\le1}\ox\F_*^{\le1}$ generated by $\zeta_1^2\ox\zeta_1$ and
$\zeta_1\ox1+1\ox\zeta_1$, without involvement of $\zeta_1\ox\zeta_1$. This
means that the image of the difference ${S_\F}_*-K_*(1+T)$ under the
restriction map $\F_*^{\le2}\onto{R_\F}_*^{\le2}$ is zero.
\end{proof}

%% file: dHascoAss_c8.tex
\chapter{Computation of the algebra of secondary cohomology operations and its
dual}\label{A}

We first describe explicit splittings of the pair algebra $\R^\FF$ of relations
in the Steenrod algebra and its dual $\R_\FF$. Then we describe in terms of
these splittings $s$ the multiplication maps $A^s$ for the Hopf pair algebra
$\aB^\FF$ of secondary cohomology operations and we describe the dual maps
$A_s$ determining the Hopf pair coalgebra $\aB_\FF$ dual to $\aB^\FF$. On the
basis of the main result in \cite{Baues} we describe systems of
equations which can be solved inductively by a computer and which yield the
multiplication maps $A^s$ and $A_s$ as a solution. It turns out that $A_s$ is
explicitly given by a formula in which only the values $A_s(\zeta_n)$,
$n\ge1$, have to be computed where $\zeta_n$ is the Milnor generator in the
dual Steenrod algebra $\A_*$.

\section{Computation of $\R^\FF$ and $\R_\FF$}\label{rcomp}

Let us fix a function $\chi:\FF\to\GG$ which splits the projection
$\GG\to\FF$, namely, take
\begin{equation}\label{chieq}
\chi(k\!\!\mod p) = k\!\!\mod p^2, \ 0\le k<p.
\end{equation}
We will use $\chi$ to define splittings of
$\R^\FF=\left(\R^\FF_1\xto\d\R^\FF_0\right)$. Here a \emph{splitting} $s$
of $\R^\FF$ is an $\FF$-linear map for which the diagram
\begin{equation}\label{s}
\alignbox{
\xymatrix{
&&\R^\FF_1\ar[d]^\d\\
R_\F\ar@{ (->}[r]\ar[urr]^s
&\F_0\ar@{=}[r]
&\R^\FF_0
}
}
\end{equation}
commutes with $R_\F=\im(\d)=\ker(q_\F:\F_0\to\A)$. We only consider the
case $p=2$.

\begin{Definition}[The right equivariant splitting of $\R^\FF$]\label{chi}
Using $\chi$, all \emph{Adem relations}
$$
[a,b]:=\Sq^a\Sq^b+\sum_{k=0}^{\left[\frac
a2\right]}\binom{b-k-1}{a-2k}\Sq^{a+b-k}\Sq^k
$$
for $a,b>0$, $a<2b$, can be lifted to elements $[a,b]_\chi\in R_\aB$ by
applying $\chi$ to all coefficients, i.~e. by interpreting $[a,b]$ as an
element of $\aB$. As shown in \cite{Baues}*{16.5.2}, $R_\F$ is a free right
$\F_0$-module with a basis consisting of \emph{preadmissible relations}. For
$p=2$ these are elements of the form
$$
\Sq^{a_1}\cdots\Sq^{a_{k-1}}[a_k,a]\in R_\F
$$
satisfying $a_1\ge2a_2$, ..., $a_{k-2}\ge2a_{k-1}$, $a_{k-1}\ge2a_k$,
$a_k<2a$. Sending such an element to
$$
\Sq^{a_1}\cdots\Sq^{a_{k-1}}[a_k,a]_\chi\in R_\B
$$
determines then a unique right $\F_0$-equivariant splitting $\phi$ in the pair
algebra $\R^\FF$; that is, we get a commutative diagram
$$
\xymatrix@!C=5em{
R_\F\ar[r]^-\phi\ar@{ (->}[d]
&**[l]R_\aB\ox\FF=\R^\FF_1\ar[d]^\d\\
\F_0\ar@{=}[r]
&\R^\FF_0.
}
$$
\end{Definition}

For a splitting $s$ of $\R^\FF$ the map $s\!\ox\!1\oplus1\!\ox\!s$
induces the map $s_\#$ in the diagram
\begin{equation}\label{rsplit}
\alignbox{
\xymatrix{
\R^\FF_1\ar[r]^-\Delta\ar@{->>}[d]^\d
&(\R^\FF\hat\ox\R^\FF)_1\ar@{->>}[d]^{\d_{\hat\ox}}
&\R^\FF_1\!\ox\!\F_0\oplus\F_0\!\ox\!\R^\FF_1\ar[l]\\
R_\F\ar@/^/@{-->}[u]^s\ar@{ (->}[d]\ar[r]^{\Delta_R}
&R^{(2)}_\F\ar@/^/@{-->}[u]^{s_\#}\ar@{ (->}[d]
&R_\F\!\ox\!\F_0\oplus\F_0\!\ox\!R_\F\ar[l]\ar@{-->}[u]_{s\ox\!1\oplus1\!\ox\!s}\\
\F_0\ar[r]^-\Delta
&\F_0\ox\F_0.
}
}
\end{equation}
Then the difference
$U=s_\#\Delta_R-\Delta s:R_\F\to(\R^\FF\hat\ox\R^\FF)_1$ satisfies $\d_{\hat\ox}U=0$ since
$$
\d_{\hat\ox} s_\#\Delta_R=\Delta_R=\Delta_R\d s=\d_{\hat\ox}\Delta s.
$$
Thus $U$ lifts to $\ker\d_{\hat\ox}\cong\A\ox\A$ and gives an $\FF$-linear
map
\begin{equation}\label{ueq}
U^s:R_\F\to\A\ox\A.
\end{equation}
If we use the splitting $s$ to identify $\R^\FF_1$ with the direct sum
$\A\oplus R_\F$, then it is clear that knowledge of the map $U^s$ determines the
diagonal $\R^\FF_1\to(\R^\FF\hat\ox\R^\FF)_1$ completely. Indeed $s_\#$
yields the identification $(\R^\FF\hat\ox\R^\FF)_1\cong\A\!\ox\!\A\oplus
R^{(2)}_\F$, and under these identifications
$\Delta:\R^\FF_1\to(\R^\FF\hat\ox\R^\FF)_1$ corresponds to a map which by
commutativity of \eqref{rsplit} must have the form
\begin{equation}\label{udia}
\A\oplus R_\F
\xto{\left(\begin{smallmatrix}\Delta_\A&U^s\\0&\Delta_R\end{smallmatrix}\right)}
\A\!\ox\!\A\oplus R^{(2)}_\F
\end{equation}
and is thus determined by $U^s$.

One readily checks that the map $U^s$ for $s=\phi$ in \bref{chi} coincides with the map $U$
defined in \cite{Baues}*{16.4.3} in terms of the algebra $\aB$.

Given the splitting $s$ and the map $U^s$, the only piece of structure
remaining to determine the $\Alg^\pair_\oo$-comonoid structure of
$\R^\FF$ completely is the $\F_0$-$\F_0$-bimodule structure on
$\R^\FF_1\cong\A\oplus R_\F$. Consider for $f\in\F_0$, $r\in R_\F$ the
difference $ s(fr)-f s(r)$. It belongs to the kernel of $\d$ since
$$
\d s(fr)=fr=f\d s(r)=\d(f s(r)).
$$
Thus we obtain the \emph{left multiplication map}
\begin{equation}\label{aeq}
a^s:\F_0\ox R_\F\to\A.
\end{equation}

Similarly we obtain the \emph{right multiplication map}
$$
b^s:R_\F\ox\F_0\to\A
$$
by the difference $s(rf)-s(r)f$.

\begin{Lemma}\label{chisp}
For $s=\phi$ in \bref{chi} the right multiplication map $b^\phi$ is trivial,
that is $\phi$ is right equivariant, and the left multiplication map factors
through $q_\F\ox1$ inducing the map
$$
a^\phi:\A\ox R_\F\to\A.
$$
\end{Lemma}

\begin{proof}
Right equivariance holds by definition. As for the factorization, $R_\F\ox
R_\F\into\F_0\ox R_\F$ is in the kernel of $a^\phi:\F_0\ox R_\F\to\A$, since
by right equivariance of $s$ and by the pair algebra property \eqref{paireq}
for $\R^\FF$ one has for any $r,r'\in R_\F$
$$
s(rr')=s(r)r'=s(r)\d s(r')=(\d s(r))s(r')=rs(r').
$$
Hence factoring the above map through $(\F_0\ox R_\F)/(R_\F\ox
R_\F)\cong\A\ox R_\F$ we obtain a map
$$
\A\ox R_\F\to\A.
$$
\end{proof}

Summarizing the above, we thus have proved

\begin{Proposition}\label{detrel}
Using the splitting $s=\phi$ of $\R^\FF$ in \bref{chi}
the comonoid $\R^\FF$ in the category $\Alg^\pair_\oo$ described
in \bref{relcom} is completely determined by the maps
$$
U^\phi:R_\F\to\A\ox\A
$$
and
$$
a^\phi:\A\ox R_\F\to\A
$$
given in \eqref{ueq} and \eqref{aeq} respectively.
\end{Proposition}\qed

We next introduce another splitting $s=\psi$ for which $U^s=0$. For this we use
the fact that $\A_*=\Hom(\A,\FF)$ and
\begin{equation}
\aB_\#=\Hom(\aB_0,\GG)
\end{equation}
with $\aB_0=T_\GG(E_\A)$ both are polynomial algebras in such a way that
generators of $\A_*$ are also (part of the) generators of $\aB_\#$.

Using $\chi$ in \eqref{chieq} we obtain the function
\begin{equation}\label{psichi}
\psi_\chi:\A_*\to\aB_\#
\end{equation}
(which is not $\FF$-linear) as follows. Each element $x$ in $\A_*$ is
uniquely an $\FF$-linear combination $x=\sum_\alpha n_\alpha\alpha$ where
$\alpha$ runs through the monomials in Milnor generators. Such a monomial
can be also considered as an element in $\aB_\#$ by \bref{eslp} so that we
can define
$$
\psi_\chi(x)=\sum_\alpha\chi(n_\alpha)\alpha\in\aB_\#.
$$

\begin{Definition}[The comultiplicative splitting of $\R^\FF$]\label{psi}
Consider the following commutative diagram with exact rows and columns
$$
\xymatrix{
\A_*\ar@{ >->}[r]\ar@{-->}[dr]_{\psi_\chi}&\F_*\ar[r]&\Hom(R_\aB,\FF)\\
&\aB_\#\ar[r]^-q\ar@{->>}[u]&\Hom(R_\aB,\GG)\ar[u]\\
\A_*\ar@{ >->}[r]^{{q_\F}_*}&\F_*\ar@{ >->}[u]_j\ar[r]
&\Hom(R_\aB,\FF)\ar@{ >->}[u]_{j_R}\ar@{->>}[r]&\A_*\ar@{-->}[dl]^{q\psi_\chi}\\
&\R_\FF^0\ar@{=}[u]&\R_\FF^1\ar@{=}[u]
}
$$
with the columns induced by the short exact sequence $\FF\into\GG\onto\FF$
and the rows induced by \eqref{drel}. In particular $q$ is induced by the
inclusion $R_\aB\subset\aB_0$. Now it is clear that $\psi_\chi$ yields a map
$q\psi_\chi$ which lifts to $\Hom(R_\aB,\FF)$ so that we get the map
$$
q\psi_\chi:\A_*\to\R_\FF^1
$$
which splits the projection $\R_\FF^1\onto\A_*$. Moreover $q\psi_\chi$ is
$\FF$-linear since for all $x,y\in\A_*$ the elements
$\psi_\chi(x)+\psi_\chi(y)-\psi_\chi(x+y)\in\aB_\#$ are in the image of the
inclusion $j{q_\F}_*:\A_*\into\aB_\#$ and thus go to zero under $q$.

The dual of $q\psi_\chi$ is thus a retraction $(q\psi_\chi)^*$ in the
short exact sequence
$$
\xymatrix{
R_\F\ar@{-->}[dr]_{(q\psi_\chi)^*_\perp}&R_\aB\ox\FF\ar@{=}[d]\ar@{->>}[l]_\pi&\A\ar@{
>->}[l]_-\iota\\
&\R^\FF_1\ar@{-->}[ur]_{(q\psi_\chi)^*}
}
$$
which induces the splitting $\psi=(q\psi_\chi)^*_\perp$ of $\R^\FF$
determined by
$$
\psi(\pi(x))=x-\iota((q\psi_\chi)^*(x)).
$$
\end{Definition}

\begin{Lemma}
For the splitting $s=\psi$ of $\R^\FF$ we have $U^\psi=0$.
\end{Lemma}

\begin{proof}
We must show that the following diagram commutes:
$$
\xymatrix{
R_\F\ar[r]^{\Delta_R}\ar[d]_\psi&R_\F^{(2)}\ar[d]^{\psi_\#}\\
R_\aB\ox\FF\ar[r]^-\Delta&(\R^\FF\hat\ox\R^\FF)_1.
}
$$
Obviously this is equivalent to commutativity of the dual diagram
$$
\xymatrix{
(\R_\FF\check\ox\R_\FF)^1\ar[r]^-{\Delta_*}\ar[d]_{(\psi_\#)_*}&\Hom(R_\aB,\FF)\ar[d]^{\psi_*}\\
{R_\F^{(2)}}_*\ar[r]^{{\Delta_R}_*}&{R_\F}_*
}
$$
which in turn is equivalent to commutativity of
\begin{equation}\label{psicom}
\alignbox{
\xymatrix{
\A_*\ox\A_*\ar[d]_{(\psi_\#)_*^\perp}\ar[r]^{\delta_*}&\A_*\ar[d]^{q\psi_\chi}\\
(\R_\FF\check\ox\R_\FF)^1\ar[r]^-{\Delta_*}&\Hom(R_\aB,\FF).
}
}
\end{equation}
On the other hand, the left hand vertical map in the latter diagram can be
included into another commutative diagram
$$
\xymatrix{
\F_*\!\ox\!\A_*\oplus\A_*\!\ox\!\F_*
\ar[d]_{1\!\ox\!q\psi_\chi\oplus q\psi_\chi\!\ox\!1}
&\A_*\ox\A_*
\ar@{ )->}[l]_-{\binom{i\ox1}{1\ox i}}\ar[d]^{(\psi_\#)_*^\perp}\\
\F_*\!\ox\!\R_\FF^1\oplus\R_\FF^1\!\ox\!\F_*&(\R_\FF\ox\R_\FF)^1\ar@{ )->}[l]
}
$$
It follows that on elements, commutativity of \eqref{psicom} means that the
equality
$$
q\psi_\chi(xy)=i(x)q\psi_\chi(y)+q\psi_\chi(x)i(y)
$$
holds for any $x,y\in\A_*$. By linearity, it is clearly enough to prove
this when $x$ and $y$ are monomials in Milnor generators.

For this observe that for any $x\in\A_*=\Hom(\A,\FF)$, the element
$q\psi_\chi(x)\in\Hom(R_\aB,\FF)$ is the unique $\FF$-linear map making the
diagram
$$
\xymatrix{
R_\aB\ar@{ >->}[r]\ar@{-->}[d]_{q\psi_\chi(x)}
&\aB_0\ar[d]_{\psi_\chi(x)}\ar@{->>}[r]
&\A\ar[d]^x\\
\FF\ar@{ >->}[r]
&\GG\ar@{->>}[r]
&\FF
}
$$
commute. This uniqueness implies the equality we need in view of the
following commutative diagram with exact columns:
$$
\xymatrix@!C=5em{
R_\aB\ar[r]\ar@{ >->}[d]
&R_\aB^{(2)}\ar@{ >->}[d]\ar@{-->}[r]
&\FF\ox\FF\ar[r]^\cong\ar@{ >->}[d]
&\FF\ar@{ >->}[d]\\
\aB_0\ar[r]^-\Delta\ar@{->>}[d]
&\aB_0\ox\aB_0\ar@{->>}[d]\ar[r]^{\psi_\chi(x)\ox\psi_\chi(y)}
&\GG\ox\GG\ar@{->>}[d]\ar[r]^\cong
&\GG\ar@{->>}[d]\\
\A\ar[r]^-\delta\ar@/_3ex/[rrr]_{xy}
&\A\ox\A\ar[r]^{x\ox y}
&\FF\ox\FF\ar[r]^\cong
&\FF,
}
$$
since when $x$ and $y$ are monomials in Milnor generators, one has
$\psi_\chi(xy)=\psi_\chi(x)\psi_\chi(y)$.
\end{proof}

Therefore we call $\psi$ the comultiplicative splitting of $\R^\FF$. We now
want to compute the left and right multiplication maps $a^\psi$ and
$b^\psi$ defined in \eqref{aeq}. The dual maps $a_\psi=(a^\psi)_*$ and
$b_\psi=(b^\psi)_*$ can be described by the diagrams
\begin{equation}\label{ml}
\alignbox{
\xymatrix{
&(R_\aB)_*\ar@{->>}[r]\ar[d]_{m_*^\l}&\A_*\ar@/_1em/@{-->}[l]_{q\psi_\chi}\ar[d]^{m_*}\\
\F_*\ox(R_\F)_*\ar@{
(->}[r]&\F_*\ox(R_\aB)_*\ar@{->>}[r]
&\A_*\ox\A_*\ar@/^1em/@{-->}[l]^{i\ox q\psi_\chi}
}
}
\end{equation}
and
\begin{equation}\label{mr}
\alignbox{
\xymatrix{
&(R_\aB)_*\ar@{->>}[r]\ar[d]_{m_*^\r}&\A_*\ar@/_1em/@{-->}[l]_{q\psi_\chi}\ar[d]^{m_*}\\
(R_\F)_*\ox\F_*\ar@{
(->}[r]&(R_\aB)_*\ox\F_*\ar@{->>}[r]
&\A_*\ox\A_*.\ar@/^1em/@{-->}[l]^{q\psi_\chi\ox i}
}
}
\end{equation}
Here $m_*$ is dual to the multiplication in $\A$ and $m_*^\l$ and $m_*^\r$
are induced by the $\F_0$-$\F_0$-bimodule structure of $R_\aB\ox\FF$. One
readily checks
\begin{align*}
a_\psi&=m_*^\l q\psi_\chi-(i\ox q\psi_\chi)m_*\\
b_\psi&=m_*^\r q\psi_\chi-(q\psi_\chi\ox i)m_*.
\end{align*}
We now consider the diagram
$$
\xymatrix{
&\aB_\#\ar[d]_{m_*^\GG}&\A_*\ar[l]_{\psi_\chi}\ar[d]^{m_*}\\
\F_*\ox\F_*\ar@{ (->}[r]&\aB_\#\ox\aB_\#&\A_*\ox\A_*\ar[l]_-{\psi_\chi^\ox}
}
$$
Here $\psi_\chi^\ox$ is defined similarly as $\psi_\chi$ in \eqref{psichi}
by the formula
$$
\psi_\chi^\ox\left(\sum_{\alpha,\beta}n_{\alpha\beta}\alpha\ox\beta\right)=
\sum_{\alpha,\beta}\chi(n_{\alpha\beta})\alpha\ox\beta
$$
where $\alpha$, $\beta$ run through the monomials in Milnor generators.
Moreover $m_*^\GG$ is the dual of the multiplication map $m^\GG$ of
$\aB_0=T_\GG(E_\A)$.

\begin{Lemma}\label{mulstar}
The difference $m^\GG_*\psi_\chi-\psi_\chi^\ox m_*$ lifts to an
$\FF$-linear map $\nabla_\chi:\A_*\to\F_*\ox\F_*$ such that one has
\begin{align*}
a_\psi&=(1\ox\pi)\nabla_\chi\\
b_\psi&=(\pi\ox1)\nabla_\chi.
\end{align*}
Here $\pi:\F_*\onto{R_\F}_*$ is induced by the inclusion $R_\F\subset\F_0$.
\end{Lemma}

\begin{proof}
We will only prove the first equality; the proof for the second one is
similar.

The following diagram
$$
\xymatrix@!C=3em{
&{R_\aB}_*\ar@{ >->}[rrd]_{j_R}\ar[ddddd]_{m^\l_*}&&&&&&\A_*\ar[llllll]_{q\psi_\chi}\ar@{=}[ddl]\ar[ddddd]^{m_*}\\
&&&{R_\aB}_\#\ar[ddd]_{m^\l_\#}\\
&&&&&\aB_\#\ar[d]_{m_*^\GG}\ar[ull]^\pi&\A_*\ar[l]_{\psi_\chi}\ar[d]^{m_*}\\
&&&&&\aB_\#\ox\aB_\#\ar[dll]_{1\ox\pi}&\A_*\ox\A_*\ar[l]_-{\psi_\chi^\ox}\\
&&&\aB_\#\ox{R_\aB}_\#&\F_*\ox\F_*\ar[lllldd]|\hole^>(.75){1\ox\pi}\ar@{ (->}[ur]\\
&\F_*\ox{R_\aB}_*\ar@{ >->}[urr]^{j\ox j_R}&&&&&&\A_*\ox\A_*\ar[llllll]^{i\ox q\psi_\chi}\ar@{=}[uul]\\
\F_*\ox{R_\F}_*\ar@{ (->}[ur]
}
$$
commutes except for the innermost square, whose deviation from
commutativity is $\nabla_\chi$ and lies in the image of
$\F_*\ox\F_*\incl\aB_\#\ox\aB_\#$, and the outermost square, whose deviation
from commutativity is $a_\psi$ and lies in the image of
$\F_*\ox{R_\F}_*\incl\F_*\ox{R_\aB}_*$. It follows that
$(1\ox\pi)\nabla_\chi$ and $a_\psi$ have the same image under $j\ox
j_R$, and since the latter map is injective we are done.
\end{proof}

Let us describe the map $\nabla_\chi$ more explicitly.

\begin{Lemma}
The map $\nabla_\chi$ factors as follows
$$
\A_*\xto{\bar\nabla\chi}\F_*\ox\A_*\xto{1\ox i}\F_*\ox\F_*.
$$
\end{Lemma}

\begin{proof}
Let $\A_\#\subset\aB_\#$ be the subring generated by the elements $M_1$,
$M_{21}$, $M_{421}$, $M_{8421}$, .... It is then clear that the image of
$\psi_\chi$ lies in $\A_\#$ and the reduction $\aB_\#\onto\F_*$ carries
$\A_\#$ to $\A_*$. Moreover obviously the image of $\psi^\ox m_*$ lies in
$\A_\#$, hence it only remains to show the inclusion
$$
m_*^\GG(\A_\#)\subset\aB_\#\ox\A_\#.
$$
Since $m_*^\GG$ is a ring homomorphism, it suffices to check this on the
generators $M_1$, $M_{21}$, $M_{421}$, $M_{8421}$, .... But this is clear from
\eqref{dizeta}.
\end{proof}

\begin{Corollary}\label{calcab}
For the comultiplicative splitting $\psi$ one has
$$
a_\psi=0.
$$
Moreover the map $b_\psi$ factors as follows
$$
\A_*\xto{\bar b_\psi}{R_\F}_*\ox\A_*\xto{1\ox i}{R_\F}_*\ox\F_*.
$$
\end{Corollary}

\begin{proof}
The first statement follows as by definition $\pi(\A_*)=0$; the second is
obvious.
\end{proof}

Using the splitting $\psi$ we get the following analogue of \bref{detrel}.

\begin{Proposition}\label{psicomon}
The comonoid $\R^\FF$ in the category $\Alg^\pair_\oo$ described in
\bref{relcom} is completely determined by the multiplication map
$$
\bar b^\psi: R_\F\ox\A\to\A
$$
dual to the map $\bar b_\psi$ from \ref{calcab}. In fact, the identification
$$
\R^\FF_1=\A\oplus R_\F
$$
induced by the splitting $s=\psi$ identifies the diagonal of $\R^\FF$ with
$\Delta_\A\oplus\Delta_R$ (see \eqref{ueq}, \eqref{udia}), and the bimodule
structure of $\R_1^\FF$ with
\begin{align*}
f(\alpha,r)&=(f\alpha,fr)\\
(\alpha,r)f&=(\alpha\qf f-\bar b^\psi(r,\qf f),rf)
\end{align*}
for $f\in\F_0$, $r\in R_\F$, $\alpha\in\A$.
\end{Proposition}

\section{Computation of the Hopf pair algebra $\aB^\FF$}\label{bcomp}

The Hopf pair algebra $\V=\aB^\FF$ in \bref{unique}, given by the algebra of
secondary cohomology operations, satisfies the following crucial condition
which we deduce from \cite{Baues}*{16.1.5}.

\begin{Theorem}\label{splitr}
There exists a right $\F_0$-equivariant splitting
$$
u:\R^\FF_1=R_\aB\ox\FF\to\aB_1\ox\FF=\aB^\FF_1
$$
of the projection $\aB^\FF_1\to\R^\FF_1$, see \eqref{hpad}, such that the
following holds. The diagram
$$
\xymatrix{
\A\oplus_\k\Sigma\A\ar@/_/@{-->}[d]_q\ar@{ >->}[r]
&\aB^\FF_1\ar@/_/@{-->}[d]_q\ar[r]
&\aB^\FF_0\ar@{=}[d]\ar@{->>}[r]
&\A\ar@{=}[d]\\
\A\ar[u]_{\bar u}\ar@{ >->}[r]
&\R^\FF_1\ar[u]_u\ar[r]
&\R^\FF_0\ar@{->>}[r]
&\A
}
$$
commutes, where $\bar u$ is the inclusion. Moreover in the diagram of
diagonals, see \eqref{diacomp},
$$
\xymatrix{
\aB^\FF_1\ar[r]^-{\Delta_\aB}
&(\aB^\FF\hat\ox\aB^\FF)_1
&\Sigma\A\ox\A\ar@{ )->}[l]\\
\R^\FF_1\ar[r]^-{\Delta_R}\ar[u]^u
&(\R^\FF\hat\ox\R^\FF)_1\ar[u]^{u\hat\ox u}
}
$$
the difference $\Delta_\aB u-(u\hat\ox u)\Delta_R$ lifts to $\Sigma\A\ox\A$
and satisfies
$$
\xi\bar\pi=\Delta_\aB u-(u\hat\ox
u)\Delta_R:\xymatrix@1{\R_\FF^1\ar@{->>}[r]^{\bar\pi}&\bar
R\ar[r]^-\xi&\Sigma\A\ox\A}
$$
where $\xi$ is dual to $\xi_*$ in \bref{cosxi}. Here $\bar\pi$ is the
projection $\R_\FF\onto R_\F\onto\bar R$. The cocycle $\xi$ is trivial if
$p$ is odd.
\end{Theorem}

\begin{Definition}\label{multop}
Using a splitting $u$ of $\aB^\FF$ as in \bref{splitr} we define a
\emph{multiplication operator}
$$
A:\A\ox R_\aB\to\Sigma\A
$$
by the equation
$$
A(\bar\alpha\ox x)=u(\alpha x)-\alpha u(x)
$$
for $\alpha\in\F_0$, $x\in R_\aB$. Thus $-A$ is a multiplication map as
studied in \cite{Baues}*{16.3.1}. Fixing a splitting $s$ of $\R^\FF$ as in
\eqref{s} we define an \emph{$s$-multiplication operator} $A^s$ to be the
composite
$$
A^s:\xymatrix@1{\A\ox R_\F\ar[r]^-{1\ox s}&\A\ox R_\aB\ar[r]^-A&\Sigma\A}.
$$
Such operators have the properties of the following $s$-multiplication
maps.
\end{Definition}

\begin{Definition}\label{mulmap}
Let $s$ be a splitting of $\R^\FF$ as in \eqref{s} and let $U^s$, $a^s$,
$b^s$ be defined as in section \ref{rcomp}. An \emph{$s$-multiplication
map}
$$
A^s:\A\ox R_\F\to\A
$$
is an $\FF$-linear map of degree $-1$ satisfying the following conditions
with $\alpha,\alpha',\beta,\beta'\in\F_0$, $x,y\in R_\F$
\begin{enumerate}
\item $A^s(\alpha,x\beta)=A^s(\alpha,x)\beta+\k(\alpha)b^s(x,\beta)$
\item
$A^s(\alpha\alpha',x)=A^s(\alpha,\alpha'x)+\k(\alpha)a^s(\alpha',x)+(-1)^{\deg(\alpha)}\alpha
A^s(\alpha',x)$
\item $\delta A^s(\alpha,x)=A^s_\ox(\alpha\ox\Delta
x)+L(\alpha,x)+\nabla_\xi(\alpha,x)+\delta\k(\alpha)U^s(x)$.
\end{enumerate}
Here $A^s_\ox:\A\ox R_\F^{(2)}\to\A\ox\A$ is defined by the equalities
\begin{align*}
A^s_\ox(\alpha\ox x\ox\beta')
&=\sum(-1)^{\deg(\alpha_\r)\deg(x)}A^s(\alpha_\l,x)\ox\alpha_\r\beta',\\
A^s_\ox(\alpha\ox\beta\ox y)
&=\sum(-1)^{\deg(\alpha_\r)\deg(\beta)+\deg(\alpha_\l)+\deg(\beta)}\alpha_\l\beta\ox
A^s(\alpha_\r,y),
\end{align*}
where as always
$$
\delta(\alpha)=\sum\alpha_\l\ox\alpha_\r\in\A\ox\A.
$$

Two $s$-multiplication maps $A^s$ and ${A^s}'$  are \emph{equivalent} if there exists an
$\FF$-linear map
$$
\gamma:R_\F\to\A
$$
of degree $-1$ such that the equality
$$
A^s(\alpha,x)-{A^s}'(\alpha,x)
=\gamma(\alpha x)-(-1)^{\deg(\alpha)}\alpha\gamma(x)
$$
holds for any $\alpha\in\A$, $x\in R_\F$ and moreover $\gamma$ is right
$\F_0$-equivariant and the diagram
$$
\xymatrix{
\A\ar[r]^-\delta&\A\ox\A\\
R_\F\ar[u]_\gamma\ar[r]_\Delta&R_\F^{(2)}\ar[u]_{\gamma_\ox}
}
$$
commutes, with $\gamma_\ox$ given by
\begin{align*}
\gamma_\ox(x\ox\beta)&=\gamma(x)\ox\beta,\\
\gamma_\ox(\alpha\ox y)&=(-1)^{\deg(\alpha)}\alpha\ox\gamma(y)
\end{align*}
for $\alpha,\beta\in\F_0$, $x,y\in R_\F$.
\end{Definition}

\begin{Theorem}\label{exmul}
There exists an $s$-multiplication map $A^s$ and any two such
$s$-multiplication maps are equivalent. Moreover each $s$-multiplication
map is an $s$-multiplication operator as in \bref{multop} and vice versa.
\end{Theorem}

\begin{proof}
We apply \cite{Baues}*{16.3.3}. In fact, we obtain by $A^s$ the
multiplication operator
$$
A:\A\ox R_\aB=\A\!\ox\!\A\oplus\A\!\ox\!R_\F\to\Sigma\A
$$
with
\begin{equation}\label{exmulf}
A(\alpha\ox x)=A^s(\alpha\ox\bar x)+\k(\alpha)\xi
\end{equation}
where $(\bar x,\xi)\in R_\F\oplus\A=R_\aB\ox\FF$ corresponds to $x$, that is $s(\bar x)+\iota(\xi)=x$ for $\iota:\A\subset R_\aB\ox\FF$.
\end{proof}

\begin{Remark}
For the splitting $s=\phi$ of $\R^\FF$ in \bref{chi} the maps
$$
A_{n,m}:\A\to\A
$$
are defined by $A_{n,m}(\alpha)=A^\phi(\alpha\ox[n,m])$, with $[n,m]$
the Adem relations in $R_\F$. Using formul\ae\ in \bref{mulmap} the maps
$A_{n,m}$ determine the $\phi$-multiplication map $A^\phi$ completely. The
maps $A_{n,m}$ coincide with the corresponding maps $A_{n,m}$ in
\cite{Baues}*{16.4.4}. In \cite{Baues}*{16.6} an algorithm for
determination of $A_{n,m}$ is described, leading to a list of values of
$A_{n,m}$ on the elements of the admissible basis of $\A$. The algorithm
for the computation of $A_{n,m}$ can be deduced from theorem \bref{exmul}
above.
\end{Remark}

\begin{Remark}
Triple Massey products $\brk{\alpha,\beta,\gamma}$ with
$\alpha,\beta,\gamma\in\A$, $\alpha\beta=0=\beta\gamma$, as in \bref{tmp}
can be computed by $A^s$ as follows. Let $\bar\beta\bar\gamma\in R_\aB$ be
given as in \bref{tmp}. Then $\bar\beta\bar\gamma\ox1\in R_\aB\ox\FF$
satisfies
$$
\bar\beta\bar\gamma\ox1=s(\bar x)+\iota(\xi)
$$
with $\bar x\in R_\F$, $\xi\in\A$ and $\brk{\alpha,\beta,\gamma}$ satisfies
$$
A^s(\alpha\ox\bar x)+\k(\alpha)\xi\in\brk{\alpha,\beta,\gamma}.
$$
Compare \cite{Baues}*{16.3.4}.
\end{Remark}

Now it is clear how to introduce via $a^s$, $b^s$, $U^s$, $\xi$, $\k$, and $A^s$ a
Hopf pair algebra structure on
\begin{equation}\label{hopfin}
\alignbox{
\xymatrix{
\A\oplus\Sigma\A\oplus R_\F\ar[r]^-q\ar@{=}[d]&\A\oplus R_\F\ar@{=}[d]\\
\aB^\FF_1&\R^\FF_1
}
}
\end{equation}
which is isomorphic to $\aB^\FF$, compare \bref{detrel}.

In the next section we describe an algorithm for the computation of a
$\psi$-multiplication map, where $\psi$ is the comultiplicative splitting
of $\R^\FF$ in \bref{psi}. For this we compute the dual map $A_\psi$ of
$A^\psi$.

\section{Computation of the Hopf pair coalgebra $\aB_\FF$}\label{cobcomp}

For the comultiplicative splitting $s=\psi$ of $\R^\FF$ in \bref{psi} we
introduce the following $\psi$-comultiplication maps which are dual to the
$\psi$-multiplication maps in \bref{mulmap}.

\begin{Definition}\label{apsi}
Let $\bar b_\psi$ be given as in \ref{calcab}. A
\emph{$\psi$-comultiplication map}
$$
A_\psi:\A_*\to\A_*\ox{R_\F}_*
$$
is an $\FF$-linear map of degree $+1$ satisfying the following conditions.
\begin{enumerate}
\item\label{mreqs}
The maps in the diagram
$$
\xymatrix{
\A_*\ox{R_\F}_*\ar[d]_{1\ox m^\r_*}&\A_*\ar[d]^{m_*}\ar[l]_-{A_\psi}\\
\A_*\ox{R_\F}_*\ox\F_*&\A_*\ox\A_*\ar[l]^-{A_\psi\ox i}
}
$$
satisfy
$$
(1\ox m^\r_*)A_\psi=(A_\psi\ox i)m_*+(\k_*\ox\bar b_\psi)m_*.
$$
Here $\k_*$ is computed in \bref{dkappa} and $m^\r_*$ is defined in
\eqref{mr}.
\item\label{mleqs}
The maps in the diagram
$$
\xymatrix{
\A_*\ox{R_\F}_*\ar[d]^{1\ox m^\l_*}&&\A_*\ar[d]^{A_\psi}\ar[ll]_{A_\psi}\\
\A_*\ox\F_*\ox{R_\F}_*&\A_*\ox\A_*\ox{R_\F}_*\ar[l]^{1\ox i\ox1}
&\A_*\ox{R_\F}_*\ar[l]^-{m_*\ox1}
}
$$
satisfy
$$
(1\ox m^\l_*)A_\psi=(1\ox i\ox1)(m_*\ox1)A_\psi-(\tau\ox i\ox1)(1\ox A_\psi)m_*.
$$
Here $m^\l_*$ is as in \eqref{ml}, and $\tau:\A_*\to\A_*$ is given
by $\tau(\alpha)=(-1)^{\deg(\alpha)}\alpha$.
\item\label{mult}
For $x,y\in\A_*$ the product $xy$ in the algebra $\A_*$ satisfies the
formula
$$
A_\psi(xy)=A_\psi(x)m_*(y)+(-1)^{\deg(x)}m_*(x)A_\psi(y)+L_*(x,y)+{\nabla_\xi}_*(x,y).
$$
Here $L_*$ and ${\nabla_\xi}_*$ are given in \ref{L} and
\ref{nablaelts} respectively, with $L_*={\nabla_\xi}_*=0$ for $p$ odd.
\end{enumerate}

Two $\psi$-comultiplication maps $A_\psi$, $A_\psi'$ are
\emph{equivalent} if there is a derivation
$$
\gamma_*:\A_*\to{R_\F}_*
$$
of degree $+1$ satisfying the equality
$$
A_\psi-A_\psi'=m^\l_*\gamma_*-(\tau\ox\gamma_*)m_*.
$$
\end{Definition}

As a dual statement to \bref{exmul} we get

\begin{Theorem}
There exists a $\psi$-comultiplication map $A_\psi$ and any two such
$\psi$-comultiplication maps are equivalent. Moreover each
$\psi$-comultiplication map $A_\psi$ is the dual of a $\psi$-multiplication
map $A^\psi$ in \bref{exmul} with $A_\psi={A^\psi}_*$.
\end{Theorem}\qed

Now dually to \eqref{hopfin}, it is clear how to introduce via
$a_\psi$, $b_\psi$, $\xi_*$, $\k_*$, and $A_\psi$ a Hopf pair
coalgebra structure on
$$
\xymatrix{
\A_*\oplus\Sigma\A_*\oplus{R_\F}_*\ar@{=}[d]&\A_*\oplus{R_\F}_*\ar[l]_-i\ar@{=}[d]\\
\aB_\FF^1&\R_\FF^1
}
$$
which is isomorphic to $\aB_\FF$, compare \bref{psicomon}.

We now embark on the simplification and solution of the equations
\ref{apsi}\eqref{mreqs} and \ref{apsi}\eqref{mleqs}. To begin with, note that
the equations \ref{apsi}\eqref{mreqs} imply that the image of the composite
map
$$
\A_*\xto{A_\psi}\A_*\ox{R_\F}_*\xto{1\ox m^\r_*}\A_*\ox{R_\F}_*\ox\F_*
$$
actually lies in
$$
\A_*\ox{R_\F}_*\ox\A_*\subset\A_*\ox{R_\F}_*\ox\F_*;
$$
similarly \ref{apsi}\eqref{mleqs} implies that the image of
$$
\A_*\xto{A_\psi}\A_*\ox{R_\F}_*\xto{1\ox m^\l_*}\A_*\ox\F_*\ox{R_\F}_*
$$
lies in
$$
\A_*\ox\A_*\ox{R_\F}_*\subset\A_*\ox\F_*\ox{R_\F}_*.
$$

\begin{Lemma}
The following conditions on an element $x\in{R_\F}_*=\Hom(R_\F,\FF)$ are
equivalent:
\begin{itemize}
\item $m^\l_*(x)\in\A_*\ox{R_\F}_*\subset\F_*\ox{R_\F}_*$;
\item $m^\r_*(x)\in{R_\F}_*\ox\A_*\subset{R_\F}_*\ox\F_*$;
\item $x\in\bar R_*\subset{R_\F}_*$.
\end{itemize}
\end{Lemma}

\begin{proof}
Recall that $\bar R=R_\F/{R_\F}^2$, i.~e. $\bar R_*$ is the space of linear forms
on $R_\F$ which vanish on ${R_\F}^2$. Then the first condition means that
$x:R_\F\to\FF$ has the property that the composite
$$
\F_0\ox R_\F\xto{m^\l} R_\F\xto x\FF
$$
vanishes on $R_\F\ox R_\F\subset\F_0\ox R_\F$; but the image of $R_\F\ox R_\F$
under $m^\l$ is precisely ${R_\F}^2$. Similarly for the second condition.
\end{proof}

We thus conclude that the image of $A_\psi$ lies in $\A_*\ox\bar R_*$.

Next note that the condition \ref{apsi}\eqref{mult} implies
\begin{equation}\label{apsisq}
A_\psi(x^2)=L_*(x,x)+\nabla_{\xi_*}(x,x)
\end{equation}
for any $x\in\A_*$. Moreover the latter formula also implies

\begin{Proposition}\label{4=0}
For any $x\in\A_*$ one has
$$
A_\psi(x^4)=0.
$$
\end{Proposition}

\begin{proof}
Since the squaring map is an algebra endomorphism, by \ref{bider} one has
$$
L_*(x,y^2)=\sum\zeta_1x_\l y_{\l'}^2\ox\tilde L_*(x_r,y_{r'}^2),
$$
with
$$
m_*(x)=\sum x_\l\ox x_r,\ \ m_*(y)=\sum y_{\l'}\ox y_{r'}.
$$
But $\tilde L_*$ vanishes on squares since it is a biderivation, so $L_*$ also
vanishes on squares. Moreover by \eqref{nablaelts}
$$
\nabla_{\xi_*}(x^2,y^2)=\sum\xi_*(x^2,y^2)_\A\ox\xi_*(x^2,y^2)_R-\sum
x_\l^2y_{\l'}^2\ox\xi_*(x_r^2,y_{r'}^2);
$$
this is zero since $\xi_*(x^2,y^2)=0$ for any $x$ and $y$ by \eqref{xifromS}.
\end{proof}

Taking the above into account, and identifying the image of $i:\A_*\into\F_*$
with $\A_*$, \ref{apsi}\eqref{mreqs} can be rewritten as follows:
$$
(1\ox m_*^r)A_\psi(\zeta_n)
=A_\psi(\zeta_n)\ox1+\left(L_*(\zeta_{n-1},\zeta_{n-1})+\nabla_{\xi_*}(\zeta_{n-1},\zeta_{n-1})\right)\ox\zeta_1
+\sum_{i=0}^n\zeta_1\zeta_{n-i}^{2^i}\ox\bar b_\psi(\zeta_i),
$$
or
$$
(1\ox\tilde m_*^r)A_\psi(\zeta_n)
=\left(L_*(\zeta_{n-1},\zeta_{n-1})+\nabla_{\xi_*}(\zeta_{n-1},\zeta_{n-1})\right)\ox\zeta_1
+\sum_{i=0}^n\zeta_1\zeta_{n-i}^{2^i}\ox\bar b_\psi(\zeta_i).
$$
Still more explicitly one has
$$
L_*(\zeta_k,\zeta_k)=\sum_{0\le i,j\le k}\zeta_1\zeta_{k-i}^{2^i}\zeta_{k-j}^{2^j}\ox\tilde
L_*(\zeta_i,\zeta_j)=\sum_{0\le i\le k}\zeta_1\zeta_{k-i}^{2^{i+1}}\ox\tilde
L_*(\zeta_i,\zeta_i)
+\sum_{0\le i<j\le k}\zeta_1\zeta_{k-i}^{2^i}\zeta_{k-j}^{2^j}\ox\tilde
L^S_*(\zeta_i,\zeta_j),
$$
where we have denoted
$$
\tilde L^S_*(\zeta_i,\zeta_j):=\tilde
L_*(\zeta_i,\zeta_j)+\tilde L_*(\zeta_j,\zeta_i);
$$
similarly
$$
\nabla_{\xi_*}(\zeta_k,\zeta_k)=
\sum_{0\le i<j\le k}\zeta_{k-i}^{2^i}\zeta_{k-j}^{2^j}\ox
S_*(\zeta_i,\zeta_j).
$$
As for $b_\psi(\zeta_i)$, by \ref{mulstar} it can be calculated by the formula
\begin{equation}\label{bpsibar}
\bar b_\psi(\zeta_i)=\sum_{0<j<i}v_{i-j}^{2^{j-1}}\ox\zeta_j,
\end{equation}
where $v_k$ are determined by the equalities
$$
M_{2^k,2^{k-1},...,2}-M_{2^{k-1},2^{k-2},...,1}^2\equiv2v_k\mod4
$$
in $\aB_\#$. For example,
\begin{align*}
v_1&=M_{11},\\
v_2&=M_{411}+M_{231}+M_{222}+M_{2121},\\
v_3&
=M_{8411}
+M_{8231}
+M_{8222}
+M_{82121}
+M_{4631}
+M_{4622}
+M_{46121}
+M_{4442}
+M_{42521}
+M_{42431}
+M_{42422}\\
&+M_{424121}
+M_{421421},
\end{align*}
etc.

Thus putting everything together we see

\begin{Lemma}\label{mrC}
The equation \ref{apsi}\eqref{mreqs} for the
value on $\zeta_n$ is equivalent to
$$
(1\ox\tilde m^r_*)A_\psi(\zeta_n)=\sum_{0<k<n}C^{(n)}_{2^n-2^k+1}\ox\zeta_k
$$
where
\begin{multline*}
C^{(n)}_{2^n-1}=
\sum_{0<i<n}\zeta_1\zeta_{n-1-i}^{2^{i+1}}\ox\left(\tilde
L_*(\zeta_i,\zeta_i)+v_i\right)
+\sum_{0<i<j<n}\zeta_1\zeta_{n-1-i}^{2^i}\zeta_{n-1-j}^{2^j}\ox\tilde
L^S_*(\zeta_i,\zeta_j)
+\sum_{0<i<j<n}\zeta_{n-1-i}^{2^i}\zeta_{n-1-j}^{2^j}\ox
S_*(\zeta_i,\zeta_j)
\end{multline*}
and, for $1<k<n$,
$$
C^{(n)}_{2^n-2^k+1}=\sum_{0<i\le n-k}\zeta_1\zeta_{n-k-i}^{2^{k+i}}\ox
v_i^{2^{k-1}}.
$$
\end{Lemma}\qed

For low values of $n$ these equations look like
\begin{align*}
(1\ox\tilde m^r_*)A_\psi(\zeta_2)&=0,\\
(1\ox\tilde m^r_*)A_\psi(\zeta_3)&
=\zeta_1\ox(\pi(M_{222})\ox\zeta_1+\pi(M_{22})\ox\zeta_2)
+\zeta_1^2\ox\pi(M_{32}+M_{23}+M_{212}+M_{122})\ox\zeta_1\\
&+\zeta_1^3\ox\pi M_{22}\ox\zeta_1,\\
(1\ox\tilde m^r_*)A_\psi(\zeta_4)&
=\zeta_1\ox\left(\pi(M_{8222}+M_{722}+M_{4622}+M_{4442}+M_{42422})\ox\zeta_1\right.\\
&\ \ \ \ \ \ \ \left.+\pi(M_{822}+M_{462}+M_{444}+M_{4242})\ox\zeta_2+\pi(M_{44})\ox\zeta_3\right)\\
&+\zeta_1^4\ox\pi(M_{632}+M_{623}+M_{6212}+M_{6122}+M_{542}+M_{452}+M_{443}+M_{4412}+M_{4142}+M_{3422}\\
&\ \ \ \ \ \ \ +M_{2522}+M_{2432}+M_{2423}+M_{24212}+M_{24122}+M_{21422}+M_{1622}+M_{1442}+M_{12422})\ox\zeta_1\\
&+\zeta_1^5\ox\pi(M_{622}+M_{442}+M_{2422})\ox\zeta_1\\
&+\zeta_2^2\ox\pi(M_{522}+M_{432}+M_{423}+M_{4212}+M_{4122}+M_{1422})\ox\zeta_1
+\zeta_1\zeta_2^2\ox\pi(M_{422})\ox\zeta_1\\
&+\zeta_1^9\ox\left(\pi(M_{222})\ox\zeta_1+\pi(M_{22})\ox\zeta_2\right)
+\zeta_1^4\zeta_2^2\ox\pi(M_{32}+M_{23}+M_{212}+M_{122})\ox\zeta_1\\
&+\zeta_1^5\zeta_2^2\ox\pi(M_{22})\ox\zeta_1,
\end{align*}
etc. (Note that $A_\psi(\zeta_1)=0$ by dimension considerations.)

As for the equations \ref{apsi}\eqref{mleqs}, they have form
$$
(1\ox\tilde m^\l_*)A_\psi(\zeta_n)
=(\tilde m_*\ox1)A_\psi(\zeta_n)+\zeta_1^{2^{n-1}}\ox A_\psi(\zeta_{n-1})+\zeta_2^{2^{n-2}}\ox
A_\psi(\zeta_{n-2})+...+\zeta_{n-2}^4\ox A_\psi(\zeta_2)+\zeta_{n-1}^2\ox
A_\psi(\zeta_1).
$$

\begin{Lemma}
Suppose given a map $A_\psi$ satisfying \ref{apsi}\eqref{mult} and those
instances of \ref{apsi}\eqref{mreqs}, \ref{apsi}\eqref{mleqs} which involve
starting value of $\A_\psi$ on the Milnor generators $i(\zeta_1)$,
$i(\zeta_2)$, ..., where $i:\A_*\to\F_*$ is the inclusion. Then $\A_\psi$
satisfies these equations for all other values too.
\end{Lemma}

Now recall that, as already mentioned in \ref{L*}, according to
\cite{Baues}*{16.5} $\bar R$ is a free right $\A$-module generated by the set
$\PAR\subset\bar R$ of preadmissible relations. More explicitly, the composite
$$
R^\pre\ox\A\xto{\textrm{inclusion}\ox1}\bar R\ox\A\xto{m^\r}\bar R
$$
is an isomorphism of right $\A$-modules, where $R^\pre$ is the $\FF$-vector
space spanned by the set $\PAR$ of preadmissible relations.
Dually it follows that the composite
$$
\Phi^\r_*:\bar R_*\xto{m^\r_*}\bar R_*\ox\A_*\xto{\ro\ox1}R_\pre\ox\A_*
$$
is an isomorphism of right $\A_*$-comodules. Here $\ro:\bar R_*\onto R_\pre$ denotes
the restriction homomorphism from the space $\bar R_*$ of $\FF$-linear forms
on $\bar R$ to the space $R_\pre$ of linear forms on its subspace
$R^\pre\subset\bar R$ spanned by $\PAR$.

It thus follows that we will obtain equations equivalent to \ref{apsi}\eqref{mreqs} if
we compose both sides of these equations with the isomorphism
$1\ox\Phi^\r_*:\A_*\ox\bar R_*\to\A_*\ox R_\pre\ox\A_*$. Let us then denote
$$
(1\ox\Phi^\r_*)A_\psi(\zeta_n)=\sum_\mu\rho_{2^n-|\mu|}(\mu)\ox\mu
$$
with some unknown elements $\rho_{j}(\mu)\in(\A_*\ox R_\pre)_j$, where
$\mu$ runs through some basis of $\A_*$.

Now freedom of the right $\A_*$-comodule $\bar R_*$ on $R_\pre$ means that the
above isomorphism $\Phi^\r_*$ fits in the commutative diagram
$$
\xymatrix{
\bar R_*\ar[r]^-{\Phi^\r_*}\ar[d]^{m^\r_*}&R_\pre\ox\A_*\ar[d]^{1\ox m_*}\\
\bar R_*\ox\A_*\ar[r]^-{\Phi^\r_*\ox1}&R_\pre\ox\A_*\ox\A_*.
}
$$
It follows that we have
$$
(1\ox1\ox m_*)(1\ox\Phi^\r_*)A_\psi(\zeta_n)
=(1\ox\Phi^\r_*\ox1)(1\ox m^r_*)A_\psi(\zeta_n).
$$

Then taking into account \ref{mrC} this gives equations
$$
\sum_\mu\rho_{2^n-|\mu|}(\mu)\ox m_*(\mu)=
\sum_\mu\rho_{2^n-|\mu|}(\mu)\ox\mu\ox1+\sum_{0<k<n}(1\ox\Phi^\r_*)(C^{(n)}_{2^n-2^k+1})\ox\zeta_k,
$$
with the constants $C^{(j)}_n$ as in \ref{mrC}. This immediately determines
the elements $\rho_j(\mu)$ for $|\mu|>0$. Indeed, the above equation implies
that $(1\ox\Phi^\r_*)A_\psi(\zeta_n)$ actually lies in the subspace $\A_*\ox
R_\pre\ox\Pi\subset\A_*\ox R_\pre\ox\A_*$ where $\Pi\subset\A_*$ is the
following subspace:
$$
\Pi=\set{x\in\A_*\ \mid\ m_*(x)\in\bigoplus_{k\ge0}\A_*\ox\FF\zeta_k}.
$$
It is easy to see that actually
$$
\Pi=\bigoplus_{k\ge0}\FF\zeta_k,
$$
so we can write
$$
(1\ox\Phi^\r_*)A_\psi(\zeta_n)=\sum_{k\ge0}\rho_{2^n-2^k+1}(\zeta_k)\ox\zeta_k
$$
where we necessarily have
$$
\rho_{2^n-2^k+1}(\zeta_k)\ox1+\rho_{2^n-2^{k+1}+1}(\zeta_{k+1})\ox\zeta_1^{2^k}+\rho_{2^n-2^{j+k}+1}(\zeta_{k+2})\ox\zeta_2^{2^k}+...
=(1\ox\Phi_\r)(C^{(n)}_{2^n-2^k+1}).
$$
for all $k\ge1$. By dimension considerations, $\rho_{2^n-2^k+1}(\zeta_k)$ can only be nonzero
for $k<n$, so the number of unknowns in these equations
strictly decreases as $k$ grows. Thus moving ``backwards'' and using
successive elimination we determine all $\rho_{2^n-2^k+1}(\zeta_k)$ for $k>0$.

It is easy to compute values of the isomorphism $1\ox\Phi^\r_*$ on all elements
involved in the constants $C^{(n)}_j$. In particular, elements of the form
$\Phi^\r_*(v_j^{2^k})$ can be given by an explicit formula. One has
$$
\Phi^\r_*(v_k)=\sum_{0\le i<k}\left(\Sq^{2^k}\Sq^{2^{k-1}}\cdots\Sq^{2^{i+2}}[2^i,2^i]\right)_*\ox\zeta_i^2
$$
and
$$
\Phi^\r_*(v_k^{2^{j-1}})
=\sum_{0\le i<k}\left(\Sq^{2^{k+j-1}}\Sq^{2^{k+j-2}}\cdots\Sq^{2^{i+j+1}}[2^{i+j-1},2^{i+j-1}]\right)_*\ox\zeta_i^{2^j},
$$
so our ``upside-down'' solving gives
\begin{align*}
\rho_{2^{n-1}+1}(\zeta_{n-1})&=\zeta_1\ox[2^{n-2},2^{n-2}]_*,\\
\rho_{2^n-2^{n-2}+1}(\zeta_{n-2})&=\zeta_1^{1+2^{n-1}}\ox[2^{n-3},2^{n-3}]_*+\zeta_1\ox\left(\Sq^{2^{n-1}}[2^{n-3},2^{n-3}]\right)_*\\
\rho_{2^n-2^{n-3}+1}(\zeta_{n-3})&=\zeta_1\zeta_2^{2^{n-2}}\ox[2^{n-4},2^{n-4}]_*+\zeta_1^{1+2^{n-1}}\ox\left(\Sq^{2^{n-2}}[2^{n-4},2^{n-4}]\right)_*+\zeta_1\ox\left(\Sq^{2^{n-1}}\Sq^{2^{n-2}}[2^{n-4},2^{n-4}]\right)_*\\
\cdots\\
\rho_{2^n-2^{n-k}+1}(\zeta_{n-k})
&=\sum_{1\le i\le k}\zeta_1\zeta_{k-i}^{2^{n-k+i}}\ox\left(\Sq^{2^{n-k+i-1}}\Sq^{2^{n-k+i-2}}\cdots\Sq^{2^{n-k+1}}[2^{n-k-1},2^{n-k-1}]\right)_*
\end{align*}
for $k<n-1$.

As for $\rho_{2^n-1}(\zeta_1)$, here we do not have a general formula, but
nevertheless it is easy to compute this value explicitly. In this way we obtain, for example,
\begin{align*}
\rho_1(\zeta_1)&=0,\\
\rho_3(\zeta_1)&=0,\\
\rho_7(\zeta_1)
&
=\zeta_1^3\ox[2,2]_*
+\zeta_1^2\ox\left([3,2]_*+[2,3]_*\right),\\
\rho_{15}(\zeta_1)
&
=\zeta_1^5\zeta_2^2\ox[2,2]_*
+\zeta_1^4\zeta_2^2\ox\left([3,2]_*+[2,3]_*\right)
+\zeta_1\zeta_2^2\ox\left(\Sq^4[2,2]\right)_*
+\zeta_2^2\ox\left((\Sq^5[2,2])_*+(\Sq^4[2,3])_*\right)\\
&
+\zeta_1^5\ox\left(\Sq^6[2,2]\right)_*
+\zeta_1^4\ox\left((\Sq^7[2,2])_*+(\Sq^6[3,2])_*+(\Sq^6[2,3])_*\right),\\
\rho_{31}(\zeta_1)
&
=\zeta_1\zeta_2^4\zeta_3^2\ox[2,2]_*
+\zeta_2^4\zeta_3^2\ox\left([3,2]_*+[2,3]_*\right)
+\zeta_1^9\zeta_3^2\ox\left(\Sq^4[2,2]\right)_*\\
&
+\zeta_1^8\zeta_3^2\ox\left((\Sq^5[2,2])_*+(\Sq^4[2,3])_*\right)
+\zeta_1^9\zeta_2^4\ox\left(\Sq^6[2,2]\right)_*\\
&
+\zeta_1^8\zeta_2^4\ox\left((\Sq^7[2,2])_*+(\Sq^6[3,2])_*+(\Sq^6[2,3])_*\right)
+\zeta_1\zeta_3^2\ox\left(\Sq^8\Sq^4[2,2]\right)_*\\
&
+\zeta_3^2\ox\left((\Sq^9\Sq^4[2,2])_*+(\Sq^8\Sq^4[2,3])_*\right)
+\zeta_1\zeta_2^4\ox\left(\Sq^{10}\Sq^4[2,2]\right)_*\\
&
+\zeta_2^4\ox\left((\Sq^{11}\Sq^4[2,2])_*+(\Sq^{10}\Sq^5[2,2])_*+(\Sq^{10}\Sq^4[2,3])_*\right)
+\zeta_1^9\ox\left(\Sq^{12}\Sq^6[2,2]\right)_*\\
&
+\zeta_1^8\ox\left((\Sq^{13}\Sq^6[2,2])_*+(\Sq^{12}\Sq^6[3,2])_*+(\Sq^{12}\Sq^6[2,3])_*\right),
\end{align*}
etc.

To summarize, let us state

\begin{Proposition}\label{mrrho}
The general solution of \ref{apsi}\eqref{mreqs} for the value on $\zeta_n$ is
given by the formula
$$
A_\psi(\zeta_n)=(1\ox\Phi^\r_*)\1\sum_{k\ge0}\rho_{2^n-2^k+1}(\zeta_k)\ox\zeta_k,
$$
where the elements $\rho_j(\zeta_k)\in(\A_*\ox R_\pre)_j$ are the ones explicitly
given above for $k>0$ while $\rho_{2^n}(1)\in(\A_*\ox R_\pre)_{2^n}$ is
arbitrary.
\end{Proposition}\qed

\

Let us now treat the equations \ref{apsi}\eqref{mleqs} in a similar way, now using the
fact that $\bar R$ is a free \emph{left} $\A$-module on an explicit basis
$\PAR'$ (see \ref{barr} again).

Then similarly to the above dualization it follows that the
composite
$$
\Phi^\l_*:\bar R_*\xto{m^\l_*}\A_*\ox\bar R_*\xto{1\ox\ro'}\A_*\ox R'_\pre
$$
is an isomorphism of left $\A_*$-comodules, where $\ro':\bar R_*\onto
R'_\pre$ denotes the restriction homomorphism from the space $\bar R_*$ of
$\FF$-linear forms on $\bar R$ to the space $R'_\pre$ of linear forms on the
subspace ${R^\pre}'$ of $\bar R$ spanned by $\PAR'$.

Thus similarly to the above the equations \ref{apsi}\eqref{mleqs} are
equivalent to ones obtained by composing them with the isomorphism
$1\ox\Phi^\l_*:\A_*\ox\bar R_*\to\A_*\ox\A_*\ox R'_\pre$. Let us then denote
$$
(1\ox\Phi^\l_*)A_\psi(\zeta_n)=\sum_{\pi\in\PAR'}\sigma_{2^n-|\pi|}(\pi)\ox\pi_*
$$
with some unknown elements $\sigma_j(\pi)\in(\A_*\ox\A_*)_j$, where $\pi_*$
denotes the corresponding element of the dual basis, i.~e. the unique linear form on
$R'_\pre$ assigning 1 to $\pi$ and 0 to all other elements of $\PAR'$.

Now again as above, freedom of the left $\A_*$-comodule $\bar R_*$ on $R'_\pre$ means that the
above isomorphism $\Phi^\l_*$ fits in the commutative diagram
$$
\xymatrix{
\bar R_*\ar[r]^-{\Phi^\l_*}\ar[d]^{m^\l_*}&\A_*\ox R'_\pre\ar[d]^{m_*\ox1}\\
\A_*\ox\bar R_*\ar[r]^-{1\ox\Phi^\l_*}&\A_*\ox\A_*\ox R'_\pre.
}
$$
In particular one has
$$
(1\ox1\ox\Phi^\l_*)(1\ox m^\l_*)A_\psi(\zeta_n)
=(1\ox m_*\ox1)(1\ox\Phi^\l_*)A_\psi(\zeta_n).
$$
Using this, we obtain that the equations \ref{apsi}\eqref{mleqs} are
equivalent to the following system of equations
$$
(1\ox m_*-m_*\ox1)(\sigma_{2^n-|\pi|}(\pi))
=1\ox\sigma_{2^n-|\pi|}(\pi)+\varSigma_{2^n-|\pi|}(\pi),
$$
where we denote
$$
\varSigma_{2^n-|\pi|}(\pi)=\zeta_1^{2^{n-1}}\ox\sigma_{2^{n-1}-|\pi|}(\pi)+\zeta_2^{2^{n-2}}\ox\sigma_{2^{n-2}-|\pi|}(\pi)+...+\zeta_{n-2}^4\ox\sigma_{4-|\pi|}(\pi)+\zeta_{n-1}^2\ox\sigma_{2-|\pi|}(\pi).
$$

We next use the following standard fact:

\begin{Proposition}\label{contra}
For any coalgebra $C$ with the diagonal $m_*:C\to C\ox C$ and counit
$\eps:C\to\FF$ there is a contractible cochain complex of the form
$$
\xymatrix{
C\ar[r]^{d_1}
&C^{\ox2}\ar@/^/[l]^{s_1}\ar[r]^{d_2}
&C^{\ox3}\ar@/^/[l]^{s_2}\ar[r]^{d_3}
&C^{\ox4}\ar@/^/[l]^{s_3}\ar[r]^{d_4}
&\cdots,\ar@/^/[l]^{s_4}
}
$$
i.~e. one has
$$
s_nd_n+d_{n-1}s_{n-1}=1_{C^{\ox n}}
$$
for all $n$. Here,
\begin{align*}
d_1&=m_*,\\
d_2&=1\ox m_*-m_*\ox1,\\
d_3&=1\ox1\ox m_*-1\ox m_*\ox1+m_*\ox1\ox1,\\
d_4&=1\ox1\ox1\ox m_*-1\ox1\ox m_*\ox1+1\ox m_*\ox1\ox1-m_*\ox1\ox1\ox1,
\end{align*}
etc., while $s_n$ can be taken to be equal to either
$$
s_n=\eps\ox 1_{C^{\ox n}}
$$
or
$$
s_n=1_{C^{\ox n}}\ox\eps.
$$
\end{Proposition}
\qed

Now suppose given the elements
$\sigma_{2^k-|\pi|}(\pi)$, $k<n$, satisfying the equations; we must then find
$\sigma_{2^n-|\pi|}(\pi)$ with
$$
d_2\sigma_{2^n-|\pi|}(\pi)=1\ox\sigma_{2^n-|\pi|}(\pi)+\varSigma_{2^n-|\pi|}(\pi),
$$
with $\varSigma_{2^n-|\pi|}(\pi)$ as above. Then since $d_3d_2=0$, it will follow
$$
d_3(1\ox\sigma_{2^n-|\pi|}(\pi)+\varSigma_{2^n-|\pi|}(\pi))=0.
$$
Then
$$
1\ox\sigma_{2^n-|\pi|}(\pi)+\varSigma_{2^n-|\pi|}(\pi)
=(s_3d_3+d_2s_2)(1\ox\sigma_{2^n-|\pi|}(\pi)+\varSigma_{2^n-|\pi|}(\pi))
=d_2s_2(1\ox\sigma_{2^n-|\pi|}(\pi)+\varSigma_{2^n-|\pi|}(\pi))
$$
Taking here $s_n$ from the second equality of \ref{contra}, we see that one
has
$$
1\ox\sigma_{2^n-|\pi|}(\pi)
=\varSigma_{2^n-|\pi|}(\pi)+d_2\left(1\ox(1\ox\eps)(\sigma_{2^n-|\pi|}(\pi))+(1\ox1\ox\eps)(\varSigma_{2^n-|\pi|}(\pi))\right).
$$
It follows that we can reconstruct the terms $\sigma_{2^n-|\pi|}(\pi)$ from
$(1\ox\eps)\sigma_{2^n-|\pi|}(\pi)$, i.~e. from their components that lie in
$\A_*\ox\FF\subset\A_*\ox\A_*$.

Then denoting
$$
\sigma_{2^n-|\pi|}(\pi)=x_{2^n-|\pi|}(\pi)\ox1+\sigma'_{2^n-|\pi|}(\pi),
$$
with
$$
\sigma'_{2^n-|\pi|}(\pi)\in\A_*\ox\tilde\A_*,
$$
the last equation gives
$$
1\ox x_{2^n-|\pi|}(\pi)\ox1+1\ox\sigma'_{2^n-|\pi|}(\pi)
=\varSigma_{2^n-|\pi|}(\pi)+(m_*\ox1+1\ox m_*)\sum_{i\ge0}\zeta_i^{2^{n-i}}\ox
x_{2^{n-i}-|\pi|}(\pi).
$$
By collecting terms of the form $1\ox...$ on both sides, we conclude that any solution for $\sigma$ satisfies
$$
\sigma_{2^n-|\pi|}(\pi)
=m_*(x_{2^n-|\pi|}(\pi))+\sum_{i\ge0}\zeta_i^{2^{n-i}}\ox
x_{2^{n-i}-|\pi|}(\pi).
$$
Thus the equation \ref{apsi}\eqref{mleqs} is equivalent to the system of
equations
$$
(1\ox m_*+m_*\ox1)\sum_{i\ge0}\zeta_i^{2^{n-i}}\ox x_{2^{n-i}-|\pi|}(\pi)
=1\ox m_*(x_{2^n-|\pi|}(\pi))+\sum_{i\ge0}1\ox\zeta_i^{2^{n-i}}\ox x_{2^{n-i}-|\pi|}(\pi)
+\varSigma_{2^n-|\pi|}(\pi)
$$
on the elements $x_j(\pi)\in\A_j$. Substituting here back the value of
$\varSigma_{2^n-|\pi|}(\pi)$ we obtain the equations
\begin{multline*}
\sum_{i\ge0}\zeta_i^{2^{n-i}}\ox m_*(x_{2^{n-i}-|\pi|}(\pi))
+\sum_{i\ge0}m_*(\zeta_i)^{2^{n-i}}\ox x_{2^{n-i}-|\pi|}(\pi)
=1\ox m_*(x_{2^n-|\pi|}(\pi))+\sum_{i\ge0}1\ox\zeta_i^{2^{n-i}}\ox
x_{2^{n-i}-|\pi|}(\pi)\\
+\sum_{i>0}\zeta_i^{2^{n-i}}\ox m_*(x_{2^{n-i}-|\pi|}(\pi))
+\sum_{i'>0,j\ge0}\zeta_{i'}^{2^{n-i'}}\ox\zeta_j^{2^{n-i'-j}}\ox x_{2^{n-i'-j}-|\pi|}(\pi).
\end{multline*}
These equations easily reduce to
$$
m_*(\zeta_i)^{2^{n-i}}=1\ox\zeta_i^{2^{n-i}}+\sum_{0\le
j<i}\zeta_{i-j}^{2^{n-(i-j)}}\ox\zeta_j^{2^{n-i}},
$$
which is identically true. We thus conclude

\begin{Proposition}\label{mlx}
The general solution $A_\psi(\zeta_n)$ of \ref{apsi}\eqref{mleqs} is determined by
$$
A_\psi(\zeta_n)=(1\ox\Phi^\l_*)\1\sum_{\pi\in\PAR'}\left(x_{2^n-|\pi|}(\pi)\ox1+\tilde
m_*(x_{2^n-|\pi|}(\pi))+\sum_{i>0}\zeta_i^{2^{n-i}}\ox
x_{2^{n-i}-|\pi|}(\pi)\right)\ox\pi_*,
$$
where $x_j(\pi)\in\A_j$ are arbitrary homogeneous elements.
\end{Proposition}\qed

Now to put together \ref{mrrho} and \ref{mlx} we must use the dual
$$
\Phi_*:R_\pre\ox\A_*\to\A_*\ox R'_\pre
$$
of the composite isomorphism
$$
\Phi:\A\ox{R^\pre}'\xto{{\Phi^\l}\1}\bar R\xto{\Phi^\r}R^\pre\ox\A.
$$

We will need

\begin{Lemma}
There is an inclusion
$$
\Phi_*\left(R_\pre\ox\FF1\right)\subset\A_*\ox{R'_\pre}^{\le2},
$$
where
$$
{R'_\pre}^{\le2}\subset R'_\pre
$$
is the subspace of those linear forms on ${R^\pre}'$ which vanish on all left
preadmissible elements $[n,m]a\in\PAR'$ with $a\in\tilde\A$.

Similarly, there is an inclusion
$$
\Phi_*\1\left(\FF1\ox R'_\pre\right)\subset{R_\pre}^{\le2}\ox\A_*,
$$
where
$$
{R_\pre}^{\le2}\subset R_\pre
$$
is the subspace of those linear forms on $R^\pre$ which vanish on all right
preadmissible elements $a[n,m]$ with $a\in\tilde\A$.
\end{Lemma}

\begin{proof}
Dualizing, for the first inclusion what we have to prove is that given any
admissible monomial $a\in\A$ and any $[n,m]b\in\PAR'$ with $b\in\tilde\A$, in
$\bar R$ one has the equality
$$
a[n,m]b=\sum_ia_i[n_i,m_i]b_i
$$
with $a_i[n_i,m_i]\in\PAR$ and admissible monomials $b_i\in\tilde\A$. Indeed,
considering $a$ as a monomial in $\F_0$ there is a unique way to write
$$
a[n,m]=\sum_ia_i[n_i,m_i]c_i
$$
in $\F_0$, with $a_i[n_i,m_i]\in\PAR$ and $c_i$ some (not necessarily
admissible or belonging to $\tilde\F_0$) monomials in the $\Sq^k$ generators
of $\F_0$. Thus in $\F_0$ we have
$$
a[n,m]b=\sum_ia_i[n_i,m_i]c_ib.
$$
In $\bar R$ we may replace each $c_ib$ with a sum of admissible monomials of
the same degree; obviously this degree is positive as $b\in\tilde\A$.

The proof for the second inclusion is exactly similar.
\end{proof}

This lemma implies that for any simultaneous solution $A_\psi(\zeta_n)$ of
\ref{apsi}\eqref{mreqs} and \ref{apsi}\eqref{mleqs}, the elements in $\A_*\ox
R_\pre\ox\A_*$ and $\A_*\ox\A_*\ox R'_\pre$ corresponding to it according to,
respectively, \ref{mrrho} and \ref{mlx}, satisfy
\begin{multline*}
\sum_{\substack{a\in\tilde\A\\{}[k,l]a\in\PAR'}}\left(x_{2^n-k-l-|a|}([k,l]a)\ox1+\tilde
m_*(x_{2^n-k-l-|a|}([k,l]a))+\sum_{i>0}\zeta_i^{2^{n-i}}\ox
x_{2^{n-i}-k-l-|a|}([k,l]a)\right)\ox([k,l]a)_*\\
=(1\ox1\ox\varrho^{>2})(1\ox\Phi_*)\left(\sum_{k>0}\rho_{2^n-2^k+1}(\zeta_k)\ox\zeta_k\right),
\end{multline*}
where
$$
\varrho^{>2}:R_\pre'\onto{R_\pre'}^{>2}
$$
is the restriction of linear forms on ${R^\pre}'$ to the subspace spanned by
the subset of $\PAR'$ consisting of the left preadmissible relations of the
form $[k,l]a$ with $a\in\tilde\A$. Indeed the remaining part of the element
from \ref{mrrho} is
$$
\rho_{2^n}(1)\ox1,
$$
and according to the lemma its image under $1\ox\Phi_*$ goes to zero under the
map $\varrho^{>2}$.

Since the elements $\rho_{2^n-2^k+1}(\zeta_k)$ are explicitly given for all $k>0$, this
allows us to explicitly determine all elements $x_j([k,l]a)$ for
$[k,l]a\in\PAR'$ with $a\in\tilde\A$. For example, in low degrees we obtain
\begin{align*}
x_2([2,3]\Sq^1)=x_2([3,2]\Sq^1)&=\zeta_1^2,\\
x_3([2,2]\Sq^1)&=\zeta_1^3,\\
x_{10}([2,3]\Sq^1)=x_{10}([3,2]\Sq^1)&=\zeta_1^4\zeta_2^2,\\
x_{11}([2,2]\Sq^1)&=\zeta_1^5\zeta_2^2,\\
x_{26}([2,3]\Sq^1)=x_{26}([3,2]\Sq^1)&=\zeta_2^4\zeta_3^2,\\
x_{27}([2,2]\Sq^1)&=\zeta_1\zeta_2^4\zeta_3^2,
\end{align*}
with all other $x_j([k,l]a)=0$ for $j<32$ and $[k,l]a\in\PAR'$ with
$a\in\tilde\A$. 

\begin{Remark}\label{conj>0}
Calculations can be performed for larger $j$ too. But in fact a pattern is
clearly apparent here. It suggests the conjecture that actually all
elements $x_j([k,l]a)$ for $[k,l]a\in\PAR'$ with $a\in\tilde\A$ can be chosen
to be
\begin{align*}
x_{2^n-6}([2,3]\Sq^1)=x_{2^n-6}([3,2]\Sq^1)&=\zeta_{n-3}^4\zeta_{n-2}^2,\\
x_{2^n-5}([2,2]\Sq^1)&=\zeta_1\zeta_{n-3}^4\zeta_{n-2}^2,
\end{align*}
for $n\ge3$, with all other $x_j([k,l]a)=0$.
\end{Remark}

It remains to deal with  the elements $x_j([k,l])$. These shall satisfy
\begin{multline*}
\sum_{k<2l}\left(x_{2^n-k-l}([k,l])\ox1+\tilde
m_*(x_{2^n-k-l}([k,l]))+\sum_{i>0}\zeta_i^{2^{n-i}}\ox
x_{2^{n-i}-k-l}([k,l])\right)\ox[k,l]_*\\
=(1\ox\Phi_*)\left(\rho_{2^n}(1)\ox1\right)+(1\ox1\ox\varrho^{\le2})(1\ox\Phi_*)\left(\sum_{k>0}\rho_{2^n-2^k+1}(\zeta_k)\ox\zeta_k\right),
\end{multline*}
where now
$$
\varrho^{\le2}:R_\pre'\onto{R_\pre'}^{\le2}
$$
is the restriction of linear forms on ${R^\pre}'$ to the subspace spanned by
the Adem relations. The last summand
$D_n=(1\ox1\ox\varrho^{\le2})(1\ox\Phi_*)\left(\sum_{k>0}\rho_{2^n-2^k+1}(\zeta_k)\ox\zeta_k\right)$
is again explicitly given; for example, in low degrees it is equal to
\begin{align*}
D_1&=0,\\
D_2&=0,\\
D_3&=\left(\zeta_1\ox\zeta_1\right)^2\ox[2,2]_*,\\
D_4&=\left(\zeta_1^2\zeta_2\ox\zeta_1+\zeta_2\ox\zeta_2+\zeta_1^2\ox\zeta_1\zeta_2\right)^2\ox[2,2]_*,\\
D_5&=\left(
 \zeta_2^2\zeta_3\ox\zeta_1
+\zeta_1^4\zeta_3\ox\zeta_2
+\zeta_1^4\zeta_2^2\ox\zeta_1\zeta_2
+\zeta_1^4\ox\zeta_2\zeta_3
+\zeta_3\ox\zeta_3
+\zeta_2^2\ox\zeta_1\zeta_3
\right)^2\ox[2,2]_*.
\end{align*}
Then finally the equations that remain to be solved can be equivalently
written as follows:
\begin{multline*}
(1\ox1\ox\tilde\eps)(1\ox\Phi_*)\1\left(\sum_{k<2l}\left(x_{2^n-k-l}([k,l])\ox1+\tilde
m_*(x_{2^n-k-l}([k,l]))+\sum_{i>0}\zeta_i^{2^{n-i}}\ox
x_{2^{n-i}-k-l}([k,l])\right)\ox[k,l]_*\right)\\
=(1\ox1\ox\tilde\eps)(1\ox\Phi_*)\1(D_n),
\end{multline*}
where
$$
\tilde\eps:\A_*\onto\tilde\A_*
$$
is the projection to the positive degree part, i.~e. maps 1 to 0 and all
homogeneous positive degree elements to themselves. Again, the right hand
sides of these equations are explicitly given constants, for example, in low
degrees they are given by

\begin{tabular}{cl}
$0$,&$n=1$;\\
$0$,&$n=2$;\\
$\zeta_1^2\ox[2,2]_*\ox\zeta_1^2$,&$n=3$;\\
$\left(\zeta_1^4\zeta_2^2\ox[2,2]_*+\zeta_2^2\ox(\Sq^4[2,2])_*+\zeta_1^4\ox(\Sq^6[2,2])_*\right)\ox\zeta_1^2$,&$n=4$;\\
$\left(\zeta_2^4\zeta_3^2\ox[2,2]_*+\zeta_1^8\zeta_3^2\ox(\Sq^4[2,2])_*+\zeta_1^8\zeta_2^4\ox(\Sq^6[2,2])_*+\zeta_3^2\ox(\Sq^8\Sq^4[2,2])_*\right.$\\
$\left.+\zeta_2^4\ox(\Sq^{10}\Sq^4[2,2])_*+\zeta_1^8\ox(\Sq^{12}\Sq^6[2,2])_*\right)\ox\zeta_1^2$,&$n=5$.
\end{tabular}

One possible set of solutions for $\zeta_k$ with $k\le5$ is given by
\begin{align*}
x_5([1,2])&=\zeta_1^2\zeta_2,\\
x_4([1,3])&=\zeta_1^4,\\
x_{13}([1,2])&=\zeta_2^2\zeta_3,\\
x_{12}([1,3])&=\zeta_2^4,\\
x_{29}([1,2])&=\zeta_3^2\zeta_4,\\
x_{28}([1,3])&=\zeta_3^4
\end{align*}
and all remaining $x_j([k,l])=0$ for $j+k+l\le32$.

Or equivalently one might give the same solution ``on the other side of $\Phi$'' by
\begin{align*}
\rho_2(1)&=0,\\
\rho_4(1)&=0,\\
\rho_8(1)&
=\zeta_1^2\zeta_2\ox[1,2]_*
+\zeta_1^4\ox[1,3]_*
+\zeta_2\ox(\Sq^2[1,2])_*
+\zeta_1^2\ox(\Sq^3[1,2])_*,\\
\rho_{16}(1)&
=\zeta_2^2\zeta_3\ox[1,2]_*
+\zeta_2^4\ox[1,3]_*
+\zeta_1^4\zeta_3\ox\left(\Sq^2[1,2]\right)_*
+\zeta_1^4\zeta_2^2\ox\left(\Sq^3[1,2]\right)_*\\
&
+\zeta_3\ox\left(\Sq^4\Sq^2[1,2]\right)_*
+\zeta_2^2\ox\left(\Sq^5\Sq^2[1,2]\right)_*
+\zeta_1^4\ox\left(\Sq^6\Sq^3[1,2]\right)_*,\\
\rho_{32}(1)&
=\zeta_3^2\zeta_4\ox[1,2]_*.
+\zeta_3^4\ox[1,3]_*
+\zeta_2^4\zeta_4\ox\left(\Sq^2[1,2]\right)_*
+\zeta_2^4\zeta_3^2\ox\left(\Sq^3[1,2]\right)_*\\
&
+\zeta_1^8\zeta_4\ox\left(\Sq^4\Sq^2[1,2]\right)_*
+\zeta_1^8\zeta_3^2\ox\left(\Sq^5\Sq^2[1,2]\right)_*
+\zeta_1^8\zeta_2^4\ox\left(\Sq^6\Sq^3[1,2]\right)_*\\
&
+\zeta_4\ox\left(\Sq^8\Sq^4\Sq^2[1,2]\right)_*
+\zeta_3^2\ox\left(\Sq^9\Sq^4\Sq^2[1,2]\right)_*
+\zeta_2^4\ox\left(\Sq^{10}\Sq^5\Sq^2[1,2]\right)_*
+\zeta_1^8\ox\left(\Sq^{12}\Sq^6\Sq^3[1,2]\right)_*
\end{align*}

\begin{Remark}\label{conj0}
As in \ref{conj>0}, here one also has a suggestive pattern which leads to a
conjecture that a simultaneous solution of \eqref{mreqs} and
\eqref{mleqs} is determined by putting
\begin{align*}
x_{2^n-3}([1,2])&=\zeta_{n-2}^2\zeta_{n-1},\\
x_{2^n-4}([1,3])&=\zeta_{n-2}^4
\end{align*}
for $n\ge3$, with all other $x_j([k,l])=0$.
\end{Remark}

This then gives the solution itself as follows:
\begin{align*}
A_\psi(\zeta_1)&=0,\\
\\
A_\psi(\zeta_2)&=0,\\
\\
A_\psi(\zeta_3)
=\zeta_1^2\zeta_2&\ox M_3\\
+\zeta_1^4&\ox\left(M_{31}+\zeta_1M_3\right)\\
+\zeta_1^3&\ox M_{221}\\
+\zeta_2&\ox\left(M_5+M_{41}+M_{32}+\zeta_1^2M_3\right)\\
+\zeta_1^2&\ox\left(M_{51}+M_{321}+M_{231}+M_{2121}+\zeta_1(M_{5}+M_{41}+M_{32}+M_{221})+\zeta_1^2M_{11}^2+(\zeta_1^3+\zeta_2)M_3\right)\\
+\zeta_1&\ox M_{2221},\\
\\
A_\psi(\zeta_4)
=\zeta_2^2\zeta_3&\ox M_3\\
+\zeta_2^4&\ox\left(M_{31}+\zeta_1M_3\right)\\
+\zeta_1^5\zeta_2^2&\ox M_{221}\\
+\zeta_1^4\zeta_3&\ox\left(M_5+M_{41}+M_{32}+\zeta_1^2M_3\right)\\
+\zeta_1^4\zeta_2^2&\ox\left(M_{51}+M_{321}+M_{231}+M_{2121}+\zeta_1(M_5+M_{41}+M_{32}+M_{221})+\zeta_1^2M_{11}^2+(\zeta_1^3+\zeta_2)M_3\right)\\
+\zeta_1^9&\ox M_{2221}\\
+\zeta_1\zeta_2^2&\ox M_{4221}\\
+\zeta_3&\ox\left(M_9+M_{72}+M_{621}+M_{54}+M_{441}+M_{432}+M_{342}+M_{2421}+\zeta_1^4M_5+\zeta_2^2M_3\right)\\
+\zeta_2^2&\ox\left(
M_{721}
+M_{451}
+M_{4321}
+M_{4231}
+M_{42121}
+M_{3421}
+(M_5
+M_{41}
+M_{32}
+M_{2111})^2\right.\\
&
\left.
+\zeta_1
(M_9
+M_{72}
+M_{621}
+M_{54}
+M_{441}
+M_{432}
+M_{4221}
+M_{342}
+M_{2421}
)
+\zeta_1^4M_3^2
+\zeta_1^5M_5
+(\zeta_1\zeta_2^2
+\zeta_3)M_3
\right)\\
+\zeta_1^5&\ox\left(M_{6221}+M_{4421}+M_{24221}\right)\\
+\zeta_1^4&\ox\left(
M_{831}
+M_{8121}
+M_{651}
+M_{6321}
+M_{6231}
+M_{62121}
+M_{4521}
+M_{4431}
+M_{44121}
+M_{41421}\right.\\
&
+M_{2721}
+M_{2451}
+M_{24321}
+M_{24231}
+M_{242121}
+M_{23421}\\
&
+\zeta_1(M_{6221}+M_{4421}+M_{24221})
+\zeta_1^2(M_5+M_{41}+M_{32}+M_{2111})^2\\
&
+\zeta_2(M_9+M_{72}+M_{621}+M_{54}+M_{441}+M_{432}+M_{342}+M_{2421})
+\zeta_1^4M_{211}^2
+\zeta_1^6M_3^2
+\zeta_3(M_5+M_{41}+M_{32})\\
&\left.
+\zeta_1^4\zeta_2M_5
+(\zeta_1^2\zeta_3+\zeta_2^3)M_3
\right)\\
+\zeta_1&\ox\left(M_{82221}+M_{44421}+M_{46221}+M_{424221}\right),
\end{align*}
\begin{align*}
A_\psi(\zeta_5)
=\zeta_3^2\zeta_4&\ox M_3\\
+\zeta_3^4&\ox\left(M_{31}+\zeta_1M_3\right)\\
+\zeta_1\zeta_2^4\zeta_3^2&\ox M_{221}\\
+\zeta_2^4\zeta_4&\ox\left(M_5+M_{41}+M_{32}+\zeta_1^2M_3\right)\\
+\zeta_2^4\zeta_3^2&\ox\left(M_{51}+M_{321}+M_{231}+M_{2121}+\zeta_1(M_5+M_{41}+M_{32}+M_{221})+\zeta_1^2M_{11}^2+(\zeta_1^3+\zeta_2)M_3\right)\\
+\zeta_1\zeta_2^8&\ox M_{2221}\\
+\zeta_1^9\zeta_3^2&\ox M_{4221}\\
+\zeta_1^8\zeta_4&\ox\left(M_9+M_{72}+M_{621}+M_{54}+M_{441}+M_{432}+M_{342}+M_{2421}+\zeta_1^4M_5+\zeta_2^2M_3\right)\\
+\zeta_1^8\zeta_3^2&\ox\left(
M_{721}
+M_{451}
+M_{4321}
+M_{4231}
+M_{42121}
+M_{3421}
+(M_5
+M_{41}
+M_{32}
+M_{2111})^2\right.\\
&
+\zeta_1
(M_9
+M_{72}
+M_{621}
+M_{54}
+M_{441}
+M_{432}
+M_{4221}
+M_{342}
+M_{2421}
)\\
&\left.
+\zeta_1^4M_3^2
+\zeta_1^5M_5
+(\zeta_1\zeta_2^2
+\zeta_3)M_3
\right)\\
+\zeta_1^9\zeta_2^4&\ox\left(M_{6221}+M_{4421}+M_{24221}\right)\\
+\zeta_1^8\zeta_2^4&\ox\left(
M_{831}
+M_{8121}
+M_{651}
+M_{6321}
+M_{6231}
+M_{62121}
+M_{4521}
+M_{4431}
+M_{44121}
+M_{41421}\right.\\
&
+M_{2721}
+M_{2451}
+M_{24321}
+M_{24231}
+M_{242121}
+M_{23421}\\
&
+\zeta_1(M_{6221}+M_{4421}+M_{24221})
+\zeta_1^2(M_5+M_{41}+M_{32}+M_{2111})^2\\
&
+\zeta_2(M_9+M_{72}+M_{621}+M_{54}+M_{441}+M_{432}+M_{342}+M_{2421})
+\zeta_1^4M_{211}^2
+\zeta_1^6M_3^2\\
&\left.
+\zeta_1^4\zeta_2M_5+\zeta_3(M_5+M_{41}+M_{32})
+(\zeta_1^2\zeta_3+\zeta_2^3)M_3
\right)\\
+\zeta_1^{17}&\ox\left(M_{82221}+M_{44421}+M_{46221}+M_{424221}\right)\\
+\zeta_4&\ox\left(
M_{\underline{17}}
+M_{\underline{13}4}
+M_{\underline{11}42}
+M_{\underline{10}421}
+M_{98}
+M_{872}
+M_{8621}
+M_{854}
+M_{8441}
+M_{8432}
+M_{8342}
+M_{82421}\right.\\
&
\left.+M_{584}
+M_{3842}
+M_{28421}
+\zeta_1^8M_9
+\zeta_2^4M_5
+\zeta_3^2M_3
\right)\\
+\zeta_1\zeta_3^2&\ox M_{84221}\\
+\zeta_3^2&\ox\left(
M_{\underline{11}421}
+M_{8721}
+M_{8451}
+M_{84321}
+M_{84231}
+M_{842121}
+M_{83421}
+M_{38421}\right.\\
&+
(M_9
+M_{72}
+M_{621}
+M_{54}
+M_{441}
+M_{432}
+M_{42111}
+M_{342}
+M_{2421}
)^2\\
&
+\zeta_1
(M_{\underline{17}}
+M_{\underline{13}4}
+M_{\underline{11}42}
+M_{\underline{10}421}
+M_{98}\\
&\quad
+M_{872}
+M_{8621}
+M_{854}
+M_{8441}
+M_{8432}
+M_{84221}
+M_{8342}
+M_{82421}
+M_{584}
+M_{3842}
+M_{28421})\\
&\left.
+\zeta_1^8M_5^2
+\zeta_1^9M_9
+\zeta_2^4M_3^2
+\zeta_1\zeta_2^4M_5
+(\zeta_1\zeta_3^2+\zeta_4)M_3
\right)\\
+\zeta_1\zeta_2^4&\ox\left(M_{\underline{10}4221}+M_{86221}+M_{84421}+M_{824221}+M_{284221}\right)\\
+\zeta_2^4&\ox\left(
M_{\underline{12}521}
+M_{\underline{12}431}
+M_{\underline{12}4121}
+M_{\underline{12}1421}
+M_{\underline{10}721}
+M_{\underline{10}451}
+M_{\underline{10}4321}
+M_{\underline{10}4231}
+M_{\underline{10}42121}
+M_{\underline{10}3421}\right.\\
&
+M_{8831}
+M_{88121}
+M_{8651}
+M_{86321}
+M_{86231}
+M_{862121}
+M_{84521}
+M_{84431}
+M_{844121}
+M_{841421}\\
&
+M_{82451}
+M_{82721}
+M_{823421}
+M_{824321}
+M_{824231}
+M_{8242121}\\
&
+M_{49421}
+M_{48521}
+M_{48431}
+M_{484121}
+M_{481421}
+M_{418421}\\
&
+M_{2\underline{11}421}
+M_{28721}
+M_{28451}
+M_{284321}
+M_{284231}
+M_{2842121}
+M_{283421}
+M_{238421}\\
&
+\zeta_1(M_{\underline{10}4221}+M_{86221}+M_{84421}+M_{824221}+M_{284221})\\
&
+\zeta_1^2(M_9+M_{72}+M_{621}+M_{54}+M_{441}+M_{432}+M_{42111}+M_{342}+M_{2421})^2\\
&
+\zeta_2
(M_{\underline{17}}
+M_{\underline{13}4}
+M_{\underline{11}42}
+M_{\underline{10}421}
+M_{98}
+M_{872}
+M_{8621}
+M_{854}
+M_{8441}
+M_{8432}
+M_{8342}
+M_{82421}\\
&\quad
+M_{584}
+M_{3842}
+M_{28421})\\
&\left.
+\zeta_1^4M_{4211}^2
+\zeta_1^{10}M_5^2
+\zeta_1^8\zeta_2M_9
+\zeta_1^2\zeta_2^4M_3^2
+\zeta_2^5M_5
+\zeta_4(M_5+M_{41}+M_{32})
+(\zeta_1^2\zeta_4+\zeta_2\zeta_3^2)M_3
\right)
\end{align*}
\begin{align*}
+\zeta_1^9&\ox\left(M_{\underline{12}6221}+M_{\underline{12}4421}+M_{\underline{12}24221}+M_{4\underline{10}4221}+M_{88421}+M_{486221}+M_{484421}+M_{4824221}+M_{4284221}\right)\\
+\zeta_1^8&\ox\left(
M_{\underline{14}631}
+M_{\underline{14}6121}
+M_{\underline{14}2521}
+M_{\underline{14}2431}
+M_{\underline{14}24121}
+M_{\underline{14}21421}\right.\\
&
+M_{\underline{12}831}
+M_{\underline{12}8121}
+M_{\underline{12}651}
+M_{\underline{12}6321}
+M_{\underline{12}6231}
+M_{\underline{12}62121}
+M_{\underline{12}4521}
+M_{\underline{12}4431}
+M_{\underline{12}44121}
+M_{\underline{12}41421}\\
&
+M_{\underline{12}2721}
+M_{\underline{12}2451}
+M_{\underline{12}24321}
+M_{\underline{12}24231}
+M_{\underline{12}242121}
+M_{\underline{12}23421}\\
&
+M_{86631}
+M_{866121}
+M_{862521}
+M_{862431}
+M_{8624121}
+M_{8621421}\\
&
+M_{844521}
+M_{844431}
+M_{8444121}
+M_{8441421}
+M_{842631}
+M_{8426121}
+M_{8423421}
+M_{84212421}\\
&
+M_{6\underline{10}521}
+M_{6\underline{10}431}
+M_{6\underline{10}4121}
+M_{6\underline{10}1421}
+M_{68631}
+M_{686121}
+M_{682521}
+M_{682431}
+M_{6824121}
+M_{6821421}\\
&
+M_{629421}
+M_{628521}
+M_{628431}
+M_{6284121}
+M_{6281421}
+M_{6218421}\\
&
+M_{4\underline{12}521}
+M_{4\underline{12}431}
+M_{4\underline{12}4121}
+M_{4\underline{12}1421}
+M_{4\underline{10}721}
+M_{4\underline{10}451}
+M_{4\underline{10}4321}
+M_{4\underline{10}4231}
+M_{4\underline{10}42121}
+M_{4\underline{10}3421}\\
&
+M_{48831}
+M_{488121}
+M_{48651}
+M_{486321}
+M_{486231}
+M_{4862121}
+M_{484521}
+M_{484431}
+M_{4844121}
+M_{4841421}\\
&
+M_{482721}
+M_{482451}
+M_{4824321}
+M_{4824231}
+M_{48242121}
+M_{4823421}\\
&
+M_{449421}
+M_{448521}
+M_{448431}
+M_{4484121}
+M_{4481421}
+M_{4418421}\\
&
+M_{42\underline{11}421}
+M_{428721}
+M_{428451}
+M_{4284321}
+M_{4284231}
+M_{42842121}
+M_{4283421}
+M_{4238421}\\
&
+(M_{831}
+M_{8121}
+M_{7311}
+M_{7221}
+M_{71211}
+M_{651}
+M_{6411}
+M_{6321}
+M_{63111}
+M_{62211}
+M_{612111}\\
&\quad
+M_{43311}
+M_{43221}
+M_{431211}
+M_{422211}
+M_{421311}
+M_{421221}
+M_{41421}\\
&\quad
+M_{35211}
+M_{34311}
+M_{34221}
+M_{341211}
+M_{314211}\\
&\quad
+M_{2721}
+M_{26211}
+M_{252111}
+M_{2451}
+M_{24411}
+M_{24321}
+M_{243111}
+M_{242211}
+M_{2412111}\\
&\quad
+M_{23421}
+M_{224211}
+M_{2142111})^2\\
&
+(M_{51}+M_{411}+M_{321})^4+M_3^8\\
&
+\zeta_1(M_{\underline{12}6221}+M_{\underline{12}4421}+M_{\underline{12}24221}+M_{88421}+M_{4\underline{10}4221}+M_{486221}+M_{484421}+M_{4824221}+M_{4284221})\\
&
+\zeta_1^4\left(M_5+M_{41}+M_{32}\right)^4\\
&
+\zeta_2^2\left(M_9+M_{72}+M_{621}+M_{54}+M_{441}+M_{432}+M_{342}+M_{2421}\right)^2\\
&
+\zeta_3(M_{\underline{17}}+M_{\underline{13}4}+M_{\underline{11}42}+M_{\underline{10}421}+M_{98}+M_{872}+M_{8621}+M_{854}+M_{8441}+M_{8432}+M_{8342}+M_{82421}\\
&\quad
+M_{584}+M_{3842}+M_{28421})\\
&
+\zeta_1^8\zeta_2^2M_5^2+\zeta_3^2(M_5+M_{41}+M_{32})^2+\zeta_1^8\zeta_3M_9\\
&
+\zeta_4(M_9+M_{72}+M_{621}+M_{54}+M_{441}+M_{432}+M_{342})\\
&
\left.+(\zeta_1^{12}+\zeta_2^4)M_3^4+(\zeta_1^4\zeta_3^2+\zeta_2^6)M_3^2+(\zeta_1^4\zeta_4+\zeta_2^4\zeta_3)M_5+(\zeta_2^2\zeta_4+\zeta_3^3)M_3
\right)\\
+\zeta_1&\ox\left(
M_{\underline{16}82221}
+M_{\underline{16}46221}
+M_{\underline{16}44421}
+M_{\underline{16}424221}\right.\\
&\left.
+M_{8\underline{12}4421}
+M_{8\underline{12}6221}
+M_{8\underline{12}24221}
+M_{888421}
+M_{84\underline{10}4221}
+M_{8486221}
+M_{8484421}
+M_{84824221}
+M_{84284221}
\right)
\end{align*}

The formul\ae\ above were obtained via computer calculations. They lead to the
general patterns in \ref{conj>0} and \ref{conj0} which would determine the map
$A_\psi$ completely.

%% file: dHascoAss_c9.tex
\chapter{The dual $d_{(2)}$ differential
}\label{diff}

In this chapter we will compute the $d_{(2)}$ differential in the $\E^2$ term
$$
\E_2^{p,q}=\Cotor^p_{\A_*}(\FF,\FF)^q\cong\Ext^p_{\A}(\FF,\FF)^q
$$
of the Adams spectral sequence. For this we will first set up algebraic formalism necessary to carry out an analog of the computations in Chapter \ref{E3} in the dual setting. First let us recall how the above isomorphism is obtained.

\section{Secondary coresolution}

One starts with a projective resolution of the $\A$-module $\FF$, e.~g. with the minimal resolution as in \eqref{minireso}. Its graded $\FF$-linear dual
\begin{equation}\label{miniresod}
\FF\to\A_*^\set{g_0^0}\to\bigoplus_{n\ge0}\A_*^\set{g_1^{2^n}}\to\bigoplus_{|i-j|\ne1}\A_*^\set{g_2^{2^i+2^j}}\to...
\end{equation}
is then an injective resolution of $\FF$ in the category of right $\A_*$-comodules. (This is not entirely trivial since we take \emph{graded} duals. However all (co)modules that we encounter will be degreewise finite, i.~e. having generating sets with finite number of elements in each degree. Obviously then graded duality is a contravariant equivalence between the categories of such (co)modules.)

There are isomorphisms
$$
\Hom_\A(M,N)\cong M_*\cote_{\A_*}N
$$
for any left $\A$-modules $M$ and $N$ of the above kind (i.~e. of graded finite type), where on the right the graded dual $M_*$ is considered as a right $\A_*$-comodule and $N$ as a left $\A_*$-comodule in the standard way. It follows that applying $\Hom_\A(-,\FF)$ to \eqref{minireso} and applying $-\cote_{\A_*}\FF$ to \eqref{miniresod} gives isomorphic cochain complexes (of $\FF$-vector spaces). But by definition cohomology of the latter complex is given by
$$
H^p(\eqref{miniresod}\cote_{\A_*}\FF)^q=\Cotor^p_{\A_*}(\FF,\FF)^q.
$$

It then follows from \eqref{d2gen} that in these terms the secondary differential
$$
d_{(2)}^{pq}:\Cotor^p_{\A_*}(\FF,\FF)^q\to\Cotor^{p+2}_{\A_*}(\FF,\FF)^{q+1}
$$
is given by
\begin{equation}\label{d2delta}
d_{(2)}^{pq}(\hat g_p^q)=\sum_{\textrm{$g_p^q$ appears in
$\delta(g_{p+2}^{q+1})$}}\hat g_{p+2}^{q+1}=\delta_*(\hat g_p^q)^0.
\end{equation}
Here,
$$
\delta_*:\bigoplus_q\Sigma\A_*^{\set{g_p^q}}\to\bigoplus_q\A_*^{\set{g_{p+2}^q}}
$$
is the dual of the map
$$
\delta:\A\brk{g_{p+2}^*}\to\Sigma\A\brk{g_p^*}
$$
determined in \ref{delta}, whereas $\hat g_*^*$ denotes the dual basis of $g_*^*$, i.~e. $\hat g_p^q\in\A_*^\set{g_*^*}$ is the vector with the $g_p^q$-th coordinate equal to 1 and all other coordinates equal to zero. Moreover by $\delta_*(\hat g_p^q)^0$ is denoted the zero degree component of $\delta_*(\hat g_p^q)$, i.~e. the result of applying to the element
$$
\delta_*(\hat g_p^q)\in\bigoplus_{j\ge0}\A_j^\set{g_{p+2}^{q+j+1}}
$$
the projection to the $(j=0)$-th component
$$
\bigoplus_{j\ge0}\A_j^\set{g_{p+2}^{q+j+1}}\to\A_0^\set{g_{p+2}^{q+1}}.
$$

Instead of directly dualizing the map $\delta$, it is more convenient from the computational point of view to dualize the conditions of \ref{delta} using \eqref{diadelta} and determine $\delta_*$ directly from these dualized conditions. In fact using \ref{exmulf} we can further detalize the diagram \eqref{diadelta} in the following way:
\begin{equation}\label{diadeltad}
\alignbox{
\xymatrix{
&\Sigma\A\ox V_{p+1}\ar[r]^-{1\ox d}&\Sigma\A\ox\A\ox V_p\ar[dr]^{m\ox1}\\
V_{p+3}\ar[ur]^{\delta^\A_{p+1}}\ar[dr]_d&\A\ox\A\ox V_p\ar[ur]^{\k\ox1\ox1}&\A\ox R_\F\ox V_p\ar[r]^{1\ox A^s}&\Sigma\A\ox V_p\\
&\A\ox V_{p+2}\ar[u]_{1\ox\ph^{\A,s}}\ar[ur]_{1\ox\ph^{R,s}}\ar[r]_-{1\ox\delta^\A_p}&\A\ox\Sigma\A\ox V_p\ar[ur]_{m\ox1}
}
}
\end{equation}
where $A^s$ is the multiplication map corresponding to a splitting $s$ of the
$\GG$-relation pair algebra used, as in \ref{rcomp}, to identify $R_\aB$ with
$\A\oplus R_\F$, and $(\ph^{\A,s},\ph^{R,s})$ are the components of the corresponding composite map
$$
V_{p+2}\xto\ph R_\aB\ox V_p=\A\!\ox\!V_p\oplus R_\F\!\ox\!V_p,
$$
with $\ph$ as defined in \eqref{dd}.

Moreover just as the map $\delta$ is completely determined by its restriction to $V_{p+2}$, its dual $\delta_*$ is determined by the composite $\delta_0$ as in
$$
\Hom(V_p,\Sigma\A_*)\xto{\delta_*}\Hom(V_{p+2},\A_*)\xto{\Hom(V_{p+2},\eps)}\Hom(V_{p+2},\FF),
$$
where graded $\Hom$ is meant, and $\eps$ is the augmentation of $\A_*$. In fact we only need this composite map $\delta_0$ as by \eqref{d2delta} above we have
\begin{equation}\label{d2delta0}
d_{(2)}^{pq}(\hat g_p^q)=\delta_0(\hat g_p^q).
\end{equation}

Now the dual to diagram \eqref{diadeltad} is easy to identify; it is
\begin{equation}\label{diadeltadu}
\alignbox{
\xymatrix{
&\Sigma\A_*\ox\hat V_{p+1}\ar[dl]_{\delta_0}&\Sigma\A_*\ox\A_*\ox\hat V_p\ar[l]_-{1\ox d_*}\ar[dl]_{\zeta_1\ox1\ox1}\\
\hat V_{p+3}&\A_*\ox\A_*\ox\hat V_p\ar[d]_{1\ox\ph_*^{\A,s}}&\A_*\ox{R_\F}_*\ox\hat V_p\ar[dl]_{1\ox\ph_*^{R,s}}&\Sigma\A_*\ox\hat V_p\ar[ul]_{m_*\ox1}\ar[l]_-{A_s\ox1}\ar[dl]^{m_*\ox1}\\
&\A_*\ox\hat V_{p+2}\ar[ul]^{d_*}&\A_*\ox\Sigma\A_*\ox\hat V_p\ar[l]^-{1\ox\delta_0}
}
}
\end{equation}
where $\hat V_p$ are the graded dual spaces of $V_p$.

It is straightforward to reformulate the above in terms of elements: the
values of the map $\delta_0$ on arbitrary elements $a\ox g\in\Sigma\A_*\ox\hat
V_p$ must satisfy
\begin{equation}\label{deltamain}
\begin{aligned}
\delta_0(\sum a_\l\ox d_*(a_\r\ox g))
&=d_*(\sum a_l\ox\delta_0(a_\r\ox g))\\
&+d_*(\sum\zeta_1a_\l\ox\ph_*^{\A,s}(a_\r\ox g))
+d_*(\sum a_\A\ox\ph_*^{R,s}(a_R\ox g)),
\end{aligned}
\end{equation}
where we have denoted by
$$
\Delta(a)=\sum a_\l\ox a_\r
$$
the value of the diagonal $\Delta:\A_*\to\A_*\ox\A_*$ and by
$$
A_s(a)=\sum a_\A\ox a_R
$$
the value of the comultiplication map $A_s:\Sigma\A_*\to\A_*\ox{R_\F}_*$ on
$a\in\A_*$.

We thus obtain

\begin{Proposition}
\label{proposi}
The $d_{(2)}$ differential of the Adams spectral sequence is given on the
cohomology classes represented by the generators $\hat g$ in the minimal resolution by the formula
$$
d_{(2)}(\hat g)=\delta_0(\Sigma1\ox\hat g),
$$
where
$$
\delta_0:\Sigma\A_*\ox\hat V_s\to\hat V_{s+2}
$$
are any maps satisfying the equations \eqref{deltamain}.
\end{Proposition}\qed

At this point the cooperation of the authors ended since the time of Jibladze's visit at the 
MPIM was over. Therefore our goal of doing computer calculations on the basis of \ref{proposi} 
is left to an interested reader.